\documentclass[a4paper,11pt]{amsbook}
\usepackage{thesis}

\usepackage[inner=4cm,outer=3.5cm,top=2.5cm,bottom=2.5cm,includefoot,includehead]{geometry}

\numberwithin{equation}{chapter}
 
\date{October 2014}
\title{Cohomological Finiteness Properties of Groups}

\newcommand\MFMor[1]{\left(\vcenter{\xymatrix@R-20pt@C-23pt{#1}}\right)}

\makeindex

\begin{document}

\frontmatter

\thispagestyle{empty}
\begin{center}
\includegraphics[width=0.4\textwidth]{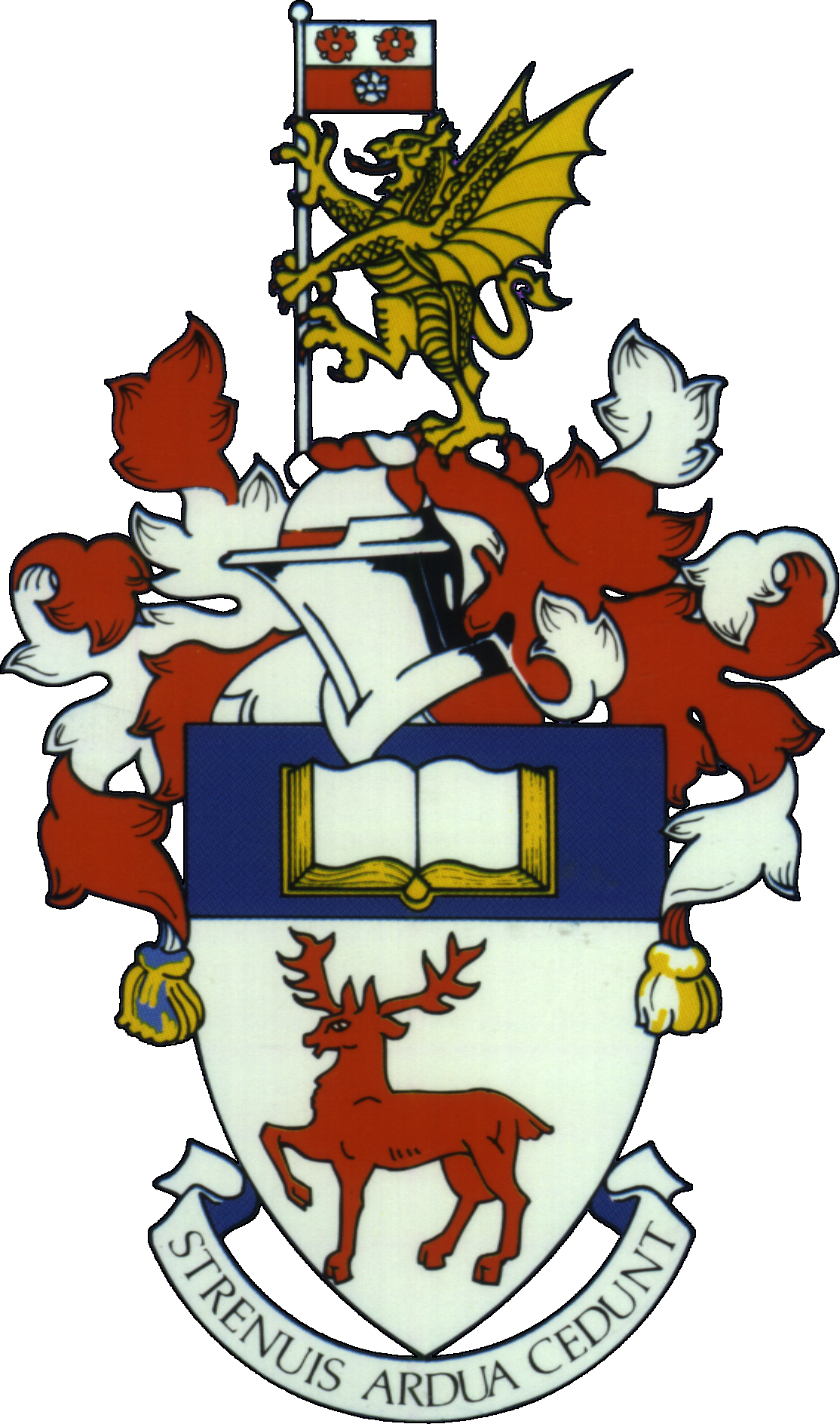}

\vspace{3ex}

\textsc{\LARGE University of Southampton}

\vspace{4ex}

\textsc{\Large Faculty of Social and Human Sciences}

\vspace{2ex}

\textsc{\Large School of Mathematics}

\vspace{7ex}

\makeatletter
{\Huge\sc \@title\par}
\makeatother

\vspace{10ex}

{\Large Simon St.~John-Green}

\vspace{5ex}

Supervisor: Prof. Brita Nucinkis

\vspace{5ex}

Thesis for the degree of Doctor of Philosophy

\vspace{1ex}

\makeatletter
\@date
\makeatother

\end{center}
\cleardoublepage











\thispagestyle{empty}
\begin{center}\large

\textsc{University of Southampton}

\textsc{\underline{Abstract}}

\textsc{Faculty of Social and Human Sciences}

{\Large School of Mathematics}

\underline{Doctor of Philosophy}

\makeatletter
\textsc{\Large \@title}
\makeatother

{by Simon St.~John-Green}
\end{center}

\vspace{2ex}

{\linespread{1}\selectfont

\noindent The main objects of study in this thesis are cohomological finiteness conditions of discrete groups.  While most of the conditions we investigate are algebraic, they are inspired by topological invariants, particularly those concerning proper actions on CW-complexes.

The first two chapters contain preliminary material necessary for the remainder of the thesis.  Chapter \ref{chapter:C} concerns modules over a category with an emphasis on finiteness conditions.  This material is well-known for (EI) categories, but we use a more general setup applicable to Mackey and cohomological Mackey functors, needed in Chapter \ref{chapter:M}.  Chapter \ref{chapter:bredon} specialises to Bredon cohomology, giving an overview of some results and detailing a few interesting examples.

In Chapter \ref{chapter:M} we study finiteness conditions associated to Bredon cohomology with coefficients restricted to Mackey functors and cohomological Mackey functors, building again on the material in Chapter \ref{chapter:C}.  In particular we characterise the corresponding $\FP_n$ conditions and prove that the Bredon cohomological dimension with coefficients restricted to cohomological Mackey functors is equal to the $\mathfrak{F}$-cohomological dimension for all groups.

We prove in Chapter \ref{chapter:G} that for groups of finite $\mathfrak{F}$-cohomological dimension, the $\mathfrak{F}$-cohomological dimension equals the Gorenstein cohomological dimension, and give an application to the behaviour of the $\mathfrak{F}$-cohomological dimension under group extensions.

If a group $G$ admits a closed manifold model for $\B G$ then $G$ is a Poincar\'e duality group, in Chapter \ref{chapter:BD} we study Bredon--Poincar\'e duality groups, a generalisation of these.  In particular if $G$ admits a cocompact manifold model $X$ for $\EFin G$ (the classifying space for proper actions) with $X^H$ a submanifold for any finite subgroup $H$ of $G$, then $G$ is a Bredon--Poincar\'e duality group.  We give several sources of examples, including using the reflection group trick of Davis to produce examples where the dimensions of the submanifolds $X^H$ are specified.  We classify Bredon--Poincar\'e duality groups in low dimensions and examine their behaviour under group extensions.

In Chapter \ref{chapter:H} we study Houghton's group $H_n$, calculating the centralisers of virtually cyclic subgroups and the Bredon cohomological dimension with respect to both the family of finite subgroups and the family of virtually cyclic subgroups.

}
\cleardoublepage

\setcounter{tocdepth}{3}
\tableofcontents
\cleardoublepage

\listoffigures
\cleardoublepage

\chapter*{Declaration of Authorship}



\noindent I, Simon St.~John-Green, declare that this thesis and the work presented in it are my own and has been generated by me as the result of my own original research.

\makeatletter
\noindent Title: \@title
\makeatother

\smallskip

\noindent I confirm that:
\begin{enumerate}
 \item This work was done wholly or mainly while in candidature for a research degree at this University;
 \item Where any part of this thesis has previously been submitted for a degree or any other qualification at this University or any other institution, this has been clearly stated;
 \item Where I have consulted the published work of others, this is always clearly attributed;
 \item Where I have quoted from the work of others, the source is always given. With the exception of such quotations, this thesis is entirely my own work;
 \item I have acknowledged all main sources of help;
 \item Where the thesis is based on work done by myself jointly with others, I have made clear exactly what was done by others and what I have contributed myself;
 \item Either none of this work has been published before submission, or parts of this work have been published as: 
 \begin{enumerate}
  \item Centralisers in Houghton's groups (2012, to appear Proc.~Edinburgh Math.~Soc.) \cite{Me-HoughtonsGroups}.
  \item Bredon--Poincar\'e duality groups (2013, to appear J.~Group Theory) \cite{Me-BredonPoincareDuality}.
  \item Finiteness conditions for Mackey and cohomological Mackey functors (J.~Algebra \textbf{411} (2014), no.~0, 225--258) \cite{Me-CohomologicalMackey}
  \item On the Gorenstein and $\mathfrak{F}$-cohomological dimensions (2013, to appear Bull.~Lond.~Math.~Soc.) \cite{Me-GorensteinAndF}.
 \end{enumerate}
\end{enumerate}

\smallskip

\noindent\begin{tabular}{@{}l l}Signed:&\raisebox{-0.7\height}{\includegraphics[width=40ex]{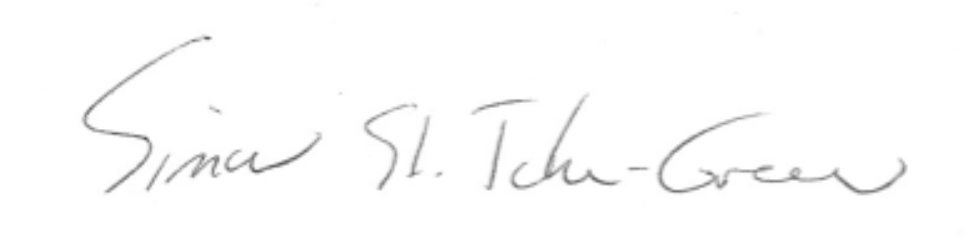}}\end{tabular}

\makeatletter
\noindent Date: \@date.
\makeatother
\chapter*{Acknowledgements}

First and foremost I'd like to thank my supervisor Brita Nucinkis, for introducing me to this beautiful area of mathematics and for being an endless source of motivation, support, and guidance throughout my studies.  

I'd like to thank my advisor Ian Leary for his support and for many useful mathematical insights, including the examples in Section \ref{section:reflection groups} and Example \ref{example:DicksLeary example PD over R not Z}.  I thank him as well for reading a first draft of this thesis. 

For making the last three years such a wonderful experience, I thank all the PhD students and post-docs, past and present, from the Southampton mathematics department, including but not limited to: Alex Bailey, Matt Burfitt, Valerio Capraro, Chris Cave, Yuyen Chien, Charles Cox, Dennis Dreesen, Michal Ferov, Robin Frankhuizen, Tom Harris, Pippa Hiscock, Ana Khukhro, Dave Matthews, Ingrid Membrillo, Raffaele Rainone, Amin Saied, Conor Smyth, Rob Snocken, Joe Tait, Edwin Tye, and Mike West.

The material in Sections \ref{subsection:houghton cdVCyc LW} and \ref{subsection:houghton cdVCyc calc} is joint work with Nansen Petrosyan and I thank him for suggesting the topic and for sharing with me some fascinating mathematics.

Thanks to Dieter Degrijse, not only for pointing out to me an error in Lemma \ref{lemma:M fg and N proj then nu iso} while reading a preprint of \cite{Me-CohomologicalMackey}, but also for many interesting conversations concerning Mackey and cohomological Mackey functors. 

I'd like to thank Conchita Mart{\'{\i}}nez-P\'{e}rez for her very helpful comments concerning a preprint of \cite{Me-HoughtonsGroups}, material now contained in Chapter \ref{chapter:H}.   Thanks to the referee of the Proceedings of the Edinburgh Mathematical Society who reviewed the same preprint for making useful comments, including suggesting the graph $\Gamma$ used in Section \ref{section:centralisersOfElements}.  Thanks to Jim Davis for showing me that some groups constructed by Block and Weinberger, appearing now in Section \ref{section:counterexample to generalised PDn}, answered a question in an early preprint of \cite{Me-BredonPoincareDuality}.
Simplifications to the proofs of Lemmas \ref{prop:mackey OFFPn implies MFFPn} and \ref{lemma:HF ind uROF is R} are due to the referee from the Journal of Algebra.

Finally but most importantly, I thank my family, especially my parents and my sister.  I would dedicate this thesis to them, only then they might feel compelled to read it, and I would not want to force that upon anyone.

\mainmatter

\chapter{Introduction}\label{section:Introduction}

This introduction contains an overview of relevant background material and details the contributions of this thesis as they arise.

\section{Free actions and group cohomology}

For any group $G$ there exists an aspherical CW-complex $X$ with fundamental group $G$, this is called a \emph{model for $\B G$} or \emph{Eilenberg--Mac Lane} space.  Such a space is unique up to homotopy equivalence, a fact observed essentially by Hurewicz \cite{Hurewicz-BeitrageZurTopologieDerDeformationen}, arguably kick-starting the field of group cohomology.  The universal cover $\tilde{X}$, a contractible CW-complex with a free $G$-action, is called \emph{a model for $\E G$} or a \emph{classifying space for free actions}.\index{Model for $\E G$}\index{Model for $\B G$}    Equivalently one could define a model for $\E G$ as the terminal object in the $G$-homotopy category of free $G$-CW-complexes.

One can use invariants of these spaces to study the groups themselves, for example defining the group cohomology $H^*(G)$ to be $H^*(\B G)$.  Alternatively there is an algebraic definition of group cohomology, replacing the space $\B G$ with a resolution of $\ZZ$ by projective $\ZZ G$-modules.  

An important invariant is the \emph{geometric dimension}  $\gd G$, the minimal dimension of a model for $\E G$.  Its algebraic counterpart is the \emph{cohomological dimension} $\cd G$, the minimal length of a resolution of $\ZZ$ by projective $\ZZ G$-modules.\index{Cohomological dimension $\cd G$}\index{Geometric dimension $\gd G$}  It's easy to see that $\gd G = 0$ if and only if $\cd G = 0$ if and only if $G$ is the trivial group.  Also, by a theorem of Stallings and Swan, $\cd G  = 1$ if and only if $\gd G = 1$ if and only if $G$ is a free group \cite{Stallings-OnTorsionFreeGroupsWithInfinitelyManyEnds,Swan-GroupsOfCohomologicalDimensionOne}.
Eilenberg and Ganea conjectured that $\cd G = \gd G$ for all groups and, along with Stallings and Swan's result for the dimension one case, proved this conjecture for all cases, except for the possibility that $\cd G = 2$ and $\gd G = 3$ \cite{EilenbergGanea-Lusternik-SchnirelmannCategory}.  That this is impossible is still an open problem, known as the Eilenberg--Ganea conjecture.

A group $G$ has \emph{type $\F_n$} if it admits a model for $\B G$ with finite $n$-skeleton, and on the algebraic side $G$ has \emph{type $\FP_n$} if $\ZZ$ admits a resolution by projective $\ZZ G$-modules, finitely generated up to dimension $n$.  A group of type $\F_n$ is necessarily of type $\FP_n$.  All groups are of type $\F_0$, since there always exists a model for $\B G$ with a single $0$-cell \cite[7.1.5]{Geoghegan}.  The conditions $\F_1$, $\FP_1$ and finitely generated are all equivalent, but the situation is more complex for larger $n$.  A group is $\F_2$ if and only if it is finitely presented, and $\FP_n$ together with $\F_2$ implies $\F_n$ \cite[7.2.1]{Geoghegan}\cite[VIII \S 7]{Brown}.  Bestvina and Brady use discrete Morse theory techniques to construct subgroups of right-angled Artin groups that are $\FP_n$ but not $\FP_{n+1}$ for all $n$, and groups of type $\FP_n$ which are not finitely presented for all $n$ \cite{BestvinaBrady-MorseTheoryAndFinitenessPropertiesOfGroups}.

We say a group is \emph{type $\F$} if it is $\F_\infty$ and $\gd G < \infty$, and \emph{type $\FP$} of it is $\FP_\infty$ and $\cd G < \infty$.

For a more detailed overview of these finiteness conditions see \cite[Chapter VIII]{Brown}, \cite{Bieri-HomDimOfDiscreteGroups} and \cite[Chapter II]{Geoghegan}.

\section{Proper actions and Bredon cohomology}

Let $\mathcal{F}$ be a family of subgroups of a group $G$, closed under conjugation and taking subgroups.\index{Family of subgroups $\mathcal{F}$}   A \emph{model for $\EF G$}, or classifying space for actions with isotropy in $\mathcal{F}$, is the terminal object in the $G$-homotopy category of $G$-CW complexes with isotropy in $\mathcal{F}$.     \index{Model for $\EF G$}
A model for $\E G$ is thus the same as a model for $\E_{\Triv} G$, where $\Triv$ denotes the family consisting of only the trivial subgroup. \index{Triv@$\Triv$, family of only the trivial subgroup}

Using the equivariant Whitehead theorem one can show that a $G$-CW-complex $X$ is a model for $\EF G$ if and only if $X$ has isotropy in $\mathcal{F}$ and $X^H$ is contractible for all $H \in \mathcal{F}$ \cite[Theorem 1.9]{Luck-SurveyOnClassifyingSpaces}.  Models for $\EF G$ always exist---there are various standard constructions including the infinite join construction of Milnor \cite{Milnor-ConstructionOfUniversalBundles}, Segals construction \cite{Segal-ClassifyingSpacesAndSpectralSequences}, and a construction where one iteratively attaches equivariant cells to build a $G$-CW-complex with contractible fixed point sets \cite[Proposition 2.3, p.35]{Lueck}.  

Let $\Fin$ denote the family of all finite subgroups of a group $G$.\index{Fin@$\Fin$, family of all finite subgroups}  There are many groups which admit natural models for $\EFin G$, for example mapping class groups, word-hyperbolic groups, and one-relator groups.  A good survey is \cite{Luck-SurveyOnClassifyingSpaces}.

There has been recent interest in models for $\EFin G$ and models for $\EVCyc G$, where $\VCyc$ denotes the family of virtually cyclic subgroups, because they appear on one side of the Baum--Connes and Farrell--Jones conjectures respectively \cite{LuckReich-BaumConnesAndFarrellJones}.\index{VCyc@$\VCyc$, family of all virtually cyclic subgroups}  These are deep conjectures with far reaching consequences in mathematics \cite{MislinValette-BaumConnes,BartelsLuckReich-FarrellJonesConjectureAndApplications}.

We denote by $\gdF G$ the minimal dimension of a model for $\EF G$, if $\mathcal{F} = \Fin$ then this is known as the \emph{proper geometric dimension} of $G$.  The cohomology theory most suited to the study of this geometric invariant is Bredon cohomology, introduced for finite groups by Bredon in \cite{Bredon-EquivariantCohomologyTheories} to study equivariant obstructions and extended to the study of infinite groups by L\"uck \cite{Lueck}.  

Fixing $G$ we consider the \emph{orbit category} $\OF$.  This is the small category whose objects are the transitive $G$-sets with stabilisers in $\mathcal{F}$ and the morphisms between two such $G$-sets is the free abelian group on the $G$-maps between them.  
\emph{Bredon modules}, or \emph{$\OF$-modules}, are contravariant additive functors from $\OF$ to the category of left $R$-modules, where $R$ is some commutative ring. 
Our definition of $\OF$ is different from the usual definition given in, for example, \cite{Lueck}, but the definitions lead to isomorphic $\OF$-module categories (Remark \ref{remark:equivalent def of orbit category}).

The category of all Bredon modules is an abelian category with frees and projectives, so one can use techniques from homological algebra to study them.  Let $\uR$ be the constant Bredon module, defined as $\uR(G/H) = R$ for all $H \in \mathcal{F}$ and $\uR(\alpha) = \id_R$ for any $G$-map $\alpha : G/H \to G/K$.\index{Constant $\OF$-module $\uR$}  As in ordinary group cohomology, using projective resolutions of $\uR$ one builds the Bredon cohomology of $G$.  Analagously to $\cd G$, we denote by $\OFcd G$ the \emph{Bredon cohomological dimension} of $G$---the minimal length of a projective resolution of $\uR$.  We denote by $\OFFP_n$ the Bredon cohomological analogue of the $\FP_n$ conditions of ordinary cohomology, so $G$ is $\OFFP_n$ if there exists a resolution of $\uR$ by projective Bredon modules which is finitely generated in all degrees $\le n$.  

That the Bredon cohomological dimension $\OFincd G$ is the correct algebraic invariant to mirror $\gdFin G$ is exemplified by the following theorem, an analogue of the classical results of Eilenberg--Ganea and Stallings--Swan.
\begin{Theorem*}\cite[Theorem 0.1]{LuckMeintrup-UniversalSpaceGrpActionsCompactIsotropy}\cite{Dunwoody-AccessabilityAndGroupsOfCohomologicalDimensionOne}
 If $R = \ZZ$ then $\OFincd G = \gdFin G$, except for the possibility that $\OFincd G = 2$ and $\gdFin G = 3$.
\end{Theorem*}
Brady, Leary and Nucinkis construct groups with $\OFincd G = 2$ and $\gdFin G = 3$ \cite{BradyLearyNucinkis-AlgAndGeoDimGroupsWithTorsion}.

If $G$ admits a model for $\EFin G$ with cocompact $n$-skeleton then $G$ is $\OFinFP_n$ over $\ZZ$.  In the other direction, if $G$ is $\OFinFP_n$ and the Weyl groups $WH = N_GH/H$ are finitely presented for all finite subgroups $H$ of $G$, then $G$ admits a model for $\EFin G$ with cocompact $n$-skeleton \cite[Theorem 0.1]{LuckMeintrup-UniversalSpaceGrpActionsCompactIsotropy}.
\begin{Prop*}\cite[Lemmas 3.1,3.2]{KMPN-CohomologicalFinitenessConditionsInBredon}
 A group $G$ is $\OFinFP_n$ if and only if $G$ has finitely many conjugacy classes of finite subgroups and the Weyl groups $WH = N_GH / H$ are $\FP_n$ for all finite subgroups $H$.
\end{Prop*}
We will discuss these conditions in more depth in Section \ref{section:bredon FPn}.

Much of this thesis is concerned with the Bredon cohomological dimension and $\OFinFP_n$ conditions, and how they interact with other cohomological finiteness conditions.  This includes those obtained by restricting the coefficients of Bredon cohomology to Mackey functors or cohomological Mackey functors (Section \ref{section:I Mackey} and Chapter \ref{chapter:M}) and those obtained via relative homological algebra, namely the Gorenstein cohomological dimension and $\mathfrak{F}$-cohomological dimension (Sections \ref{section:I Fcohomology} and \ref{section:I Gorenstein}, and Chapter \ref{chapter:G}).

There are already many results giving bounds for the Bredon cohomological dimension in terms of other algebraic invariants.  In \cite{MartinezPerez-EulerClassesAndBredonForRestrictedFamilies,MartinezPerez-SubgroupPosets}, Mart\'{\i}nez-P\'erez uses the poset of finite subgroups of a group to provide bounds for $\OFincd G$ and in \cite{KMN-CohomologicalFinitenessConditionsForElementaryAmenable} Kropholler, Mart\'{\i}nez-P\'erez, and Nucinkis show that any elementary amenable group of type $\FP_\infty$ over $\ZZ$ satisfies 
\[
 \OFincd_\ZZ G = hG = \cd_\QQ G,
\]
where $hG$ denotes the \emph{Hirsch length} of $G$.  See \cite{Hillman-ElementaryAmenableGroupsAnd4Manifolds} for a definition of Hirsch length for elementary amenable groups.

Finiteness conditions in Bredon cohomology are not well-behaved with respect to group extensions.  This is exemplified by the constructions of Leary and Nucinkis \cite{LearyNucinkis-SomeGroupsOfTypeVF} of 
\begin{enumerate}
 \item groups which are virtually-$\F$ (there exists a finite index subgroup of type $\F$) and satisify $\vcd G < \OFincd G$, and
 \item groups which are virtually-$\F$ with infinitely many conjugacy classes of finite subgroups (and hence not of type $\OFinFP_0$ by Proposition \ref{prop:uFP_0 iff finitely many conj classes of finite subgroups}).
\end{enumerate}
Interestingly, a virtually-$\F$ group cannot contain infinitely many conjugacy classes of subgroups of prime power order \cite[IX.13.2]{Brown}, but may contain infinitely many conjugacy classes of subgroups isomorphic to some finite group $H$ as long as $H$ does not have prime power order \cite{Leary-FiniteSubgroupsOfGroupsOfTypeVF}.

In Chapter \ref{chapter:bredon} we look in detail at $\OF$-modules and at some results concerning finiteness conditions in Bredon cohomology which will be needed later on in the thesis.  We also give some interesting examples of groups whose Bredon cohomological dimension is not preserved under change of rings.  Apart from these examples, this chapter contains mainly background material and straightforward extensions of known results.

\section{Modules over a category}

An {\bf Ab} category (also called a pre-additive category) is a category $\C$ enriched over the category of abelian groups---for any two objects $x$ and $y$ in $\C$ the morphisms from $x$ to $y$ form an abelian group and morphism composition distributes over addition \cite[p.28]{SaundersMaclane}.  So for any $w, x, y, z \in \C$ and morphisms
\[
 w \stackrel{f}{\rightarrow} x \stackrel[h]{g}{\rightrightarrows} y \stackrel{k}{\rightarrow} z,
\]
we have
\[
 k \circ (g + h) \circ f = k \circ g \circ f + k \circ h \circ f.
\]

If $\C$ is a small {\bf Ab} category then a \emph{$\C$-module} is a contravariant additive functor from $\C$ to the category of left $R$-modules.  The theory of modules over a category specialises to Bredon cohomology by setting $\C = \OF$.  In Chapter \ref{chapter:C} we study modules over an \textbf{Ab} category $\C$ with the property that for all objects $x$ and $y$ in $\C$, the set of morphisms from $x$ to $y$ forms a free abelian group.  We describe standard constructions including tensor products; projective, injective, and flat modules; restriction, induction, and coinduction; the $\Tor^{\C}_*$ and $\Ext_{\C}^*$ functors; projective dimension and $\CFP_n$ conditions; and the Bieri--Eckmann criterion.

Let $[x,y]_{\C}$ denote the morphisms from $x$ to $y$ in $\C$.  The category $\C$ is said to be an (EI) category if (see Remark \ref{remark:C not assumed EI}):
\begin{enumerate}
 \item[(EI)] For every $x \in \C$ there is a distinguished basis of $[x,x]_{\C}$, the elements of which are isomorphisms.
\end{enumerate}
The material in this chapter is well-known in the case that $\C$ is an (EI) category.  We study a more general case which can be specialised, not just to Bredon modules in Chapter \ref{chapter:bredon}, but also to Mackey and cohomological Mackey functors in Chapter \ref{chapter:M}---Mackey and cohomological Mackey functors may be described as modules over categories $\MF$ and $\HeckeF$ respectively, categories which do not have (EI).

\section{\texorpdfstring{\nG}{nG} and \texorpdfstring{$\mathfrak{F}$}{F}-cohomological dimension}\label{section:I Fcohomology}

Let $n_G$ denote the minimal dimension of a contractible proper $G$-CW complex.\index{nG@$n_G$, minimal dimension of a proper contractible $G$-CW complex}  Nucinkis suggested $\mathfrak{F}$-cohomology in \cite{Nucinkis-CohomologyRelativeGSet} as an algebraic analogue of $n_G$, it is a special case of the relative homology of Mac Lane \cite{MacLane-Homology} and Eilenberg--Moore \cite{EilenbergMoore-FoundationsOfRelativeHomologicalAlgebra}.  Fix a subfamily $\mathcal{F}$ of the family of finite subgroups, closed under conjugation and taking subgroups.  Let $\Delta$ be the $G$-set $\coprod_{H \in {\mathcal{F}}} G/H$ and say that a module is \emph{$\mathcal{F}$-projective} if it is a direct summand of a module of the form $N \otimes R \Delta$ where $N$ is any $R G$-module.\index{F-projective@$\mathcal{F}$-projective}  Short exact sequences are replaced with $\mathcal{F}$-split short exact sequences---short exact sequences which split when restricted to any subgroup in $\mathcal{F}$, or equivalently which split when tensored with $R \Delta$.\index{F-split@$\mathcal{F}$-split} The class of $\mathcal{F}$-split short exact sequences is allowable in the sense of Mac Lane, and the projective modules with respect to these sequences are exactly the $\mathcal{F}$-projectives.  This means an $RG$-module $P$ is $\mathcal{F}$-projective if and only if given any $\mathcal{F}$-split short exact sequence 
\[
 0 \longrightarrow A \longrightarrow B \longrightarrow C \longrightarrow 0
\]
of $RG$-modules, applying $\Hom_{RG}(P, -)$ gives a short exact sequence 
\[
 0 \longrightarrow \Hom_{RG}(P, A) \longrightarrow \Hom_{RG}(P, B) \longrightarrow \Hom_{RG}(P, C) \longrightarrow 0.
\]

There are enough $\mathcal{F}$-projectives and one can define Ext and Tor functors, denoted $\FExt^*_{RG}$ and $\FTor_*^{RG}$, for any $RG$-modules $M$ and $N$, \index{F-Ext@$\mathcal{F}$-Ext functor $\FExt^*_{RG}(M, -)$}\index{F-Tor@$\mathcal{F}$-Tor functor $\FTor_*^{RG}(M, -)$}
\[  
\FExt^*_{RG}(M, N ) = H^*\Hom_{RG}(P_*, N)  
\]
\[  
\FTor_*^{RG}(M, N ) = H_*(P_* \otimes_{RG} N)  
\]
where $P_*$ is a $\mathcal{F}$-split resolution of $M$ by $\mathcal{F}$-projective modules.  We define the $\mathcal{F}$-cohomology and $\mathcal{F}$-homology \index{F-cohomology@$\mathcal{F}$-cohomology $\FH^*(G, -)$}\index{F-homology@$\mathcal{F}$-homology $\FH_*(G, -)$}
\[ 
\FH^*(G, M) = \FExt^*_{RG}(R, M) ,
\]
\[ 
\FH_*(G, M) = \FTor_*^{RG}(R, M) .
\]

The \emph{$\mathcal{F}$-cohomological dimension}, denoted $\Fcd G$, is the shortest length of an $\mathcal{F}$-split $\mathcal{F}$-projective resolution of $R$. \index{F-cohomological dimension@$\mathcal{F}$-cohomological dimension $\Fcd$}  A group $G$ is $\FFP_n$ if there exists an $\mathcal{F}$-split $\mathcal{F}$-projective resolution of $R$, finitely generated in all degrees $\le n$. \index{FFPn@$\FFP_n$ condition} 
 
\begin{Notation}
 When $\mathcal{F} = \Fin$, we use the standard notation in the literature, writing $\Fincd$ instead of $\Fin\negthinspace\cd$ and referring to the $\mathfrak{F}$-cohomological dimension.  
\end{Notation}

A result of Bouc and Kropholler--Wall implies $\Fincd G \le n_G$ \cite{Bouc-LeComplexeDeChainesDunGComplexe,KrophollerWall-GroupActionsOnAlgebraicCellComplexes}, but it is unknown if $\Fincd G < \infty$ implies $n_G < \infty$, or if there are any groups for which the two invariants differ.  Unfortunately $\mathcal{F}$-cohomology can be very difficult to deal with, in particular it lacks some useful features such as free modules.

Since every model for $\EFin G$ is a proper contractible $G$-CW complex, it is clear that $n_G \le \gdFin G$.
\begin{Conj*}[Kropholler--Mislin Conjecture {\cite[Conjecture 43.1]{GuidoBookOfConjectures}}]
 If $n_G < \infty$ then $\gdFin G < \infty$.
\end{Conj*}

Kropholler and Mislin verified their conjecture for groups of type $\FP_\infty$ \cite{KrophollerMislin-GroupsOnFinDimSpacesWithFinStab} and later L\"uck verified the conjecture for groups with $l(G)< \infty$ \cite{Luck-TypeOfTheClassifyingSpace}.  Here $l(G)$ is the length of the longest chain 
\[
 1 = H_0 \lneq H_1 \lneq H_2 \lneq \cdots \lneq H_n \lneq G
\]
of finite subgroups in $G$. \index{lG@$l(G)$, length of a group}
Nucinkis posed an algebraic version of the conjecture, asking whether the finiteness of $\Fincd G$ and $\OFincd G$ are equivalent and verifying this for groups with $l(G) < \infty$ \cite[Conjecture on p.337, Corollary 4.5]{Nucinkis-EasyAlgebraicCharacterisationOfUniversalProperGSpaces}.

The class $\HF$ of hierarchically decomposable groups was introduced by Krop\-holler as the smallest class of groups such that if there exists a finite-dimensional contractible $G$-CW complex with stabilisers in $\HF$ then $G \in \HF$ \cite{Kropholler-OnGroupsOfTypeFP_infty}, he proves that every torsion-free group of type $\FP_\infty$ in $\HF$ has finite cohomological dimension.\index{Kropholler's class of hierarchically decomposable groups $\HF$}  The class $\HF$ is extremely large, containing for example all countable elementary amenable groups and all countable linear groups.  The first known example of a group not in Kropholler's class $\HF$ was Thompson's group $F$, since $F$ is torsion-free of type $\FP_\infty$ but with infinite cohomological dimension \cite{BrownGeoghegan-AnInfiniteDimensionTFFPinftyGroup}.  Other examples have since been found \cite{Gandini-CohomologicalInvariants,ABJLMS-InfiniteGroupsWithFixedPointProperties}.  Gandini and Nucinkis have verified the Kropholler--Mislin conjecture for a class of groups containing many groups of unbounded torsion \cite{GandiniNucinkis-SomeH1FGroupsAndAConjectureOfKrophollerAndMislin}.

In \cite[Example 3.6]{MartinezPerez-EulerClassesAndBredonForRestrictedFamilies} Mart\'{\i}nez-Per\'ez modifies the Leary--Nucinkis construction \cite{LearyNucinkis-SomeGroupsOfTypeVF} to produce an extension $G$ of a torsion-free group by a cyclic group of order $p$, with $\Fincd G = 3$ but $\OFincd G = 4$.  Taking direct products of these groups and using \cite[Theorem C]{DegrijsePetrosyan-GeometricDimensionForVCYC} gives a family of virtually torsion-free groups $G_n$ with $\OFincd G_n = \Fincd G_n + n$ for all natural numbers $n$ \cite[Remark 3.6]{Degrijse-ProperActionsAndMackeyFunctors}.  However one should note that in these examples $\Fincd G_n$ is growing linearly with $n$.  

Interestingly, it is still unknown if $\OFincd G = \Fincd G$ when $G$ is of type $\OFinFP_\infty$, although Degrijse and Mart\'{\i}nez-P\'erez have obtained some results pertaining to this question in \cite{DegrijseMartinezPerez-DimensionInvariants}.  They investigate groups admitting a cocompact model $X$ for $\EFin G$, and find a description of $\OFincd G$ as largest $n$ for which $H^n_c(X^K, X^K_{\text{sing}})$ is non-trivial, where $K$ runs over the finite subgroups of $G$, $X^K_{\text{sing}}$ denotes the subspace of cells in $X^K$ with isotropy containing $K$ but not equal to $K$, and $H^n_c$ denotes the cohomology with compact supports  \cite[Corollary 2.5]{DegrijseMartinezPerez-DimensionInvariants}.  Using this they prove that if $G$ acts properly and chamber-transitively on a building of type $(W,S)$, where $(W,S)$ is a finite Coxeter group, then $\OFincd G = \Fincd G$ \cite[Theorem 5.4]{DegrijseMartinezPerez-DimensionInvariants}.

\section{Mackey and cohomological Mackey functors}\label{section:I Mackey}
 
In \cite{MartinezPerezNucinkis-MackeyFunctorsForInfiniteGroups}, Mart\'{\i}nez-Per\'ez and Nucinkis studied cohomological finiteness conditions arising from taking the Bredon cohomology of a group $G$ but restricting to Mackey functor coefficents.  They showed that the associated Mackey cohomological dimension $\operatorname{\mathcal{M}_{\Fin}cd} G$ is always equal to both the virtual cohomological dimension and the $\mathfrak{F}$-cohomological dimension when $G$ is virtually torsion-free.  
One can view Mackey functors as contravariant functors from a small category $\mathcal{M}_{\Fin}$ into the category of left $R$-modules, and a crucial result in the paper of   Mart\'{\i}nez-Per\'ez and Nucinkis is that the Bredon cohomology with coefficients in a Mackey functor may be calculated using a projective resolution of Mackey functors.  Specifically they prove that you can induce a resolution of $\uR$ by projective Bredon modules to a resolution of the Burnside functor $B^G$ by projective Mackey functors.  This is explained in more detail in Section \ref{section:prelim mackey}.

Degrijse showed that for groups with $l(G) < \infty$ the Mackey cohomological dimension is equal to the $\mathfrak{F}$-cohomological dimension \cite[Theorem A]{Degrijse-ProperActionsAndMackeyFunctors}.  He proves this via the study of Bredon cohomology with cohomological Mackey functor coefficents and the associated notion of cohomological dimension $\operatorname{\mathcal{H}_{\Fin}cd} G$.  

The main ingredient of Chapter \ref{chapter:M} is a similar result to that of Mart\'{\i}nez-Per\'ez and Nucinkis for Bredon cohomology with cohomological Mackey functor coefficients.  Yoshida observed that a cohomological Mackey functor may be described as a contravariant functor from a small category $\mathcal{H}_{\Fin}$ to the category of left $R$-modules \cite{Yoshida-GFunctors2}.  We use Yoshida's result to prove in Section \ref{section:Homology and Cohomology of Cohomological Mackey Functors} that the Bredon cohomology with coefficients in a cohomological Mackey functor may be calculated with a projective resolution of cohomological Mackey functors, by showing that a resolution of $\uR$ by projective Bredon modules can be induced to a resolution of the fixed point functor $R^-$ by projective cohomological Mackey functors.  

Degrijse also proves in \cite{Degrijse-ProperActionsAndMackeyFunctors} that $\operatorname{\mathcal{H}_{\Fin}cd} G = \Fincd G$ when $\operatorname{\mathcal{H}_{\Fin}cd} G < \infty$, and asks if they are always equal, we can verify this:

\theoremstyle{plain}\newtheorem*{CustomThmA}{Theorem \ref{theorem:Fcd=HFcd}}
\begin{CustomThmA}
 $\operatorname{\mathcal{H}_{\Fin}cd} G = \Fincd G$ for all groups $G$.
\end{CustomThmA}

Thus for an arbitrary group $G$ we have a chain of inequalities
\[ \Fincd G = \operatorname{\mathcal{H}_{\Fin}cd} G \le \operatorname{\mathcal{M}_{\Fin}cd} G \le \operatorname{\mathcal{O}_{\Fin}cd} G. \]

Since the new invariants $\operatorname{\mathcal{M}_{\Fin}cd}$ and $\operatorname{\mathcal{H}_{\Fin}cd}$ interpolate between $\operatorname{\mathcal{O}_{\Fin}cd} G$ and $\Fincd G$, one might hope to use them to gain information about how the Kropholler--Mislin conjecture might fail.  However, few of the inequalities above are well understood.  The inequality $\Fincd G \le \OFincd G$ has already been discussed in Section \ref{section:I Fcohomology}.  For groups with $l(G)$ finite, $\HFincd G = \MFincd G$ \cite[Theorem 4.10]{Degrijse-ProperActionsAndMackeyFunctors}, we don't know of any examples where they differ.

\begin{Question}
\begin{enumerate}
 \item For an arbitrary group $G$, does the finiteness of $\operatorname{\mathcal{M}_{\Fin}cd} G$ imply the finiteness of $\operatorname{\mathcal{O}_{\Fin}cd} G$?
 \item Is there any relation between $\operatorname{\mathcal{M}_{\Fin}cd} G$ and $n_G$?
\end{enumerate}
\end{Question}

The $\OFinFP_n$ conditions are well understood---see Section \ref{section:bredon FPn}.  We study the $\mathcal{M}_{\Fin}\negthinspace\FP_n$ conditions corresponding to Mackey functors, the $\mathcal{H}_{\Fin}\negthinspace\FP_n$ conditions corresponding to cohomological Mackey functors,  and the $\FinFP_n$ conditions defined in the previous section. 

\theoremstyle{plain}\newtheorem*{CustomThmB}{Corollary \ref{cor:OFFPn iff MFFPn}}
\begin{CustomThmB}
Over any ring $R$, a group is $\mathcal{M}_{\Fin}\negthinspace\FP_n$ if and only if it is $\mathcal{O}_{\Fin}\negthinspace\FP_n$.
\end{CustomThmB}

\theoremstyle{plain}\newtheorem*{CustomThmC}{Theorem \ref{theorem:HFFPn iff FFPn}}
\begin{CustomThmC}
If $R$ is a commutative Noetherian ring, a group is $\mathcal{H}_{\Fin}\negthinspace\FP_n$ if and only if it is $\FinFP_n$. 
\end{CustomThmC}

In Section \ref{section:family of p subgroups} we prove a result similar to that shown for $\mathfrak{F}$-cohomology in \cite[\S 4]{LearyNucinkis-GroupsActingPrimePowerOrder}, that depending on the coefficient ring, $\operatorname{\mathcal{H}_{\Fin}cd}$ may be calculated using a subfamily of the family of finite subgroups.  For example when working over $\ZZ$ we need consider only the family of finite subgroups of prime power order, and over either the finite field $\FF_p$ or over $\ZZ_{(p)}$ (the integers localised at $p$), we need consider only the family of subgroups of order a power of $p$. 

\theoremstyle{plain}\newtheorem*{CustomThmD}{Theorem \ref{theorem:HF HeckeFcd = HeckePcd and HeckeFFPn = HeckePFPn}}
\begin{CustomThmD}
  Let $R$ be either $\ZZ$, $\FF_p$, or $\ZZ_{(p)}$.  If $R = \ZZ$ then denote by $\mathcal{P}$ the family of subgroups of prime-power order.  If $R = \FF_p$ or $\ZZ_{(p)}$ then let $\mathcal{P}$ denote the family of subgroups of order a power of $p$.

  For all $n \in \NN \cup \{ \infty \}$, the conditions $\operatorname{\mathcal{H}_{\Fin}cd} G = n$ and $\operatorname{\mathcal{H}_{\mathcal{P}}cd} G = n$ are equivalent, as are the conditions $\mathcal{H}_{\Fin}\negthinspace\FP_n$ and $\mathcal{H}_{\mathcal{P}}\negthinspace\FP_n$. 
\end{CustomThmD}

We also give a complete description of the condition $\mathcal{H}_{\Fin}\negthinspace\FP_n$ over $\FF_p$.

\theoremstyle{plain}\newtheorem*{CustomThmE}{Corollary \ref{cor:mackey HFFPn complete description}}
\begin{CustomThmE}
$G$ is $\mathcal{H}_{\Fin}\negthinspace\FP_n$ over $\FF_p$ if and only if $G$ has finitely many conjugacy classes of $p$-subgroups, and $WH = N_GH/H$ is $\FP_n$ over $\FF_p$ for all finite $p$-subgroups $H$.
\end{CustomThmE}
 
 \section{Gorenstein cohomological dimension}\label{section:I Gorenstein}
 
An $RG$-module is \emph{Gorenstein projective} if it is a cokernel in a strong complete resolution of $RG$-modules, these were first defined over an arbitrary ring by Enochs and Jenda \cite{EnochsJenda-GorensteinInjectiveAndProjectiveModules}.\index{Gorenstein projective}  We give a full explanation of complete resolutions in Section \ref{subsection:complete resolutions}.  The \emph{Gorenstein projective dimension} $\Gpd M$ of an $RG$-module $M$ is the minimal length of a resolution of $M$ by Gorenstein projective modules. \index{Gorenstein projective dimension $\Gpd$} Equivalently, $\Gpd M \le n$ if and only if $M$ admits a complete resolution of coincidence index $n$ \cite[p.864]{BahlekehDembegiotiTalelli-GorensteinDimensionAndProperActions}.

The \emph{Gorenstein cohomological dimension} of a group $G$, denoted $\Gcd G$, is the Gorenstein projective dimension of $R$.\index{Gorenstein cohomological dimension $\Gcd$}  If $G$ is virtually torsion-free then $\Gcd G = \vcd G$ \cite[Remark 2.9(1)]{BahlekehDembegiotiTalelli-GorensteinDimensionAndProperActions}, indeed the Gorenstein cohomology can be seen as a generalisation of the virtual cohomological dimension.  

$\Gcd G$ is closely related to the $\silp RG$ and $\spli RG$ invariants studied by Gedrich and Gruenberg \cite{GedrichGruenberg-CompleteCohomologicalFunctors} and recently shown to be equal by Emmanouil when $R = \ZZ$ \cite{Emmanouil-OnCertainCohomologicalInvariantsOfGroups}.  The invariants $\silp RG$ and $\spli RG$ are defined as the supremum of the injective lengths (injective dimensions) of the projective $R G$-modules and the supremum of the projective lengths (projective dimensions) of the injective $RG$-modules respectively.  It is known that 
\[
\Gcd G \le \spli R G \le \Gcd G + 1,
\]
and conjectured that $\Gcd G = \spli R G$ \cite[Conjecture A]{DembegiotiTalelli-RelationBetweenCertainCohomologicalInvariants}.  In fact, Dembegioti and Talelli phrase this conjecture with the generalised cohomological dimension of Ikenaga \cite{Ikenaga-HomologicalDimensionAndFarrellCohomology}, but this is always equal to the Gorenstein cohomological dimension \cite[Theorem 2.5]{BahlekehDembegiotiTalelli-GorensteinDimensionAndProperActions}.
 
By \cite[Lemma 2.21]{AsadollahiBahlekehSalarian-HierachyCohomologicalDimensions}, every permutation $R G$-module with finite stabilisers is Gorenstein projective, so combining with \cite[Lemma 3.4]{Gandini-CohomologicalInvariants} gives that $\Gcd G \le \Fincd G$.  

\theoremstyle{plain}\newtheorem*{CustomThmGA*}{Theorem \ref{thm:Fcd finite then Fcd=Gcd}}
\begin{CustomThmGA*}
If $\Fincd G < \infty$ then $\Fincd G = \Gcd G$.
\end{CustomThmGA*}

We don't know if $\Gcd G < \infty$ implies $\Fincd G < \infty$, although if $\Gcd G = 0$ or $1$ then $\Gcd G = \Fincd G = \OFincd G$ \cite[Proposition 2.19]{AsadollahiBahlekehSalarian-HierachyCohomologicalDimensions}\cite[Theorem 3.6]{BahlekehDembegiotiTalelli-GorensteinDimensionAndProperActions}.  Additionally if $G$ is in Kropholler's class $\HF$ and has a bound on the orders of its finite subgroups then $\Fincd G = \Gcd G$ (see Example \ref{example:HF}).

It is conjectured by Talelli that $\Gcd G < \infty$ if and only if $\OFincd G < \infty$ (see for example, \cite[Conjecture A]{Talelli-OnGroupsOfTypePhi}).  This is a stronger version of the conjecture of Nucinkis mentioned in Section \ref{section:I Fcohomology}, that $\Fincd G < \infty$ if and only if $\OFincd G < \infty$.  If one could strengthen Theorem \ref{thm:Fcd finite then Fcd=Gcd} to show that $\Gcd G = \Fincd G$ for all groups $G$, then the two conjectures would be equivalent.  Using \cite{MartinezPerez-ABoundForTheBredonCohomologicalDimension}, Bahlekeh, Dembegioti, and Talelli show that for groups with $\OFincd G < \infty$, there is a bound $\OFincd G < l(G) + \Gcd G $ \cite[Theorem C]{BahlekehDembegiotiTalelli-GorensteinDimensionAndProperActions}.  

Generalising a construction of Avramov--Martsinkovsky, it was shown by Asadollahi, Bahlekeh, and Salarian that if $\Gcd G < \infty$ then there is a long exact sequence of cohomological functors relating group cohomology, complete cohomology and Gorenstein cohomology \cite[\S 7]{AvramovMartsinkovsky-AbsoluteRelativeAndTateCohomology-FiniteGorenstein}\cite[\S 3]{AsadollahiBahlekehSalarian-HierachyCohomologicalDimensions}.  Theorem \ref{thm:Fcd finite then Fcd=Gcd} follows from constructing a similar long exact sequence relating $\mathfrak{F}$-cohomology, complete $\mathfrak{F}$-cohomology (defined in Section \ref{subsection:complete F cohomology}), and a new cohomology theory we call $\mathfrak{F}_G$-cohomology (defined in Section \ref{section:FGcohomology}).  

When they both exist, these two long exact sequences fit into the commutative diagram below, see Proposition \ref{prop:AM seqs diagram},
\[
\xymatrix{
 \cdots \ar[r] & \FhatH^{n-1} \ar^{\gamma_{n-1}}[d] \ar[r] & \FGH^n \ar[r] \ar^{\alpha_n}[d] & \FinH^n \ar[r] \ar^{\beta_n}[d] & \FhatH^n \ar^{\gamma_n}[d] \ar[r] & \FGH^{n+1} \ar[r] \ar^{\alpha_{n+1}}[d] & \cdots \\
 \cdots \ar[r] &  \widehat{H}^{n-1} \ar[r] & \GH^n \ar[r] \ar_{\eta_n}[ru] & H^n \ar[r] & \widehat{H}^n \ar[r]   & \GH^{n+1} \ar[r] \ar[ru] & \cdots \\
}
\]
where for conciseness we have written $H^n$ for $H^n(G, -)$ etc.  In the commutative diagram above, $\widehat{H}^n(G, -)$ is the complete cohomology, $\GH^n(G, -)$ is the Gorenstein cohomology, $\FhatH^{n}(G, -)$ is complete $\mathfrak{F}$-cohomology, and $\FGH^n(G, -)$ is the $\mathfrak{F}_G$-cohomology.
 
Since Theorem \ref{thm:Fcd finite then Fcd=Gcd} is proved via this commutative diagram, it appears that the requirement that $\Fincd G < \infty$ will be difficult to circumvent---without it we do not know how to construct the long exact sequence appearing on the top row.

In Section \ref{section:Group Extensions} we use that the Gorenstein cohomological dimension is subadditive to improve upon a result of Degrijse on the behaviour of the $\mathfrak{F}$-cohomological dimension under group extensions \cite[Theorem B]{Degrijse-ProperActionsAndMackeyFunctors}.  Degrijse phrased his result in terms of Bredon cohomological dimension of $G$ with coefficients restricted to cohomological Mackey functors, but this invariant is equal to $\Fincd G$ by Theorem \ref{theorem:Fcd=HFcd} (see previous section).

\theoremstyle{plain}\newtheorem*{CustomThmGC*}{Corollary \ref{cor:Fcd subadditive if finite}}
\begin{CustomThmGC*}
Given a short exact sequence of groups
\[
1 \longrightarrow N \longrightarrow G \longrightarrow Q \longrightarrow 1,
\]
if $\Fincd G < \infty$ then $ \Fincd G \le \Fincd N + \Fincd Q $.
\end{CustomThmGC*}

\begin{Question}
 Is the $\mathfrak{F}$-cohomological dimension subadditive under group extensions?
\end{Question}

In Section \ref{section:cdQ} we use the Avramov--Martsinkovsky long exact sequence to prove the following.
\theoremstyle{plain}\newtheorem*{CustomThmGD*}{Proposition \ref{prop:cdQG finite then cdQG le GcdG}}
\begin{CustomThmGD*}
 If $\cd_{\QQ} G < \infty$ then $\cd_{\QQ} G \le \Gcd_{\ZZ} G$.
\end{CustomThmGD*}

We know of no groups for which $\cd_{\QQ} G \le \Gcd_{\ZZ} G$ fails.  If $\cd_{\ZZ} G < \infty$ then necessarily $\cd_{\ZZ}G = \Gcd G$ \cite[Corollary 2.25]{AsadollahiBahlekehSalarian-HierachyCohomologicalDimensions}, but we cannot rule out the possibility that there exists a torsion-free group $G$ with $\cd_{\ZZ} G = \infty$ but $\Gcd G < \infty$.  In fact, the question below is still open even for torsion-free groups.

\begin{Question}
 Do there exist groups $G$ with $\cd_{\QQ} G = \infty$ but $\Gcd_{\ZZ} G < \infty$?
\end{Question}

\section{Bredon duality groups}
 
A \emph{duality group} is a group $G$ of type $\FP$ for which 
\[
H^i(G, \ZZ G) \cong \left\{ \begin{array}{l l} \text{$\ZZ$-flat} & \text{if $i = n$,} \\ 0 & \text{else,} \end{array}\right.
\]
where $n$ is necessarily the cohomological dimension of $G$.\index{Duality group}  The name duality comes from the fact that this condition is equivalent to existence of a $\ZZ G$-module $D$, giving an isomorphism
\[
H^i(G, M) \cong H_{n-i}(G, D \otimes_\ZZ M)
\]
for all $i$ and all $\ZZ G$-modules $M$.  It can be proven that given such an isomorphism, the module $D$ is necessarily $H^n(G, \ZZ G)$.  A duality group $G$ is called a \emph{Poincar\'e duality group} if in addition \index{Poincar\'e duality group}
\[H^i(G, \ZZ G) \cong \left\{ \begin{array}{l l} \ZZ & \text{if $i = n$,} \\ 0 & \text{else.} \end{array}\right.\]
These groups were first defined by Bieri \cite{Bieri-GruppenMitPoincareDualitat}, and independently by Johnson--Wall \cite{JohnsonWall-OnGroupsSatisfyingPoincareDuality}.  Duality groups were first studied by Bieri and Eckmann in \cite{BieriEckmann-HomologicalDualityGeneralizingPoincareDuality}.  See \cite{Davis-PoincareDualityGroups} and \cite[\S III]{Bieri-HomDimOfDiscreteGroups} for a thorough introduction.

If a group $G$ has an oriented manifold model for $\B G$ then $G$ is a Poincar\'e duality group \cite[p.1]{Davis-PoincareDualityGroups}.  Wall asked if the converse is true \cite{Wall-HomologicalGroupTheory}, the answer is no as Poincar\'e duality groups can be built which are not finitely presented \cite[Theorem C]{Davis-CohomologyCoxeterGroupRingCoeff}.  However the question remains a significant open problem if we include the requirement that $G$ be finitely presented.  The conjecture is known to hold only in dimensions at most $2$ \cite{Eckmann-PoincareDualityGroupsOfDimensionTwoAreSurfaceGroups}.

Let $R$ be a commutative ring.  A group $G$ is a \emph{duality group over $R$} if $G$ is $\FP$ over $R$ and
\[
H^i(G, R G) \cong \left\{ \begin{array}{l l} \text{$R$-flat} & \text{if $i = n$,} \\ 0 & \text{else.} \end{array}\right.
\]
$G$ is \emph{Poincar\'e duality over $R$} if
\[
H^i(G, R G) \cong \left\{ \begin{array}{l l} R & \text{if $i = n$,} \\ 0 & \text{else.} \end{array}\right.
\]
An analogue of Wall's conjecture is whether every torsion-free finitely presented Poincar\'e duality group over $R$ is the fundamental group of an aspherical closed $R$-homology manifold \cite[Question 3.5]{Davis-PoincareDualityGroups}.  This is answered in the negative by Fowler for $R = \QQ$ \cite{Fowler-FinitenessForRationalPoincareDuality}, but remains open for $R = \ZZ$.

We study a generalisation of Poincar\'e duality groups, looking at the algebraic analogue of the condition that $G$ admits a manifold model $M$ for $\EFin G$ such that for any finite subgroup $H$ the fixed point set $M^H$ is a submanifold.  

If $G$ admits a cocompact manifold model $M$ for $\EFin G$ then $G$ is $\OFinFP$.  Also if for any finite subgroup $H$ the fixed point set $M^H$ is a submanifold, we have the following condition on the cohomology of the Weyl groups $WH = N_GH/H$:
\[ 
H^i(WH, \ZZ [WH]) =  \left\{ \begin{array}{c c} \ZZ & \text{ if } i = \dim M^H,  \\ 0 & \text{ else,} \end{array} \right. 
\]
see \cite[p.3]{DavisLeary-DiscreteGroupActionsOnAsphericalManifolds} for a proof of the above.  
Building on this, in \cite{DavisLeary-DiscreteGroupActionsOnAsphericalManifolds} and also in \cite[Definition 5.1]{MartinezPerez-EulerClassesAndBredonForRestrictedFamilies} a \emph{Bredon duality group over $R$} is defined as a group $G$ of type $\OFinFP$ such that for every finite subgroup $H$ of $G$ there is an integer $d_H$ with \index{Bredon duality group}\index{dH@$d_H$, integer such that $H^{d_H}(WH, R[WH]) \neq 0$}
\[ 
H^i(WH, R [WH]) = \left\{ \begin{array}{l l} \text{$R$-flat} & \text{if $i = d_H$,} \\ 0 & \text{else.} \end{array} \right. 
\]
Furthermore, $G$ is said to be \emph{Bredon--Poincar\'e duality over $R$} if for all finite subgroups $H$, \index{Bredon--Poincar\'e duality group}
\[
H^{d_H}(WH, R [WH]) = R.
\]
We say that a Bredon duality group $G$ is \emph{dimension $n$} if $\OFincd G = n$.  Note that for torsion-free groups these definitions reduce to the usual definitions of duality and Poincar\'e duality groups.  

One might generalise Wall's conjecture: Let $G$ be Bredon--Poincar\'e duality over $\ZZ$, such that $WH$ is finitely presented for all finite subgroups $H$.  Does $G$ admit a cocompact manifold model $M$ for $\EFin G$?  This is false by an example of Jonathan Block and Schmuel Weinberger, suggested to us by Jim Davis.

\theoremstyle{plain}\newtheorem*{CustomThmBPDNoMfd}{Theorem \ref{theorem:BPD groups with no manifold model}}
\begin{CustomThmBPDNoMfd}
  There exist examples of Bredon--Poincar\'e duality groups over $\ZZ$, such that $WH$ is finitely presented for all finite subgroups $H$ but $G$ doesn't admit a cocompact manifold model $M$ for $\EFin G$.
\end{CustomThmBPDNoMfd}

If $G$ is Bredon--Poincar\'e duality and virtually torsion-free then $G$ is virtually Poincar\'e duality.  Thus an obvious question is whether all virtually Poincar\'e duality groups are Bredon--Poincar\'e duality, in \cite{DavisLeary-DiscreteGroupActionsOnAsphericalManifolds} it is shown that this is not the case for $R = \ZZ$.  An example is also given in \cite[\S 6]{MartinezPerez-EulerClassesAndBredonForRestrictedFamilies} which fails for both $R = \ZZ$ and for $R = \FF_p$, the finite field of $p$ elements.  One might also ask if every Bredon--Poincar\'e duality group is virtually torsion-free but this is also not the case, see for instance Examples \ref{example:FarbWeinberger} and \ref{example:disc sub of lie not vtf}. 

In \cite[Theorems D,E]{Hamilton-WhenIsGroupCohomologyFinitary} Hamilton shows that, over a field $\FF$ of characteristic $p$, given an extension $\Gamma$ of a torsion-free group $G$ of type $\FP_\infty$ by a finite $p$-group, the group $\Gamma$ will be of type $\OFinFP_\infty$ (by examples of Leary and Nucinkis, an extension by an arbitrary finite group may not even be $\OFinFP_0$ \cite{LearyNucinkis-SomeGroupsOfTypeVF}).  Mart\'{\i}nez-P\'erez builds on this result to show that if $G$ is assumed Poincar\'e duality then $\Gamma$ is Bredon--Poincar\'e duality over $\FF$ with $\OFincd_{\FF} \Gamma  =\cd_{\FF} G$ \cite[Theorem C]{MartinezPerez-EulerClassesAndBredonForRestrictedFamilies}.  However, her results do not extend to Bredon duality groups.

Given a Bredon duality group $G$ we write $\mathcal{V}(G)$ for the set \index{VG@$\mathcal{V}(G)$}
\[
\mathcal{V}(G) = \{ d_F : F \text{ a non-trivial finite subgroup of }G\} \subseteq \{0, \ldots, n \}.
\]
In Example \ref{example:duality arbitrary V} we will build Bredon duality groups with arbitrary $\mathcal{V}(G)$.  If $G$ has a manifold model, or homology manifold model, for $\EFin G$ then there are some restrictions on $\mathcal{V}(G)$---see Section \ref{subsection:actions on homology manifolds} for this.  In Section \ref{section:reflection groups} we build Bredon--Poincar\'e duality groups for many choices of $\mathcal{V}(G)$, however the following question remains open:

\begin{Question}\label{question:duality prescribed VG}
 Is it possible to construct Bredon--Poincar\'e duality groups with prescribed $\mathcal{V}(G)$?
\end{Question}

It follows from Proposition \ref{prop:FP and Ofincd finite then Gcd=Fcd=sup} that for a Bredon--Poincar\'e duality group, $d_{1} \le \OFincd G$ (recall $d_1$ is the integer for which $H^{d_1}(G, RG) \cong R$) and also, if we are working over $\ZZ$, then $d_{1} = \cd_\QQ G$ (Lemma \ref{Lemma:observations}).  Thus the following question is of interest.

\begin{Question}\label{question:duality n=nid}
 Do there exist Bredon duality groups with $\OFincd G \neq d_{1}$?
\end{Question}

Examples of groups for which $\cd_\QQ G \neq \OFincd_\ZZ G$ are known \cite{LearyNucinkis-SomeGroupsOfTypeVF}, but there are no known examples of type $\OFinFP_\infty$.  This question is also related to \cite[Question 5.8]{MartinezPerez-EulerClassesAndBredonForRestrictedFamilies} where it is asked whether a virtually torsion-free Bredon duality group satisfies $\OFincd G = \vcd G$.

One might hope to give a definition of Bredon--Poincar\'e duality groups in terms of Bredon cohomology only, we do not know if this is possible but we show in Section \ref{section:wrong notion of duality} that the na\"ive idea of asking that a group be $\OFinFP$ with 
\[
H^i_{\OFin}(G, R[?,-]) \cong \left\{ \begin{array}{l l} \uR & \text{if $i = n$,} \\ 0 & \text{else,} \end{array} \right.
\]
is not the correct definition, where in the above $H^i_{\OFin}$ denotes the Bredon cohomology and $\uR$ is the constant covariant Bredon module.  Namely we show in Theorem \ref{theorem:bredon wrong duality} that any such group is necessarily a torsion-free Poincar\'e duality group over $R$.
 
\section{Houghton's groups}
 
Houghton's group $H_n$ was introduced in \cite{Houghton-FirstCohomologyOfAGroupWithPermModCoeff}, as an example of a group acting on a set $S$ with $H^1(H_n, A \otimes \ZZ[S]) = A^{n-1}$ for any abelian group $A$.   

In \cite{Brown-FinitenessPropertiesOfGroups}, Brown used an important new technique to show that the groups $F_{n,r}, T_{n,r}$, and $V_{n,r}$ of Thompson and Higman were $\FP_\infty$.  In the same paper he showed that Houghton's group $H_n$ is interesting from the viewpoint of cohomological finiteness conditions, namely $H_n$ is $\FP_{n-1}$ but not $\FP_n$.  Brown proves this by studying the action of $H_n$ on the geometric realisation $\vmv$ of a certain poset $\mathcal{M}$.  More recently, Johnson gave a finite presentation for $H_3$ \cite{Johnson-EmbeddingSomeRecursivelyPresentedGroups}, and later Lee did the same for $H_n$ where $n \ge 3$ \cite{Lee-GeometryOfHoughtonsGroups}.  

Interestingly, $H_n$ embeds in Thompson's group $V = V_{2,1}$ for all $n \ge 0$ \cite{Rover-Thesis}.  Antol\'{\i}n, Burillo, and Martino have shown that for $n \ge 2$, the group $H_n$ has solvable conjugacy problem \cite{AntolinBurilloMartino-ConjugacyInHoughtonsGroups} and Burillo, Cleary, Martino, and R\"over have calculated the automorphism groups and abstract commensurators of $H_n$ \cite{BurilloClearyMartinoRover-CommensuratorsOfHn}. 

There has been recent interest in the structure of the centralisers of Thompson's groups and their generalisations \cite{MartinezNucinkis-GeneralizedThompsonGroups,BleakEtAl-CentralizersInVn,MartinezPerezMatucciNucinkis-CentralisersInGeneralisationsOfV}.  The results obtained here are similar to \cite[4.10,4.11]{MartinezPerezMatucciNucinkis-CentralisersInGeneralisationsOfV} where it is shown that in the groups $V_r(\Sigma)$, generalisations of Thompson's $V$, the centralisers of finite subgroups are of type $\FP_\infty$ whenever the groups $V_r(\Sigma)$ are of type $\FP_\infty$.

In Section \ref{section:centralisersOfFiniteSubgroups} we completely describe centralisers of finite subgroups and prove the following.
\theoremstyle{plain}\newtheorem*{CustomThmF}{Corollary \ref{cor:stabOfH is FPn-1 not FPn}}
\begin{CustomThmF}
If $Q$ is a finite subgroup of $H_n$ then the centraliser $C_{H_n}Q$ is $\FP_{n-1}$ but not $\FP_n$.
\end{CustomThmF}

This contrasts with \cite{KochloukovaDessislavaMartinezPerez-CentralisersOfFiniteSubgroupsInSolubleGroups} where examples are given of soluble groups of type $\FP_n$ with centralisers of finite subgroups that are not $\FP_n$, although it is shown in \cite{MartinezPerezNucinkis-VirtuallySolubleGroupsOfTypeFPinfty} that in virtually soluble groups of type $\FP_\infty$ the centralisers of all finite subgroups are of type $\FP_\infty$.

In Section \ref{section:centralisersOfElements} our analysis is extended to arbitrary elements and virtually cyclic subgroups.  Using this information elements in $H_n$ are constructed whose centralisers are $\FP_i$  for any $0 \le i \le n-3$.  

In Section \ref{section:Browns model} the space $\vmv$ mentioned previously is shown to be a model for $\EFin H_n$.

Finally Section \ref{section:finitenessConditionsSatisfied} contains a discussion of Bredon (co)homological finiteness conditions that are satisfied by Houghton's group.  
In particular we calculate the Bredon cohomological dimension with respect to the family of finite subgroups, and use a construction of L\"uck and Weiermann \cite{LueckWeiermann-ClassifyingspaceForVCYC} to calculate the Bredon cohomological dimension with respect to the family of virtually cyclic subgroups.

\theoremstyle{plain}\newtheorem*{CustomThmG}{Proposition \ref{prop:cd=gd=n} and Theorem \ref{theorem:houghton cdVCyc Hn = n}}
\begin{CustomThmG}
 $\OFincd H_n = \OVCyccd H_n = n$.
\end{CustomThmG}
\chapter{Modules over a category}\label{chapter:C}

Much of this chapter is based on \cite{Lueck}.  Although we consider a slightly more general situation, as explained in Remark \ref{remark:C not assumed EI}, the idea is the same.  The material in this chapter is used in much of this thesis, especially in Chapters \ref{chapter:bredon} and \ref{chapter:M}.

Let $R$ be a commutative ring with unit and ${\C}$ a small {\bf Ab} category (sometimes called a preadditive category) with the condition below.
\begin{enumerate}\index{Property (A)}
 \item[$(A)$] For any two objects $x$ and $y$ in ${\C}$, the set of morphisms, denoted $[x,y]_{{\C}}$, between $x$ and $y$ is a free abelian group.
\end{enumerate}
Recall that an {\bf Ab} category is one where the morphisms between any two objects form an abelian group and where morphism composition distributes over this addition \cite[A.4]{Weibel}.

\begin{Remark}\label{remark:C not assumed EI}\index{Property (EI)}
  In \cite[9.2]{Lueck}, categories $\mathfrak{X}$ are considered with the property that every endomorphism in $\mathfrak{X}$ is an isomorphism.  
  However the approach to defining modules over a category in \cite[9.2]{Lueck} is different from that used here (see also Remark \ref{remark:equivalent def of orbit category}).  One can translate between the different viewpoints in the following way:
  \[
  [x,y]_{\mathfrak{X}} = \ZZ[\{\text{Morphisms $x \to y$ in the sense of L\"uck}\}],
  \]
  where $\ZZ[X]$ denotes the free abelian group on a set $X$.
  
 The correct analogue of L\"uck's property with our definitions is the following:
 \begin{enumerate}
  \item[(EI)] For every $x \in \C$, there is a distinguished basis of $[x,x]_{\C}$, the elements of which are isomorphisms.
 \end{enumerate}
 
 The main advantage of the (EI) property is that it allows objects in $\C$ to be given a partial order: setting $x \le y$ if $[x,y]_{\C}$ is non-empty.  
 We choose not to ask for this property in this section, since we want everything discussed here to be relevant to the Mackey and Hecke categories, discussed in Chapter \ref{chapter:M}, which do not have (EI).  The motivating example of a category with (EI) is the orbit category, see Example \ref{example:orbit category}.
\end{Remark}

Throughout, the fraktur letters ${\C}$, $\mathfrak{D}$, $\mathfrak{E}$ etc.\ will always denote small {\bf Ab} categories with (A).

Define the category of covariant ${\C}$-modules over $R$ to be the category of additive covariant functors from ${\C}$ to $\RMod$, the category of left $R$-modules.\index{C-module@$\C$-module}  Similarly the category of contravariant ${\C}$-modules over $R$ is the category of additive contravariant functors from ${\C}$ to $\RMod$.  

If neither ``covariant'' or ``contravariant'' is specified in a statement about $\C$-modules, the reader should assume the statement holds for both covariant and contravariant $\C$-modules.

Since $\C$-modules form a functor category and $\RMod$ is an abelian category, the category of $\C$-modules is an abelian category \cite[44]{Murfet-AbelianCategories}.  In fact, it inherits all of Grothendieck's axioms for an abelian category which are satisfied by $\RMod$ \cite[44,55]{Murfet-AbelianCategories}, namely:
\begin{enumerate}
\item AB3 and AB4---Every small colimit exists and products of exact sequences are exact.
\item AB3* and AB4*---Every small limit exists and coproducts of exact sequences are exact.
\item AB5---Filtered colimits of exact sequences are exact.
\end{enumerate}
Again because we are working in a functor category, a sequence of $\C$-modules 
\[
0 \longrightarrow A \longrightarrow B \longrightarrow C \longrightarrow 0 
\]
is exact if and only if it is exact when evaluated at every $x \in {\C}$.  Note that $0$ denotes the zero functor, sending every object to the zero module. Similarly, using the fact that the category of $\C$-modules is a functor category and the category of abelian groups is complete, limits and colimits are computed pointwise \cite[p.8]{Murfet-AbelianCategories}. 

Since $[x,y]_{\C}$ is abelian for all $x$ and $y$ in ${\C}$, for any $y \in {\C}$ we can form a contravariant module $R[-, y]_{\C}$ by \index{Free $\C$-module $R[-,x]_{\C}$}
\[
R[-,y]_{\C} (x) = R \otimes_\ZZ [x,y]_{\C}.
\]
The analogous construction for covariant modules gives us
\[
R[y,-]_{\C} (x) = R \otimes_\ZZ [y,x]_{\C}.
\]
In Section \ref{section:C frees projectives} we will show that these modules are analogues of free modules in the category of $\C$-modules.  Since $R[x,y]_{\C}$ is a free $R$-module we write $r\alpha$ instead of $r \otimes \alpha$, for $r \in R$ and $\alpha \in [x,y]_{\C}$.

If $f \in R[x, y]_{\C}$, where $f = \sum_i r_i f_i$ for some $f_i \in [x,y]_{\C}$, and $Q$ is a $\C$-module, then we will write $Q(f)$ for the $R$-module homomorphism given by $\sum_i r_i Q(f_i)$.

Let $A$ and $B$ be any two covariant $\C$-modules, or any two contravariant $\C$-modules, then we denote by $\Hom_{\C}(A, B)$ the $\C$-module morphisms between $A$ and $B$, i.e.~the natural transformations from $A$ to $B$.\index{HomC@$\Hom_{\C}$}

\begin{Lemma}[The Yoneda-type lemma]\label{lemma:C yoneda-type}\index{Yoneda-type lemma}
For any covariant functor $A$ and $x \in {\C}$, there is an isomorphism, natural in $A$:
 \begin{align*}
 \Hom_{{\C}} \left( R[x,-]_{\C} , A\right) &\cong A(x)  \\
 f &\mapsto f(x)(\id_x)
 \end{align*}

Similarly for any contravariant functor $M$ and $x \in {\C}$, there is an isomorphism, natural in $M$:
 \begin{align*}
\Hom_{{\C}} \left( R[-,x]_{\C} , M\right) &\cong M(x) \\
 f &\mapsto f(x)(\id_x)
 \end{align*}
\end{Lemma}
The proof is a generalisation of \cite[p.9]{MislinValette-BaumConnes} into the setting of $\C$-modules. 
\begin{proof}
 We provide a proof only for covariant modules, that for contravariant modules is similar.  

Let $f$ be a morphism $f:R[x,-]_{\C} \to A$, we claim $f$ is completely determined by $f(x)$.  If $\alpha \in R[x,y]_{\C}$ then we can view 
 $\alpha $ as $\alpha = R[x,\alpha]_{\C} ( \id_x )$, thus
 \begin{align*}
 f(y)(\alpha) &= f(y) \circ R[x,\alpha]_{\C} ( \id_x ) \\ 
 &= A(\alpha) \circ f(x) ( \id_x) 
 \end{align*}
 where we use that $f$ is $R$-additive and that $f$ is a morphism in the category of $\C$-modules---so a natural transformation of functors---meaning the diagram below commutes.
 \[\xymatrix{
 R[x,x]_{\C} \ar_{R[x,\alpha]_{\C}}[d] \ar^-{f(x)}[r] & A(x) \ar^{A(\alpha)}[d]\\
 R[x,y]_{\C} \ar^-{f(y)}[r] & A(y)
 }\]
 
 Conversely, given an element $a \in A(x)$ we can define a morphism $f$, with $f(x)(\id_x) = a$, by
 \begin{align*}
 f(y):R[x,y]_{\C} &\longrightarrow A(y) \\
 \alpha &\longmapsto A(\alpha)( a ).
 \end{align*}
 \end{proof}

The endomorphisms $[x,x]_{\C}$ of an object $x \in {\C}$ form an associative ring.  This ring will appear often, so we write $\End(x)$ instead of $[x,x]_{\C}$, and write $R\End(x)$ instead of $R \otimes_\ZZ \End(x)$.   

\begin{Remark}\label{remark:C eval at x gives Endx structure}
Given a covariant module $A$, evaluating $A$ at $x$ gives a left $R\End(x)$-module, using the action
\begin{align*}
 R\End(x) \times A(x) &\longrightarrow A(x) \\
 (f , a) &\longmapsto A(f)(a).
\end{align*}
This is a left $R\End(x)$-module structure since given any two elements $g, f \in R\End(x)$,
\[
(g \circ f)\cdot x  = A(g \circ f) (x) = A(g) \circ A(f) (x) = g \cdot (f \cdot x).
\]
Similarly, for a contravariant module $M$, $M(x)$ has a right $R\End(x)$-module structure.
\end{Remark}
\begin{Remark}\label{remark:C Ax has Endx structure}
Let $\widehat{\End(x)}$ denote the category with one object and with morphisms the free abelian group $\End(x)$.  Clearly $\widehat{\End(x)}$ has property (A) and it's possible to identify covariant $\widehat{\End(x)}$-modules and left $R\End(x)$-modules, similarly contravariant $\widehat{\End(x)}$-modules and right $R\End(x)$-modules.
\end{Remark}

There is often a need to consider bi-modules.  A $\C$-$\mathfrak{D}$ bi-module (can be covariant or contravariant in either variable, although most of the bi-modules we shall use will be covariant in one variable in contravariant in the other), is a functor
\[
Q(-, ?) : \C \times \mathfrak{D} \to \RMod.
\]

\begin{Example}\label{example:C free bifunctor}
The $\C$-$\C$ bi-module $R[-,?]_{\C}$ is defined as
\[
R[-,?]_{\C} : (x,y) \mapsto R[x,y]_{\C}.
\]
\end{Example}

\begin{Example}[The orbit category]\label{example:orbit category}\index{Orbit category $\OF$}
The orbit category, denoted $\OF$, is the prototypical example of a category with property (A), and will be studied properly in Chapter \ref{chapter:bredon}.  It was introduced for finite groups by Bredon \cite{Bredon-EquivariantCohomologyTheories}, who used the associated cohomology theory, Bredon cohomology, to develop equivariant obstruction theories.  It was later generalised to arbitrary groups by L\"uck \cite{Lueck}. 

Fix a discrete group $G$ and family $\mathcal{F}$ of subgroups of $G$, closed under taking subgroups and conjugation.  Commonly studied families are the family $\Fin$ of all finite subgroups, and the family $\VCyc$ of all virtually cyclic subgroups.  The objects of the orbit category $\OF$ are all transitive $G$-sets with stabilisers in $\mathcal{F}$, ie. the $G$-sets $G/H$ where $H$ is a subgroup in $\mathcal{F}$.  The morphism set $[G/H, G/K]_{\OF}$ is the free abelian group on the set of $G$-maps $G/H \to G/K$.  A $G$-map
\begin{align*}
\alpha: G/H &\longrightarrow G/K  \\
H &\longmapsto gK
\end{align*}
is completely determined by the element $\alpha(H) = gK$, and such an element $gK \in G/K$ determines a $G$-map if and only if $HgK = gK$, usually written as $gK \in (G/K)^H$.  Equivalently $gK$ determines a $G$-map if and only if $g^{-1}Hg \le K$.  In particular if $\mathcal{F} \subseteq \Fin$ then the orbit category has (EI), since any $G$-map $\alpha: G/K \to G/K$ is automatically an automorphism.  The isomorphism classes of elements in $\OF$, denoted $\Iso \OF$, are exactly the conjugacy classes of subgroups in $\mathcal{F}$.

\end{Example}
\begin{Remark}\label{remark:equivalent def of orbit category}
 The morphisms from $G/H$ to $G/K$ in the orbit category are usually defined as just the $G$-maps $G/H \to G/K$.  We show that this definition gives an isomorphic module category. 
 
 For this remark, let $\OF^\prime$ denote the category with the same objects as $\OF$ but with morphisms from $G/H$ to $G/K$ just the $G$-maps $G/H \to G/K$.  Let $\iota : \OF^\prime \to \OF$ be the faithful inclusion and given an $\OF$-module $M$ define an $\OF^\prime$ module $M^\prime = M\circ \iota$.  We claim that the functor $ M \mapsto M^\prime$ gives an isomorphism of categories between $\OF$-modules and $\OF^\prime$ modules.  
 
 Any $\OF^\prime$-module $M^\prime$ extends uniquely to an $\OF$-module $M$ by setting:
 \begin{enumerate}
  \item $M(G/H) = M^\prime(G/H)$ for all $H \in \mathcal{F}$.
  \item $M(\sum_i z_i \alpha_i ) = \sum_i z_i M^\prime(\alpha_i)$ for any $\OF$-morphism $\sum_i z_i \alpha_i$ written as the sum of $G$-maps $\alpha_i$.
 \end{enumerate}
 This gives an inverse to the functor $M\mapsto M^\prime$ described above.
\end{Remark}

\section{Tensor products}\label{section:C tensor products}

\subsection{Tensor product over \texorpdfstring{${\C}$}{C}}\label{subsection:C tensor product over C}

We describe a construction, due to L\"uck \cite[9.12]{Lueck}, of the categorical tensor product of \cite[16.7]{Schubert-KategorienII}\cite{Fisher-TensorProductOfFunctorsSatellitesAndDerivedFunctors} for the categories of $\C$-modules over $R$.  

For $M$ contravariant and $A$ covariant, the tensor product over $\C$ of $M$ and $A$ is \index{Tensor product over a category, $\otimes_{\C}$}
\[
M \otimes_{{\C}} A = \left. \bigoplus_{x \in {\C}} M(x) \otimes_R A(x) \right/ \sim
\]
where $M(\alpha)(m) \otimes a \sim m \otimes A(\alpha)(a)$ for all morphisms $\alpha \in [x, y]$ in ${\C}$, elements $m \in M(y)$ and $a \in A(x)$, and objects $x,y \in {\C}$.  Since $R$ is commutative, this construction yields an $R$-module.  The tensor product is associative \cite[Lemma 3.1]{MartinezPerez-SpectralSequenceBredon}, and commutes with direct sums.

\begin{Example}
 If $A$ is a left $R\End(x)$-module and $M$ is a right $R\End(x)$-module then, by Remark \ref{remark:C Ax has Endx structure}, $A$ and $M$ can be regarded as covariant and contravariant $\widehat{\End(x)}$-modules.  It's easy to check that
\[
M \otimes_{\widehat{\End(x)}} A \cong M \otimes_{R\End(x)} A.
\]
\end{Example}

\begin{Prop}\label{prop:C adjoint isomorphism of tensor product}\cite[p.166]{Lueck}\cite{MartinezPerez-SpectralSequenceBredon}
There are adjoint natural isomorphisms of $R$-modules:
\[
\Hom_{\mathfrak{D}}(M( {?}) \otimes_{{\C}} Q( {?}, {-}), N( {-})) \cong \Hom_{\C}(M( {?}), \Hom_{\mathfrak{D}}(Q( {?}, {-}), N( {-}))) 
\]
\[
\Hom_{\C}(Q( {?}, {-} ) \otimes_{\mathfrak{D}} A( {-}), B( {?})) \cong \Hom_{\mathfrak{D}}(A( {-}), \Hom_{{\C}}(Q( {?}, {-}), B( {?}))).
\]
Here $M$ and $N$ are contravariant modules, $A$ and $B$ are covariant modules, and $Q(?,-)$ is an $\mathfrak{D}$-${\C}$-bi-module---a contravariant $\mathfrak{D}$-module in ``\/$-$'' and a covariant $\C$-module in ``\/$?$''.
\end{Prop}

\begin{Lemma}\label{lemma:C yoneda type isos in tensor product}\cite[p.14]{MislinValette-BaumConnes}
There are natural isomorphisms of $R$-modules for any contravariant module $M$ and covariant module $A$:
\[
M \otimes_{{\C}} R[x, -]_{\C} \cong M(x)
\]
\[
R[-,x]_{\C} \otimes_{{\C}} A \cong A(x).
\]
\end{Lemma}

\subsection{Tensor product over \texorpdfstring{$R$}{R}}\label{subsection:C tensor product over R}

We describe the tensor product over $R$ as in \cite[9.13]{Lueck}.  If $A$ and $B$ are $\C$-modules, either both covariant or both contravariant, then the tensor product over $R$ is the $\C$-module
\[
(A \otimes_R B)(x) = A(x) \otimes_R B(x).
\]
If $\alpha: x \to y$ is a morphism in $\C$, then
\[
(A \otimes_R B) (\alpha) = A(\alpha) \otimes_R B(\alpha).
\]

\section{Frees, projectives, injectives and flats}\label{section:C frees projectives}

Free objects in a category are usually defined as left adjoint to some forgetful functor, often with codomain $\Set$.  For $\C$-modules the necessary forgetful functor is 
 \begin{align*}
U:\{\text{ $\C$-modules }\} &\longrightarrow [\text{Ob}(\C), \Set] \\
UA : x &\longmapsto A(x).
 \end{align*}
Here $[\text{Ob}(\C), \Set]$ denotes the category of functors $\text{Ob}(\C) \to \Set$, where $\text{Ob}(\C)$ is the category whose objects are the objects of $\C$ but with only the identity morphisms at each object.  The functor $F$ left adjoint to $U$ is, for $X \in [\text{Ob}(\C), \Set]$, 
\[
FX = \bigoplus_{x \in \C} \bigoplus_{X(x)} R[x,-]_{\C}.
\]
Analagously, if we are working with contravariant functors,
\[
FX = \bigoplus_{x \in \C} \bigoplus_{X(x)} R[-,x]_{\C}. 
\]

That $(F, U)$ form an adjoint pair is a consequence of the Yoneda-type Lemma \ref{lemma:C yoneda-type}---for any covariant module $A$,
\begin{align*}
 \Hom_{\C} (FX, A) &= \Hom_{\C} \left(\bigoplus_{x \in {\C}} \bigoplus_{X(x)} R\left[x,- \right]_{\C}, A\right) \\
 &\cong \prod_{x \in \C} \prod_{X(x)} \Hom_{\C} \left( R\left[x,- \right]_{\C}, A \right) \\
 &\cong \prod_{x \in {\C}} \prod_{X(x)} A(x) \\
 &\cong \Hom_{[\text{Ob}({\C}), \Set]} (X, UA).
\end{align*}
The proof for contravariant functors is analogous.

Projective and injective modules are defined as in any abelian category---a $\C$-module $P$ is projective if $\Hom_{\C}(P, -)$ is exact and a $\C$-module $I$ is injective if $\Hom_{\C}(-,I)$ is exact \cite[\S 2.2]{Weibel}.  Free modules are projective since if 
\[
0 \longrightarrow A \longrightarrow B \longrightarrow C \longrightarrow 0 
\]
is an exact sequence of $\C$-modules then, by the Yoneda-type Lemma (\ref{lemma:C yoneda-type}), applying $\Hom_{\C}(R[x,?]_{\C}, -)$ gives the exact sequence
\[
0 \longrightarrow A(x) \longrightarrow B(x) \longrightarrow C(x) \longrightarrow 0.
\]
Since direct sums of projectives are projective in any abelian category, this is enough to show the category of $\C$-modules has enough projectives, in fact the counit of the adjunction between $F$ and $U$,
\[
\eta: (FUA) \longrightarrow A,
\]
is always an epimorphism:  By construction, 
\[
FUA = \bigoplus_{x \in \C} \bigoplus_{a \in A(x)} F_a(x, -)  
\]
where $F_a (x, -) \cong R[x, -]_{\C}$.  The counit is the map defined on $F_a(x, -)$, via the Yoneda-type Lemma \ref{lemma:C yoneda-type}, by $\id_x \mapsto a$.  It's clear that every $a \in A(x)$ is in the image of $\eta(x)$, and thus $\eta$ is an epimorphism.

The category of $\C$-modules also has enough injectives, see Remark \ref{remark:C enough injectives} for a proof using coinduction.

A covariant (respectively contravariant) $\C$-module $F$ is flat if the functor $F \otimes_{\C} \dagger$ (respectively $\dagger \otimes_{\C} F$) is flat.  Lemma \ref{lemma:C yoneda type isos in tensor product} shows free modules are flat, and since the tensor product commutes with direct sums, projectives are flat also.

A covariant $\C$-module $M$ is said to be \emph{finitely generated} if there exists an epimorphism
\[
 \bigoplus_{x \in I} R[x,- ]_{\C} \longtwoheadrightarrow M,
\]
for some finite indexing set $I$ of objects in $\C$.  There is an analogous definition for contravariant $\C$-modules.  

\section{Restriction, induction and coinduction}\label{section:C restriction induction and coinduction}
Given a functor $\iota : \C \to \mathfrak{D}$, we define restriction, induction, and coinduction functors.  Induction and restriction can be found in \cite[\S 9.8]{Lueck} but with the names extension and restriction, he also defines an adjoint pair of functors called ``splitting'' and ``inclusion''.  We don't include these here as the adjointness of these functors relies on the (EI) property which we are not assuming holds.  

Restriction and induction are, for covariant modules:\index{Restriction $\Res_\iota$}\index{Induction $\Ind_\iota$}\index{Coinducton $\CoInd_\iota$}
\begin{align*}
 \Res_\iota : \{\text{Covariant $\mathfrak{D}$-modules} \} &\longrightarrow \{\text{Covariant $\C$-modules} \} \\
 \Res_\iota : A &\longmapsto A \circ \iota
\end{align*}
\begin{align*}
 \Ind_\iota : \{\text{Covariant $\C$-modules} \} &\longrightarrow \{\text{Covariant $\mathfrak{D}$-modules} \} \\
 \Ind_\iota : A &\longmapsto R[\iota(?), -]_{\mathfrak{D}} \otimes_{\C} A(?).
\end{align*}
Where the notation $R[\iota(?), -]_{\mathfrak{D}}$ means that in the variable ``?'', this functor should be regarded as a $\C$-module using $\iota$.  Coinduction is, for covariant modules:
\begin{align*}
 \CoInd_\iota : \{\text{Covariant $\C$-modules} \} &\longrightarrow \{\text{Covariant $\mathfrak{D}$-modules} \} \\
 \CoInd_\iota : A &\longmapsto \Hom_{\C} ( R[-, \iota(?)]_{\mathfrak{D}}, A(?) ).
\end{align*}

For contravariant functors, the definition of restriction is identical, and for induction and coinduction is nearly identical:
\begin{align*}
 \Ind_\iota : \{\text{Contravariant $\C$-modules} \} &\longrightarrow \{\text{Contravariant $\mathfrak{D}$-modules} \} \\
 \Ind_\iota : M &\longmapsto M(?) \otimes_{\C} R[-, \iota(?)]_{\mathfrak{D}} 
\end{align*}
\begin{align*}
 \CoInd_\iota : \{\text{Contravariant $\C$-modules} \} &\longrightarrow \{\text{Contravariant $\mathfrak{D}$-modules} \} \\
 \CoInd_\iota : M(-) &\longmapsto \Hom_{\C} ( R[\iota( ?), -]_{\mathfrak{D}}, M( ?) ).
\end{align*}

Usually the functor $\iota$ will be implicit, and we will use the notation $\Res_{\C}^\mathfrak{D}$ for $\Res_{\iota}$, and similarly for induction and coinduction.  We will also write $\Res_x^{\C}$ instead of $\Res_{\widehat{\End(x)}}^{\C}$ and similarly for induction and coinduction.

Note that for any left $R\End(x)$-module $A$, 
\[
\Ind_x^{\C}A(x) =  R[x, x] \otimes_{R\End(x)} A \cong A 
\]
\[
\CoInd_x^{\C} A(x) = \Hom_{R\End(x)} ( R[x,x], A ) \cong A. 
\]
Similarly for contravariant induction and coinduction.

\begin{Prop}\label{prop:ind coind and res adjoint isos}\cite[\S 2]{MartinezPerezNucinkis-MackeyFunctorsForInfiniteGroups}
 Induction is left adjoint to restriction and coinduction is right adjoint to restriction.
\end{Prop}

The following proposition is almost entirely a consequence of this adjointness.  

\begin{Prop}\label{prop:C properties of res ind coind}~
 \begin{enumerate}
  \item Restriction is exact.
  \item Induction is right exact and preserves frees, projectives, flats and ``finitely generated''.
  \item Coinduction preserves injectives.
  \item Induction and restriction preserve colimits and coinduction and restriction preserve limits.
 \end{enumerate}
\end{Prop}
\begin{proof}
\begin{enumerate}
 \item Since a short exact sequence of modules over $\C$ is exact if and only if it's exact when evaluated at every element of $\C$, restriction is always exact.  
 \item Since induction has an exact right adjoint it preserves projectives and is right-exact \cite[2.3.10, 2.6.1]{Weibel}. 
 
 That induction takes frees to frees is a consequence of Lemma \ref{lemma:C yoneda type isos in tensor product}, 
\[
\Ind_{\C}^{\mathfrak{D}} R[\dagger,-]_{\C} =  R[\iota( {?}), -]_{\mathfrak{D}} \otimes_{\C} R[\dagger, ?]_{\C} \cong R[\dagger,-]_{\mathfrak{D}}  
\]
and similarly for contravariant modules.  

 That induction takes flats to flats is a consequence of Lemma \ref{lemma:C nat iso ind res and tensor} below.  In the covariant case, this implies the functor $? \otimes_{\mathfrak{D}} \Ind^{\mathfrak{D}}_{\C} F$ is naturally isomorphic to the functor $ ( \Res^{\mathfrak{D}}_{\C} ? ) \otimes_{\C} F$.  Thus if $F$ is assumed flat then $? \otimes_{\mathfrak{D}} \Ind^{\mathfrak{D}}_{\C} F$ is exact.  An analogous proof holds for contravariant $F$.
 
 If $A$ is a finitely generated $\C$-module then there is an epimorphism $F \longtwoheadrightarrow A$ for some finitely generated free $F$.  Induction is right exact so there is an epimorphism 
 \[
 \Ind_{\C}^\mathfrak{D} F \longtwoheadrightarrow \Ind_{\C}^\mathfrak{D} A.
 \]
 Induction takes finitely generated frees to finitely generated frees so $\Ind_{\C}^\mathfrak{D} A$ is finitely generated.
 \item Since coinduction has an exact left adjoint it preserves injectives \cite[2.3.10]{Weibel} and is left-exact \cite[2.6.1]{Weibel}
 \item This is another consequence of adjointness \cite[p.118]{SaundersMaclane}.
\end{enumerate}
\end{proof}

\begin{Lemma}\label{lemma:C nat iso ind res and tensor}
 There exist natural isomorphisms for any contravariant $\C$-module $M$ and covariant $\C$-module $A$:
\[
M \otimes_{\mathfrak{D}}  \Ind_{\C}^\mathfrak{D} A \cong \Res_{\C}^\mathfrak{D} M \otimes_{\C} A 
\]
\[
\Ind_{\C}^\mathfrak{D} M \otimes_{\mathfrak{D}}   A \cong  M \otimes_{\C} \Res_{\C}^\mathfrak{D} A.
\]
\end{Lemma}
\begin{proof}
We prove the first natural isomorphism, the second is analogous:
  \begin{align*}
M \otimes_{\mathfrak{D}}  \Ind_{\C}^\mathfrak{D} A
&\cong M( -)  \otimes_{\mathfrak{D}} \big(  R[ ?,  -]_{\mathfrak{D}} \otimes_{\C} A( ?) \big) \\
&\cong \big( M( -)  \otimes_{\mathfrak{D}}  R[ ?,  -]_{\mathfrak{D}} \big) \otimes_{\C} A( ?) \\
&\cong \Res_{\C}^\mathfrak{D} M \otimes_{\C} A.
\end{align*}
\end{proof}

\begin{Remark}[The category of $\C$-modules has enough injectives]\label{remark:C enough injectives}
A consequence of Proposition \ref{prop:C properties of res ind coind}(3) is that the category of $\C$-modules has enough injectives.  For any ring $S$ and module $M$ over $S$ there always exists an injective module $I$ and injection $M \longhookrightarrow I$ \cite[2.3.11]{Weibel}.  Given a $\C$-module $M$ choose injective $R\End(x)$-modules $I_x$ such that $M(x)$ injects into $I_x$ for all $x \in \C$, and consider the map
 \[
 \prod_{x \in \C} \eta_x : M \longrightarrow \prod_{x \in \C} \CoInd_{R\End(x)}^{\C} I_x 
 \]
where $\eta_x$ is chosen, via the adjointness of coinduction and restriction, such that $\eta_x(x)$ is the inclusion of $M(x)$ into $\CoInd_x^{\C}I_x(x) = I_x$.
 
Clearly the product of the $\eta_x$ maps is an injection.  The module on the right is injective by Proposition \ref{prop:C properties of res ind coind}(3) and the fact that in any abelian category, products of injective modules are injective.
\end{Remark}

\begin{Example}
 If $A$ and $B$ are covariant $\C$-modules, we define a $\C$-$\C$ bi-module:
\[
A(?) \otimes_R B(-) : (x, y)  \mapsto A(x) \otimes A(y).
\]
 Denote by $\Delta : \C \to \C \times \C$ the diagonal functor $ \Delta : x \to (x,x)$.  The tensor product over $R$ defined in Section \ref{subsection:C tensor product over R} could be defined as
\[
A \otimes_R B = \Res_\Delta (A(?) \otimes_R B(-)).
\]
\end{Example}

\section{Tor and Ext}\label{section:C Tor and Ext}
Since the categories of $\C$-modules are abelian and have enough projectives, it is possible to use techniques from homological algebra to study them.  For $M$ a $\C$-module, a projective resolution $P_*$ of $M$ is an exact chain complex of $\C$-modules,
\[
 \cdots \longrightarrow P_i \longrightarrow P_{i-1} \longrightarrow \cdots \longrightarrow P_0 \longrightarrow M \longrightarrow 0
\]
where each $P_i$ is projective.

If $A$ is a covariant $\C$-module and $P_*$ a projective resolution of $A$ then for any covariant module $B$ and contravariant module $M$, we define $\Ext^*_{\C}$ and $\Tor_*^{\C}$ as \index{TorC@$\Tor^{\C}_*$}\index{ExtC@$\Ext^*_{\C}$}
\[
\Ext_{\C}^k (A, B ) = H^k \Hom_{\C} \big( P_* , B \big) 
\]
\[ 
\Tor^{\C}_k (M, A ) = H_k \big( M \otimes_{\C} P_* \big).
\]
We make the same definitions for contravariant modules, if $M$ is a contravariant module, $Q_*$ a projective resolution of $M$, $A$ a covariant module and $N$ a contravariant module then
\[ 
\Ext_{\C}^k (M, N ) = H^k \Hom_{\C} \big( Q_* , N \big) 
\]
\[ 
\Tor^{\C}_k (M, A ) = H_k \big( Q_* \otimes_{\C} A \big).
\]
A priori $\Tor^{\C}_*$ has just been given two definitions, these are equivalent by Proposition \ref{prop:C balancing Tor} below, an analogue of the classical result that $\Tor$ for modules over a ring can be computed using a resolution in either variable.

\begin{Prop}\label{prop:C balancing Tor}
If $A$ is any covariant module, $M$ is any contravariant module, $P_*$ is a projective covariant resolution of $A$, and $Q_*$ is a projective contravariant resolution of $M$ then for all $k$,
\[
H_k \big( M \otimes_{\C} P_* \big) \cong H_k \big( Q_* \otimes_{\C} A \big).
\]
\end{Prop}
We need some notation for the proof:  If $(C_*, \partial_*)$ is an arbitrary chain complex of $\C$-modules then we write $C_{*+j}$ for the chain complex whose degree $i$ term is $C_{i+j}$, and differential $(-1)^j\partial_{i+j}$.  This change in the differential doesn't affect exactness, as the homology groups of the new complex are simply $H_n(C_{*+j}) = H_{n+j}(C_*)$.
\begin{proof}
 The proof is a generalisation of \cite[Theorem 2.7.2, p.58]{Weibel} into the setting of $\C$-modules.  Form three double complexes, $M \otimes_{\C} P_*$, $Q_* \otimes_{\C} P_* $ and $Q_* \otimes_{\C} A$.  The augmentation maps $\varepsilon: P_* \longrightarrow A$ and $\eta: Q_* \longrightarrow M$ induce maps between the total complexes,
 \[
  \text{Tot}\left(Q_* \otimes_{\C} P_*\right) \longrightarrow \text{Tot} \left( M \otimes_{\C} P \right) \cong M \otimes_{\C} P_* 
 \]
 \[
 \text{Tot}\left(Q_* \otimes_{\C} P_*\right) \longrightarrow \text{Tot} \left( Q_* \otimes_{\C} A \right) \cong Q_* \otimes_{\C} A_* 
 \]
 where $\text{Tot}$ denotes the total complex of a bicomplex of $R$-modules (see \cite[1.2.6]{Weibel} for the definition of total complex).  We claim that these maps are weak equivalences.  Define a new double complex $C_{**}$, by adding $A_* \otimes_{\C} Q_{*-1}$ in the $(-1)$ column of $P_* \otimes_{\C} Q_*$, giving the following complex.  Note that we need to shift $Q_*$ so that the resulting complex is a bi-complex, without the shift the horizontal and vertical differentials would not anti-commute.
\[
\xymatrix{
\cdots \ar[d] & \cdots \ar[d] & \cdots \ar[d] \\
A \otimes_{\C} Q_2 \ar[d] & \ar[l] \ar[d] P_0 \otimes_{\C} Q_2 & P_1 \otimes_{\C} Q_2 \ar[l] \ar[d] & \cdots \ar[l] \\
A \otimes_{\C} Q_1 \ar[d] & \ar[l] \ar[d] P_0 \otimes_{\C} Q_1 & P_1 \otimes_{\C} Q_1 \ar[l] \ar[d] & \cdots \ar[l] \\
A \otimes_{\C} Q_0 \ar[d] & \ar[l] \ar[d] P_0 \otimes_{\C} Q_0 & P_1 \otimes_{\C} Q_0 \ar[l] \ar[d] & \cdots \ar[l] \\
0 & 0 & 0
}
\]
The complex $\text{Tot}(C_{**})_{*+1}$ is the mapping cone of $\varepsilon \otimes_{\C} \id_Q$, so it suffices to show that it is acyclic (see \cite[\S 1.5]{Weibel}).  But this follows from the Acyclic Assembly Lemma \cite[2.7.3]{Weibel}, since the flatness of $Q_i$ means the functor $\dagger \otimes_{\C} Q_i$ is exact for all $i$ and hence the rows of $C_{**}$ are exact.  

  Similarly, the mapping cone of $\id_P \otimes_{\C} \eta$ is the complex $\text{Tot}(D_{**})_{*+1}$, where $D_{**}$ is the double complex obtained by adding $P_{*-1} \otimes_{\C} B$ in row $(-1)$ to the complex $P_* \otimes_{\C} Q_*$.  Since $P_i$ is flat for all $i$, $P_i  \otimes_{\C} \dagger$ is exact, and the columns of $D_{**}$ are exact.  Thus $\text{Tot}(D_{**})_{*+1}$ is acyclic, again by the Acyclic Assembly Lemma \cite[2.7.3]{Weibel}, showing $\id_P \otimes_{\C} \eta$ is a weak equivalence.
 
\end{proof}

$\Tor_*^{\C}$ could also be calculated using flat resolutions instead of projective resolutions.  The standard proof of this in the case of modules over a ring goes through with almost no modification, see for example \cite[3.2.8]{Weibel}.  Similarly, we could calculate $\Ext_{\C}^*$ using injective resolutions, again the proof is the standard one. 

\section{Finiteness conditions}\label{section:C finiteness conditions}
 We define projective and flat dimensions as one would expect, the \emph{projective dimension} $\Cpd A$ of a contravariant $\C$-module $A$ is the minimal length of a projective resolution of $A$ and the \emph{flat dimension} $\Cfd A$ is the minimal length of a flat resolution.  These can be characterised as the vanishing of the $\Ext^*_{\C}$ and $\Tor_*^{\C}$ groups as in ordinary homological algebra. \index{Projective dimension $\Cpd$}
 
 Recall that a $\C$-module is finitely generated if it admits an epimorphism from a finite direct sum of modules of the form $R[x, -]_{\C}$ for some $x \in \C$.  We say a $\C$-module $A$ is $\CFP_n$ if there is a projective resolution of $A$ which is finitely generated up to degree $n$.\index{CFPn@$\CFP_n$ condition}  Additionally we call $\CFP_0$ modules \emph{finitely generated} and $\CFP_1$ modules \emph{finitely presented}.  There is an analogue of the Bieri-Eckmann criterion \cite{BieriEckmann-FinitenessPropertiesOfDualityGroups}, see also \cite[Theorem 1.3]{Bieri-HomDimOfDiscreteGroups}.  A proof in the case that $\C = \OF$ appears in \cite[Theorem 5.3]{MartinezNucinkis-GeneralizedThompsonGroups} and no substantial change is required to prove for $\C$-modules.
 
 \begin{Theorem}[Bieri--Eckmann Criterion]\label{theorem:C bieri-eckmann criterion}\index{Bieri--Eckmann criterion}
 The following conditions on any contravariant $\C$-module $A$ are equivalent:
 \begin{enumerate}
  \item $A$ is $\CFP_n$.
  \item If $B_\lambda$, for $\lambda \in \Lambda$, is a filtered system of $\C$-modules then the natural map 
  \[
  \varinjlim_\Lambda \Ext_{\C}^k ( A, B_{\lambda} ) \longrightarrow \Ext_{\C}^k(A, \varinjlim_\Lambda B_\lambda)    
  \]
  is an isomorphism for $k \le n-1$ and a monomorphism for $k = n$.
  \item For any filtered system $B_\lambda$, for $\lambda \in \Lambda$, such that $\varinjlim_\Lambda B_{\lambda} = 0$, 
  \[
  \varinjlim_\Lambda \Ext_{\C}^k ( A, B_{\lambda} ) = 0 
  \]
  for all $k \le n$.
   \item For any collection of indexing sets $\Lambda_x$, for $x \in \C$, the natural map
\[\Tor_k^{\C} \left(  M , \prod_{x \in \operatorname{Ob}\C} \prod_{\Lambda_x} R[x, -]_{\C} \right) \longrightarrow \prod_{x \in \operatorname{Ob}\C} \prod_{\Lambda_x} \Tor_k^{\C} \left( M, R[x, -]_{\C} \right) \]
 is an isomorphism for $k < n$ and an epimorphism for $k = n$.
 \end{enumerate}
 There is a similar result for covariant modules.
 \end{Theorem}

 \begin{Lemma}\label{lemma:C CFPn conditions for SES}
  If 
  \[0 \longrightarrow A \longrightarrow B \longrightarrow C \longrightarrow 0\]
  is a short exact sequence of $\C$-modules then
  \begin{enumerate}
   \item If $A$ and $B$ are $\CFP_n$ then $C$ is $\CFP_n$.
   \item If $A$ and $C$ are $\CFP_n$ then $B$ is $\CFP_n$.
   \item If $B$ and $C$ are $\CFP_n$ then $A$ is $\CFP_{n-1}$.
  \end{enumerate}
 \end{Lemma}
\begin{proof}
 This follows from the long exact sequence associated to $\Ext_{\C}^*$ and the Bieri--Eckmann criterion (Theorem \ref{theorem:C bieri-eckmann criterion}).
\end{proof}
\chapter{Bredon modules}\label{chapter:bredon}

Fix a family $\mathcal{F}$ of subgroups of $G$, closed under subgroups and conjugation, and recall from Example \ref{example:orbit category} that the orbit category $\OF$ is the category whose objects are all transitive $G$-sets with stabilisers in $\mathcal{F}$ and whose morphism set $[G/H, G/K]_{\OF}$ is the free abelian group on the set of $G$-maps $G/H \longrightarrow G/K$.  Common families to study are the family $\Fin$ of all finite subgroups and the family $\VCyc$ of all virtually cyclic subgroups.

Contravariant $\OFin$-modules and their associated finiteness conditions provide a good algebraic reflection of the geometric world of proper actions.  This background has already been discussed in the introduction and we will discuss connections with geometry in Sections \ref{section:bredon cohomology of spaces} and \ref{section:bredon cd} also.

Sections \ref{section:bredon free modules} and \ref{section:bredon res ind coind} specialise information from Chapter \ref{chapter:C} to modules over the orbit category, and the later sections discuss finiteness conditions for contravariant $\OF$-modules.

Recall that a $G$-map $\alpha: G/H \longrightarrow G/K $ is completely determined by the element $\alpha(H) = gK$, and such an element $gK \in G/K$ determines a $G$-map if and only if $HgK = gK$, equivalently $g^{-1}Hg \le K$.  

\section{Free modules}\label{section:bredon free modules}

For this section we require that $\mathcal{F} \subseteq \Fin$.  In this section we describe the structure of free $\OF$-modules. Throughout this section $H$ and $K$ will denote subgroups in $\mathcal{F}$ and $\alpha_g$ will denote a $G$-map $\alpha_g : G/H \to G/K$ sending $H \mapsto gK$ for any $H$ and $K$. 

\begin{Remark}[Structure of $\End(G/H)$]\label{remark:cov left AutG/H is right WH}
If $\alpha_g : G/H \longrightarrow G/H $ is the $G$-map sending $H \longmapsto gH$ then necessarily $g \in N_GH$ and two such $g$ determine the same $G$-map if they are in the same left $H$-coset.  Furthermore $\alpha_h \circ \alpha_g = \alpha_{gh}$ so, denoting by $WH$ the Weyl group $N_GH/H$, \index{Weyl group $WH$}
\[
\End(G/H) = \widehat{\ZZ[WH]}^\text{op}.
\]
Here $\ZZ[\widehat{WH}]$ denotes the category of one element and morphisms given by $\ZZ[WH]$, and $\widehat{\ZZ[WH]}^\text{op}$ is the opposite of that category.
As described in Remark \ref{remark:C eval at x gives Endx structure}, if $A$ is a covariant $\C$-module then evaluating at $x$ gives $A(x)$ a left $R\End(x)$-structure.  Thus evaluating a covariant $\OF$-module at $G/H$ gives a left $\widehat{R[WH]}^\text{op}$-structure, equivalently a right $R[WH]$-structure.  

Similarly, if $M$ is a contravariant $\OF$-module then evaluating at $G/H$ gives $M(G/H)$ a left $R[\End(G/H)]$ structure, equivalently a left $R[WH]$-module structure. 

Note that this description may fail when $\mathcal{F} \not\subseteq \Fin$, as it is possible to have an infinite cyclic subgroup $H$ of a group $G$ along with an element $g \in G$ such that $g^{-1}Kg$ is a proper subgroup of $K$.  This occurs for example in the Baumslag--Solitar group $\operatorname{BS}(1,2)$.
\end{Remark}

\begin{Example}[Right action of {$R[WK]$} on {$R[G/H, -]_{\OF}(G/K)$}]\label{example:right action on R[G/H G/K]}
The action of $WK$ on $ R[G/H,G/K]_{\OF}$ is as follows:  If $f :G/H \to G/K$ with $f(H) = gK$ and $w \in WK$ then 
\[
f \cdot w = R[G/H, \alpha_w]_{\OF} (f) = \alpha_w \circ f.
\]
Since $(\alpha_w \circ f)(1) = gwK$, under the identification 
\[
R[G/H, G/K]_{\OF} \cong R[(G/K)^H], 
\]
the action is given by $gK\cdot w = gwK$.
\end{Example}

\begin{Lemma}\label{lemma:R G/HG/K as an RWK module}
 There is an isomorphism of right $R[WK]$-modules
 \[
R[G/H, -]_{\OF}(G/K) = R[G/H, G/K]_{\OF} \cong \bigoplus_{\substack{g N_GK \in G/N_GK \\ g^{-1}Hg \le K}} R[WK].
 \]
\end{Lemma}
\begin{proof}
 Firstly, $R[G/H, G/K]_{\OF} \cong R[(G/K)^H]$ is a free $WK$-module, since if $n \in N_GK$ is such that $gnK = gK$ then $nK = K$ and hence $n \in K$.  Now, $gK$ and $g^\prime K$ lie in the same $WK$ orbit if and only if $g(WK) K = g^\prime (WK) K$, equivalently $gN_GK = g^\prime N_GK$, and $gK$ determines an element of $R[(G/K)^K]$ if and only if $g^{-1}Hg \le K$.  Thus there is one $R[WK]$ orbit for each element in the set 
\[
\{ gN_GK \in G / N_G K \: : \: g^{-1} H g \le K\}.
\]  
\end{proof}

For contravariant modules the situation is more complex, evaluating at $G/H$ doesn't always give a free $R[WH]$-module, although it does always give a $R[WH]$-module of type $\FP_\infty$.  This is proved in the case $R=\ZZ$ in \cite[Proof of 3.2]{KMN-CohomologicalFinitenessConditionsForElementaryAmenable}, the proof for general rings $R$ requires no substantial change, and is given in Corollary \ref{cor:Z[G/H,G/K] as finite sum of FPinfty WH modules}.

\begin{Example}[Left action of {$R[WH]$} on {$R[-, G/K]_{\OF}(G/H)$}]\label{example:left action on R[G/H G/K]}
A similar argument to the previous example shows that under the identification 
\[R[G/H, G/K]_{\OF} \cong R[(G/K)^H]\]
the action of $R[WH]$ is given by $w \cdot gK = wgK$.
\end{Example}

\begin{Lemma}\label{lemma:Z[G/H,G/K] as sum of WH submodules}
 There is an isomorphism of left $R[WH]$-modules 
\[
R[-, G/K]_{\OF}(G/H) = R[G/H, G/K]_{\OF} = \bigoplus_{x} R[WH / WH_{xK}] 
\]
where $x$ runs over a set of coset representatives of the subset of the set of $N_G H$-$K$ double cosets
\[
\{ x \in N_G H \char92 G / K \: : \: x^{-1}Hx \le K\},
\]
and the stabilisers are given by
\[
WH_{xK} = \left(N_G H \cap xKx^{-1}\right) / H.
\]
\end{Lemma}
\begin{proof}
Recall the identification $R[G/H, G/K]_{\OF} = R[(G/K)^H]$.  The elements $xK$ and $yK$ are in the same $WH$-orbit if there exists some $nH \in WH$ (where $n \in N_G H$) with 
\[
nHxK = yK \Leftrightarrow nxK = yK \Leftrightarrow (N_G H) x K = (N_G H) y K.
\]
Combining this with the fact that $xK \in (G/K)^H$ if and only if $x^{-1}H x \le K$ means there is a $WH$-orbit for each $N_G H$-$K$ double coset $N_G HxK$ such that $x^{-1} H x \le K$, i.e.~coset representatives for
\[
\{ x \in N_G H \char92 G / K \: : \: {x^{-1}}Hx \le K\}
\]
are orbit representatives for the $R[WH]$-orbits in $R[G/H, G/K]_{\OF}$.  

The $N_G(H)$-stabiliser of the point $xK \in (G/K)^H$ is the set
\[
\{ g \in N_G(H) \: : \: gxK = xK \} = \{g \in N_G(H) \: : \: g \in xKx^{-1}\} = N_G(H) \cap xKx^{-1}.
\]
So the $WH$-stabiliser of $xK \in (G/K)^H$ is $WH_{xK} = (N_G(H) \cap xKx^{-1})/H$.
\end{proof}

\begin{Cor}\label{cor:Z[G/H,G/K] as finite sum of FPinfty WH modules}
The $\OF$-module $R[-, G/K]_{\OF}(G/H) = R[G/H, G/K]_{\OF}$ is a finite direct sum of projective $R[WH]$-permutation modules of type  $\FP_\infty$ with stabilisers in $\mathcal{F}$.  In particular $R[G/H, G/K]_{\OF}$ is $\FP_\infty$.
\end{Cor}
\begin{proof}
Since $K$ is finite, the set $\{ x \in N_G H \backslash G / K \: : \: {x^{-1}}Hx \le K\} $ is finite and $R[G/H,G/K]_{\OF}$ can be written as a finite direct sum
\[
R[G/H, G/K]_{\OF} = \bigoplus_{x} R[WH/WH_{xK}] 
\]
where the $WH_{xK}$ are finite groups.  Since $R$ is $\FP_\infty$ as a $R[WH_{xK}]$-module and 
\[
R[WH/WH_{xK}] = \Ind_{R[WH_{xK}]}^{R[WH]} R,
\]
we can apply Lemma \ref{lemma:FPinfty as a F-module then as a G module} below and deduce that $ R[WH/WH_{xK}]$ is $\FP_\infty$ as an $RG$-module.  Finally, any finite direct sum of $\FP_\infty$ modules is $\FP_\infty$.
\end{proof}

\begin{Lemma}\label{lemma:FPinfty as a F-module then as a G module}
If $M$ is $\FP_\infty$ as an $R F$-module for some subgroup $F \le G$, then $\Ind_{RF}^{RG} M = R G \otimes_{R F} M$ is $\FP_\infty$ as an $R G$-module.
\end{Lemma}
\begin{proof}
Let $\prod_i N_i$ be an arbitrary direct product of $R G$-modules, then
\begin{align*}
 \Tor_*^{RG} \left( \Ind_{RF}^{RG} M, \prod_i N_i \right) &= \Tor_*^{RF} \left( M, \prod_i N_i \right) \\
&= \prod_i \Tor_*^{RF} \left( M, N_i \right) \\
&= \prod_i \Tor_*^{RG} \left( \Ind_{RF}^{RG} M, N_i \right)
\end{align*}
where the first and third equalities come from Shapiro's Lemma.  This finishes the proof as $\Ind_{RF}^{RG} M$ is $\FP_\infty$ if and only if $\Tor_*^{R G}(\Ind_{RF}^{RG} M, -)$ commutes with direct products \cite[Theorem VIII.4.8]{Brown}.
\end{proof}

\section{Restriction, induction and coinduction}\label{section:bredon res ind coind}

In this section we require that $\mathcal{F} \subseteq \Fin$.  We specialise the constructions of Section \ref{section:C restriction induction and coinduction} to the categories of covariant and contravariant $\OF$-modules.  In order to match the literature, we write $\Ind_{WH}^{\OF} A $ instead of $\Ind_{G/H}^{\OF} A $ for induction with the inclusion functor 
\[
\iota: \widehat{\End(G/H)} \longhookrightarrow \OF
\]
and similarly for restriction and coinduction.  Note that the notation for covariant and contravariant induction is the same, if neither covariant or contravariant is specified then contravariant should be assumed.

\begin{Example}\label{example:cov E_1R is uR}
If $R$ is the trivial $RG$ module then inducing to a covariant $\OF$-module gives
\[
\Ind^{\OF}_{RG}R : G/H \longmapsto R \otimes_{RG} R[G/H] = R .
\]
Checking the morphisms as well, $\Ind^{\OF}_{RG}R = \uR$, the constant covariant functor on $R$---sending every object to $R$ and every $G$-map to the identity.
\end{Example}

A group is said to \emph{contain no $R$-torsion} if for every finite subgroup $F \le G$, $\vert F \vert$ is invertible in $R$.  For example every group has no $\QQ$-torsion.  If \index{R-torsion@$R$-torsion}\label{Rtorsionfree}
\[
\vert F \vert = p_1^{n_1} \cdots p_m^{n_m} 
\]
is a prime factorisation of $\vert F \vert$ then for each $p_i$ there is an element of order $p_i$ by Cauchy's Theorem \cite[1.6.17]{Robinson}.  Since the invertible elements $R^*$ form a group, if all the $p_i$ are invertible in $R$ then so is $\vert F \vert$.  Hence a group has no $R$-torsion if and only if the order of every finite-order element is invertible in $R$.

Recall from Proposition \ref{prop:C properties of res ind coind} that covariant and contravariant restriction are exact, in addition we have the following:
\begin{Prop}\label{prop:co and contra res properties}~
\begin{enumerate}
 \item Covariant restriction preserves projectives and flats.
 \item Contravariant restriction preserves finite generation.
 \item 
 Contravariant restriction at $H$ preserves projectives and flats if $WH$ is $R$-torsion-free, if not then contravariant restriction takes projectives to $\FP_\infty$-modules.
\end{enumerate}
\end{Prop}
\begin{proof}
\begin{enumerate}
 \item If $P$ is a projective covariant $\OF$-module and $F$ a free covariant $\OF$-module with a split epimorphism $F \longtwoheadrightarrow P$ then restricting at $G/H$ yields a split epimorphism $F(G/H) \longtwoheadrightarrow P(G/H)$, by Lemma \ref{lemma:R G/HG/K as an RWK module} $F(G/H)$ is free and thus $P(G/H)$ is projective.

If $F$ is a flat covariant module and $M$ any left $R[WH]$-module then,  
\begin{align*}
F(G/H) \otimes_{R[WH]} M &\cong \left( R[-,G/H]_{\OF} \otimes_{\OF}F  \right) \otimes_{R[WH]} M \\
&\cong \left( R[-,G/H]_{\OF} \otimes_{R[WH]} M \right) \otimes_{\OF} F 
\end{align*}

Thus for any short exact sequence of left $R[WH]$-modules
\[
0 \longrightarrow M^\prime \longrightarrow M \longrightarrow M^{\prime\prime} \longrightarrow 0 
\]
applying $ F(G/H) \otimes_{R[WH]} -$ is equivalent to applying first the contravariant induction functor and then $\dagger \otimes_{\OF} F$.  Since contravariant induction is exact (Proposition \ref{prop:co and contra ind properties}(2)) and $F$ is assumed flat, exactness is preserved, and thus $F(G/H)$ is flat as required.

 \item Use the argument of the previous part, noting that Lemma \ref{lemma:Z[G/H,G/K] as sum of WH submodules} implies that for contravariant frees restricting at $G/H$ preserves finite generation.
 \item If $WH$ is $R$-torsion-free then, using Lemma \ref{lemma:Z[G/H,G/K] as sum of WH submodules}, restricting any free at $G/H$ gives a projective module, and the result follows.  To see that in this case, restriction preserves flats, let $F$ be a contravariant flat module and consider a short exact sequence
\[
0 \longrightarrow A \longrightarrow B \longrightarrow C \longrightarrow 0  
\]
of left $R[WH]$-modules, thus by Proposition \ref{prop:co and contra ind properties} below,
\[
0 \longrightarrow \Ind_{WH}^{\OF} A \longrightarrow \Ind_{WH}^{\OF} B \longrightarrow \Ind_{WH}^{\OF} C \longrightarrow 0
\]
is a short exact sequence of covariant modules.  Since $F$ is flat, the functor $\dagger \otimes_{\OF} F$ is exact, applying this to the above and using Lemma \ref{lemma:C nat iso ind res and tensor} gives a short exact sequence
\[
  0 \longrightarrow A \otimes_R F(G/H) \longrightarrow B \otimes_R F(G/H) \longrightarrow C \otimes_R F(G/H) \longrightarrow 0 
\]
showing $F(G/H)$ is flat.

If $WH$ is not $R$-torsion free then the result is just Corollary \ref{cor:Z[G/H,G/K] as finite sum of FPinfty WH modules}.
\end{enumerate}
\end{proof}

\begin{Example}\label{example:bredon cov res doesnt preserve fg}
 Unlike in the contravariant case, the covariant restriction functor does not preserve ``finitely generated'' in general:  Take for example the infinite dihedral group $D_\infty = (\ZZ/2\ZZ) \ast (\ZZ/2\ZZ)$ generated by the two elements $a$ and $b$ of order $2$.  The finite subgroup $\langle a \rangle$ is self-normalising, thus $R[W\langle a \rangle] = R$ and Lemma \ref{lemma:R G/HG/K as an RWK module} implies that as $R$-modules,
\[
R[D_\infty/1, D_\infty/\langle a \rangle]_{\OFin} = \bigoplus_{g\langle a \rangle \in D_\infty / \langle a \rangle} R.
\]
\end{Example}

\begin{Remark}\label{remark:cov restriction at 1 preserves fg}
The covariant restriction functor $\Res_{G}^{\OF}$ preserves ``finitely generated''.  Recall that
\[
R[G/K, G/1]_{\OF} \cong \left\{ \begin{array}{l l} RG & \text{ if } K = 1 \\ 0 & \text{ else.} \end{array} \right.  
\]
So if $A$ is an arbitrary finitely generated covariant $\OF$-module and $F$ a free covariant $\OF$-module with an epimorphism onto $A$ then $F(G/1)$ is finitely generated as an $RG$-module and since $\Res^{\OF}_{RG}$ is exact there is a surjection $F(G/1) \longtwoheadrightarrow A(G/1)$.
\end{Remark}

Recall from Proposition \ref{prop:C properties of res ind coind} that contravariant and covariant induction both preserve projectives, flats and finitely generation.  In addition we have the following facts.
\begin{Prop}\label{prop:co and contra ind properties}~
\begin{enumerate}
 \item 
 If $WH$ has no $R$-torsion the covariant induction functor $\Ind_{WH}^{\OF}$ is exact.  
 \item Contravariant induction is always exact.
\end{enumerate}
\end{Prop}
\begin{proof}
\begin{enumerate}
 \item Assume that $WH$ has no $R$-torsion, we must check that the functor 
\[ A \longmapsto A \otimes_{R[WH]} R[G/H,-]_{\OF} \]
is exact, where $A$ is an $R[WH]$-module.  Equivalently that for any subgroup $K$ in $\mathcal{F}$, the functor
\[ - \otimes_{R[WH]} R[G/H,G/K]_{\OF}  \]
is exact, but by Lemma \ref{lemma:Z[G/H,G/K] as sum of WH submodules}
\[R[G/H,G/K]_{\OF} = \bigoplus_{x \in I} R\left[WH/WH_{x}\right]  \]
for some finite indexing set $I$ and $WH_x$ finite subgroups of $WH$.  By Maschke's Theorem, $R\left[WH/WH_{x}\right]$ is projective, and hence flat, as an $R[WH]$-module.  Hence $ - \otimes_{R[WH]} R[G/H,G/K]_{\OF}  $ is indeed exact.

\item Similarly to the above, we must check the functor 
\[R[G/K, G/H]_{\OF} \otimes_{R[WH]} - \]
is exact, but by Lemma \ref{lemma:R G/HG/K as an RWK module}, $R[G/K, G/H]_{\OF}$ is free as an $R[WH]$-module, so this is automatic.
\end{enumerate}
\end{proof}

\section{Bredon homology and cohomology of spaces}\label{section:bredon cohomology of spaces}

Recall that a space $X$ is a \emph{$G$-CW complex} \cite[\S II.1]{Dieck-TransformationGroups} if there exists a filtration $\{X_i\}_{i \in \ZZ}$ of $X$ such that \index{G-CW-complex@$G$-CW-complex}

\begin{enumerate}
 \item $X$ has the colimit topology with respect to the filtration. 
 \item $X_{-1} = \emptyset$.
 \item $X_n$ is obtained from $X_{n-1}$ via a pushout of $G$-spaces:
\[
 \xymatrix{
   \displaystyle\coprod_{j \in \Delta_n} G/H_j \times S^{n-1} \ar[r] \ar[d] & X_{n-1} \ar[d] \\
   \displaystyle\coprod_{j \in \Delta_n} G/H_j \times D^{n} \ar[r] & X_n
 }
\]
\end{enumerate}

For any $j \in \Delta_n$ the image of $G/H_j \times D^{n}$ in $X$ is called an \emph{equivariant $n$-cell} with isotropy $H_j$.  We say $X$ has \emph{isotropy in $\mathcal{F}$} if the subgroups $H_j$ are elements of $\mathcal{F}$.  For example a $G$-CW complex is proper if and only if it has isotropy in $\Fin$ and is free if and only if it has isotropy in $\Triv$ (the family consisting of only the trivial subgroup).

\begin{Remark}\cite[II.(1.15)]{Dieck-TransformationGroups}
 If $X$ is a CW-complex with a $G$-action such that 
 \begin{enumerate}
  \item For all $g \in G$, the map $x \longmapsto gx$ takes cells to cells.
  \item If $g \in G$ fixes a cell $\sigma$ setwise then $g$ fixes $\sigma$ pointwise.
 \end{enumerate}
Then $X$ is a $G$-CW-complex.  Such an action is often called cellular or rigid.
\end{Remark}

$X$ is \emph{finite-dimensional} if $X = X_n$ for some integer $n$ and the minimal such $n$ is called the dimension, and $X$ is \emph{finite-type} if for all $n$, $X_n$ is obtained from $X_{n-1}$ by attaching finitely many equivariant $n$-cells (ie. the set $\Delta_n$ is finite).  $X$ is \emph{finite} (equivalently cocompact) if it is both finite-dimensional and finite-type. 

Let $\mathcal{F}$ be a family of subgroups such that $X$ has isotropy in $\mathcal{F}$.  Define the contravariant $\OF$-module \index{Cellular chain complex $C_*^{\OF}$}
\[ C_n^{\OF}(X) = \bigoplus_{j \in \Delta_n} \ZZ [ -,G/H_j]_{\OF}.\] 
Denoting by $C_n(X)$ the ordinary cellular chain complex of $X$ (see for example \cite[p.139]{Hatcher})
\[ C_n^{\OF}(X)(G/K) = C_n(X^K). \]
The boundary maps from the cellular chain complexes $C_n(X^K)$ give boundary maps for $C_n^{\OF}(X)$, including an augmentation map,
\[\varepsilon : C_0^{\OF}(X) \longrightarrow \uZZ, \]
which maps every $0$-cell in $C_0(X^K)$ to $1 \in \uZZ(G/K) \cong \ZZ$.  We obtain a free chain complex of contravariant $\OF$-modules 
\[ \cdots \longrightarrow C_n^{\OF}(X) \longrightarrow C_{n-1}^{\OF}(X) \longrightarrow \cdots \longrightarrow C_0^{\OF}(X) \longrightarrow \uZZ \longrightarrow 0. \]
The \emph{Bredon homology of $X$} with coefficients in some covariant module $A$ is
\[H_*^{\OF}(X, A) = H_* ( C_*^{\OF}(X) \otimes_{\OF} A) \]
and similarly the \emph{Bredon cohomology of $X$} with coefficients in some contravariant module $M$ is 
\[H^*_{\OF}(X, M) = H^* \Hom_{\OF}(C_*^{\OF}(X), M). \]

\section{Homology and cohomology of groups}\label{section:bredon homology and cohomology of groups}
Recall the definitions of $\Tor^{\OF}_*$ and $\Ext^*_{\OF}$ from Section \ref{section:C Tor and Ext}.  For a group $G$, covariant module $A$, and contravariant module $M$ we define \index{Bredon (co)homology $H^*_{\OF}(G, -)$ and $H_*^{\OF}(G, -)$}
\[ H^*_{\OF} (G, M) = \Ext^*_{\OF}(\uR, M) \]
\[ H_*^{\OF} (G, A) = \Tor_*^{\OF}(\uR, A). \]
Note that in both statements above, $\uR$ denotes the contravariant constant functor on $R$.

If $X$ is a model for $\EF G$ then the chain complex $C_*^{\OF}$ is exact, since evaluating at $G/H$ for any $H \in \mathcal{F}$ gives the cellular chain complex of the contractible space $X^H$.  Thus there are isomorphisms for any covariant $\OF$-module $A$ and contravariant module $M$:
\[ H^*_{\OF} (G, M) = H^*_{\OF}(X, M) \]
\[ H_*^{\OF} (G, A) = H_*^{\OF}(X, A). \] 

\section{Cohomological dimension}\label{section:bredon cd}

Recall from Section \ref{section:C finiteness conditions} that the projective dimension of an $\OF$-module $M$, denoted $\OFpd M$, is the minimal length of a projective $\OF$-module resolution of $M$.  We say that $G$ has \emph{Bredon cohomological dimension} $n$, written $\OFcd G = n$, if $\OFpd \uR = n$ where $\uR$ is the constant contravariant $\OF$-module.  If we want to emphasize the ring $R$ we will write $\OFcd_R$ instead of $\OFcd$. \index{Bredon cohomological dimension $\OFcd G$}

As mentioned in the previous section, if $X$ is a model for $\EF G$ then $C_*^{\OF}(X)$ is an exact resolution of $\uZZ$ by free $\OF$-modules, hence $\OFcd_{\ZZ} G \le \gdF G$ (recall $\gd_\mathcal{F}$ is the minimal dimension of a model for $\EF G$).  A theorem of L\"uck and Meintrup in the high dimensional case and of Dunwoody in the dimension $1$ case provides an inequality in the other direction when $\mathcal{F} = \Fin$.

\begin{Theorem}\cite[Theorem 0.1]{LuckMeintrup-UniversalSpaceGrpActionsCompactIsotropy}\cite{Dunwoody-AccessabilityAndGroupsOfCohomologicalDimensionOne}
 Except for the possibility that $\OFincd_{\ZZ} G = 2$ and $\gd_{\Fin} G = 3$, $\OFincd_{\ZZ} G = \gd_{\Fin} G$.
\end{Theorem}

In \cite{BradyLearyNucinkis-AlgAndGeoDimGroupsWithTorsion}, Brady, Leary and Nucinkis construct groups $G$ with $\OFincd_{\ZZ} G = 2$, but $\gd_{\Fin} G = 3$ , showing the bound is sharp.  

\subsection{Low dimensions}\label{subsection:Low Dimensions}
Recall that $\cd_R G = 0$ if and only if $G$ is finite with no $R$-torsion \cite[Proposition 4.12]{Bieri-HomDimOfDiscreteGroups}.

\begin{Prop}\label{prop:ucdG=0 over R iff G finite}
 $\OFcd_R G = 0$ if and only if there exists a subgroup $H \in \mathcal{F}$ with $\lvert G/H \rvert$ invertible in $R$ and every $K \in \mathcal{F}$ is subconjugate to $H$.  In particular, $\OFincd_R G = 0$ if and only if $G$ is finite and and $\OVCyccd_R G  =0$ if and only if $G$ is virtually cyclic.
\end{Prop}

A proof of this when $R = \ZZ$ is available in \cite[Prop 3.20]{Fluch}, there are some minor modifications needed to generalise to arbitrary rings $R$.

\begin{proof}
 Using a more general definition of family of subgroups $\mathcal{F}$ than we use here, Symonds proves that $\uR$ is projective if and only if every component of $\mathcal{F}$ has a unique maximal element $M$ and $\lvert N_GM : M \rvert$ is finite and invertible in $R$, where he views $\mathcal{F}$ as a poset with inclusion \cite[Lemma 2.5]{Symonds-BredonCohomologyOfSubgroupComplexes}.  Since we assume $\mathcal{F}$ is closed under intersection, for us $\mathcal{F}$ may have only one component.  Also, since we assume $\mathcal{F}$ is closed under conjugation we must have $N_GM = G$---if $g \in G \setminus N_GM$ then since $M$ is maximal $M^g \lneq M$ and thus $M \lneq M^{g^{-1}}$ contradicting maximality of $M$.  The proposition now follows immediately from Symonds' result and the fact that $\OFcd G = 0$ if and only if $\uR$ is projective.
\end{proof}

Combining \cite[Proposition 4.12]{Bieri-HomDimOfDiscreteGroups} and Proposition \ref{prop:ucdG=0 over R iff G finite}, $\OFincd_\ZZ G = 0$ if and only if $\cd_\QQ G = 0$ if and only if $G$ is finite.

Recall that $\cd_\ZZ G = 1$ if and only if $G$ is a free group \cite{Stallings-OnTorsionFreeGroupsWithInfinitelyManyEnds,Swan-GroupsOfCohomologicalDimensionOne}, $\cd_R G = 1$ if and only if $G$ is $R$-torsion-free and acts properly on a tree, and $\cd_\QQ G = 1$ if and only if $G$ acts properly on a tree or equivalently $G$ is virtually-free \cite{Dunwoody-AccessabilityAndGroupsOfCohomologicalDimensionOne}.

\begin{Lemma}\label{lemma:ucdZG = 1 iff cdQG = 1}
For any group $G$, $\OFincd_\ZZ G = 1$ if and only if $\cd_\QQ G = 1$.  
\end{Lemma}
\begin{proof} 
If $\OFcd_\ZZ G = 1$ then Lemma \ref{lemma:ucdR less ucdZZ and OFFPn over ZZ implies OFFPn over R} implies $\OFcd_\QQ G \le 1$ and Lemma \ref{lemma:no R torsion then cdR less ucd R etc} implies $\cd_\QQ G \le 1$.  Since $G$ is not finite, $\cd_\QQ G = 1$.

If $\cd_\QQ G = 1$ then by \cite[Theorem 1.1]{Dunwoody-AccessabilityAndGroupsOfCohomologicalDimensionOne}, $G$ acts properly and with finite stabilisers on a tree $T$. For any finite subgroup $H \le G$, $H$ acts on $T$, $T^H \neq \emptyset$ and in particular $T^H$ is a sub-tree of $T$ \cite[6.1, 6.3.1]{Serre}.  $T$ is thus a model for $\text{E}_{\Fin} G$ and $\OFcd_\ZZ G = 1$.
\end{proof}

\begin{Cor}\label{cor:tfae cdRG less 1}
The following are equivalent for an infinite group $G$, and any ring $R$:
\begin{enumerate}
 \item $\cd_R G = 1$.
 \item $G$ has no $R$-torsion and $\OFincd_\ZZ G = 1$.
 \item $G$ has no $R$-torsion and $\OFincd_R G = 1$.
\end{enumerate}
\end{Cor}
\begin{proof}
\begin{itemize}
\item[$1\Rightarrow2$] If $\cd_R G = 1$ then $G$ has no $R$-torsion \cite[Proposition 4.11]{Bieri-HomDimOfDiscreteGroups} and $G$ acts properly on a tree \cite[Theorem 1.1]{Dunwoody-AccessabilityAndGroupsOfCohomologicalDimensionOne}.  By the argument of Lemma \ref{lemma:ucdZG = 1 iff cdQG = 1} the tree is a model for $\text{E}_{\Fin} G$ and hence $\OFincd_\ZZ G = 1$.  
\item[$2\Rightarrow3$]  Lemma \ref{lemma:ucdR less ucdZZ and OFFPn over ZZ implies OFFPn over R}.
\item[$3\Rightarrow1$]  Lemma \ref{lemma:no R torsion then cdR less ucd R etc}.
\end{itemize}
\end{proof}

\begin{Question}\label{question:bredon OFcdR G = 1}
What does the condition $\OFcd_R G = 1$ represent?  Is it equivalent to $\OFcd_\ZZ G = 1$?
\end{Question}

\section{\texorpdfstring{$\FPn$}{FPn} conditions}\label{section:bredon FPn}

Recall from Section \ref{section:C finiteness conditions} that an $\OF$-module $M$ is $\OFFP_n$ if there is a resolution of $M$ by projective $\OF$-modules, finitely generated up to dimension $n$.  We say $G$ is $\OFFP_n$ if $\uR$ is $\OFFP_n$, if $G$ is $\OFFP_\infty$ with finite Bredon cohomological dimension then we say $G$ is $\OFFP$.  \index{OFFPn@$\OFFP_n$ condition}

If $G$ admits a model $X$ for $\EF G$ with cocompact $n$-skeleton then the chain complex $C_*^{\OF}(X)$ is finitely generated up to dimension $n$ and so $G$ is of type $\OFFP_n$ over $\ZZ$.  Conversely, if $G$ is $\OFFP_n$ over $\ZZ$ and $WH$ is finitely presented for every finite subgroup then $G$ has a model for $\EF G$ with cocompact $n$-skeleton \cite[Theorem 0.1]{LuckMeintrup-UniversalSpaceGrpActionsCompactIsotropy}.

\begin{Prop}\label{prop:uFP_0 iff finitely many conj classes of finite subgroups}
$G$ is $\OFFP_0$ over $R$ if and only if there exists a finite set $H_1, \ldots H_m \in \mathcal{F}$ such that every $K \in \mathcal{F}$ is subconjugate to some $H_i$.  In particular, if $\mathcal{F} \subseteq \Fin$, $G$ is $\OFFP_0$ over $R$ if and only if there are finitely many conjugacy classes of subgroups in $\mathcal{F}$.
\end{Prop}
The case $\mathcal{F} = \Fin$ appears in \cite[Lemma 3.1]{KMN-CohomologicalFinitenessConditionsForElementaryAmenable}.
\begin{proof}
If $G$ is $\OFFP_0$ then there is an epimorphism,
\[\bigoplus_{i=1}^m R[-,G/H_i]_{\OF} \longtwoheadrightarrow \uR,\]
where the indexing set $I$ is finite.  Let $K$ be a subgroup in $\mathcal{F}$, evaluating at $G/K$ gives a surjection
\[\bigoplus_{i=1}^m R[G/K,G/H_i]_{\OF} \longtwoheadrightarrow R,\]
so for some $i$ we have $R[G/K,G/H_i]_{\OF} \neq 0$ and hence $K$ is subconjugate to one of the $H_i$.  

For the converse, one checks that the augmentation map 
\[ \bigoplus_{i =1}^m R[-,G/H_i]_{\OF} \longrightarrow \uR\] 
is a surjection.

If $\mathcal{F} \subseteq \Fin$ then observe that each $H_i$ has at most finitely many subconjugate subgroups, so the existence of such a collection $H_1, \ldots, H_m$ is equivalent to $\mathcal{F}$ having finitely many conjugacy classes.
\end{proof}

\begin{Prop}\label{prop:M is uFPn iff uFP0 and M(G/K) FPn over WK}
  Let $G$ be $\OFFP_0$ and $\mathcal{F} \subseteq \Fin$, then a contravariant module $M$ is $\OFFP_n$ ($n \ge 1$) over $R$ if and only if $M(G/K)$ is of type $\FP_n$ over $R[WK]$ for all subgroups $K$ in $\mathcal{F}$.
\end{Prop}
The proof in the case $R = \ZZ$ appears as \cite[Lemma 3.2]{KMN-CohomologicalFinitenessConditionsForElementaryAmenable} and requires no substantial alteration to generalise to arbitrary rings $R$.

\begin{Question}\label{question:characterise OFFPn for any F}
 Is there an easy characterisation of the condition $\OFFP_n$ for arbitrary $\mathcal{F}$, or for $\mathcal{F} = \VCyc$?
\end{Question}

\begin{Cor}\label{cor:uFPn equivalent conditions}
The following are equivalent for a group $G$ and $\mathcal{F} \subseteq \Fin$,
\begin{enumerate}
 \item $G$ is $\OFFP_n$ over $R$.
 \item $G$ is $\OFFP_0$ and the Weyl groups $WK$ are $\FP_n$ over $R$ for all $K \in \mathcal{F}$.
 \item $G$ is $\OFFP_0$ and the centralisers $C_G K$ are $\FP_n$ over $R$ for all $K \in \mathcal{F}$.
\end{enumerate}
\end{Cor}
\begin{proof}
 By the previous Proposition (1) and (2) are equivalent.  To see the equivalence of (2) and (3) consider the short exact sequence
\[0 \longrightarrow K \longrightarrow N_G K  \longrightarrow WK \longrightarrow 0. \]
$K$ is finite and hence $\FP_\infty$, so $WK$ is $\FP_n$ over $R$ if and only if $N_G K$ is $\FP_n$ over $R$ \cite[Proposition 2.7]{Bieri-HomDimOfDiscreteGroups}.  Since $K$ is finite, so $C_G K$ is finite index in $N_G K $ \cite[1.6.13]{Robinson} and so $C_G K$ is $\FP_n$ over $R$ if and only if $N_G K$ is $\FP_n$ over $R$.  Combining the last two results gives $WK$ is $\FP_n$ over $R$ if and only if $C_G K$ is $\FP_n$ over $R$.
\end{proof}

\begin{Example}\label{example:abels group}
 In \cite{BieriStrebel-ValuationsAndFinitelyPresentedMetabelianGroups}, it's shown that Abels' group is $\FP_2$ over $\QQ$ but not over $\ZZ$.  The Bestvina Brady groups also provide examples of groups which are $\FP_n$ over some rings but not others \cite{BestvinaBrady-MorseTheoryAndFinitenessPropertiesOfGroups}.
\end{Example}

\subsection{Quasi-\texorpdfstring{$\OFFP_n$}{OFFPn} conditions}\label{subsection:bredon quasi-FPn}\index{Quasi-$\OFFP_n$ condition}

In \cite[\S 6]{MartinezNucinkis-GeneralizedThompsonGroups}, Mart\'{\i}nez-P\'erez and Nucinkis define the \emph{quasi-$\OFFP_n$} condition, a weakening of $\OFFP_n$, these are defined for all families $\mathcal{F} \subseteq \Fin$.  We will need these conditions in Chapter \ref{chapter:H}.  A group $G$ is quasi-$\OFFP_n$ if $WK$ is $\FP_n$ for all $K \in \mathcal{F}$ and $G$ has finitely many conjugacy classes of subgroups in $\mathcal{F}$ isomorphic to a given finite subgroup. 

For any positive integer $k$, define the module 
\[ \uR_k (G/H) = \left\{ \begin{array}{l l} R & \text{if }\lvert H \rvert \le k, \\0 & \text{otherwise.} \end{array} \right. \]
Then $G$ is quasi-$\OFFP_n$ if and only if $\uR_k$ is $\OFFP_n$ for all positive integers $k$ \cite[Proposition 6.5]{MartinezNucinkis-GeneralizedThompsonGroups}.  

Say $G$ is quasi-$\OF\negthinspace\F_n$ if for all positive integers $k$, $G$ has a finite type model for $\mathop{E_{\mathcal{F}_k}G}$, where $\mathcal{F}_k$ is the subfamily of $\mathcal{F}$ containing all subgroups of order less than $k$.  If $G$ is quasi-$\OFFP_n$ then $G$ is quasi-$\OF\negthinspace\F_n$ if and only if the centralisers $C_GK$ are finitely presented for all $K \in \mathcal{F}$ \cite[Proposition 6.10]{MartinezNucinkis-GeneralizedThompsonGroups}.

We can give a geometric meaning to the quasi-$\OF\negthinspace\F_n$ conditions---a group $G$ is quasi-$\OF\negthinspace\F_n$ if and only if $G$ admits a model for $\EF G$ which is a mapping telescope of models for $\mathop{E_{\mathcal{F}_k}G}$ with cocompact $n$-skeleta \cite[Theorem 6.11]{MartinezNucinkis-GeneralizedThompsonGroups}.

\section{Change of rings}\label{section:bredon change of rings}

If $\varphi: R_1 \to R_2$ is a ring homomorphism then we define the change of rings functor $\varphi^*$ from $\OF$-modules over $R_2$ to $\OF$-modules over $R_1$ as follows,
\[\varphi^*A : G/H \mapsto A(G/H), \]
where we are viewing $A(G/H)$ as an $R_1$-module via $\varphi$.  On morphisms:
\[ \varphi^*A\left(\sum_i r_i \alpha_i \right) = \sum_i \varphi(r_i)A(\alpha_i) \]
where $r_i \in R_1$ and the $\alpha_i$ are morphisms $G/H \to G/K$ for some $G/H, G/K \in \OF$.

We also define a functor $R_2 \otimes_{R_1}-$ from $\OF$-modules over $R_1$ to $\OF$-modules over $R_2$ by 
\[
R_2 \otimes_{R_1} A : G/H \mapsto R_2 \otimes_{R_1}A(G/H)
\]
where we are using $\varphi$ to view $R_2$ as an $R_1$-module.  Applying this to a free module gives
\[ 
R_2  \otimes_{R_1} R_1[-, G/H]_{\OF}  \cong R_2[-, G/H]_{\OF}.
\]
Hence if $P$ is a projective $\OF$-module over $R_1$ then $R_2 \otimes_{R_1} P$ is a projective $\OF$-module over $R_2$.

\begin{Lemma}\label{lemma:ucdR less ucdZZ and OFFPn over ZZ implies OFFPn over R}
If $\OFcd_\ZZ G \le n$ then $\OFcd_R G \le n$ for all rings $R$.  Similarly if $G$ is $\OFFP_n$ over $\ZZ$ then $G$ is $\OFFP_n$ over $R$ for all rings $R$.
\end{Lemma}
\begin{proof}
 For the first part, take a projective resolution of $\uZZ$ by contravariant $\OF$-modules of length $n$ and define a new resolution by $Q_n(G/H) = R \otimes_{\ZZ} P_n(G/H) $ for all $n \in \NN$ and $G/H \in \OF$.  Since for any $H \in \mathcal{F}$ the complex $P_*(G/H)$ is $\ZZ$-acyclic and hence $\ZZ$-split so $Q_*(G/H)$ is acyclic also.  Finally each $Q_n$ is projective.  
 
 The second part is similar---choose the projective $\OF$-module resolution of $\uZZ$ to be finitely generated in all degrees $i \le n$ and use that if $P_i$ is finitely generated then so is $Q_i$.
\end{proof}

\begin{Lemma}\label{lemma:no R torsion then cdR less ucd R etc}
If $G$ has no $R$-torsion then $\cd_R G \le \OFcd_R G$.
\end{Lemma}
\begin{proof}
Take a projective resolution of $\uR$ by contravariant modules of length $n$ and evaluate at $G/1$, since $G$ is $R$-torsion-free, Proposition \ref{prop:co and contra res properties}(3) implies that $P_*(G/1)$ is a length $n$ projective resolution of $R$ by $R G$-modules.
\end{proof}

\begin{Prop}\label{prop:bredon change of rings}
If $\varphi:R_1 \to R_2$ is a ring homomorphism and $A$ is a $\OF$-module over $R_2$ then
\[\Tor^{R_1, \OF}_*(\uR_1, \varphi^* A) \cong \Tor^{R_2, \OF}_*(\uR_2, A). \]
There are similar isomorphisms for contravariant modules and for $\Ext_{\OF}^*$.
\end{Prop}
\begin{proof}
Firstly, consider the case $\varphi: \ZZ \to R$ for some ring $R$, we prove
\[ \Tor^{\ZZ, \OF}_*(\uZZ, \varphi^* A) = \Tor^{R, \OF}_*(\uR, A).\]
 Choose a resolution $P_*$ of $\uZZ$ by contravariant projective $\OF$-modules over $\ZZ$.  For any $G/H$ in $\OF$, $P_*(G/H)$ is a $\ZZ$-split resolution, so applying the functor $R \otimes_\ZZ -$ to $P_*$ yields a projective resolution of $\uR$ by projective $\OF$-modules over $R$.  Observing that 
\[ P_* \otimes_{\OF, \ZZ} \varphi^*A \cong (P_* \otimes_\ZZ R ) \otimes_{\OF, R} A \]
completes the proof.  

For the general case, let $\varphi_1: \ZZ \to R_1$ and $\varphi_2:\ZZ \to R_2$ be (unique) ring homomorphisms, then $\varphi \circ \varphi_1 = \varphi_2$ and $\varphi_1^* \circ \varphi^* = \varphi_2^* $.  Applying the previous part twice gives
\begin{align*}
 \Tor^{R_1, \OF}_*(\uR_1, \varphi^* A) &\cong \Tor^{\ZZ, \OF}_*(\uZZ, \varphi_1^*\circ \varphi^* A) \\
&\cong \Tor^{R_2, \OF}_*(\uR_2, A).
\end{align*}
\end{proof}

The next result is essentially \cite[1.4.3]{Hamilton-Thesis}, where it is proved for rings of prime characteristic in the setting of ordinary group cohomology.
\begin{Prop}
 Given some integer $m > 0$ and ring $R$ with characteristic $m$, then $G$ is $\OFFP_n$ over $R$ if and only if $G$ is $\OFFP_n$ over $\ZZ/m\ZZ$.
\end{Prop}
\begin{proof}
The proof below is for contravariant modules, the proof for covariant modules is analogous.

Assume that $G$ is $\OFFP_n$ over $\ZZ/m\ZZ$.  If $M_*$ is any directed system of contravariant $\OF$-modules over $R$ with $\varinjlim M_* = 0$, we necessarily have $\varinjlim \varphi^* M_* = 0$.  By the Bieri--Eckmann criterion (Theorem \ref{theorem:C bieri-eckmann criterion}), and the fact that $\underline{\ZZ/m\ZZ}$ is assumed $\OFFP_n$ over $\ZZ/m\ZZ$, we have that for all $i \le n$,
\[ \varinjlim \Ext^i_{\OF, \ZZ/m\ZZ}(\underline{\ZZ/m\ZZ}, \varphi^*M_*) = 0. \]
Thus by Proposition \ref{prop:bredon change of rings} applied to the canonical map $\ZZ/m\ZZ \to R$,
\[\varinjlim \Ext^i_{\OF, R}(\uR, M_*) =  0. \]
The Bieri--Eckmann criterion (Theorem \ref{theorem:C bieri-eckmann criterion}) gives that $\uR $ is $\OFFP_n$ over $R$.

For the ``only if'' direction, suppose $M_*$ is a directed system of $\OF$-modules over $\ZZ/m\ZZ$, with $\varinjlim M_* = 0$ thus $\varinjlim M_* \otimes_{\ZZ / m\ZZ} R = 0$ and by Theorem \ref{theorem:C bieri-eckmann criterion} for all $i \le n$, 
\[\varinjlim \Ext_{\OF, R}^i(\uR, M_* \otimes_{\ZZ/ m\ZZ} R) = 0. \]
Combining with Proposition \ref{prop:bredon change of rings}
\begin{align*} 
\varinjlim \Ext_{\ZZ, \OF}^i(\underline{\ZZ/m\ZZ}, M_* \otimes_{\ZZ/ m\ZZ} R) &= \varinjlim \Ext_{R, \OF}^i(\uR, M_* \otimes_{\ZZ/ m\ZZ} R) \\
&=0.
\end{align*}

Since $\ZZ/ m\ZZ$ is self-injective \cite[Cor 3.13]{Lam}, $R$ splits as a $\ZZ/m\ZZ$ module into $ R \cong \ZZ / m\ZZ \oplus N$ where $N$ is some $\ZZ/m\ZZ$ module.  Thus we have
\begin{align*}
\varinjlim \Big(  \Ext^i_{\ZZ/m\ZZ, \OF}&(\underline{\ZZ/m\ZZ}, M_*) & \\
&\oplus \Ext^i_{\ZZ/m\ZZ, \OF}(\underline{\ZZ/m\ZZ}, M_* \otimes_{\ZZ/m\ZZ}N )  \Big) = 0.
\end{align*}

In particular 
\[\varinjlim \Ext_{\ZZ/m\ZZ, \OF}^i(\underline{\ZZ/m\ZZ}, M_*) = 0\]
so by the Bieri--Eckmann criterion (Theorem \ref{theorem:C bieri-eckmann criterion}) $\underline{\ZZ/m\ZZ}$ is $\OFFP_n$ over $\ZZ/m\ZZ$, i.e.~$G$ is $\OFFP_n$ over $\ZZ/m\ZZ$.
\end{proof}

\begin{Remark}
 This proposition fails in characteristic zero as the ring $\ZZ$ is not self-injective.  For example $\QQ$ is not isomorphic, as a $\ZZ$-module, to $N \oplus \ZZ$ for any $\ZZ$-module $N$.
\end{Remark}

\section{Some interesting examples}\label{section:bredon interesting examples}

 By the right-angled Coxeter group $(W,S)$ corresponding to some flag complex $L$ we mean the group $W$ generated by a set $S$ of involutions where $S$ is in bijection with the vertices of $L$ and two involutions commute if and only if they are adjacent in $L$.  Given such a $(W,S)$ we let $\mathcal{S}$ be the poset of spherical subsets of $S$ (subsets generating a finite subgroup of $W$) and form the geometric realisations $K = \lvert \mathcal{S} \rvert$ and $\partial K = \lvert \mathcal{S}_{> \emptyset} \rvert$.  
 
 We form the simplicial complexes $\mathcal{U}(W, \partial K)$ and $\mathcal{U}(W, K)$ as in \cite[5.1.2]{Davis}, both admit $W$-actions and $\mathcal{U}(W, \partial K)$ is the singular set of $\mathcal{U}(W, K)$ (sub-complex with non-trivial isotropy).  The complex $\mathcal{U}(W, K)$, often called the Davis complex, is known to be a model for $\text{E}_{\Fin} W$ \cite[Theorem 12.3.4(ii)]{Davis}. \index{Davis complex $\mathcal{U}$}

\begin{Lemma}\label{lemma:bredon UWK acyclic}
\cite[8.2.8]{Davis} $\mathcal{U}(W, \partial K)$ is $R$-acyclic if and only if $(\partial K)_T$ is $R$-acyclic for all spherical subsets $T \in \mathcal{S}$, where
\[
K_T = \bigcap_{s \in T} \lvert \mathcal{S}_{\ge s} \rvert.
\]
\end{Lemma}
 
 The example below first appeared in \cite{Bestvina-VirtualCohomologicalDimensionOfCoxeterGroups}, see also \cite[8.5.8]{Davis}, and much of the following argument appears in \cite[proof of Theorem 2]{DicksLeary-SubgroupsOfCoxeterGroups}.  
 
\begin{Example}[A group $W$ with {$\OFincd_{\FF_3} W = 2$} but {$\OFincd_\ZZ W = 3$} which is not torsion-free]
Consider the right-angled Coxeter group $(W, S)$ corresponding to the barycentric subdivision $L$ of the ordinary triangulation of $\RR \PP^2$.

\emph{Claim: $\OFincd_{\ZZ}W = 3$.}
 Since $\mathcal{U}(W, K)$ is a model for $\text{E}_{\Fin} W$ and one can calculate that it is $3$-dimensional, we conclude $\OFincd_\ZZ W \le 3$.  To see that $\OFincd_\ZZ W = 3$ we calculate $H^n_{\OFin}(W, \ZZ[-,W/1]_{\OFin})$ as in \cite[p.147]{LearyNucinkis-SomeGroupsOfTypeVF}, using Lemma \ref{lemma:bredon tech lem for cdQ neq OFcdQ} at the end of this section,
 \begin{align*}
 H^n_{\OFin}(W, \ZZ[-,W/1]_{\OFin}) &\cong H^3_W( \mathcal{U}(W, K), \mathcal{U}(W, K)^{\text{sing}}; \ZZ W) \\
  &\cong H^3_W( \mathcal{U}(W, K), \mathcal{U}(W, \partial K); \ZZ W).
 \end{align*}

 Recall that $\mathcal{U}(W, K) = W \times K / \sim$ where the identification is only on $W \times \partial K$, that $K$ is a fundamental domain for the $W$-action on $\mathcal{U}(W, K)$, and that $(K, \partial K) \simeq (\mathcal{C} \RR \PP^2, \RR \PP^2)$.  Here $\mathcal{C}X$ denotes the cone on a space $X$.  The action of $W$ on $C_*(\mathcal{U}(W, K),\mathcal{U}(W, \partial K))$ is free, so 
 \begin{align*}
   H^*_W(\mathcal{U}(W, K), \mathcal{U}(W, \partial K); \ZZ W) &\cong H^* (K, \partial K; \ZZ) \otimes_{\ZZ} \ZZ W \\
   &\cong H^*(\mathcal{C} \RR \PP^2 , \RR \PP^2;\ZZ) \otimes_\ZZ \ZZ W.
 \end{align*}
In particular, in dimension $3$, 
\[
H^3(\mathcal{C} \RR \PP^2 , \RR \PP^2;\ZZ) \otimes_\ZZ \ZZ W \cong H^2(\RR \PP^2 ;\ZZ) \otimes_\ZZ \ZZ W \cong \FF_2 W.
\]
We conclude $\OFincd_\ZZ W = 3$.
 
\emph{Claim: $\OFincd_{\FF_3}W = 2$.}
 $\mathcal{U}(W, \partial K)$ is the singular set of $\mathcal{U}(W, K)$, so in particular the fixed point sets of finite subgroups (except for the trivial subgroup) agree.  They are contractible and hence $\FF_3$-acyclic.  We claim $\mathcal{U}(W, \partial K)$ is also $\FF_3$-acyclic.  We use Lemma \ref{lemma:bredon UWK acyclic}---if $T \neq \emptyset$ then $(\partial K)_T = K_T$ which is contractible and hence $\FF_3$-acyclic and if $T = \emptyset$ then $(\partial K)_T = \partial K$ which is the barycentric subdivision of $L = \RR \PP^2$ and hence $\FF_3$-acyclic.  Taking the Bredon chain complex 
 \[P_* = C_*^{\OFin}(\mathcal{U}(W, \partial K)) \otimes_\ZZ \FF_3\]
 gives that $\OFincd_{\FF_3}W \le 2$.
 
 $W$ is a right-angled Coxeter group so every finite subgroup of $W$ has order a power of $2$, in particular $W$ has no $\FF_3$-torsion.  By Corollary \ref{cor:tfae cdRG less 1}, $\OFincd_{\FF_3} W = 1$ if and only if $\OFincd_\ZZ W = 1$ but we have already shown that $\OFincd_\ZZ W = 3$ proving that $\OFincd_{\FF_3} W \neq 1$ and in fact $\OFincd_{\FF_3} W = 2$.
\end{Example}

\begin{Example}[A group with {$\cd_\QQ G \neq \OFincd_\QQ G$}]
In \cite{LearyNucinkis-SomeGroupsOfTypeVF}, Leary and Nucinkis construct examples of groups with $\vcd_\ZZ G = nm$ and $\OFincd_\ZZ G = m(n+1)$  for various integers $n$ and $m$.  We show that these groups have $\OFincd_\QQ G = m(n+1)$ as well, so since $\cd_\QQ G \le \vcd_\ZZ G$ this provides examples of groups with $\cd_\QQ G \neq \OFincd_\QQ G$.  

We prove that groups $G$ satisfying the assumptions of \cite[Theorem 6]{LearyNucinkis-SomeGroupsOfTypeVF} satisfy $\OFincd_\QQ G \ge m(n+1)$ also, since combining this with the inequality $\OFincd_\QQ G \le \OFincd_\ZZ G$ gives $\OFincd_\QQ G = m(n+1)$ as required.  

Leary and Nucinkis show there exists a model $X$ for $\text{E}_{\Fin} G$ such that the cellular chain complex $C_*(X^{m(n+1)},(X^{m(n+1)})^{\text{sing}} )$ contains a copy of $\ZZ G$ in dimension $m(n+1)$ as a direct summand.  Here $X^i$ denotes the $i$ skeleton of some CW-complex $X$.  Using Lemma \ref{lemma:bredon tech lem for cdQ neq OFcdQ} below, 
\begin{align*}
 H^{m(n+1)}_{\OFin} (G, \QQ[ - G/1]_{\OFin}) &\cong H^{m(n+1)}_G(C_*(X,X^{\text{sing}} ); \QQ G)  \\
&\cong H^{m(n+1)}_G(C_*(X^{m(n+1)},(X^{m(n+1)})^{\text{sing}} ); \QQ G) \\
&\neq 0
\end{align*}
showing $\OFincd_\QQ G \ge m(n+1)$. 
\end{Example}

The examples constructed with the method above can never be of type $\OFinFP_\infty$ \cite[Question 2, p.154]{LearyNucinkis-SomeGroupsOfTypeVF}, so a natural question is:
\begin{Question}\label{question:bredon OFFPinfty with cdQQ neq OFcdQQ}
 Are there groups $G$ with $\cd_\QQ G \neq \OFincd_\QQ G$ and type $\OFinFP_\infty$?
\end{Question}

\begin{Lemma}\label{lemma:bredon tech lem for cdQ neq OFcdQ}
For any group $G$ and model $X$ for $\EFin G$
\[H^*_{\OFin}(G, R[-,G/1]_{\OFin}) \cong H^*_G(C_*(X, X^\text{sing}); RG) \]
where $C_*(X, X^\text{sing})$ denotes the cellular chain complex of the pair $(X, X^\text{sing})$.
\end{Lemma}
\begin{proof}
 Firstly, 
\[ H^* \left( \Hom_{\OFin} \left( C_*^{\OFin}(X^\text{sing}), R[-,G/1]_{\OFin} \right) \right) = 0\]
since the $G$-orbits of cells in $X^\text{sing}$ all give rise to contravariant modules of the form $R[-,G/H]_{\OFin}$ for $H \neq 1$, and by the Yoneda-type Lemma \ref{lemma:C yoneda-type},
\[\Hom_{\OFin}( R[-, G/H]_{\OFin}, R[-,G/1]_{\OFin}) \cong R[G/H, G/1]_{\OFin} = 0. \]
Using the long exact sequence in homology associated to the pair $(X, X^\text{sing})$,
\begin{align*}
H^*_{\OFin}(G, R[-, G/1]_{\OFin}) 
&\cong H^* \Hom_{\OFin}( C_*^{\OFin}(X), R[-,G/1]_{\OFin}) \\
&\cong H^* \Hom_{\OFin} (C_*^{\OFin}(X, X^\text{sing}), R[-,G/1]_{\OFin}  ). \tag{$\star$}
\end{align*}
Via the Yoneda-type Lemma \ref{lemma:C yoneda-type}, there is a chain of natural isomorphisms:
\begin{align*} 
\Hom_{\OFin} (C_*^{\OFin}(X &, X^\text{sing}), R[-, G/1]_{\OFin}) \\
&\cong \Hom_{\OFin} \left( \bigoplus_{\substack{\text{$G$-orbits of $i$-cells}\\ \text{with trivial isotropy}}} R[-, G/1]_{\OFin}, R[-, G/1]_{\OFin} \right) \\
&\cong \prod \Hom_{\OFin} \left( R[-, G/1]_{\OFin}, R[-, G/1]_{\OFin} \right) \\
&\cong \prod \Hom_{RG} (RG, RG) \\
&\cong \Hom_{RG} \left(\bigoplus RG, RG \right) \\
&\cong \Hom_{RG} ( C_*^{\OFin}(X, X^\text{sing}), RG ).
\end{align*}
Thus,
\[
H^*\Hom_{\OFin} (C_*^{\OFin}(X, X^\text{sing}) \cong H^* \Hom_{RG} ( C_*(X, X^\text{sing}), RG ),
\]
and combining this with the isomorphism $(\star)$ completes the proof.
\end{proof}

\section{Finitely generated projectives and duality}\label{section:bredon fg projectives and duality}

In this section we require $\mathcal{F} \subseteq \Fin$.  This section contains a number of technical results concerning dual $\OF$-modules, they are all analogs of results for modules over group rings that can be found in \cite{Bieri-HomDimOfDiscreteGroups}.  The results in this section are built on in Section \ref{section:wrong notion of duality} and utilised in Section \ref{section:Homology and Cohomology of Cohomological Mackey Functors}.

For $M$ a contravariant module, denote by $M^D$ the dual module \index{Dual module $M^D$}
\[M^D = \Hom_{\OF} \left( M(-), R[ -, ?]_{\OF} \right). \]
Similarly for $A$ a covariant module,
\[A^D = \Hom_{\OF} \left( A(-), R[?, -]_{\OF} \right). \]
This definition should be compared with that of the dual of an $RG$-module $M$, namely $M^D = \Hom_{RG}(M, RG)$ \cite[\S 3.1]{Bieri-HomDimOfDiscreteGroups}.

\begin{Example}\label{example:dual of constant functor}
 If $G$ is an infinite group and $\uR$ is the covariant constant functor on $R$ then $\uR^D = 0$, as
\begin{align*}
\uR^D &= \Hom_{\OF}( \uR(?), R[-, ?]_{\OF} ) \\
 &\cong \Hom_{\OF}( \Ind_G^{\OF} R(?), R[-,?]_{\OF} ) \\
 &\cong \Hom_{RG}( R, R[-,G/1]_{\OF} ),
\end{align*}
using Example \ref{example:cov E_1R is uR} and the adjointness of induction and restriction.  Finally, $\Hom_{RG}( R, R[-,G/1]_{\OF} )$ is the zero module since $G$ is infinite.
\end{Example}

\begin{Lemma}\label{lemma:double dual is nat iso}
 The dual functor takes projectives to projectives and the double-dual functor $-^{DD}: \{\text{$\OF$-modules}\} \to \{\text{$\OF$-modules}\}$ is a natural isomorphism when restricted to the subcategory of finitely generated projective $\OF$-modules.
\end{Lemma}
\begin{proof}
By the Yoneda-type Lemma \ref{lemma:C yoneda-type},
\[
 R[-,G/H]_{\OF}^D \cong \Hom_{\OF} (R[?,G/H]_{\OF}, R[?, -]_{\OF}) \cong R[G/H, -]_{\OF}.
\]
The proof for covariant frees is identical.

For any module $M$, there is a map $\zeta: M \longrightarrow M^{DD}$, given by $\zeta (m)(f) = f(m)$.  If $M = R[-,G/H]_{\OF}$ then applying the Yoneda-type lemma twice shows $M^{DD} = M$.  This generalises to projectives since the duality functor represents direct sums.

Naturality follows from naturality of the map $\zeta$.
\end{proof}

For $M$ and $N$ contravariant $\OF$-modules, we construct an $R$-module homomorphism
\[\nu : N \otimes_{\OF} M^D  \longrightarrow \Hom_{\OF} \left( M, N \right). \]
The main result of this section will be Lemma \ref{lemma:M fg proj then nu iso}, that $\nu$ is an isomorphism when $M$ is finitely generated projective and Proposition \ref{prop:G FPn then nat iso for bredon cohomology}, that $\nu$ induces an isomorphism
\[ N(?) \otimes_{\OF} H^i_{\OF}\left(G, R[-,?]_{\OF}\right) \cong H^i_{\OF} (G, N)\]
for all $i \le n$ when $G$ is $\OFFP_n$.  

Recall that elements of $N \otimes_{\OF} M^D$ are equivalence classes of finite sums of elements of the form
\[ n_H \otimes \varphi_H  \in \bigoplus_{G/H \in \OF}  N(G/H) \otimes_R \Hom_{\OF} \left(M, R[ -,G/H]_{\OF}\right). \]
For any $G/L \in \OF$ and $m \in M(G/L)$ we define 
\begin{align*}
\nu \left(n_H \otimes_R \varphi_H  \right)(G/L) : M(G/L) &\longrightarrow N(G/L) \\
 m &\longmapsto N\left( \varphi_H(G/L)(m) \right)(n_H).
\end{align*}

This makes sense because $\varphi_H(G/L)(m) \in R[G/L, G/H]_{\OF}$ and $N$ is a contravariant module so 
\[N\big(\varphi_H(G/L)(m)\big): N(G/H) \longrightarrow N(G/L).\]

We must check that $\nu(n_H \otimes_R \varphi_H)$ is a natural transformation, it's well defined including that it doesn't depend on the choice of equivalence class in $N({?}) \otimes_{\OF} \Hom_{\OF} \left(M({-}), R[{-},{?}]_{\OF}\right) $, and that it is an $R$-module homomorphism.

\noindent\emph{$\nu(n_H \otimes_R \varphi_H)$ is a natural transformation:} 

Let $\alpha: G/L_1 \longmapsto G/L_2$ be a $G$-map and $G/L_i \in \OF$, we must check the following diagram commutes.
\[
\xymatrix@C+150pt{
M(G/L_1) \ar^{\nu(n_H \otimes_R \varphi_H)(G/L_1) : m \longmapsto N\big(\varphi_H(G/L_1)(m)\big)(n_H)}[r] & N(G/L_1) \\
M(G/L_2) \ar^{M(\alpha)}[u] \ar^{\nu(n_H \otimes_R \varphi_H)(G/L_2) : m \longmapsto N\big(\varphi_H(G/L_2)(m)\big)(n_H)}[r] & N(G/L_2) \ar^{N(\alpha)}[u]
}
\]

\begin{align*}
 N(\alpha) \circ \big( \nu(n_H \otimes_R \varphi_H )(G/L_2) \big) (m) &= N(\alpha) \circ  N\left(\varphi_H(G/L_2)(m)\right)(n_H)  \\
 &= N\big( \varphi_H(G/L_2)(m) \circ \alpha \big)(n_H) \\
 &= N \big( \left( R[\alpha, G/H]_{\OF} \circ \varphi_H(G/L_2)\right) (m) \big) (n_H) \\
 &= N \big(  (\varphi_H(G/L_1) \circ M(\alpha))(m) \big) (n_H) \\
 &= \big( \nu(n_H \otimes_R \varphi_H)(G/L_2)\circ M(\alpha) \big) (m)
\end{align*}
Where the second equality is because $N$ is a contravariant functor, the third is because by definition $\varphi_H(G/L_2)(m) \circ \alpha = \left( R[\alpha, G/H]_{\OF} \circ \varphi_H(G/L_2)\right) (m) $, and the fourth is because $\varphi_H$ is itself a natural transformation and hence the following diagram commutes.
\begin{equation}
\xymatrix@C+25pt{
M(G/L_1) \ar^-{\varphi_H(G/L_1)}[r] & R[G/L_1, G/H]_{\OF} \\
M(G/L_2) \ar^-{\varphi_H(G/L_2)}[r] \ar^{M(\alpha)}[u] & R[G/L_2, G/H]_{\OF} \ar^{R[\alpha, G/H]_{\OF}}[u]
} \tag{$\dagger$}
\end{equation}

\noindent\emph{$\nu$ is well-defined:}
Firstly,
\[\nu(r n_H \otimes \varphi_H ) = \nu(n_H \otimes r \varphi_H) \]
this is because 
\begin{align*}
\nu \left(n_H\cdot r \otimes_R \varphi_H  \right)(G/L) (m) &= N\left( \varphi_H(G/L)(m) \right)(r n_H)  \\
&= r N\left( \varphi_H(G/L)(m) \right)(n_H)\\
&= N\left( r \varphi_H(G/L)(m) \right)(n_H)\\
&= \nu(n_H \otimes r \varphi_H).
\end{align*}
Secondly, we show $\nu$ doesn't depend on the choice of equivalence class in
$$N({?})  \otimes_{\OF} \Hom_{\OF} \left(M({-}), R[{-},{?}]_{\OF}\right). $$  
Choose $n_H \in N(G/H)$, $\varphi_M \in \Hom_{\OF} \left(M({-}), R[{-},G/M]_{\OF}\right)$, $\alpha: G/H \to G/M$ a $G$-map and $G/H, G/M \in \OF$, we must show that 
$$\nu \big( N(\alpha)( n_H) \otimes_R \varphi_M\big) = \nu\left( n_H \otimes_R  \Big( \Hom_{\OF}\left(M({-}), R[{-}, \alpha]_{\OF}\right) \Big) ( \varphi_M)\right). $$
Let $G/L \in \OF$, then 
\begin{align*}
\nu\big(N(\alpha)( n_H) &\otimes_R \varphi_M\big)(G/L)(m) = N\left(\varphi_H(G/L_1)(m) \right) \left( N(\alpha)( n_H) \right) \\ 
&= N\left(\alpha \circ \varphi_H(G/L)(m) \right) ( n_H) \\
&= N\left(R[G/L, \alpha]_{\OF} (\varphi_H (G/L_1)(m)) \right) ( n_H) \\
&= N\left(  \Hom_{\OF}\left(M({-}), R[{-}, \alpha]_{\OF}\right)(\varphi_H) (G/L_1)(m) \right) ( n_H) \\
&= \nu\big( n_H \otimes_R  \Hom_{\OF}\left(M({-}), R[{-}, \alpha]_{\OF}\right)(\varphi_M) \big) (G/L)(m).
\end{align*}

\noindent\emph{$\nu$ is a map of $R$-modules:}
It's clear that $\nu$ is additive, and 
\[\nu(r n_H \otimes \varphi_H ) = r \nu(n_H \otimes \varphi_H) \]
since $N$ being a module over $R$ implies that $N(\varphi_H(G/L)(m))$ is an $R$-module homomorphism.

\begin{Lemma}\label{lemma:nu is natural}
 $\nu$ is natural in $N$ in $M$.
\end{Lemma}
\begin{proof}
We only prove naturality in $N$, the proof for $M$ is similar.  Let $F$ be a morphism of contravariant modules $N \to N^\prime$, we must show that the following diagram of $R$-modules commutes.
\[
\xymatrix@C+25pt{
N({?}) \otimes_{\OF} \Hom_{\OF}(M({-}),R[{-},{?}]_{\OF}) \ar^-{\nu_N}[r] \ar^{F({?}) \otimes_{\OF} \Hom_{\OF}(M({-}),R[{-},{?}]_{\OF})}[d] & \Hom_{\OF}(M({-}),N({-})) \ar^{\Hom_{\OF}(M({-}), F({-}))}[d]\\
N^\prime({?}) \otimes_{\OF} \Hom_{\OF}(M({-}),R[{-},{?}]_{\OF}) \ar^-{\nu_{N^\prime}}[r] & \Hom_{\OF}(M({-}),N^\prime({-}))
}
\]
Let $n_H \otimes \varphi_H \in N(G/H) \otimes_{\OF} \Hom_{\OF}(M({-}), R[{-},G/H]_{\OF})$ then moving about the top right of the diagram yields
\[
\big( \Hom_{\OF}(M({-}),F({-})) \circ \nu_N ( n_H \otimes \varphi_H ) \big) (G/L)(m) 
\]
\[
= F(G/L) \circ N(\varphi_H(G/L)(m))(n_H).
\]
Moving around the bottom left yields
\[ \big( \nu_{N^\prime} \circ F({?}) \otimes \Hom_{\OF}(M({-}),R[{-},{?}]_{\OF}) (n_H \otimes \varphi_H )\big) (G/L)(m) \]
\[ = \nu_{N^\prime}\big(F(G/H)(n_H) \otimes \varphi_H)\big) (G/L)(m) \]
\[ = N^\prime \big( \varphi_H(G/L)(m)\big) \big( F(G/H)(n_H) \big).\]
That these two are equivalent is because $F$ is a natural transformation, so the diagram below commutes.
\[
\xymatrix{
N(G/L) \ar^{F(G/L)}[r] & N^\prime(G/L) \\
N(G/H) \ar^{F(G/H)}[r] \ar^{N(\varphi_H(G/L)(m))}[u] & N^\prime(G/H) \ar_{N^\prime(\varphi_H(G/L)(m))}[u]
}
\]
\end{proof}

The next lemma is an $\OF$ module version of \cite[Proposition 3.1]{Bieri-HomDimOfDiscreteGroups}.

\begin{Lemma}\label{lemma:M fg proj then nu iso}
 If $M$ is finitely generated projective then $\nu$ is an isomorphism.
\end{Lemma}
\begin{proof}
 Consider first the case $M = R[-,G/H]_{\OF}$, then the map $\nu$ becomes
 \[
\nu : N({?}) \otimes_{\OF} \Hom \left(R[{-},G/H]_{\OF}, R[{-},{?}]_{\OF}\right)  \longrightarrow \Hom_{\OF} \left( R[{-},G/H]_{\OF}, N({-}) \right).  
 \]
 But, using Lemmas \ref{lemma:C yoneda-type} and \ref{lemma:C yoneda type isos in tensor product}, the left hand side collapses to
\begin{align*}
N({?}) \otimes_{\OF} \Hom \left(R[{-},G/H]_{\OF}, R[{-},{?}]_{\OF}\right) &\cong N({?}) \otimes_R R[G/H, {?}]_{\OF} \\
&\cong N(G/H). \tag{$\star$}
\end{align*}
Under these isomorphisms $n_H \in N(G/H)$ maps to 
\[n_H \otimes \id_H \in N({?}) \otimes_R R[G/H,{?}]_{\OF}\]
and then to $n_H \otimes \varphi$ where $\varphi$ is the unique natural transformation $\varphi$ such that $\varphi(G/H)(\id_H) = \id_H$.

The right hand side collapses to
\begin{align*}
\Hom_{\OF} \left( R[{-},G/H]_{\OF}, N({-}) \right) \cong N(G/H) \tag{$\dagger$} 
\end{align*}
again by the Yoneda-type Lemma \ref{lemma:C yoneda-type}, where $n_H$ maps to the unique natural transformation $\psi$ with $\psi(G/H)(\id) = n_H$.
\[
\nu(n_H \otimes \varphi)(G/H)(\id_H) = N(\varphi(G/H)(\id_H))(n_H) = N(\id_H)(n_H) = n_H 
\]
Precomposing $\nu$ with the isomorphism from $(\star)$ and postcomposing with the isomorphism from $(\dagger)$ gives the identity map $N(G/H) \to N(G/H)$ and hence $\nu$ is an isomorphism.

The case for finitely generated free modules follows as all the necessary functors commute with finite direct sums, and for projectives from naturality of $\nu$ proved in Lemma \ref{lemma:nu is natural}.
\end{proof}

The following result is an analog of \cite[5.2(a,c)]{Bieri-HomDimOfDiscreteGroups}.  
\begin{Lemma}\label{lemma:M fg and N proj then nu iso}~
\begin{enumerate}
 \item If $M$ is finitely presented and $N$ is flat then $\nu$ is an isomorphism.
 \item If $M$ is finitely generated and $N$ is projective then $\nu$ is an isomorphism.
\end{enumerate}
\end{Lemma}
\begin{proof}
\begin{enumerate}
 \item 
If 
\[
F_1 \longrightarrow F_0 \longrightarrow M \longrightarrow 0
\]
is an exact sequence with $F_i$ finitely generated free for $i = 0,1$ then by the naturality of $\nu$ and flatness of $N$ we have the following commutative diagram with exact rows (for brevity we write $\Hom$ for $\Hom_{\OF}$ and $\otimes$ for $\otimes_{\OF}$).
\[ \xymatrix@C-9pt{
0 \ar[r] & N(?) \otimes M^D(?) \ar[r] \ar[d] & N(?) \otimes F_0^D(?) \ar[r] \ar[d]  & N(?) \otimes  F_1^D(?) \ar[d]  \\
0 \ar[r] & \Hom(M, N) \ar[r] & \Hom(F_0, N) \ar[r] & \Hom(F_1, N) 
} \]
The right hand and middle vertical maps are isomorphisms by Lemma \ref{lemma:M fg proj then nu iso}, the result follows from the $5$-Lemma.

\item 
If $F(?)$ is free then by Lemma \ref{lemma:C yoneda type isos in tensor product} there is an isomorphism
\[F(?)\otimes_{\OF} \Hom(M(-), R[-,?]_{\OF}) \cong \Hom(M, F). \]
Checking the definition of this isomorphism shows it's induced by $\nu$.  If $N(?)$ is projective and $i:N(?) \longhookrightarrow F(?)$ is a split injection then by naturality of $\nu$, the following diagram commutes:
\[\xymatrix{
N(?)\otimes_{\OF} \Hom(M(-), R[-,?]_{\OF}) \ar[d]\ar[r] & \Hom(M, N) \ar[d]\\
F(?)\otimes_{\OF} \Hom(M(-), R[-,?]_{\OF}) \ar^-{\cong}[r] & \Hom(M, F)
}\]
Since $i$ is a split injection, the left hand map is an injection and the top map must be an injection.  Consider the commutative diagram in the proof of part 1, only $F_0$ is known to be projective so the middle vertical map is an isomorphism.  Since $N$ is projective the left and right hand vertical maps are monomorphisms and the Four Lemma completes the proof, implying that the left hand vertical map is an isomorphism.
\end{enumerate}
\end{proof}

\begin{Lemma}\label{lemma:alpha is welldef and natural}
If $P_*$ is any chain complex of contravariant $\OF$-modules and $N$ is any contravariant $\OF$-module, the following morphism is both well defined and natural in $P_*$ and $N$:
\begin{align*}
\xi^i : N(?) \otimes_{\OF} H^iP_*(?)^D &\rightarrow 
H^i \big( N(?) \otimes_{\OF} P_*(?)^D \big)\\
\xi^i : N(?) \otimes_{\OF} H^i(\Hom(P_*(-), R[-,?]) &\rightarrow 
H^i \big( N(?) \otimes_{\OF} \Hom(P_*(-), R[-, ?]) \big) \\
n_H \otimes [\varphi_H] &\mapsto [n_H \otimes \varphi_H],
\end{align*}
where $H^iP_*(?)^D: G/H \mapsto H^iP_*(G/H)^D$.
\end{Lemma}
\begin{proof}
 If $\varphi_H$ is a cocycle, $n_H \otimes \varphi_H $ is also a cocycle and similarly if $\varphi_H$ is a coboundary then $n_H \otimes \varphi_H$ is a coboundary.  

 If $\alpha:G/L \to G/H$ is a $G$-map then by definition $\alpha_*[\varphi_H ] = [\alpha_*\varphi_H]$ and
\begin{align*}
\xi^i\big(  \alpha^*n_H \otimes [\varphi_H] - n_H \otimes \alpha_*[\varphi_H ] \big) &= \xi^i\big(  \alpha^*n_H \otimes [\varphi_H] - n_H \otimes [\alpha_*\varphi_H ] \big) \\
&=[\alpha^*n_H \otimes \varphi_H - n_H \otimes \alpha_*\varphi_H]  \\
&= 0.
\end{align*}
Finally naturality follows because the the functors $H^i$ and $\Hom_{\OF}(-,?)$ are natural, and so is the process of taking tensor products.
\end{proof}

Since $\nu$ is natural (Lemma \ref{lemma:nu is natural}), if $P_*$ is a projective resolution of $\uR$ by contravariant modules then $\nu$ induces chain homomorphisms 
\[
N(?) \otimes_{\OF} P_*(?)^D \longrightarrow \Hom_{\OF}(P_*,N)
\]
which in turn induce maps on cohomology
\[
H^i \big( N(?) \otimes_{\OF} P_*(?)^D \big) \longrightarrow H^i_{\OF} (G, N). 
\]
Precomposing this with $\xi^i$ gives a map 
\[ 
\nu^i: N(?) \otimes_{\OF}  H^i_{\OF}(G, R[-,?]_{\OF}) \longrightarrow H^i_{\OF}(G, N).  
\]

\begin{Prop}\label{prop:G FPn and N proj then nu iso}
 If $G$ is $\OFFP_n$ over $R$ and $N$ is projective then $\nu^i$ is an isomorphism for all $i \le n$.
\end{Prop}
\begin{proof}
 Choose a projective resolution $P_* \longtwoheadrightarrow \uR$, finitely generated up to dimension $n$ and write $K_i$ for the $i^\text{th}$ syzygy of $P_*$.  Since $N$ is projective it is also flat and we have the following commutative diagram with exact rows, where we omit the $_{\OF}$ on $\otimes$, $\Hom$, and $H^i$.
 \[
 \xymatrix@C-13pt{
 N(?) \otimes P_{i-1}^D(?) \ar[r] \ar^{\nu}[d] &
 N(?) \otimes K_{i-1}^D(?) \ar[r] \ar^{\nu}[d] &
 N(?) \otimes H^i(G, R[-, ?]_{\OF}) \ar[r] \ar^{\nu^i}[d] & 0 \\
 \Hom(P_{i-1}, N) \ar[r] &
 \Hom(K_{i-1}, N) \ar[r] &
 H^i(G, N) \ar[r] & 0
 }
 \]
 Since $G$ is $\OFFP_n$, $K_{i-1}$ and $P_{i-1}$ are finitely generated, Lemma \ref{lemma:M fg and N proj then nu iso}(2) implies the middle and left hand vertical maps are isomorphisms.  The $5$-Lemma completes the proof.
\end{proof}

The following result is an analog of \cite[9.1]{Bieri-HomDimOfDiscreteGroups}.
\begin{Prop}\label{prop:G FPn then nat iso for bredon cohomology}If $G$ is $\OFFP$ over $R$, with $\OFcd_R G = n$, and $N$ is any contravariant module then there is a natural isomorphism
\[\nu^n :  N(?) \otimes_{\OF} H^n_{\OF}\left(G, R[-,?]_{\OF}\right) \cong H^n_{\OF} (G, N).\]
\end{Prop}
\begin{proof}
 Let 
\[
0 \longrightarrow K \longrightarrow F \longrightarrow N \longrightarrow  0 
\]
be a short exact sequence of contravariant $\OF$-modules over $R$ with $F$ free.  By the naturality of $\nu^n$ we have the following commutative diagram with exact rows, we omit the ${\OF}$ decorations on $\otimes$, $H^*$, and $R[-,?]$ for brevity.
\[ \xymatrix@C-18pt{
K(?)\otimes H^n(G, R[-,?]) \ar[r] \ar[d] & F(?) \otimes H^n(G, R[-,?]) \ar[r]\ar[d]  & N(?) \otimes H^n(G, R[-,?]) \ar[r] \ar[d]  & 0 \\
H^n(G, K) \ar[r] & H^n(G, F) \ar[r] &  H^n(G, N) \ar[r]  & 0
}\]
The middle vertical map is an isomorphism by Proposition \ref{prop:G FPn and N proj then nu iso}, thus by the Four Lemma, the right hand vertical map is an epimorphism.  Since there are no restrictions on $N$, we conclude that the left hand vertical map is an epimorphism and by the 5-Lemma that the right hand map is an isomorphism.
\end{proof}
\chapter{Mackey and cohomological Mackey functors}\label{chapter:M}

This chapter contains material that has appeared in:
\begin{itemize}
  \item Finiteness conditions for Mackey and cohomological Mackey functors (J.~Algebra \textbf{411} (2014), no.~0, 225--258) \cite{Me-CohomologicalMackey}
\end{itemize}

Throughout this section we will work over an arbitrary subfamily $\mathcal{F}$ of $\Fin$, closed under conjugation and taking subgroups.  One could also work over larger families of subgroups such as $\VCyc$ \cite[p.101]{Degrijse-Thesis}, however this necessitates a change in the construction of Mackey and cohomological Mackey functors and we shall not consider it.

In Section \ref{section:mackey introduction} we give an overview of Mackey functors and cohomological Mackey functors including the description due to Yoshida of cohomological Mackey functors as modules over the category $\HeckeF$ \cite{Yoshida-GFunctors2}.

Section \ref{section:mackey FPn for mackey} contains a complete description of the condition $\MFFP_n$, the Mackey functor analogue of the $\OFFP_n$ conditions.

\theoremstyle{plain}\newtheorem*{CustomThmJ}{Corollary \ref{cor:OFFPn iff MFFPn}}
\begin{CustomThmJ}
Over any ring $R$, a group is $\MFFP_n$ if and only if it is $\OFFP_n$.
\end{CustomThmJ}

The main result of Section \ref{section:Homology and Cohomology of Cohomological Mackey Functors} is that the Bredon cohomology with coefficients in a cohomological Mackey functor may be calculated with a projective resolution of cohomological Mackey functors.  We show in Proposition \ref{prop:HF sigma ind preserves proj res of R} that a projective resolution of $\uR$ by Bredon modules can be induced to a projective resolution of the fixed point functor $R^-$ by cohomological Mackey functors, this is an analogue of \cite[Theorem 3.8]{MartinezPerezNucinkis-MackeyFunctorsForInfiniteGroups}---that one can induce a projective resolution of $\uR$ by Bredon modules to a projective resolution of the Burnside functor $B^G$ by Mackey functors.

Building on this, in Section \ref{section:HF FPn} we study the $\HFFP_n$ conditions, the cohomological Mackey functor analogue of the $\OFFP_n$ conditions, relating them to the $\FFP_n$ conditions defined in Section \ref{section:I Fcohomology}.
\theoremstyle{plain}\newtheorem*{CustomThmK}{Theorem \ref{theorem:HFFPn iff FFPn}}
\begin{CustomThmK}
If $R$ is a commutative Noetherian ring, a group is $\HFFP_n$ if and only if it is $\FFP_n$. 
\end{CustomThmK}

In Section \ref{section:HF cohomological dimension} the main result is the following.
\theoremstyle{plain}\newtheorem*{CustomThmL}{Theorem \ref{theorem:Fcd=HFcd}}
\begin{CustomThmL}
 $\HFcd G = \Fcd G$ for all groups $G$.
\end{CustomThmL}

In Section \ref{section:family of p subgroups} we prove that, depending on the coefficient ring, $\HFcd G$ may be calculated using a proper subfamily of $\mathcal{F}$.  When working over $\ZZ$ we need consider only the family $\mathcal{P}$ of subgroups in $\mathcal{F}$ with prime power order, and over either the finite field $\FF_p$ or over $\ZZ_{(p)}$ (the integers localised at $p$), we need consider only the family $\mathcal{P}$ of subgroups of $\mathcal{F}$ with order a power of $p$.  This is similar to a result of Leary and Nucinkis for $\mathcal{F}$-cohomology \cite[\S 4]{LearyNucinkis-GroupsActingPrimePowerOrder}.

\theoremstyle{plain}\newtheorem*{CustomThmM}{Theorem \ref{theorem:HF HeckeFcd = HeckePcd and HeckeFFPn = HeckePFPn}}
\begin{CustomThmM}
For all $n \in \NN \cup \{ \infty \}$, the conditions $\HFcd G = n$ and $\HPcd G = n$ are equivalent, as are the conditions $\HFFP_n$ and $\HPFP_n$. 
\end{CustomThmM}

Over the finite field $\FF_p$ we can be even more precise.

\theoremstyle{plain}\newtheorem*{CustomThmN}{Corollary \ref{cor:mackey HFFPn complete description}}
\begin{CustomThmN}
$G$ is $\HFFP_n$ over $\FF_p$ if and only if $\mathcal{P}$ has finitely many conjugacy classes and $WH$ is $\FP_n$ over $\FF_p$ for all $H \in \mathcal{P}$.
\end{CustomThmN}

\section{Introduction}\label{section:mackey introduction}

\subsection{Mackey functors}\label{section:prelim mackey}

There are many constructions of Mackey functors, we use the construction coming from modules over a category, an approach due to Linder \cite{Lindner-RemarkOnMackeyFunctors}.\index{Mackey category $\MF$}  Another construction is mentioned in Remark \ref{remark:mackey green description}.  We begin by building a small category $\MF$ then Mackey functors will be contravariant $\MF$-modules.  As in $\OF$, the objects of $\MF$ are the transitive $G$-sets with stabilisers in $\mathcal{F}$, the morphism set however is much larger.  A \emph{basic morphism} from $G/H$ to $G/K$, where $H$ and $K$ are in $\mathcal{F}$, is an equivalence class of diagrams of the form 
\[ G/H \stackrel{\alpha}{\longleftarrow} G/L \stackrel{\beta}{\longrightarrow} G/K \]
where the maps are $G$-maps, and $L \in \mathcal{F}$.  This basic morphism is equivalent to 
\[ G/H \stackrel{\alpha^\prime}{\longleftarrow} G/L^\prime \stackrel{\beta^\prime}{\longrightarrow} G/K \]
if there is a bijective $G$-map $\sigma:G/L \to G/L^\prime $, fitting into the commutative diagram below:
\[\xymatrix@-15pt{
& G/L \ar_{\alpha}[ld] \ar^{\beta}[rd] \ar_\cong^\sigma[dd]  & \\ 
G/S & & G/K \\
& G/L^\prime \ar^{\alpha^\prime}[ul] \ar_{\beta^\prime}[ur] & 
}\]
Form the free abelian monoid on these basic morphisms, and complete this free abelian monoid to a group, denoted $[G/H, G/K]_{\MF}$.  This is the set of morphisms in $\MF$ from $G/H$ to $G/K$.  

\begin{Remark}
When building the Mackey category, we could instead have started with equivalence classes of diagrams 
\[ G/H \leftarrow \Delta \rightarrow G/K \]
where $\Delta$ is any finitely generated $G$-set with stabilisers in $\mathcal{F}$ and the maps are $G$-maps.  This can be shown to be the free abelian monoid on the basic morphisms \cite[Proposition 2.2]{ThevenazWebb-StructureOfMackeyFunctors}.  Because of this alternative construction, we will pass freely between writing 
\[ \left(  G/H \leftarrow G/L \rightarrow G/K \right) + \left(  G/H \leftarrow G/L^\prime \rightarrow G/K \right) \]
and 
\[ \left(  G/H \leftarrow G/L \coprod G/L^\prime \rightarrow G/K \right ). \]
\end{Remark}

To complete the description of $\MF$, we must describe composition of morphisms.  It's sufficient to describe composition of basic morphisms, and then use distributivity to extend this to all morphisms.  If
\[ G/H \leftarrow G/L \rightarrow G/K \]
and
\[ G/K \leftarrow G/S \rightarrow G/Q \]
are two basic morphisms then their composition is the pullback diagram below in the category of $G$-sets.
\[
\xymatrix@-15pt{
& & \Delta \ar[dl] \ar[dr] & & \\
& G/L\ar[ld] \ar[rd] & & G/S \ar[ld] \ar[rd] & \\
G/H & & G/K & & G/Q 
}
\]

\begin{Lemma}[Composition of morphisms in $\MF$]\label{lemma:mackey morphism pullback}\cite[\S 3]{MartinezPerezNucinkis-MackeyFunctorsForInfiniteGroups}
The diagram below is a pullback in the category of $G$-sets.
\begin{equation*} \sum_{x \in L^g \backslash K / S^{g^\prime} }\left( \vcenter{\xymatrix@-15pt{
 & G/ \big( L^g \cap S^{g^\prime x^{-1}} \big) \ar^{\alpha_{g^{-1}}}[ld] \ar_{\alpha_{x (g^\prime)^{-1}}}[rd] & \\
G/L \ar_{\alpha_g}[rd] & & G/S \ar^{\alpha_{g^\prime}}[ld] \\
& G/K & 
}} \right) \end{equation*}
Notice that the subgroup $ L^g \cap S^{g^\prime x^{-1}}$ is both a subgroup of $K$ via the maps on the left and subconjugated to $K$ via the map $\alpha_x$, which is the composition of the maps on the right.  
\end{Lemma}
If $H$ is a subgroup of $G$ the notation $H^g$ means the conjugate $g^{-1}Hg$.

\begin{Lemma}[Standard form for morphisms in $\MF$]\label{lemma:mackey morphism standard form}\cite[Lemma 2.1]{ThevenazWebb-StructureOfMackeyFunctors}
Any basic morphism is equivalent to one in the standard form:
\[\xymatrix@-15pt{
& G/L \ar^{\alpha_g}[rd] \ar_{\id}[ld] & \\ 
G/K & & G/S
}\]
\end{Lemma}
Recall that two such basic morphisms are equivalent if there is a commutative diagram of the form:

\[\xymatrix@-15pt{
& G/L \ar^{\alpha_g}[rd] \ar_{\id}[ld] \ar^{\alpha_x}_{\cong}[dd] & \\ 
G/K & & G/S \\
& G/L^x \ar_{\alpha_{g^\prime}}[ru] \ar^{\id}[lu] & 
}\]
The commutativity of the left hand triangle ensures that $x \in K$, and that of the right hand diagram gives $\alpha_g = \alpha_{g^\prime} \circ \alpha_x$, or more concisely $gS = xg^\prime S$. This means $KgS = Kg^\prime S$ and $x = gS(g^{\prime})^{-1} \cap K = gSg^{-1} \cap K$.  Thus a basic morphism is determined by both an element of $K \backslash G / S$ and a subgroup $L$, subconjugate to $K$, unique up to conjugation by an element $x \in gSg^{-1} \cap K$.  In summary,
\begin{equation}\label{eq:mackey morphisms}
[ G/K, G/S ]_{\MF} = \bigoplus_{g \in K \backslash G/S } \bigoplus_{\substack{L \le gSg^{-1} \cap K \\ \text{Up to $gSg^{-1} \cap K$-conjugacy}}} \ZZ_{L, g},
\end{equation}
where $\ZZ_{L, g} \cong \ZZ$ for all $L$ and $g$.

\begin{Example}\label{example:mackey morphisms for S=1}
 If $S = 1$ then \eqref{eq:mackey morphisms} becomes 
\[
[ G/K, G/1 ]_{\MF} = \bigoplus_{g \in K \backslash G } \ZZ_{g} \cong \ZZ[K\backslash G].
\]
\end{Example}

\begin{Remark}\label{remark:mackey MF has A not EI}
 The category $\MF$ has property (A) by construction, but it does not have property (EI).  For example, given any non-trivial $H \in \mathcal{F}$, the endomorphism 
\[e = \left( G/H \stackrel{\alpha_1}{\longleftarrow} G/1 \stackrel{\alpha_1}{\longrightarrow} G/H \right)\]
is not an isomorphism.  If 
\[ m = \left( G/H \stackrel{\alpha_1}{\longleftarrow} G/K \stackrel{\alpha_g}{\longrightarrow} G/H \right) \]
is some other basic morphism then their composition is
\[ 
m \circ e = \sum_{x \in H/K} \left( G/H \stackrel{\alpha_1}{\longleftarrow} G/1 \stackrel{\alpha_{xg}}{\longrightarrow} G/H\right).
\]
So it's clear that $e$ cannot be a sum of automorphisms of $G/H$.
\end{Remark}

Following \cite{MartinezPerezNucinkis-MackeyFunctorsForInfiniteGroups}, we will mostly consider contravariant Mackey functors.  From here on, whenever we write $\MF$-module, we mean contravariant $\MF$-module.  

\begin{Remark}[Green's alternative description of Mackey functors]\label{remark:mackey green description}
There is an alternative description of Mackey functors, due to Green \cite{Green-AxiomaticRepresentationTheoryForFiniteGroups}, which we include here in full because when we later study cohomological Mackey functors we will need some of the language.

Green defined a Mackey functor $M$ as a mapping,
\[ M : \{ G/H \: : \: \text{$H \in \mathcal{F}$}\}  \longrightarrow \RMod \]
with morphisms for any finite subgroups $K \le H$ in $\mathcal{F}$,
\begin{align*}
 M(I_K^H) &: M(G/K) \longrightarrow M(G/H) \\
 M(R_K^H) &: M(G/H) \longrightarrow M(G/K) \\
 M(c_g) &: M(G/H) \longrightarrow M(G/H^{g^{-1}}) 
\end{align*}
 called \emph{induction}, \emph{restriction} and \emph{conjugation} respectively.  Induction is sometimes also called transfer.  In the literature, $M(I_K^H)$, $M(R_K^H)$ and $M(c_g)$ are often written as just $I_K^H $, $R_K^H$ and $c_g$, omitting the $M$ entirely.  We choose to use different notation so that we can identify $I_K^H $, $R_K^H$ and $c_g$ with specific morphisms in $\MF$ (see the end of this remark).

This mapping $M$ must satisfy the following axioms,
\begin{enumerate}
 \item[(0)] $M(I_H^H)$, $M(R^H_H)$ and $M(c_h)$ are the identity morphism for all $h \in H$.
 \item[(1)] $M(R_J^K) \circ M(R^H_K) = M(R^H_J)$, where $J \le K \le H$ and $J,K,H \in \mathcal{F}$.
 \item[(2)] $M(I^H_K) \circ M(I^K_J) = M(I^H_J)$, where $J \le K \le H$ and $J,K,H \in \mathcal{F}$.
 \item[(3)] $M(c_g) \circ M(c_h) = M(c_{gh})$ for all $g, h \in G$.
 \item[(4)] $M(R^{H^{g^{-1}}}_{K^{g^{-1}}}) \circ M(c_g )= M(c_g) \circ M(R^H_K)$, where $K \le H$ and $K,H \in \mathcal{F}$ and $g \in G$.
 \item[(5)] $M(I^{H^{g^{-1}}}_{K^{g^{-1}}}) \circ M(c_g )= M(c_g) \circ M(I^H_K)$, where $K \le H$ and $K,H \in \mathcal{F}$ and $g \in G$.
 \item[(6)] $M(R^H_J) \circ M(I^H_K) = \sum_{x \in J \backslash H / K} M(I^J_{J \cap K^{x^{-1}}}) \circ M(c_x) \circ M( R^K_{J^x \cap K}) $, where $J, K \le H$ and $J,K,H \in \mathcal{F}$.
\end{enumerate}
Axiom (6) is often called the Mackey axiom.  Converting between this description and our previous description is done by rewriting induction, restriction and conjugation in terms of morphisms of $\MF$. 
\begin{align*} 
M(I_K^H) &\longleftrightarrow M\big( G/H \stackrel{\alpha_1}{\longleftarrow} G/K \stackrel{\id}{\longrightarrow} G/K \big) \\
M(R_K^H) &\longleftrightarrow M\big( G/K \stackrel{\id}{\longleftarrow} G/K \stackrel{\alpha_1}{\longrightarrow} G/H \big) \\ 
M(c_g) &\longleftrightarrow   M\big( G/H^{g^{-1}} \stackrel{\id}\longleftarrow G/H^{g^{-1}} \stackrel{\alpha_g}{\longrightarrow} G/H \big) 
\end{align*}
Because of the above, we make the following definitions, \index{I, Induction@$I_K^H$, Induction morphism in $\MF$}\index{R, Restriction@$R_K^H$, Restriction morphism in $\MF$}\index{cg, conjugation@$c_g$, Conjugation morphism in $\MF$}
\begin{align*}
 I_K^H &= \big( G/H \stackrel{\alpha_1}{\longleftarrow} G/K \stackrel{\id}{\longrightarrow} G/K \big) \\
R_K^H &= \big( G/K \stackrel{\id}{\longleftarrow} G/K \stackrel{\alpha_1}{\longrightarrow} G/H \big) \\ 
c_g &=   \big( G/H^{g^{-1}} \stackrel{\id}\longleftarrow G/H^{g^{-1}} \stackrel{\alpha_g}{\longrightarrow} G/H \big).
\end{align*}
It is possible to write any morphism in $\MF$ as a composition of the three types of morphisms above.

One can check that Green's axioms all follow from the description of the composition of morphisms in $\MF$ as pullbacks (Lemma \ref{lemma:mackey morphism pullback}), and vice versa.  Complete proofs of the equivalence of this definition with our previous one can be found in \cite[\S 2]{ThevenazWebb-StructureOfMackeyFunctors}.
\end{Remark}

\subsubsection{Free modules}
In this section we describe the structure of $\End(G/H)$ and study free $\MF$-modules.

\begin{Remark}[Structure of $\End(G/H)$]\label{remark:mackey structure of Aut}
As mentioned in Remark \ref{remark:mackey MF has A not EI}, $\MF$ doesn't have property (EI).  Consider the endomorphisms of an object given by the diagrams of the form
\[
a_g = \big( G/H \stackrel{\alpha_1}{\longleftarrow} G/H \stackrel{\alpha_g}{\longrightarrow} G/H \big).
\]
Every $g \in WH$ uniquely determines a $G$-map $\alpha_g :G/H \to G/H$ and every $G$-map comes from such a $g$.  Finally, since $a_g \circ a_h = a_{hg}$, we determine that such endomorphisms give a copy of $\ZZ[WH]^\text{op}$ inside $\End(G/H)$.  This is similar to the situation over the orbit category, where $\End_{\OF}(G/H) \cong \ZZ[WH^\text{op}]$.  Thus, as with $\OF$-modules, if $M$ is a Mackey functor, then $M(G/H)$ is a right $R[WH^{\text{op}}]$ module, equivalently a left $R[WH]$-module.

A basic morphism in $\End(G/H)$ is determined by a morphism in standard form
\[ e_{L, g} = \big( G/H \stackrel{\alpha_1}{\longleftarrow} G/L \stackrel{\alpha_g}{\longrightarrow} G/H \big)  \]
where $L$ is some subgroup of $G$.  As such we can filter $\End(G/H)$ via the poset $\mathcal{F}/G$ of conjugacy classes of subgroups in $\mathcal{F}$.  If $L$ is a finite subgroup of $G$ then we write $ \End(G/H)_L $ for the basic morphisms $e_{L, g}$ for all $g \in G$.  Note that in particular, $\End(G/H)_H \cong R[WH]$ by the paragraph above.  Addition gives $\End(G/H)_L $ an abelian group structure.  Composing two elements of $\End(G/H)_L $ doesn't necessarily give an element of $\End(G/H)_L $, but pre-composing an element of $\End(G/H)_L $ by some $a_w$ does, since 
\[ e_{L,g} \circ a_{w}  \cong e_{L, wg}. \]
Thus $ R\End(G/H)_L $ is a left $R[WH]$-module.  In summary, there is an $R[WH]$-module isomorphism
\[ R\End(G/H) \cong \bigoplus_{L \in \mathcal{F}/G} R\End(G/H)_L \]
where $R\End(G/H)_H \cong R[WH]$.
\end{Remark}

\begin{Example}
 Using \eqref{eq:mackey morphisms}, 
\[
R\End(G/H)_1 \cong \bigoplus_{H \backslash G / H} R_{1, g},
\]
with left action of $w \in W_G H$ taking $g \mapsto wg$.  In other words, $R\End(G/H)_1 \cong R[H \backslash G / H]$ with the canonical action of $W_G H$.  This is not in general finitely generated---take for example $G = D_\infty$, the infinite dihedral group generated by the involutions $a$ and $b$, and $H = \langle a \rangle$.  Then $W_G H$ is the trivial group but $H \backslash G / H$ is an infinite set so $R[H \backslash G / H]$ is not a finitely generated $R$-module.
\end{Example}

\begin{Lemma}\label{lemma:restriction of free mackey at 1 is FPinfty}
 As a left $R[W_GS]$-module, $R[ G/S, G/K ]_{\MF}$ is an $R[W_GS]$-permutation module with finite stabilisers.  In addition, $R[ G/1, G/K ]_{\MF}$ is $\FP_\infty$ over $RG$.
\end{Lemma}
\begin{proof}
The left action of $w \in W_GS$ on $[ G/S, G/K ]_{\MF} $ is the action given by pre-composing any basic morphism $G/S \stackrel{\id}{\leftarrow} G/L \stackrel{\alpha_g}{\rightarrow} G/K$ with the morphism $G/S \stackrel{\id}{\leftarrow} G/S \stackrel{\alpha_w}{\rightarrow} G/S$ to yield the morphism 
\[ 
G/S \stackrel{\alpha_1}{\leftarrow} G/L \stackrel{\alpha_{wg}}{\rightarrow} G/K.
\]
To show this we calculate the pullback:
\[ \xymatrix@-15pt{
& & G/L \ar_{\id}[dl] \ar^{\alpha_w}[rd] & & \\
& G/S \ar_{\alpha_w}[dr] \ar_{\id}[dl]& & G/L \ar^{\id}[dl] \ar^{\alpha_g}[dr] &  \\
G/S & & G/S & & G/K
} \]
Under the identification \eqref{eq:mackey morphisms}, $w$ maps $R_{L, g}$ onto $R_{L, wg}$, so the stabiliser of this action is the stabiliser of the action of $R[W_GS
]$ on $R[S \backslash G / K]$, which is finite.  In particular $R[ G/S, G/K ]_{\MF} $ is an $R[W_GS]$-permutation module with finite stabilisers.  If $S = 1$ then, using \eqref{eq:mackey morphisms}, $R[G/1, G/K]_{\MF} \cong R[G/K]$ with $RG$ acting by multiplication on the left, thus $R[G/1, G/K]_{\MF}$ is $\FP_\infty$ as a left $RG$-module.
\end{proof}

\begin{Remark}
 $R[ G/S, G/K ]_{\MF}$ is not in general finitely generated as a left $R[W_GS]$-module.  For an example of this let $\mathcal{F}$ be all finite subgroups and choose a group $G$ with a finite subgroup $S$ such that $S \backslash G $ has infinitely many $W_GS $-orbits.  Then, by Example \ref{example:mackey morphisms for S=1},
 \[ R[G/S, G/1]_{\MF} \cong R[S \backslash G] \]
 which is not finitely generated as a left $R[W_GS]$ module.
\end{Remark}

\subsubsection{Induction}

Let $\sigma: \OF \to \MF$ be the covariant functor sending
\[ \sigma(G/H) = G/H \]
\[ \sigma(G/H \stackrel{\alpha}{\longrightarrow} G/K) = (G/H \stackrel{\id}{\longleftarrow} G/H \stackrel{\alpha}{\longrightarrow} G/K  ). \]
Thus $\sigma$ induces restriction, induction, and coinduction between $\OF$-modules and $\MF$-modules.

\begin{Lemma}\label{lemma:mackey struc of cov free as OF}\cite[Proposition 3.6]{MartinezPerezNucinkis-MackeyFunctorsForInfiniteGroups}
There is an $\OF$-module isomorphism:
\[ \Res_\sigma R[G/H, -]_{\MF} \cong \bigoplus_{L \le H} R \otimes_{W_HL} R[G/L, -]_{\OF}. \]
\end{Lemma}

Let $B^G$ denote the Burnside functor $B^G$ which, by an abuse of notation since $G/G$ is not an object of $\MF$, can be defined as \index{Burnside functor $B^G$}
\[
B^G = R[-,G/G]_{\MF}. 
\]
Upon evaluation at $G/K$ for some $K \in \mathcal{F}$,
\[
B^G(G/K) = \bigoplus_{\substack{L \le K \\ \text{Up to $K$-conjugacy}}} R_L. 
\]
This is not so dissimilar from the case of the orbit category $\OF$ where, using a similar abuse of notation, one could view $\uR$ as $R[-,G/G]_{\OF}$. 

\begin{Example}
 If $\uR$ is the constant contravariant $\OF$-module then using Lemma \ref{lemma:mackey struc of cov free as OF},
\begin{align*} 
\Ind_\sigma \uR(G/H) &\cong R[G/H, \sigma(-)]_{\MF}  \otimes_{\OF} \uR  \\
 &\cong \bigoplus_{L \le H} R \otimes_{W_HL} R[G/L, -]_{\OF} \otimes_{\OF} \uR \\
&\cong \bigoplus_{L \le H} R.
\end{align*}
Checking the morphisms as well, one sees that
\[
\Ind_\sigma \uR \cong B^G. 
\]
\end{Example}

\begin{Prop}\label{prop:mackey indOF to MF take pr of R to pr of BG}\cite[Theorem 3.8]{MartinezPerezNucinkis-MackeyFunctorsForInfiniteGroups}
 Although induction with $\sigma$ is not exact in general, induction with $\sigma$ takes contravariant resolutions of $\uR$ by projective $\OF$-modules to resolutions of $B^G$ by projective $\MF$-modules.
\end{Prop}

\subsubsection{Homology and cohomology}
We define the Mackey cohomology and Mackey homology for any contravariant $\MF$-module $M$ and covariant $\MF$-module $A$ as \index{Mackey (co)homology $H^*_{\MF}(G,-)$ and $H_*^{\MF}(G,-)$}
\[  
H^*_{\MF} ( G, M) = \Ext_{\MF}^*(B^G, M)  
\]
\[ 
H_*^{\MF} ( G, A) = \Tor_{\MF}^*(B^G, A)  .
\]

A corollary of Proposition \ref{prop:mackey indOF to MF take pr of R to pr of BG} is the following.
\begin{Cor}\label{cor:mackey H^n_MF is H^n_OF}\cite[Theorem 3.8]{MartinezPerezNucinkis-MackeyFunctorsForInfiniteGroups}
 \[ H^n_{\MF}(G, M) \cong H^n_{\OF}(G, \Res_{\sigma}M). \]
\end{Cor}

$G$ is said to be $\MFFP_n$ if there is a projective resolution of $B^G$, finitely generated up to degree $n$, and $G$ has $\MFcd G \le n$ if there is a length $n$ projective resolution of $B^G$ by $\MF$-modules. \index{Mackey cohomological dimension $\MFcd$}\index{MFFPn@$\MFFP_n$ condition}

\subsection{Cohomological Mackey functors}

A Mackey functor $M$ is called cohomological if, using the language of Remark \ref{remark:mackey green description}, it satisfies
\begin{equation}\label{eq:cohom mackey axiom}
M(I^H_K) \circ M(R^H_K) = \left( m \mapsto \vert H : K \vert m \right)
 \end{equation}
for all subgroups $K \le H$ in $\mathcal{F}$.  Recall from Remark \ref{remark:mackey green description} that to describe a Mackey functor $M$ it is sufficient to describe it on objects and on the induction, restriction and conjugation morphisms in $\MF$ ($I^H_K$, $R^H_K$ and $c_g$), we use this in the examples below.
\begin{Example}[Group cohomology]
The group cohomology functor is cohomological Mackey, more precisely the functor
\begin{align*}
H^n(-, R) : G/H &\longmapsto H^n(H, R).
\end{align*}
Where $H^n(-, R) (c_g)$ is induced by conjugation, $H^n(-, R) (R_K^H)$ is the usual restriction map and $H^n(-, R) (I_K^H)$ is the transfer (see for example  \cite[\S III.9]{Brown}).  That the group cohomology functor satisfies \eqref{eq:cohom mackey axiom} is \cite[III.9.5(ii)]{Brown}. 
\end{Example}

\begin{Example}[Fixed point and fixed quotient functors]\label{example:mackey FP and FQ}\index{Fixed point functor $M^-$}
If $M$ is a $R G$-module then we write $M^-$ for the fixed point functor
\[  M^- : G/H \longmapsto M^H \] 
where $M^H = \Hom_{RH}(R, M)$.  For any $K \le H$ in $\mathcal{F}$, $M^- (R_K^H)$ is the inclusion, $M^-(I_K^H)$ is the trace $m \mapsto \sum_{h \in H/K} hm$, and $M^-(c_g)$ is the map $m \mapsto gm$.

We write $M_-$ for the fixed quotient functor\index{Fixed quotient functor $M_-$}
\[ M_- : G/H \longmapsto M_H \]
where $M_H = R \otimes_{RH} M$.  For any $K \le H$ in $\mathcal{F}$, $M_- (R_K^H)$ is the trace $1 \otimes m \mapsto 1 \otimes \sum_{h \in H/K} hm$, $M_-(I_K^H)$ is the inclusion, and $M_-(c_g)$ is the map $m \mapsto gm$.
\end{Example}

\begin{Lemma}\label{lemma:HF FP and FQ adjunction}\cite[Lemma 4.2]{MartinezPerezNucinkis-MackeyFunctorsForInfiniteGroups}\cite[6.1]{ThevenazWebb-SimpleMackeyFunctors} 
There are Mackey functor isomorphisms for any $RG$-module $M$, 
\[ \CoInd_{RG}^{\MF} M \cong M^- \]
\[ \Ind_{RG}^{\MF} M \cong M_- \]
where induction and coinduction are with the functor $\widehat{\ZZ G} \to \MF$ given by composition of the usual inclusion functor $\widehat{\ZZ G} \to \OF$ and the functor $\sigma : \OF \to \MF$.  Thus there are also adjoint isomorphisms, for any Mackey functor $N$.
\[ \Hom_{RG}(N(G/1), M) \cong \Hom_{\MF}(N, M^-) \]
\[ \Hom_{RG}(M, N(G/1)) \cong \Hom_{\MF}(M^-, N) \]
\end{Lemma}

As observed by Th\'evenaz and Webb in \cite[\S 16]{ThevenazWebb-StructureOfMackeyFunctors}, in \cite{Yoshida-GFunctors2} Yoshida proves that the category of cohomological Mackey modules is isomorphic to the category of modules over the Hecke category $\HeckeF$, which we shall describe below.\index{Hecke category $\HeckeF$}  Yoshida concentrates mainly on finite groups but observes in \cite[\S 5, Theorem 4.3$^\prime$]{Yoshida-GFunctors2} that this isomorphism will hold for $\MF$-modules, where $\mathcal{F}$ is any subfamily of the family of finite groups.  

The Hecke category $\HeckeF$ has for objects the transitive $G$-sets with stabilisers in $\mathcal{F}$.  The morphisms between the objects $G/H$ and $G/K$ are exactly the $\ZZ G$-module homomorphisms, $\Hom_{\ZZ G}(\ZZ[G/H], \ZZ[G/K])$.  

\begin{Remark}
 In \cite{Yoshida-GFunctors2}, Yoshida actually uses the category $\HeckeF^\prime$, this category has the same objects, but the morphisms between $G/H$ and $G/K$ are the $RG$-module homomorphisms $\Hom_{R G}(R [G/H], R[G/K])$.  He then studies $R$-additive functors from $\HeckeF^\prime$ into the category of left $R$-modules and proves these are exactly the cohomological Mackey functors.  
 We claim that the category of $R$-additive functors $\HeckeF^\prime \to \RMod$ and the category of additive functors $\HeckeF \to \RMod$ are isomorphic, where the isomorphism preserves the values the functors take on objects.
 
 Since $\ZZ[G/H]$ is finitely presented as a $\ZZ G$-module and $R$ is flat as a $\ZZ$-module there is an isomorphism for all $H, K \in \mathcal{F}$ (the proof is essentially the proof of \cite[3.3.8]{Weibel})
 \[
 \Hom_{\ZZ G}(\ZZ[G/H], \ZZ[G/K]) \otimes_{\ZZ} R \cong \Hom_{R G}(R[G/H], R[G/K]).
 \]
 Using the above and that $\Hom_{\ZZ G}(\ZZ[G/H], \ZZ[G/K])$ is free as a $\ZZ$-module, there is a natural isomorphism for any $R$-module $A$
 \begin{align*}
 \Hom_{\ZZ} ( \Hom_{\ZZ G}(\ZZ[G/H], &\ZZ[G/K]), A) \\
 &\cong \Hom_{R} ( \Hom_{\ZZ G}(\ZZ[G/H], \ZZ[G/K]) \otimes_{\ZZ} R, A) \\
 &\cong \Hom_{R} ( \Hom_{R G}(R[G/H], R[G/K]) , A).
 \end{align*}
The claim follows from this isomorphism.
 \end{Remark}

\begin{Remark}
 In \cite{Degrijse-ProperActionsAndMackeyFunctors} Degrijse considers the categories $\text{Mack}_\mathcal{F} G$ and $\text{coMack}_\mathcal{F}G$.  In the notation used here $\text{Mack}_\mathcal{F} G$ is the category of $\MF$-modules and $\text{coMack}_\mathcal{F}G$ is the subcategory of cohomological Mackey functors, Degrijse doesn't study modules over $\HeckeF$.
\end{Remark}

\begin{Lemma}[Free and projective $\HeckeF$-modules]\label{lemma:HF frees are fp}\cite[Theorem 16.5(ii)]{ThevenazWebb-StructureOfMackeyFunctors}\par
 The free $\HeckeF$-modules are exactly the fixed point functors of permutation modules with stabilisers in $\mathcal{F}$, and the projective $\HeckeF$-modules are exactly the fixed point functors of direct summands of permutation modules with stabilisers in $\mathcal{F}$.
\end{Lemma}

Th\'evenaz and Webb describe a map $\pi: \MF \to \HeckeF$ (they call this map $\alpha$), taking objects $G/H$ in $\MF$ to $G/H$ in $\HeckeF$ and on morphisms as follows, for any $K \le H$,
\begin{itemize}
 \item $\pi (R_K^H) $ is the natural projection map $\ZZ[G/K] \to \ZZ[G/H]$.
 \item $\pi (I_K^H) $ takes $gH \mapsto \sum_{h \in H/K} ghK $.
 \item $\pi (c_x) $ takes $gH \mapsto gx H^x $.
\end{itemize}
If $M$ is an $\HeckeF$-module then it is straightforward to check that $M \circ \pi $ is a $\MF$-module, see for example \cite[p.809]{Tambara-HomologicalPropertiesOfTheEndomorphismRings} for a proof.  Moreover, every cohomological Mackey functor $M : \MF \to \RMod$ factors through the map $\pi$, this is the main result in \cite{Yoshida-GFunctors2}, see also \cite[\S 7]{Webb-GuideToMackeyFunctors}.  Thus we may pass freely between cohomological Mackey functors and modules over $\HeckeF$.

\begin{Lemma}\label{lemma:mackey yoshida CF struc}\cite[Lemma 3.1$^\prime$]{Yoshida-GFunctors2}  There is an isomorphism for any finite subgroups $H$ and $K$ of $G$,
\[
 R[H \backslash G / K] \cong R[G/H, G/K]_{\HeckeF} .
\]
 Under this identification, morphism composition is given by
 \[ 
 (HxK) \cdot (KyL) = \sum_{z \in H \backslash G / L} \vert (HxK \cap zLy^{-1}K) / K\vert \, ( HzL).
 \]
\end{Lemma}
\begin{Remark}\label{remark:mackey yoshida CF struc}
The identification in the lemma above relates to the usual definition of $R[G/H, G/K]_{\HeckeF}$ as $\Hom_{RG}(R[G/H], R[G/H])$ with the isomorphism
\begin{align*}
 \psi: R[H \backslash G / K] &\stackrel{\cong}{\longrightarrow} \Hom_{RG} (R[G/H], R[G/K]) \\
 HxK &\longmapsto \left( gH \longmapsto \sum_{u \in H/(H \cap xKx^{-1})} guxK \right)
 \end{align*} 
Notice that $\psi$ satisfies
\[ 
\psi( (HxK) \cdot (KxL) ) = \psi(KxL) \circ \psi(HxK).
\]
\end{Remark}

\begin{Lemma}\label{lemma:induced from alpha on yoshida description}
 If $\alpha:G/L \longrightarrow G/K$ is the $G$-map $L \mapsto xK$ then the induced map $\alpha_*$ on $R[G/H, -]_{\HeckeF}$ can be written as
 \begin{align*}
 \alpha_* : R[H \backslash G / L] &\longrightarrow R[H \backslash G / K] \\
 (HzL) &\longmapsto \sum_{y \in H \cap K^{(zx)^{-1}} / H \cap L^{z^{-1}}} (HyzxK).
 \end{align*}
\end{Lemma}
\begin{proof}
 We calculate 
 \begin{align*}
  \psi &\left( \sum_{y \in H\cap K^{(zx)^{-1}}/H \cap L^{z^{-1}}} HyzxK \right) \\
  &= \left( H \longmapsto \sum_{y \in H\cap K^{(zx)^{-1}} / H \cap L^{z^{-1}}} \sum_{u \in H / H \cap K^{(yzx)^{-1}}} uyzxK \right) \\
  &= \left( H \longmapsto \sum_{u \in H / H \cap L^{z^{-1}}} uzxL \right)
 \end{align*}
 which is exactly $\alpha_* (\psi (HzL))$.  The final equality comes from the fact that $y \in K^{(zx)^{-1}}$ so $K^{(yzx)^{-1}} = K^{(zx)^{-1}}$.
\end{proof}

\subsubsection{Explicit description of \texorpdfstring{$\pi$}{the projection map}} 
 Using the identification of Lemma \ref{lemma:mackey yoshida CF struc}, for any $K \le H$, we can describe $\pi$ as follows.
\begin{itemize}
 \item $\pi (R_K^H) = KH $, since according to Lemma \ref{lemma:mackey yoshida CF struc}, $KH $ corresponds to the morphism $gK \mapsto gH$, which is exactly Th\'evenaz and Webb's description of $\pi(R_K^H)$.
 \item $\pi (I_K^H) = HK$, as according to Lemma \ref{lemma:mackey yoshida CF struc}, $HK$  corresponds to the morphism $gH \mapsto \sum_{u \in H/K} uK$, which is Th\'evenaz and Webb's description of $\pi(I_K^H)$.
 \item $\pi (c_x) = HxH^x$, similarly to the above because $HxH^x$ corresponds to the morphism $g H \mapsto g x H^x$.
\end{itemize}

It is interesting to write down the effect of $\pi$ on a basic morphism
\[ m = \left( \vcenter{\xymatrix@-15pt{
& G/L \ar_{\alpha_1}[ld] \ar^{\alpha_x}[rd] & \\
G/H & & G/K
}} \right) \]
This morphism may be rewritten as
\begin{align*} \MFMor{
& G/L \ar_{\alpha_1}[ld] \ar^{\alpha_x}[rd] & \\
G/H & & G/L
} \circ \MFMor{
& G/L \ar_{\alpha_1}[ld] \ar^{\alpha_x}[rd] & \\
G/L & & G/L^x
} \circ \MFMor{
& G/L^x \ar_{\alpha_1}[ld] \ar^{\alpha_1}[rd] & \\
G/L^x & & G/K
} 
\end{align*}
So, 
\[
m= R^K_{L^x} \circ c_x \circ I_L^H. 
\]
Using the definition of $\pi$ in \cite[\S 16]{ThevenazWebb-StructureOfMackeyFunctors}, $\pi(m)$ maps 
\begin{align*}
 \pi(m) &= \pi(R^K_{L^x} \circ c_x \circ I_L^H) \\
 &= \left( H \mapsto \sum_{h \in H/L} hx K \right) \\
 &= \left( H \mapsto \sum_{h \in H / H \cap K^{x^{-1}}} \sum_{y \in H \cap K^{x^{-1}} / L} hyx K \right) \\
 &= \left( H \mapsto \sum_{y \in H \cap K^{x^{-1}} / L} \sum_{h \in H / H \cap K^{(yx)^{-1}}} hyxK \right) \\
 &= \sum_{y \in H \cap K^{x^{-1}} / L} (HyxK) \\
 &= \vert H \cap K^{x^{-1}} : L \vert ( HxK ).
\end{align*}

In summary, 
\begin{equation*}
 \pi \left( \vcenter{\xymatrix@-15pt{
& G/L \ar_{\alpha_1}[ld] \ar^{\alpha_x}[rd] & \\
G/H & & G/K
}} \right) = \vert H \cap K^{x^{-1}}  : L \vert (HxK).
\end{equation*}

\subsubsection{Homology and cohomology}\label{subsubsection:mackey cohom mackey homology and cohomology}

In Section \ref{section:Homology and Cohomology of Cohomological Mackey Functors} we will prove results similar to Proposition \ref{prop:mackey indOF to MF take pr of R to pr of BG} and Corollary \ref{cor:mackey H^n_MF is H^n_OF}, showing that inducing a projective resolution of $\uR$ by projective $\OF$-modules yields a projective resolution of $R^-$ by projective $\HeckeF$-modules.  For any group $G$, we define the cohomology and homology functors $H^*_{\HeckeF}(G, -)$ and $H_*^{\HeckeF}(G, -)$ as \index{Cohomological Mackey (co)homology $H^*_{\HeckeF}(G, -)$ and $H_*^{\HeckeF}(G, -)$}
\[ H^*_{\HeckeF} ( G, M) = \Ext_{\HeckeF}^*(R^-, M)\]
\[H_*^{\HeckeF} ( G, A) = \Tor_{\HeckeF}^*(R^-, A)\]
where $M$ is any contravariant $\HeckeF$-module and $A$ is any covariant $\HeckeF$-module.  In Proposition \ref{prop:HF HF cohomology is OF cohomology} we show that there is an isomorphism
\[H^n_{\HeckeF}(G, M) \cong H^n_{\OF}(G, \Res_{\pi \circ \sigma} M). \]
The $\HeckeF$ cohomological dimension of a group $G$, denoted $\HFcd G$, is defined to be the length of the shortest projective resolution of $R^-$ by $\HeckeF$-modules, or equivalently\index{Cohomological Mackey cohomological dimension $\HFcd G$}
\[
\HFcd G = \sup \{ n \: : \: H^n_{\HeckeF}(G, M) \neq 0) , \: \text{$M$ some $\HeckeF$-module.}\} 
\]
Note that in  \cite{Degrijse-ProperActionsAndMackeyFunctors} the $\HeckeF$ cohomological dimension is defined as
\[
\HFcd G = \sup \{ n \: : \: H^n_{\OF}(G, \Res_{\pi \circ \sigma} M )\neq 0 , \: \text{$M$ some $\HeckeF$-module.}\} 
\]
These two definitions are equivalent by the isomorphism of Proposition \ref{prop:HF HF cohomology is OF cohomology} mentioned above.

We say $G$ is $\HFFP_n$ if there exists a projective $\HeckeF$-module resolution of $R^-$, finitely generated up to degree $n$.\index{HFFPn@$\HFFP_n$ condition}

\section{\texorpdfstring{$\FP_n$}{FPn} conditions for Mackey functors}\label{section:mackey FPn for mackey}

As far as we are aware, there are no results in the literature on the conditions $\MFFP_n$.  We show in this section that the conditions $\MFFP_n$ and $\OFFP_n$ are equivalent.  From this point on, unless otherwise stated, all results are valid over any commutative ring $R$.

\begin{Prop}\label{prop:mackey OFFPn implies MFFPn}
 If $G$ is $\OFFP_n$ then $G$ is $\MFFP_n$.
\end{Prop}
\begin{proof}
Combine Proposition \ref{prop:mackey indOF to MF take pr of R to pr of BG} with the fact that induction preserves finite generation (Proposition \ref{prop:C properties of res ind coind}(2)).
\end{proof}

Recall that $G$ is $\OFFP_0$ if and only if $\mathcal{F}$ has finitely many conjugacy classes
(Corollary \ref{cor:uFPn equivalent conditions}).  In the lemmas below $\mathcal{F}/G$ denotes the poset of conjugacy classes in $\mathcal{F}$, ordered by subconjugation.  We write $H \le_G K$ if $H$ is subconjugate to $K$. 

\begin{Lemma}\label{lemma:mackey MFFP0}
 $G$ is $\MFFP_0$ if and only if $G$ is $\OFFP_0$.
\end{Lemma}

\begin{proof}
 We prove first that if $G$ is $\MFFP_0$ then $\mathcal{F}/G$ has a finite cofinal subset, since $\mathcal{F}$ is a subfamily of the family of finite subgroups this implies that $\mathcal{F}/G$ is finite.

 Let $f$ be an $\MF$-module morphism 
\[ f: R[-,G/K]_{\MF} \longrightarrow B^G \cong R[-,G/G]_{\MF}. \]
 Firstly, we claim that the element $m$ of $R[G/S,G/G]_{\MF}$ given by 
\[ m = \big( G/S \stackrel{\id}{\longleftarrow} G/S \longrightarrow G/G \big) \]
 cannot be in the image of $f(G/S)$ unless $S$ is subconjugate to $K$.  Assume for a contradiction that $S$ is not subconjugate to $K$ and assume $m$ is in the image of $f(G/S)$.  Thus $m = f(G/S)\varphi$ for some $\varphi \in [G/S, G/K]_{\MF}$.   Thinking of $f$ as a natural transformation gives the commutative diagram below
\[
\xymatrix@C+20pt{
R[G/S, G/K]_{\MF} \ar^{f(G/S)}[r] & R[G/S, G/G]_{\MF} \\
R[G/K, G/K]_{\MF} \ar^{f(G/K)}[r] \ar_{\varphi^*}[u] & R [G/K, G/G]_{\MF} \ar_{\varphi^*}[u]
}
\]
where
\begin{align*}
 m &= f(G/S)\varphi \\
&= f(G/S) \circ \varphi^* \id_{[G/K,G/K]_{\MF}} \\
&= \big( \varphi^* \circ f(G/K) \big) ( \id_{[G/K,G/K]_{\MF}} ).
\end{align*}
Let $f(G/K) ( \id_{[G/K,G/K]_{\MF}} ) = \sum_i r_i x_i$, where $r_i \in R$ and the $x_i$ are basic morphisms in $R[G/K, G/G]_{\MF}$.  Similarly, let $\varphi = \sum_j s_j y_j$ for $s_j \in R$ and where the $y_j$ are basic morphisms in $R[G/S, G/K]_{\MF}$.  By assumption we have that
\begin{align*}
 m &= \varphi^* \sum_i r_i x_i \\
&=  \sum_i r_i x_i \circ \sum_j s_j y_j \\
&= \sum_{i,j} (r_is_j) x_i \circ y_j .
\end{align*}
There must exist some $i$ and $j$ for which $x_i \circ y_j $ is a morphism which, when written as a sum of basic morphisms, has one component some multiple of $m$.  
We calculate $x_i \circ y_j$ for this $i$ and $j$.  Write $x_i$ and $y_j$ in their standard forms as below,
\[x_i = \Big( G/K \longleftarrow G/L_i \longrightarrow G/G \Big)\]
\[y_j = \Big( G/S \longleftarrow G/J_j \longrightarrow G/K \Big).\]
Their composition is (see Lemma \ref{lemma:mackey morphism pullback})
\[
x_i \circ y_j = \sum_k \left( \vcenter{\xymatrix@-15pt{
&  &  G/X_k\ar[rd] \ar[ld] & & \\
& G/J_j \ar[rd] \ar[ld] & & G/L_i \ar[rd] \ar[ld]  & \\
G/S & & G/K & & G/G
}}\right)
\]
where $X_k$ is some finite subgroup of $G$ which is subconjugate to both $L_i$ and $J_j$.  

We claim $\lvert J_j \rvert $ is strictly smaller than $\lvert S \rvert$.  Since $J_j$ is subconjugate to $S$ we have $\lvert J_j \rvert \le \lvert S \rvert$.  If the cardinalities were equal then $S$ and $J_j$ would be conjugate, but $J_j$ is subconjugate to $K$ whereas by assumption $S$ is not subconjugate to $K$.

Since $\lvert X_k \rvert \le \lvert J_j \rvert \lneq \lvert S \rvert$, the subgroup $X_k$ cannot be conjugate to $S$.  This contradicts our earlier assertion that $x_i \circ y_j $ when written as a sum of basic morphisms, has one component some multiple of $m$.  Thus, for $m$ to be in the image of $f(G/S)$, $S$ must be subconjugate to $K$.

Now, if $G$ is $\MFFP_0$ then $B^G$ admits an epimorphism from some finitely generated free
\[ \bigoplus_{i \in I} R[ -, G/K_i]_{\MF} \longtwoheadrightarrow B^G .\]
As this set $I$ is finite, the argument above implies that all the subgroups in $\mathcal{F}$ are subconjugate to one of a finite collection of subgroups in $\mathcal{F}$.  Thus there is a finite cofinal subset of $\mathcal{F}/G$ and $\mathcal{F}/G$ is finite.

For the converse, use Proposition \ref{prop:mackey OFFPn implies MFFPn}.
\end{proof}

This remainder of this section is devoted to a proof that for any $n$, $\MFFP_n$ implies $\OFFP_n$.  We will assume $G$ is $\MFFP_0$, equivalently $\mathcal{F}$ contains finitely many conjugacy classes.

In \cite[4.9, 4.10]{HambletonPamukYalcin-EquivariantCWComplexesAndTheOrbitCategory}, there are the following definitions, for $M$ an $\OF$-module 
\[ D_H M  = \CoInd_{R[WH]}^{\OF} M(G/H) \]
\[ j_H : M \longrightarrow D_H M \]
where $\CoInd_{R[WH]}^{\OF}$ denotes coinduction (see Section \ref{section:C restriction induction and coinduction} for the definition of coinduction) with the functor $\iota : \ZZ[WH] \longrightarrow \OF$.  Here we view $\ZZ[WH]$ as a category with one object and morphisms elements of $\ZZ[WH]$, then $\ZZ[WH]$ has property $(A)$ and $\iota$ maps the one object to $G/H$ and morphisms to the free abelian group on the automorphisms of $G/H$ in $\OF$.  Equivalently, 
\[ \CoInd_{R[WH]}^{\OF} M(G/H) \cong \Hom_{R[WH]}(R[G/H, -]_{\OF}, M(G/H)). \]
The map $j_H$ is the counit of the adjunction between coinduction and restriction.  Since evaluation, coinduction, and counits are all natural constructions, $D_H$ and $j_H$ are natural.  Crucially the $\OF$-module $D_HM$ extends to a Mackey functor \cite[Example 4.8]{HambletonPamukYalcin-EquivariantCWComplexesAndTheOrbitCategory}.  Also defined are:
\[ D M = \prod_{H \in \mathcal{F} /G} D_H M \]
\[ CM  = \text{CoKer} \left( C  \stackrel{\prod j_H}{\longrightarrow} DM \right). \]
Again all the constructions are natural and $DM$ extends to a Mackey functor.  Naturality means that if $M_\lambda$, for $\lambda \in \Lambda$, is a directed system of $\OF$-modules then $D M_\lambda$ and $CM_\lambda$ form directed systems also.

\begin{Lemma}\label{lemma:colim M is 0 then colim DM is 0}
If $M_\lambda$ is a directed system of $\OF$-modules with $\varinjlim M_\lambda  = 0$ then $\varinjlim DM_\lambda  = 0$. 
\end{Lemma}
\begin{proof}
Since the colimit of $M_{\lambda}$ is zero, so is the colimit of $M_\lambda (G/H)$, and for any $K \in \mathcal{F}$,
\begin{align*} 
\varinjlim D_H M_{\lambda}(G/K) &= \varinjlim \CoInd_{R[WH]}^{\OF} M_{\lambda}(G/H) (G/K) \\
&= \varinjlim \Hom_{R[WH]} ( R[G/H, G/K]_{\OF}, M_\lambda(G/H) ) \\
&= \varinjlim \Hom_{R[WH]} \left( \bigoplus_{i \in I} R[WH/WH_i], M_\lambda (G/H) \right) 
\end{align*}
where the last line is Lemma \ref{lemma:Z[G/H,G/K] as sum of WH submodules}, the indexing set $I$ is finite and $WH_i$ is a finite subgroup of $WH$.  Hence
\begin{align*}
\varinjlim \Hom_{R[WH]} &\left( \bigoplus_{i \in I} R[WH/WH_i], M_\lambda (G/H) \right) \\
&\cong \bigoplus_{i \in I} \varinjlim \Hom_{R[WH]} \left( R[WH/WH_i] , M_\lambda(G/H)\right) \\
&\cong \bigoplus_{i \in I} \varinjlim \Hom_{R[WH_i]} ( R, M_\lambda (G/H) ) \\
&= 0
\end{align*}
where the final zero is by the Bieri--Eckmann criterion (Theorem \ref{theorem:C bieri-eckmann criterion}), since $R$ is $R[WH_i]$-finitely generated.  Thus 
\begin{align*}
\varinjlim D M_{\lambda}(G/K)  &= \varinjlim \prod_{H \in \mathcal{F}/G} D_H M_\lambda (G/K) \\
&= \prod_{H \in \mathcal{F}/G} \varinjlim D_HM_\lambda(G/K) \\
&= 0
\end{align*}
where the commuting of the product and the colimit is because the product is finite ($\mathcal{F}/G$ is assumed finite).
\end{proof}

\begin{Lemma}
 If $M_\lambda$ is a directed system of $\OF$-modules with $\varinjlim M_\lambda = 0$ then $\varinjlim CM_\lambda = 0$.
\end{Lemma}
\begin{proof}
There is a natural short exact sequence for each $\lambda$
\[ 0 \longrightarrow M_{\lambda} \longrightarrow DM_{\lambda} \longrightarrow CM_{\lambda} \longrightarrow 0 . \]
Since the colimit of the left hand and centre term are zero (Lemma \ref{lemma:colim M is 0 then colim DM is 0}), and colimits are exact in the category of $\OF$-modules, so $\varinjlim CM_\lambda = 0$ also.
\end{proof}

\begin{Prop}\label{prop:MFFPn implies OFFPn}
 If $G$ is $\MFFP_n$ then $G$ is $\OFFP_n$.
\end{Prop}
\begin{proof}
Let $G$ be of type $\MFFP_n$ and let $M_\lambda $, for $\lambda \in \Lambda$, be a directed system of $\OF$-modules with colimit zero.  Following the notation in \cite{Degrijse-ProperActionsAndMackeyFunctors}, we define 
\[ C^0M_\lambda  = M_\lambda \]
\[ C^iM_\lambda = C C^{i-1}M_\lambda \]
for all natural numbers $i \ge 0$ and all $\lambda \in \Lambda$.  There are short exact sequences of directed systems, 
\[ 0 \longrightarrow C^iM_\lambda \longrightarrow DC^i M_\lambda \longrightarrow C^{i+1}M_\lambda \longrightarrow 0 \]
all the terms of which have colimit zero. 

As $G$ is assumed $\MFFP_n$ and $DC^iM_\lambda$ extends to a Mackey functor for all $i$, the Bieri--Eckmann criterion (Theorem \ref{theorem:C bieri-eckmann criterion}) gives that for all $m \le n$,
\[ \varinjlim H^m_{\OF} (G, DC^iM_\lambda) = 0 \]
and thus using exactness of colimits and the long exact sequence associated to cohomology gives that for all non-negative integers $m$ and $i$,
\[ \varinjlim H^m_{\OF} (G, C^{i+1}M_{\lambda}) = \varinjlim H^{m+1}_{\OF} (G, C^{i}M_{\lambda}). \]
So,
\begin{align*}
 \varinjlim H^{m}_{\OF} (G,M_{\lambda}) &= \varinjlim H^{m-1}_{\OF}(G, C^1 M_\lambda) \\
&\cong \cdots \\
&\cong \varinjlim H^0_{\OF}(G, C^mM_\lambda) \\
&\cong 0
\end{align*}
where the zero is from the Bieri--Eckmann criterion (Theorem \ref{theorem:C bieri-eckmann criterion}), because $G$ is assumed $\MFFP_0$ hence $\OFFP_0$ by Lemma \ref{lemma:mackey MFFP0}.  Using the Bieri--Eckmann criterion again, $G$ is $\OFFP_n$. 
\end{proof}

\begin{Cor}\label{cor:OFFPn iff MFFPn}
 The conditions $\OFFP_n$ and $\MFFP_n$ are equivalent.
\end{Cor}
\begin{proof}
 Combine Propositions \ref{prop:mackey OFFPn implies MFFPn} and \ref{prop:MFFPn implies OFFPn}.
\end{proof}

\section{Homology and cohomology of cohomological Mackey functors}\label{section:Homology and Cohomology of Cohomological Mackey Functors}

The main result of this section is Proposition \ref{prop:HF sigma ind preserves proj res of R}, that we may induce projective  $\OF$-module resolutions of $\uR$ to projective $\HeckeF$-module resolutions of $R^-$.  The following diagram shows the relationship between the different induction functors we will be using (for a small category $\C$, we denote by $\C$-Mod the category of contravariant $\C$-modules).

\[
 \xymatrix{
   & \text{$\HeckeF$-Mod} \\
   \text{$\OF$-Mod} \ar^-{\Ind_\sigma}[r] \ar^-{\Ind_{\pi \circ \sigma}}[ur] & \text{$\MF$-Mod} \ar_-{\Ind_{\pi}}[u]
 }
\]

\begin{Lemma}\label{lemma:HF struc of Res cov frees}
For any $L \in \mathcal{F}$, there is an isomorphism of covariant $\OF$-modules
\[ \Res_{\pi \circ \sigma} R[G/L, -]_{\HeckeF} \cong \Hom_{RL} (R, R[G/1, -]_{\OF}). \]
\end{Lemma}
\begin{proof}
If $H$ is a subgroup in $\mathcal{F}$, then evaluating the left hand side at $G/H$ yields $R[G/L, G/H]_{\HeckeF}$ while evaluating the right hand side at $G/H$ yields 
\begin{align*}
\Hom_{RL}(R, R[G/H]) &\cong \Hom_{RG}(RG \otimes_{RL} R, R[G/H]) \\
&\cong \Hom_{RG} (R[G/L], R[G/H]) \\
&\cong R[G/L, G/H]_{\HeckeF} 
\end{align*}
where the first isomorphism is \cite[p.63 (3.3)]{Brown}.  If $ \alpha_x:G/H \to G/K$ is the $G$-map $H \mapsto xK$ then, looking at the left hand side,
\begin{align*}
 \Res_{\pi \circ \sigma} R[G/L, -]_{\HeckeF} (\alpha_x) &= R[G/H, -]_{\HeckeF}( c_x \circ R_H^{K^{x^{-1}}}  ) \\
 &\cong R[G/H, -]_{\HeckeF}( c_x ) \circ R[G/H, -]_{\HeckeF}( R_H^{K^{x^{-1}}} ).
\end{align*}
But $R[G/H, -]_{\HeckeF}( R_H^{K^{x^{-1}}} )$ is post-composition with the $G$-map 
\[ \alpha_1: G/H \to G/K^{x^{-1}} \]
and $R[G/H, -]_{\HeckeF}( c_x )$ is post-composition with the $G$-map 
\[ \alpha_x: G/K^{x^{-1}} \to G/K. \]
In summary, $\Res_{\pi \circ \sigma} R[G/L, -]_{\HeckeF} (\alpha_x)$ is the map:
\begin{align*}
\Hom_{RG}(R[G/L], R[G/H]) &\longrightarrow \Hom_{RG}(R[G/L], R[G/K]) \\
f &\longmapsto  \alpha_x \circ f 
\end{align*}

Now, the right hand side, recall that 
\[ R[G/L, -]_{\OF} (\alpha_x) : f \longmapsto \alpha_x \circ f \]
so $\Hom_{RL} (R, R[G/1, -]_{\OF})(\alpha_x)$ is the map:
\begin{align*} 
 \Hom_{RL} (R, R[G/H]) &\longrightarrow \Hom_{RL} (R, R[G/K]) \\
 f &\longmapsto \alpha_x \circ f
 \end{align*}
Showing the left and right hand sides agree on morphisms.
\end{proof}

Recall that $\Ind_{RG}^{\OF}$ denotes induction with the functor $\iota:\ZZ G \longrightarrow \OF$, where we view $\ZZ G$ as the single object category whose morphisms are elements of $\ZZ G$ and $\iota$ maps the single object to $G/1$.  Equivalently for an $RG$-module $M$, 
\[ \Ind_{R G}^{\OF} M \cong R[-, G/1]_{\OF} \otimes_{R G} M. \]

\begin{Lemma}\label{lemma:contravariant induction is exact}
 The functor $\Ind_{R G}^{\OF}$ is exact.
\end{Lemma}
\begin{proof}
 This is because for any $H \in \mathcal{F}$,
\[ \Ind_{R G}^{\OF} M (G/H) = \left\{ \begin{array}{l l} M & \text{ if $H = 1$} \\ 0 & \text{ else.} \end{array} \right. \]
\end{proof}

\begin{Lemma}\label{lemma:mackey ind RH to OFG of R is OFFPinfty}
 For any finite subgroup $H$ of $G$, the $\OF$-module $\Ind_{R G}^{\OF} \Ind_{R H}^{R G} R$ is of type $\OFFP_\infty$.
\end{Lemma}
\begin{proof}
 Since $R$ is $\FP_\infty$ as a $R H$ module, $\Ind_{R H}^{R G}R $ is of type $\FP_\infty$ over $R G$.  Choose a finite type free resolution $F_*$ of $\Ind_{R H}^{R G}R$ by $RG$-modules, by Lemma \ref{lemma:contravariant induction is exact}, $\Ind_{R G}^{\OF}F_*$ is a finite type free resolution of $\Ind_{R G}^{\OF}\Ind_{R H}^{R G}R$ by $\OF$-modules.
\end{proof}

\begin{Lemma}\label{lemma:mackey tech lamma for ind}
 If $N$ is a projective $\OF$-module and $H \in \mathcal{F}$, there is an isomorphism
\[ N \otimes_{\OF} \Res_{\pi \circ \sigma} R[G/H, -]_{\HeckeF} \cong \Hom_{RH}(R, N(G/1)). \]
\end{Lemma}

\begin{proof}
The adjointness of induction and restriction gives an isomorphism of $\OF$-modules, for any $\OF$-module $N$,
\begin{align*}
 \Hom_{RH}(R, N(G/1)) &\cong \Hom_{RG}(\Ind_{R H}^{R G} R, N(G/1)) \\
 &\cong \Hom_{\OF} ( \Ind_{R G}^{\OF} \Ind_{R H}^{R G} R, N).
\end{align*}
There is a chain of isomorphisms,
 \begin{align*}
N &\otimes_{\OF} \Res_{\pi \circ \sigma} R[G/H, -]_{\HeckeF} \\
 &\cong N \otimes_{\OF} \Hom_{RH}(R, R[G/1, -]_{\OF}) \\
&\cong N(-) \otimes_{\OF} \Hom_{\OF} ( \Ind_{R G}^{\OF} \Ind_{R H}^{R G} R(?), R[?, -]_{\OF}) \\
&\cong \Hom_{\OF} ( \Ind_{R G}^{\OF} \Ind_{R H}^{R G} R(?), N(?)) \\
&\cong \Hom_{RH}(R, N(G/1))
 \end{align*}
where the first isomorphism is Lemma \ref{lemma:HF struc of Res cov frees} and the second and fourth are the adjoint isomorphism mentioned above.  The third isomorphism is from Lemma \ref{lemma:M fg and N proj then nu iso} for which we need that  $\Ind_{R G}^{\OF}\Ind_{R H}^{R G} R$ is finitely generated, but this is implied by Lemma \ref{lemma:mackey ind RH to OFG of R is OFFPinfty}.
\end{proof}

Recall from Example \ref{example:mackey FP and FQ} the definition of the fixed point functor.  For the constant $RG$-module $R$ the fixed point functor $R^-$ can be described explicitly as $R^H = R$ for all $H \in \mathcal{F}$, and on morphisms,
\begin{align*}
R^-(R^H_K) &= \id_R  \\
R^-(I^H_K) &= \left( r \mapsto \vert H : K \vert r\right)  \\
R^-(c_g) &= \id_R.
\end{align*}

\begin{Lemma}\label{lemma:HF ind uROF is R}
$\Ind_{\pi \circ \sigma} \uR  \cong R^-$.
\end{Lemma}

\begin{proof}
The proof is split into two parts, first we check that the two functors agree on objects, then we check they agree on morphisms.  Throughout the proof $H$, $K$ and $L$ are elements of $\mathcal{F}$.  If $\alpha : G/L \to G/K$ is a $G$-map then we will write $\alpha_*$ for the induced map 
\[ \alpha_* : \Hom_{RG}(R[G/H], R[G/L]) \longrightarrow \Hom_{RG}(R[G/H], R[G/K]) \]
and also for the induced map
\[ \alpha_* : R[H\backslash G / L] \longrightarrow R[H \backslash G / K] \]
where $R[H\backslash G / L]$ is identified with $\Hom_{RG}(R[G/H], R[G/L])$ using the isomorphism $\psi$ of Remark \ref{remark:mackey yoshida CF struc}. 

\noindent\emph{The functors $\Ind_{\pi \circ \sigma} \uR $ and $R^-$ agree on objects:}

For any subgroup $H \in \mathcal{F}$, 
\begin{align*}
 \Ind_{\pi \circ \sigma} &\uR(G/H) = \uR \otimes_{\OF} \Res_{\pi \circ \sigma} R[G/H, -]_{\HeckeF} \\
 &\cong \uR \otimes_{\OF} \Hom_{RG}(R[G/H], R[G/1, -]_{\OF}) \\
 &\cong \left. \bigoplus_{K \in \mathcal{F}} \Hom_{RG}(R[G/H], R[G/K]) \right/ \substack{ x_K \sim \alpha_* x_L \\ \text{$\alpha:G/L \to G/K$ any $G$ map} \\ x_K \in \Hom_{RG}(R[G/H], R[G/K]) \\ x_L \in \Hom_{RG}(R[G/H], R[G/L]) } \\
&\cong \left.  \bigoplus_{K \in \mathcal{F}} R[H \backslash G / K] \right/ \substack{ (HxK) \sim \alpha_* (HxL) \\ \text{$\alpha:G/L \to G/K$ any $G$ map} }
\end{align*}
where the first isomorphism is Lemma \ref{lemma:HF struc of Res cov frees} and the last is Lemma \ref{lemma:mackey yoshida CF struc}.  Let $(HxK) \in R[H\backslash G / K]$ be an arbitrary element, and consider the $G$-map 
\begin{align*}
\alpha_x : G/(H\cap K^{x^{-1}}) &\longrightarrow G/K \\
( H \cap K^{x^{-1}} ) &\longmapsto xK.
\end{align*}
Then, using Lemma \ref{lemma:induced from alpha on yoshida description} we calculate
\[
(\alpha_x)_* \left(H1(H \cap K^{x^{-1}})\right) = (HxK)
\]
so in $\Ind_{\pi \circ \sigma} \uR(G/H)$, the elements $[H \cdot x \cdot K]$ and $[H \cdot 1 \cdot H \cap K^{x^{-1}}]$ are equal, where $[-]$ denotes an equivalence class of elements under the relation $\sim$.  Similarly if $K \le H$ then $[H\cdot 1 \cdot K] = \lvert H : K \rvert [H \cdot 1 \cdot H]$ since if $\alpha_1 : G/K \to G/H$ is the projection, then using Lemma \ref{lemma:induced from alpha on yoshida description} again,
\[
   {\alpha_1}_* (H1K) = \lvert H : K \rvert (H1H) .
\]
Combining the two facts proved above, 
\begin{align*}[H \cdot x \cdot K] = \lvert H : H \cap K^{x^{-1}} \rvert [H \cdot 1 \cdot H].\tag{$\star$}\end{align*}
In particular, any element $[H \cdot x \cdot K]$ is equal to some multiple of $[H \cdot 1 \cdot H]$, so 
\[ \Ind_{\pi \circ \sigma} \uR(G/H) \cong R. \]

\noindent\emph{The functors $\Ind_{\pi \circ \sigma} \uR $ and $R^-$ agree on morphisms:}

Recall from Remark \ref{remark:mackey green description} that we must only check this for the morphisms $R^H_K$, $I^H_K$ and $c_x$. 

Following the generator $[H \cdot 1 \cdot H]$ up the chain of isomorphisms at the beginning of the proof shows the element 
\[ 1 \otimes \id_{R[G/H]} \in R \otimes_{\OF} \Res_{\pi \circ \sigma} R[G/H, - ]_{\HeckeF}\]
generates $\Ind_{\pi \circ \sigma} \uR(G/H) \cong R$, where 
\[\id_{R[G/H]} \in \Hom_{RG}(R[G/H], R[G/H]) \cong R[G/H, G/H]_{\HeckeF}.\]
For some subgroup $K \in \mathcal{F}$ with $K \le H$, 
\[\Ind_{\pi \circ \sigma} \uR(R_K^H) : 1 \otimes \id_{R[G/H]} \mapsto 1 \otimes \pi\]
where $\pi : R[G/K] \mapsto R[G/H]$ is the projection map.  Following this back down the chain of isomorphisms at the beginning of the proof, gives the element $[K \cdot 1 \cdot H]$.  Using $(\star)$, $[K \cdot 1 \cdot H] = [K \cdot 1 \cdot K]$, so $\Ind_{\pi \circ \sigma} \uR(R_K^H)$ is the identity on $R$, as required.

Similarly, for some $L \in \mathcal{F}$ with $H \le L$, we calculate 
\[\Ind_{\pi \circ \sigma} \uR(I_H^L): 1 \otimes \id_{R[G/H]} \longmapsto 1 \otimes t_{L/H}\]
where $t_{L/H} \in \Hom_{RG}(R[G/L], R[G/H])$ denotes the map $L \mapsto \sum_{l \in L/H} lH $.  Following this element back down the chain of isomorphisms we get the element $[L \cdot 1 \cdot H]$, which by $(\star)$ is equal to $\vert L : H \vert  [H \cdot 1 \cdot H]$.  Thus $\Ind_{\pi \circ \sigma} \uR (I_H^L)$ acts as multiplication by $\vert L : H\vert $ on $R$, as required.

For any element $x \in G$, we calculate
\[\Ind_{\pi \circ \sigma} \uR(c_x): 1 \otimes \id_{R[G/H]} \longmapsto 1 \otimes \gamma_x\]
where $\gamma_x \in \Hom_{RG}(R[G/H^{x^{-1}}], R[G/H])$ is the map $H^{x^{-1}} \mapsto xH$.  Following this down the chain of isomorphisms we get the element $[H^{x^{-1}} \cdot x \cdot H]$, which by $(\star)$ is equal to $[ H^{x^{-1}} \cdot 1 \cdot H^{x^{-1}} ]$.  Thus $\Ind_{\pi \circ \sigma} \uR(c_x)$ acts as the identity on $R$, as required.
\end{proof}

The next proposition should be compared with Proposition \ref{prop:mackey indOF to MF take pr of R to pr of BG}.  Recall that a chain complex is $\mathcal{F}$-split if it splits when restricted to $R H$ for all $H \in \mathcal{F}$.

\begin{Prop}\label{prop:HF sigma ind preserves proj res of R}
 Induction with $\pi \circ \sigma$ takes projective $\OF$-module resolutions of $\uR$ to projective $\HeckeF$-module resolutions of $R^-$.
\end{Prop}
\begin{proof}
Let $P_*$ be a projective resolution of $\uR$ by $\OF$-modules, then by Lemma \ref{lemma:mackey tech lamma for ind},
\begin{align*}
\Ind_{\pi \circ \sigma} P_*(G/H) &=  P_* \otimes_{\OF} \Res_{\pi \circ \sigma} R[G/H, -]_{\HeckeF}  \\
&\cong \Hom_{RH}(R ,P_*(G/1)).
\end{align*}
So inducing $P_* \longtwoheadrightarrow \uR$ with ${\pi \circ \sigma} $ and using Lemma \ref{lemma:HF ind uROF is R} gives the chain complex
\[ \Ind_{\pi \circ \sigma} P_* \longtwoheadrightarrow R^-. \]
Induction preserves projectives, so we must show only that the above is exact.  Since induction is right exact, it is necessarily exact at degree $-1$ and degree $0$.
Evaluating at $G/H$ gives the resolution
\[
\Hom_{RH}(R ,P_*(G/1)) \longrightarrow R .
\]
By \cite[Theorem 3.2]{Nucinkis-EasyAlgebraicCharacterisationOfUniversalProperGSpaces}, the resolution $P_*(G/1)$ is $\mathcal{F}$-split.  Since $\Hom_{RH}(R, -)$ preserves the exactness of $RH$-split complexes, $\Hom_{RH}(R, P_*(G/1))$ is exact at position $i$ for all $i \ge 1$, completing the proof.
\end{proof}

\begin{Remark}
 The proposition above may not hold with $\uR$ replaced by an arbitrary $\OF$-module $M$, as a resolution of $M$ by projective $\OF$-modules will not in general split when evaluated at $G/1$.
\end{Remark}

Recall that in Section \ref{subsubsection:mackey cohom mackey homology and cohomology} we defined, for any $\HeckeF$-module $M$,
\[H^*_{\HeckeF}(G, M) \cong \Ext^*_{\HeckeF}(R^{-}, M). \]
There is an analogue of Corollary \ref{cor:mackey H^n_MF is H^n_OF}:
\begin{Prop}\label{prop:HF HF cohomology is OF cohomology}
 For any $\HeckeF$-module $M$ and any natural number $n$, 
\[ H^n_{\HeckeF} (G, M) = H^n_{\OF} (G, \Res_{\pi \circ \sigma} M). \]
\end{Prop}
\begin{proof}
 Let $P_*$ be a projective $\OF$-module resolution of $\uR$, then 
\begin{align*}
 H^n_{\OF} (G, \Res_{\pi \circ \sigma} M) &= H^n \Hom_{\OF} \left( P_*, \Res_{\pi \circ \sigma}  M \right) \\
 &\cong H^n \Hom_{\HeckeF} \left( \Ind_{\pi \circ \sigma} P_*, M \right)\\
 &= H^n_{\HeckeF} (G, M)
\end{align*}
where the isomorphism is adjoint isomorphism between induction and restriction and $\Ind_{\pi \circ \sigma} P_*$ is a projective $\HeckeF$-module resolution of $R^-$ by Proposition \ref{prop:HF sigma ind preserves proj res of R}.
\end{proof}

\section{\texorpdfstring{$\FPn$}{FPn} conditions for cohomological Mackey functors}\label{section:HF FPn}

The main result of this section is Theorem \ref{theorem:HFFPn iff FFPn} below.  For a detailed construction of $\mathcal{F}$-cohomology and the condition $\FFP_n$ see \cite{Nucinkis-CohomologyRelativeGSet}, for an overview see Section \ref{section:mackey introduction}.

\begin{Theorem}\label{theorem:HFFPn iff FFPn}
 For any ring $R$, if $G$ is $\HFFP_n$ then $G$ is $\FFP_n$.  If $R$ is Noetherian and $G$ is $\FFP_n$ then $G$ is $\HFFP_n$.
\end{Theorem}

The proof is contained in Sections \ref{subsection:HFFPn implies FFPn} and \ref{subsection:FFPn implies HFFPn}.  

\begin{Prop}\label{prop:MFFPn implies HFFPn}
 If $G$ is $\MFFP_n$ then $G$ is $\HFFP_n$.
\end{Prop}
\begin{proof}
 Combining Corollary \ref{cor:mackey H^n_MF is H^n_OF} and Proposition \ref{prop:HF HF cohomology is OF cohomology} shows that for all groups $G$, and all non-negative integers $i$,
 \[
  H^i_{\HeckeF}(G, -) \cong H^i_{\MF}(G, \Res_\pi -). 
 \]
 Let $G$ be a group of type $\MFFP_n$ and let $M_\lambda$, for $\lambda \in \Lambda$, be a directed system of $\HeckeF$-modules with colimit zero.  Then the colimit of $\Res_\pi M_\lambda$ is zero also and by the Bieri--Eckmann criterion (Theorem \ref{theorem:C bieri-eckmann criterion}), for any $i \le n$,
 \begin{align*}
  \varinjlim_\Lambda H^i_{\HeckeF}(G, M_\lambda) &\cong \varinjlim_{\Lambda} H^i_{\MF}(G, \Res_\pi M_\lambda) \\ 
  &\cong 0.
 \end{align*}
 Applying the Bieri--Eckmann criterion again shows $G$ is of type $\HFFP_n$.
\end{proof}

\begin{Prop}
 If $G$ is $\HFFP_n$ then $G$ is $\FP_n$.
\end{Prop}
\begin{proof}
Let $P_* \longtwoheadrightarrow R^-$ be a resolution of $R^-$ by free $\HeckeF$-modules, finitely generated up to degree $n$.  Since the finitely generated free $\HeckeF$-modules are fixed point functors of finitely generated permutation modules with stabilisers in $\mathcal{F}$, evaluating at $G/1$ gives a resolution of $R$ by $RG$-modules of type $\FP_\infty$ and a standard dimension shifting argument completes the proof.
\end{proof}

So there is a chain of implications:
\[ \OFFP_n \Leftrightarrow \MFFP_n \Rightarrow \HFFP_n \Rightarrow \FP_n + \left\{ \substack{ \text{$G$ has finitely many} \\\text{conjugacy classes of}\\\text{finite $p$-subgroups in $\mathcal{F}$ }}\right\} \]
Where the final implication is \cite[Proposition 4.2]{LearyNucinkis-GroupsActingPrimePowerOrder}, where it is proved that $G$ is $\FFP_0$ if and only if $G$ has finitely many conjugacy classes of finite $p$-subgroups in $\mathcal{F}$, for all primes $p$.  It is conjectured in the same paper that a group $G$ of type $\FP_\infty$ with finitely many conjugacy classes of finite $p$-subgroups in $\mathcal{F}$ is $\FFP_\infty$ \cite[Conjecture 4.3]{LearyNucinkis-GroupsActingPrimePowerOrder}.

Since $G$ is $\MFFP_0$ if and only if $G$ has finitely many conjugacy classes in $\mathcal{F}$ (Lemma \ref{lemma:mackey MFFP0}), the implication $\MFFP_n \Rightarrow \HFFP_n$ is not reversible.

There are examples due to Leary and Nucinkis of groups which act properly and cocompactly on contractible $G$-CW-complexes but which are not of type $\OFFP_0$ \cite[Example 3, p.149]{LearyNucinkis-SomeGroupsOfTypeVF}.  By Remark \ref{remark:HF proper cocompact action then HFFPinfty}, these groups are of type $\HFFP_\infty$ showing that $\HFFP_\infty \not\Rightarrow \OFFP_0$.  Leary and Nucinkis also give examples of groups which act properly and cocompactly on contractible $G$-CW-complexes, are of type $\OFFP_0$ but which are not $\OFFP_\infty$ \cite[Example 4, p.150]{LearyNucinkis-SomeGroupsOfTypeVF}.  Hence there can be no implication $\HFFP_n + \OFFP_0 \Rightarrow \OFFP_n$.

\subsection{\texorpdfstring{$\HFFP_n$ implies $\FFP_n$}{HF FPn implies FFPn}}\label{subsection:HFFPn implies FFPn}
This section comprises a series of lemmas, building to the proof of Proposition \ref{prop:HF HFFPn implies FFPn}, that for any commutative ring $R$ the condition $\HFFP_n$ implies the condition $\FFP_n$.  Throughout, $H$ and $K$ are arbitrary subgroups in $\mathcal{F}$.

  We say a short exact sequence of $RG$-modules
\begin{align*}
0 \longrightarrow A \longrightarrow B \longrightarrow C \longrightarrow 0 \tag{$\star$} 
\end{align*}
is \emph{$H$-good} if \index{H-good@$H$-good}
\[ 
0 \longrightarrow A^H \longrightarrow B^H \longrightarrow C^H \longrightarrow 0 
\]
is exact.  Similarly an exact chain complex $C_*$ is $H$-good if $C_*^H$ is exact.  If an exact chain complex is $H$-good for all $H \in \mathcal{F}$ we say it is  \emph{$\mathcal{F}$-good}.\index{F-good@$\mathcal{F}$-good}  Note that an $\mathcal{F}$-split exact chain complex is automatically $\mathcal{F}$-good.
\begin{Remark}\label{remark:fp res of R is Fgood}
 If $C_*^-$ is an exact chain complex of fixed point functors then $C_*$ is $\mathcal{F}$-good.
\end{Remark}

\begin{Remark}
In general being $H$-good is a weaker property than being $RH$-split:  Applying $\Hom_{RH}(R,-)$ to $(\star)$ gives 
\[ 0 \rightarrow \Hom_{RH}(R, A) \rightarrow \Hom_{RH}(R, B) \rightarrow \Hom_{RH}(R, C) \rightarrow H^1(H, A) \rightarrow \cdots \]
So to find an example of an $H$-good short exact sequence which is not $RH$-split it is sufficient to find modules $C$ and $A$ such that $H^1(H, A) = 0$ and $\Ext^1_{RH}(C, A) \neq 0$.  For example if $H$ is any finite group we may set $R = \ZZ$, $A = \ZZ H$ and $C = (\ZZ/2\ZZ) H$.
\end{Remark}

Additionally, we say that an $R H$-module $M$ has property $(P_H)$ if for any $\mathcal{F}$-good short exact sequence $(\star)$, $\Hom_{R H}(M, -)$ preserves the exactness of $(\star)$.\index{Property $(P_H)$}  Since $\Hom_{R H}(M, -)$ is always left exact, having $(P_H)$ is equivalent to asking that for any $\mathcal{F}$-good short exact sequence $(\star)$ and any $R H$-module homomorphism $f:M \to C$, there is a $R H$-module homomorphism $l:M \to B$ such that the diagram below commutes.
\[ \xymatrix{
 & & M \ar^f[rd] \ar_l@{-->}[d] & & \\
 0 \ar[r] & A \ar[r]&  B \ar^g[r] & C \ar[r] & 0
 } \]
Note that the trivial $RG$-module $R$ has property $(P_H)$.
 
\begin{Lemma}\label{lemma:FH direct summands of PH have PH}
 If $M$ has $(P_H)$ then any direct summand of $M$, as an $R H$-module, has $(P_H)$.
\end{Lemma}
\begin{proof}
 This is, with a minor alteration, the proof of \cite[Theorem 3.5(ii)]{Rotman-HomologicalAlgebra}.  Let $N$ be a direct summand of $M$ and consider the diagram with exact bottom row.  Assume the bottom row is $\mathcal{F}$-good.
 \[ \xymatrix{
 & & M \ar^\pi@/^/[r] \ar_l@{-->}[d] & N \ar^\iota@/^/[l] \ar^f[d] & \\
 0 \ar[r] & A \ar[r]&  B \ar^g[r] & C \ar[r] & 0
 } \]
 Here $f$ is some arbitrary homomorphism, and $\pi$ and $\iota$ are the projection and inclusion maps respectively.  Since $M$ has $P_H$, there is a map $l : M \to B$ such that $g \circ l = f \circ \pi$, the composition $l \circ \iota$ is the required map.
\end{proof}

\begin{Lemma}\label{lemma:HF perm module has PH}
For any $K \in \mathcal{F}$, the permutation module $R[G/K]$ has $(P_H)$.
\end{Lemma}
\begin{proof}
Let $L$ be any subgroup of $H$, then using the natural isomorphism
\[
 \Hom_{R H}(R[H/L], -) \cong \Hom_{R L}(R, -) 
\]
we see that $R[H/L]$ has $(P_H)$.  Now use \cite[Proof of \S III.9.5(ii) on p.83]{Brown} to rewrite $R[G/K]$ (as an $RH$-module), as
 \[ R[G/K] \cong \bigoplus_{g \in H \backslash G / K} R[H/K_g] \]
where $K_g = \{h \in H : g^{-1}hg \le K \}$.  Thus:
\[ \Hom_{R H}(R[G/K], -) \cong \prod_{g \in H \backslash G / K} \Hom_{RH} (R[H/K_g], -).\]
Now use that $R[H/L]$ has $(P_H)$ and that direct products of exact sequences are exact.
\end{proof}

\begin{Lemma}\label{lemma:HF PH then splits}
If
\[
  0 \longrightarrow A \longrightarrow B \longrightarrow C \longrightarrow 0 
\]
is an $H$-good short exact sequence and both $B$ and $C$ have $(P_H)$ then the short exact sequence is $H$-split and $A$ has $(P_H)$.
\end{Lemma}
\begin{proof}
 Apply $\Hom_{RH}(C, -)$ to see that the short exact sequence is $H$-split.  Then since, as $RH$-modules, $B$ is the direct sum of $C$ and $A$, $A$ necessarily has $(P_H)$ by Lemma \ref{lemma:FH direct summands of PH have PH}.
\end{proof}

\begin{Lemma}\label{lemma:HF Fgood is Fsplit}
 If $P_*$ is an $\mathcal{F}$-good resolution of $R$ by permutation $RG$-modules with stabilisers in $\mathcal{F}$, then $P_*$ is $\mathcal{F}$-split.
\end{Lemma}
\begin{proof}
Use induction with Lemmas \ref{lemma:HF perm module has PH} and \ref{lemma:HF PH then splits}.
\end{proof}

\begin{Remark}
 Similarly to Proposition \ref{prop:HF sigma ind preserves proj res of R}, the above lemma may fail for $\mathcal{F}$-good resolutions of arbitrary modules.
\end{Remark}

\begin{Prop}\label{prop:HF HFFPn implies FFPn}
 If $G$ is $\HFFP_n$ then $G$ is $\FFP_n$.
\end{Prop}
\begin{proof}
 Find a free $\HeckeF$-module resolution $P_*$ of $R^-$, finitely generated up to dimension $n$.  By Remark \ref{remark:fp res of R is Fgood}, $P_*(G/1)$ is an $\mathcal{F}$-good resolution of $R$ by permutation $RG$-modules with stabilisers in $\mathcal{F}$.  By Lemma \ref{lemma:HF Fgood is Fsplit} $P_*$ is $\mathcal{F}$-split, and by \cite[Definition 2.2]{Nucinkis-CohomologyRelativeGSet} permutation $RG$-modules with stabilisers in $\mathcal{F}$ are $\mathcal{F}$-projective.
\end{proof}

\subsection{\texorpdfstring{$\FFP_n$ implies $\HFFP_n$}{FFPn implies HF FPn}}\label{subsection:FFPn implies HFFPn}

This section comprises a series of lemmas, building to the proof of Proposition \ref{prop:FFP implies HFFPn}, that if $R$ is commutative Noetherian and $G$ is $\FFP_n$ then $G$ is $\HFFP_n$.

\begin{Lemma}
 For any $H \in \mathcal{F}$, inducing $R[G/H]$ to a \emph{covariant} $\HeckeF$-module gives the free module $R[G/H, -]_{\HeckeF}$.
\end{Lemma}
\begin{proof}
On objects the two functors are equal:
 \begin{align*}
\Ind^{\HeckeF}_{R G} R[G/H] (G/K) &= R[G/H] \otimes_{R G} R[G/1, -]_{\HeckeF} (G/K) \\  
&= R[G/H] \otimes_{RG} \Hom_{RG} (RG, R[G/K]) \\
&= R[G/H] \otimes_{RG} R[G/K] \\
&= R[H \backslash G / K] \\
&= \Hom_{RG} (R[G/H], R[G/K]).
 \end{align*}
If $L \le K$ are in $\mathcal{F}$, and $\sum_I g_i L \in R[G/L]^H$ then
\begin{align*}
R[G/1, -]_{\HeckeF} (R_L^K) : R[G/1, G/L]_{\HeckeF} &\longrightarrow R[G/1, G/K]_{\HeckeF} \\
 \sum_I g_i L &\longmapsto \sum_I g_i K.
\end{align*}
Following this down the chain of isomorphisms, then 
\begin{align*}
\Ind_{RG}^{\HeckeF} R[G/H] (R_L^K) : \Hom_{RG}(R[G/H], R[G/L]) &\longrightarrow \Hom_{RG}(R[G/H], R[G/K]) \\
\sum_I g_i L &\longmapsto \sum_I g_i K
\end{align*}
as required.  Similarly, if $\sum_I g_i K \in R[G/K]^H$ then 
\begin{align*}
R[G/1, -]_{\HeckeF} (I_L^K) :  R[G/1, G/K]_{\HeckeF} &\longrightarrow R[G/1, G/L]_{\HeckeF}   \\
  \sum_I g_i K &\longmapsto \sum_{k \in K / L} \sum_I g_i k  L .
\end{align*}
Following this down the chain of isomorphisms, 
\begin{align*}
 \Ind_{RG}^{\HeckeF} R[G/H] ( I_L^K ): \Hom_{RG} (R[G/H], R[G/K]) &\longrightarrow \Hom_{RG}(R[G/H], R[G/L])\\
\sum_I g_i K &\longmapsto \sum_{k \in K / L} \sum_I g_i k  L
\end{align*}
again as required.

The proof for the conjugation morphisms $c_g$ is similar to the above.
 \end{proof}

\begin{Lemma}\label{lemma:inducing prods of perms}
\[ \Ind_{RG}^{\HeckeF} \prod_{H \in \mathcal{F}/G} \prod_{\Lambda_H} R[G / H] = \prod_{H \in \mathcal{F}/G} \prod_{\Lambda_H} R[G/H, -]_{\HeckeF} \]
where for each $H \in \mathcal{F}/G$, $\Lambda_H$ is any indexing set and we are using covariant induction. 
\end{Lemma}
\begin{proof}
In this proof, we use $\prod$ as a shorthand for $\prod_{H \in \mathcal{F}/G} \prod_{\Lambda_H}$.  On objects, the two functors are equal:  
 \begin{align*}
\Ind^{\HeckeF}_{R G} \prod R[G/H] (G/K) &= \left(\prod R[G/H]\right) \otimes_{R G} R[G/1, -]_{\HeckeF} (G/K) \\  
&= \left(\prod R[G/H]\right) \otimes_{RG} \Hom_{RG} (RG, R[G/K]) \\
&= \left(\prod R[G/H]\right) \otimes_{RG} R[G/K] \\
&\cong \prod \left( R[G/H] \otimes_{RG} R[G/K] \right) \tag{$\star$} \\
&= \prod R[H \backslash G/K]\\
&= \prod \Hom_{RG}(R[G/H], R[G/K]).
 \end{align*}
 Where the isomorphism marked $(\star)$ is the Bieri--Eckmann criterion \cite[Theorem 1.3]{Bieri-HomDimOfDiscreteGroups}, which is valid because $R[G/K]$ is $\FP_\infty$.  That the morphisms are equal can be checked as in the previous lemma.
 \end{proof}

\begin{Lemma}\label{lemma:FHn with products of perm module coeff}
 \[ \mathcal{F}H_*\left(G, \prod_{H \in \mathcal{F}/G} \prod_{\Lambda_H} R[G/H]\right) = H_*^{\HeckeF}\left(G, \prod_{H \in \mathcal{F}/G} \prod_{\Lambda_H} R[G/H, -]_{\HeckeF}\right) \]
where for each $H \in \mathcal{F}/G$, $\Lambda_H$ is any indexing set.
\end{Lemma}
\begin{proof}
 Again we use $\prod$ to stand for  $\prod_{H \in \mathcal{F}/G} \prod_{\Lambda_H}$.
 Let $P_*$ be a free $\HeckeF$-module resolution of $R^-$, then $P_*(G/1)$ is an $\mathcal{F}$-split resolution of $R$ by $\mathcal{F}$-projective modules by Lemma \ref{lemma:HF Fgood is Fsplit}, so
 \begin{align*}
\mathcal{F}H_*(G, R[G/H]) &\cong H_*\left(P_*(G/1) \otimes_{R G} \prod R[G/H] \right) \\
&\cong H_* \left( P_* \otimes_{\HeckeF} \Ind^{\HeckeF}_{R G} \prod R[G/H] \right) \\
&\cong H_* (P_* \otimes_{\HeckeF} \prod R[G/H, -]_{\HeckeF} )\\
&\cong H_*^{\HeckeF} (G,\prod R[G/H, -]_{\HeckeF})
 \end{align*}
where the second isomorphism is the adjoint isomorphism between induction and restriction and the third is Lemma \ref{lemma:inducing prods of perms}
\end{proof}

\begin{Lemma}\label{lemma:FTor commutes if FFPn}
 For any group $G$, any commutative Noetherian ring $R$, any $RG$-module $A$ of type $\FFP_n$, and any exact limit, the natural map 
\[
\mathcal{F}\Tor_i^{RG} (A, \varprojlim_{\lambda \in \Lambda} M_\lambda ) \longrightarrow \varprojlim_{\lambda \in \Lambda} \mathcal{F}\Tor_i^{RG} (A,  M_\lambda ) 
\]
is an isomorphism for $i < n$ and an epimorphism for $i = n$.
\end{Lemma}
\begin{proof}

The proof is analogous to \cite[Theorem 1.3]{Bieri-HomDimOfDiscreteGroups} and \cite[Theorem 7.1]{Nucinkis-CohomologyRelativeGSet}, using \cite[Proposition 6.3]{Nucinkis-CohomologyRelativeGSet} which states that for $R$ commutative Noetherian, finitely generated $\mathcal{F}$-projective modules over $RG$ are of type $\FP_\infty$.

\end{proof}

Specialising the previous lemma to $M = R$:
\begin{Cor}\label{cor:FFP then FHn is cts}
 If $R$ is commutative Noetherian and $G$ is $\FFP_n$ over $R$, then for any exact limit, the natural map
\[
\mathcal{F}H_i (G, \varprojlim_{\lambda \in \Lambda} M_\lambda ) \longrightarrow \varprojlim_{\lambda \in \Lambda} \mathcal{F}H_i (G,  M_\lambda ) 
\]
 is an isomorphism for $i < n$ and an epimorphism for $i = n$.
\end{Cor}

\begin{Prop}\label{prop:FFP implies HFFPn}
 If $R$ is commutative Noetherian and $G$ is $\FFP_n$ over $R$ then $G$ is $\HFFP_n$ over $R$.
\end{Prop}
\begin{proof}
In this proof, we write $\prod$ for $\prod_{H \in \mathcal{F}/G} \prod_{\Lambda_H} $ where $\Lambda_H$ is any indexing set.  Using Lemmas \ref{lemma:FHn with products of perm module coeff} and \ref{cor:FFP then FHn is cts}, for any $i < n$:
 \begin{align*}
H_i^{\HeckeF}\left(G, \prod R[G/H, -]_{\HeckeF}\right) 
&= \mathcal{F}H_i\left(G, \prod R[G/H]\right) \\
&= \prod \mathcal{F}H_i\left(G, R[G/H]\right)\\
&= \prod H_i^{\HeckeF}\left(G, R[G/H, -]_{\HeckeF}\right).  \\
 \end{align*}
Thus $G$ is $\HFFP_n$ by the Bieri--Eckmann criterion (Theorem \ref{theorem:C bieri-eckmann criterion}).
 \end{proof}

\begin{Remark}
The requirement that $R$ be Noetherian was needed only for Lemma \ref{lemma:FTor commutes if FFPn}, where we need that finitely generated $\mathcal{F}$-projectives are $\FP_\infty$.  Nucinkis has given an example of a finitely generated $\mathcal{F}$-projective module which is not $\FP_\infty$ \cite[Remark on p.167]{Nucinkis-CohomologyRelativeGSet}, but the following question is still open. 
\end{Remark}

\begin{Question}
Does Proposition \ref{prop:FFP implies HFFPn} remain true if $R$ is not Noetherian? 
\end{Question}
\section{Cohomological dimension for cohomological Mackey functors}\label{section:HF cohomological dimension}

In \cite{Degrijse-ProperActionsAndMackeyFunctors}, Degrijse shows that for all groups $G$ with $\HFcd G < \infty$, 
\[
\Fcd G = \HFcd G.
\]
We can improve this.
\begin{Theorem}\label{theorem:Fcd=HFcd}
For all groups $G$,
 \[
 \Fcd G = \HFcd G.
 \]
\end{Theorem}
\begin{proof}
Remark \ref{remark:fp res of R is Fgood} and Lemma \ref{lemma:HF Fgood is Fsplit} imply $\Fcd G \le \HFcd G$.

For the opposite inequality, we first use \cite[Lemma 3.4]{Gandini-CohomologicalInvariants} which states that for a group $G$ with $\Fcd G \le n$ there is an $\mathcal{F}$-projective resolution $P_*$ of $R$ of length $n$, where each $P_i$ is a permutation module with stabilisers in $\mathcal{F}$.  Given such a $P_*$, we take fixed points of $P_*$ to get the $\HeckeF$ resolution $P_*^-$.  Since $P_*$ is $\mathcal{F}$-split, $P_*^-$ is exact.
\end{proof}

Recall that $\Fin$ denotes the family of finite subgroups of $G$ and $n_G$ denotes the minimal dimension of a proper contractible $G$-CW complex.

\begin{Prop}\label{prop:HF proper action on contr complex dim n then HFcd G less n}
For all groups $G$, 
\[
\mathcal{H}_{\Fin}\negthinspace\cd G \le n_G.
\]
\end{Prop}
This fact is well-known for $\Fincd$ instead of $\mathcal{H}_{\Fin}\negthinspace\cd$, but since a direct proof for $\mathcal{H}_{\Fin}\negthinspace\cd$ is both interesting and short we provide one.
\begin{proof}
 Let $P_*$ denote the cellular chain complex for a contractible $G$-CW-complex $X$ of dimension $n$ and take fixed points to get the complex $P_*^- \longrightarrow R^-$ of $\mathcal{H}_{\Fin}$-modules.  Since the action of $G$ on $X$ is proper the modules comprising $P_*$ are permutation modules with finite stabilisers and so $P_*^-$ is a chain complex of free $\mathcal{H}_{\Fin}$-modules.  By a result of Bouc \cite{Bouc-LeComplexeDeChainesDunGComplexe} and Kropholler--Wall \cite{KrophollerWall-GroupActionsOnAlgebraicCellComplexes} this chain complex splits when restricted to a complex of $RH$-modules for any finite subgroup $H$ of $G$.  In other words, $P_*$ is $\mathcal{F}$-good, thus $P_*^H \longrightarrow R$ is exact for any finite subgroup $H$ by Remark \ref{remark:fp res of R is Fgood}. 
\end{proof}

This leads naturally to the question:

\begin{Question}\label{question:HF HFcd finite iff n_G finite}
 Does $\mathcal{H}_{\Fin}\negthinspace\cd G < \infty$ imply $n_G < \infty$?
\end{Question}
We know of no group for which $n_G$ and $\HFcd G$ differ.
Brown has asked the following:
\begin{Question}\label{question:HF brown vcd = nG}\cite[VIII.11 p.226]{Brown} If $G$ is virtually torsion-free with $\vcd G < \infty$, then is $n_G = \vcd G$?
\end{Question}
If $G$ is virtually torsion free then $\vcd G = \mathcal{H}_{\Fin}\negthinspace\cd G$ \cite{MartinezPerezNucinkis-MackeyFunctorsForInfiniteGroups}, so a constructive answer to Question \ref{question:HF HFcd finite iff n_G finite} would give information about Question \ref{question:HF brown vcd = nG} as well. 

Related to this is the following question, posed using $\Fincd$ instead of $\mathcal{H}_{\Fin}\negthinspace\cd$ by Nucinkis.

\begin{Question}\label{question:mackey HFincd finite implies OFincd finite}\cite[p.337]{Nucinkis-EasyAlgebraicCharacterisationOfUniversalProperGSpaces}
 Does $\mathcal{H}_{\Fin}\negthinspace\cd G < \infty$ imply that $\mathcal{O}_{\Fin}\negthinspace\cd G < \infty$?
\end{Question}

\begin{Remark}\label{remark:HF proper cocompact action then HFFPinfty}
If $G$ acts properly and cocompactly on a finite-dimensional contractible $G$-CW-complex then, by a modification of the argument of the proof of Lemma \ref{prop:HF proper action on contr complex dim n then HFcd G less n}, $G$ is $\mathcal{H}_{\Fin}\negthinspace\FP_\infty$ also.  However, if $G$ acts properly on a finite type but infinite dimensional contractible complex $X$, then the theorem of Bouc and Kropholler--Wall doesn't apply, and the cellular chain complex of $X$ may not be $\mathfrak{F}$-split, thus we cannot deduce $G$ is $\mathcal{H}_{\Fin}\negthinspace\FP_\infty$.  

For an example of a group $G$ acting properly on a finite-type but infinite dimensional contractible CW-complex with a cellular chain complex which is not $\mathfrak{F}$-split, take the cyclic group $K \cong C_2$ acting antipodally on the infinite sphere $S_\infty$, with the usual CW structure of $2$ cells in each dimension.  One calculates that $C_*(S^\infty)^K$ is not exact and hence that $C_*(S_\infty)$ is not $\ZZ K$-split.
\end{Remark}

\begin{Question}\label{question:HF action on fin type implies HFFPinfty}
 If $G$ acts properly on a contractible $G$-CW-complex of finite type, but not necessarily finite dimension, then is $G$ of type $\mathcal{H}_{\Fin}\negthinspace\FP_\infty$?
\end{Question}

\subsection{Closure properties}\label{subsection:mackey HF closure properties}

The class of groups $G$ with $\HFcd G < \infty$ is closed under subgroups, free products with amalgamation, HNN extensions \cite[Corollary 2.7]{Nucinkis-EasyAlgebraicCharacterisationOfUniversalProperGSpaces}, direct products \cite[Corollary 3.9]{Gandini-CohomologicalInvariants} and extensions of finite groups by groups with $\HFcd$ finite \cite[Lemma 5.1]{Degrijse-ProperActionsAndMackeyFunctors}.

Section \ref{section:Group Extensions} contains a proof, via the Gorenstein cohomological dimension, that for a group extension
\[
 1 \longrightarrow N \longrightarrow G \longrightarrow Q \longrightarrow 1
\]
where $\HFincd G < \infty$, we have $\HFincd N + \HFincd Q \le \HFincd G $.

\begin{Prop}\cite[3.8,3.10]{Gandini-CohomologicalInvariants}
Let 
\[ 1 \longrightarrow N \longrightarrow G \longrightarrow Q \longrightarrow 1 \]
be a group extension such that for any finite extension $H$ of $N$ where $H/N$ has prime power order, $\HFincd H \le m$, then $\HFincd G \le n + m$.
\end{Prop}

\begin{Lemma}
 Let $N$ be any group and $p$ any prime.  If for any extension
\[ 1 \longrightarrow N \longrightarrow G \longrightarrow Q \longrightarrow 1 \]
we have that $\HFincd G = \HFincd N$ where $Q$ is the cyclic group of order $p$, then $\HFincd G = \HFincd N$, where $Q$ is any finite $p$-group.
\end{Lemma}

\begin{proof}
We prove by induction on the order of $Q$, the case $\lvert Q \rvert = p$ is by assumption.  Let $Q^\prime$ be a normal subgroup of index $p$ in $Q$ (such a subgroup exists by \cite[Theorem 4.6(ii)]{Rotman-Groups}) and consider the diagram below.
 \[ \xymatrix{
    1 \ar[r] & N \ar_=[d] \ar[r] & \pi^{-1}(Q^\prime) \ar[d] \ar^\pi[r] & Q^\prime \ar^{\trianglelefteq}[d] \ar[r] & 1 \\
    1 \ar[r] & N \ar[r] & G \ar^\pi[r] & Q \ar[r] & 1
   }
 \]
Since $Q^\prime$ is normal in $Q$, the preimage $\pi^{-1}(Q^\prime)$ is normal in $G$, with quotient group $G/\pi^{-1}(Q^\prime)$ of order $p$ so $\HFincd G = \HFincd \pi^{-1}(Q^\prime)$.  Finally by the induction assumption $\HFincd \pi^{-1}(Q^\prime) = \HFincd N$.  
\end{proof}

Combining the results above, if $\HFincd G$ fails to be subadditive there must exist a finite cyclic group $Q$, group $N$ with $\HFincd N < \infty$, and an extension $G$ of $Q$ by $N$ with $\HFincd G = \infty$.

\begin{Question}\label{question:mackey extensions by finite cyclic preserve HFcd}
 If $N$ is a group with $\HFincd N < \infty$ then does every extension $G$ of a cyclic group of order $p$ by $N$ satisfy $\HFincd G < \infty$?
\end{Question}

Any counterexample cannot be virtually torsion-free, since $\HFincd G = \vcd G$ for all virtually torsion-free groups \cite{MartinezPerezNucinkis-MackeyFunctorsForInfiniteGroups}, and neither can it be elementary amenable \cite[Proposition 3.13]{Gandini-CohomologicalInvariants}.

\section{The family of \texorpdfstring{\p}{p}-subgroups}\label{section:family of p subgroups}

Throughout this section $q$ is an arbitrary fixed prime and $R$ will denote one of the following rings:  the integers $\ZZ$, the finite field $\FF_q$, or the integers localised at $q$ denoted $\ZZ_{(q)}$.  If $R = \FF_q$ or $\ZZ_{(q)}$ then let $\mathcal{P}$ denote the subfamily of $\mathcal{F}$ consisting of all finite $q$-subgroups of groups in $\mathcal{F}$.  If $R = \ZZ$ then let $\mathcal{P}$ denote the subfamily of finite $p$-subgroups of groups in $\mathcal{F}$ for all primes $p$. 

We will always treat the cases $R = \FF_q$ and $R = \ZZ_{(q)}$ together, in fact the only property of these rings that we use is that for any integer $i$ coprime to $q$, the image of $i$ under the map $\ZZ \to R$ is invertible in $R$.  Hence the arguments in this section generalise to any other rings with this property, for example any ring with characteristic $q$.  The argument used for $R = \ZZ$ will go through for any ring $R$.

For $R = \ZZ$ and $\mathcal{F} = \Fin$, Leary and Nucinkis prove that the conditions $\FFP_n$ and $\mathcal{P}\negthinspace\FP_n$ are equivalent, and that $\mathcal{F}\negthinspace\cd G = \mathcal{P}\negthinspace\cd G $ \cite[Theorem 4.1]{LearyNucinkis-GroupsActingPrimePowerOrder}.  We use an averaging method similar to theirs to show that, for $R= \ZZ$, $\FF_q$, or $\ZZ_{(q)}$:

\begin{Theorem}\label{theorem:HF HeckeFcd = HeckePcd and HeckeFFPn = HeckePFPn}
  For $n \in \NN \cup \{ \infty \}$, the conditions $\HFcd G = n$ and $\HeckeP\negthinspace\cd G = n$ are equivalent, as are the conditions $\HFFP_n$ and $\HeckeP\negthinspace\FP_n$. 
\end{Theorem}

\begin{Remark}
 If $R = \ZZ_{(q)}$ or $\FF_q$ and all subgroups of $G$ have order coprime to $q$ then $\mathcal{P}$ contains only the trivial subgroup.  Thus $\HFcd_R G = \cd_R G$ and the conditions $\HFFP_n$ and $\FP_n$ are equivalent.
\end{Remark}

At the end of the section we will look at the case that $R = \KK$ is a field of characteristic zero, and prove that in this case $\HFcd G = \cd G$ and that the conditions $\HFFP_n$ and $\FP_n$ are equivalent.

The argument relies on two maps $\iota_H$ and $\rho_H$ defined for any subgroup $H$ in $\mathcal{F} \setminus \mathcal{P}$.  These maps have different definitions depending on the ring $R$.

We treat the case that $R = \FF_q$ or $R = \ZZ_{(q)}$ first.  Let $H \in \mathcal{F} \setminus \mathcal{P}$ and let $Q$ be a Sylow $q$-subgroup of $H$, define
\[
\rho_H = R^H_Q \in R[G/Q, G/H]_{\HeckeF}
\]
\[
\iota_H = (1 / \lvert H : Q \rvert) I^H_Q \in R[G/H, G/Q]_{\HeckeF}. 
\]
The map $\iota_H$ is well defined since $\lvert H : Q \rvert$ contains no powers of $q$ and hence is invertible in $R$.  

If $R = \ZZ$ and $H \in \mathcal{F} \setminus \mathcal{P}$ then let $\{P_i\}_{i \in I}$ run over the non-trivial Sylow $p$-subgroups of $H$ (choosing one subgroup for each $p$).  We necessarily have that
$ \gcd \{ \lvert H : P_i \rvert \: : \: i \in I \} = 1$ so, by B\'ezout's identity, we may choose integers $z_i$ such that $\sum_{i \in I} z_i\lvert H : P_i \rvert = 1  $.  Define, with a slight abuse of notation,
\[
\rho_H = \bigoplus_{i \in I} R^H_{P_i}.
\]
By which we mean that for any $\HeckeF$-module $M$,
\begin{align*}
  M(\rho_H) : M(G/H) &\longrightarrow \bigoplus_{i \in I} M(G/P_i) \\
  m &\longmapsto \bigoplus_{i \in I} M(R^H_{P_i})(m). \\
\end{align*}
With a similar abuse of notation we define
\[
\iota_H = \sum_{i \in I} z_i I^H_{P_i}.
\]
By which we mean that for any $\HeckeF$-module $M$,
\begin{align*}
  M(\iota_H) : \bigoplus_{i \in I} M(G/P_i) &\longrightarrow M(G/H) \\
   (m_i)_{i \in I} &\longmapsto \sum_{i \in I} z_i M(I^H_{P_i})(m_i). \\
\end{align*} 

The next couple of lemmas catalogue properties of the maps $\iota_H$ and $\rho_H$ which are needed for the proof of Theorem \ref{theorem:HF HeckeFcd = HeckePcd and HeckeFFPn = HeckePFPn}.

\begin{Lemma}\label{lemma:iota circ rho is id}
For any $\HeckeF$-module $M$ and subgroup $H \in \mathcal{F} \setminus \mathcal{P}$,
 \[
 M(\iota_H) \circ M(\rho_H) = \id_{M(G/H)}.
 \]
\end{Lemma}
\begin{proof}
 In the case $R = \FF_q$ or $R = \ZZ_{(q)}$, this follows from the fact that $M(R^H_Q \circ I^H_Q)$ is multiplication by $\lvert H : Q \rvert$.  For $R = \ZZ$,
\[
M(\iota_H) \circ M(\rho_H) = \sum_i z_i M( R^H_{P_i} \circ I^H_{P_i} ) = \sum_i z_i \lvert H : P_i \rvert = 1.
\]
\end{proof}

\begin{Lemma}\label{lemma:restriction of HeckeF proj is HeckePproj}
 If $H \in \mathcal{F}$ then $\Res^{\HeckeF}_{\HeckeP} R[-, G/H]_{\HeckeF}$ is a finitely generated projective $\HeckeP$-module.
\end{Lemma}
\begin{proof}
If $H$ is an element of $\mathcal{P}$ then this is obvious so assume that $H \not\in\mathcal{P}$.
 First, the case $R = \FF_q$ or $R = \ZZ_{(q)}$.  The projection
 \[
 s: R[-,G/Q]_{\HeckeF} \longtwoheadrightarrow R[-, G/H]_{\HeckeF} 
 \]
corresponding to $\iota_H$ under the Yoneda-type lemma (\ref{lemma:C yoneda-type}) is split by the map 
 \[ 
 i: R[-,G/H]_{\HeckeF} \longtwoheadrightarrow R[-, G/Q]_{\HeckeF} 
 \]
corresponding to $\rho_H$ under the Yoneda-type lemma:  It is sufficient to calculate 
\[
s \circ i(G/H)(\id_H) = \rho_H \circ \iota_H = \id_H.
\]
Applying $\Res_{\HeckeP}^{\HeckeF}$ gives a split surjection
\[
\Res_{\HeckeP}^{\HeckeF}s: \Res_{\HeckeP}^{\HeckeF} R[-,G/Q]_{\HeckeF} \longtwoheadrightarrow \Res_{\HeckeP}^{\HeckeF} R[-, G/H]_{\HeckeF}.
\]
Since $\Res_{\HeckeP}^{\HeckeF}R[-,G/Q]_{\HeckeF} = R[-, G/Q]_{\HeckeP}$ this completes the proof.

Now the case $R = \ZZ$, this time we construct a split surjection
\[ s: \bigoplus_{i \in I} R[-,G/P_i]_{\HeckeF} \longtwoheadrightarrow R[-, G/H]_{\HeckeF} \]
using the maps corresponding to $\iota_H$ and $\rho_H$ under the Yoneda-type lemma.  The rest of the proof is identical to the case $R = \FF_q$ or $R = \ZZ_{(q)}$.
\end{proof}

\begin{Lemma}\label{lemma:HF ch complex HFexact iff HPexact}
 A chain complex $C_*$ of $\HeckeF$-modules is exact if and only if it is exact at $G/P$ for all subgroups $P \in \mathcal{P}$.
\end{Lemma}
\begin{proof}
The ``only if'' direction is obvious so assume $C_*$ is a chain complex of $\HeckeF$-modules, exact at all $P \in \mathcal{P}$ and let $H \in \mathcal{F} \setminus \mathcal{P}$.  

We claim that the maps $C_*(\iota_H) $ and $C_*(\rho_H)$ are chain complex maps, we show this below for $R = \FF_q$ or $R = \ZZ_{(q)}$, the proof for $R= \ZZ$ is analogous.  The only non-obvious part of this claim is that the maps commute with the boundary maps $\partial_i$ of $C_*$, in other words the diagrams below commute:
\[
\xymatrix{
C_i(G/H) \ar^{\partial_i(G/H)}[r] \ar_{C_i(\rho_H)}[d] & C_{i-1}(G/H) \ar^{C_{i-1}(\rho_H)}[d] \\
C_i(G/Q) \ar^{\partial_i(G/Q)}[r] & C_{i-1}(G/Q)  \\
}
\]
\[
\xymatrix{
C_i(G/Q) \ar^{\partial_i(G/Q)}[r] \ar_{C_i(\iota_H)}[d] & C_{i-1}(G/Q) \ar^{C_{i-1}(\iota_H)}[d] \\
C_i(G/H) \ar^{\partial_i(G/H)}[r]  & C_{i-1}(G/H)  \\
}
\]
This follows from the fact that $\partial_i$ is an $\HeckeF$-module map.

Lemma \ref{lemma:iota circ rho is id} gives that $ C_*( \iota_H ) \circ C_*( \rho_H )$ is the identity on the chain complex $C_*(G/H)$.  The induced maps $\iota_H^*$ and $ \rho_H^*$ on homology satisfy 
\[\iota_H^* \circ  \rho_H^* = \id : H_*(C_*(G/H)) \longrightarrow H_*(C_*(G/H)) \]
so $\rho_H^*$ is injective.  The image of $\rho_H^*$ lies in $H_*(C_*(G/Q)) = 0$ if $R = \FF_q$ or $R = \ZZ_{(q)}$, or $\oplus_i H_*(C_*(G/P_i)) = 0$ if $R = \ZZ$, hence $H_*(C_*(G/H))$ is zero.
\end{proof}

\begin{Lemma}\label{lemma:HF proj HeckeP module extension}
If $P$ is a projective (respectively finitely generated projective) $\HeckeP$-module then there exists a $\HeckeF$-module $Q$ such that 
\[ \Res^{\HeckeF}_{\HeckeP} Q = P \]
and $Q$ is projective (resp. finitely generated projective).
\end{Lemma}
\begin{proof}
Recall from Lemma \ref{lemma:HF frees are fp} that the projective $\HeckeP$-modules are exactly those of the form $V^-$ for $V$ some direct summand of a permutation $RG$-module whose stabilisers lie in $\mathcal{P}$.
The required module is just $ V^-$ regarded as a $\HeckeF$-module.  
\end{proof}

\begin{proof}[Proof of Theorem \ref{theorem:HF HeckeFcd = HeckePcd and HeckeFFPn = HeckePFPn}]
 Assume that $\HFcd G \le n$ and let $P_*$ be a length $n$ projective resolution of $R^-$ by $\HeckeF$-modules, then restricting to the family $\mathcal{P}$ and using Lemma \ref{lemma:restriction of HeckeF proj is HeckePproj} gives a length $n$ projective resolution by $\HeckeP$-modules.  A similar argument shows that $\HFFP_n$ implies $\HeckeP\negthinspace\FP_n$.

For the converse, let $P_*$ be a length $n$ projective resolution of $R^-$ by $\HeckeP$-modules.  Lemma \ref{lemma:HF proj HeckeP module extension} gives projective $\HeckeF$-modules $Q_i$ such that $\Res^{\HeckeF}_{\HeckeP}Q_i = P_i$ for each $i$.  Denoting by $d_i$ the boundary maps in $P_*$, define boundary maps of $Q_*$ as
$\partial_i (G/P) = d_i(G/P)$ if $P \in \mathcal{P}$ and if $H \not\in \mathcal{P}$ then
\[
\partial_i(G/H) = P_{i-1}(\iota_H) \circ d_i(G/H) \circ P_i(\rho_H).
\]

One can check that these maps are indeed $\HeckeF$-module maps and that this makes $Q_*$ a chain complex:
\begin{align*}
& \partial_i (G/H) \circ \partial_{i+1} (G/H) \\
&= P_{i-1}(\iota_H) \circ d_i(G/H) \circ P_i(\rho_H) \circ P_{i}(\iota_H) \circ d_{i+1}(G/H) \circ P_{i+1}(\rho_H)  \\
&= P_{i-1}(\iota_H) \circ d_i(G/H) \circ d_{i+1}(G/H) \circ P_{i+1}(\rho_H) \\
&= 0.
\end{align*}
Finally $P_*$ is exact by Lemma \ref{lemma:HF ch complex HFexact iff HPexact}.

Since at all stages of the argument above finite generation is preserved, we get that $\HeckeP\negthinspace\FP_n $ implies $\HFFP_n$ too.
\end{proof}

For the remainder of this section $R = \KK$ is a field of characteristic zero, in this case we can reduce to the family $\Triv$ containing only the trivial subgroup.  For any $H \in \mathcal{F}$, let
\[
\rho_H = R^H_1
\]
\[
\iota_H = ( 1 / \lvert H \rvert ) I^H_1.
\]
All the arguments of the section go through with no alteration, showing:

\begin{Prop}\label{prop:HFcd and HFFPn for R a field} $\HFcd_{\KK} G = \mathcal{H}_{\Triv}\cd_{\KK} G$ 
and the conditions $\HFFP_n$ over $\KK$ and $\mathcal{H}_{\Triv}\FP_n$ over $\KK$ are equivalent for any $n \in \NN \cup \{ \infty\}$. 
\end{Prop}

\begin{Cor} $\HFcd_{\KK} G = \cd_{\KK} G $ and the conditions $\HFFP_n$ over $\KK$ and $\FP_n$ over $\KK G$ are equivalent for any $n \in \NN \cup \{ \infty\}$. 
\end{Cor}
\begin{proof}
 The category of $\mathcal{H}_{\Triv}$-modules is isomorphic to the category of $\KK G$-modules.
\end{proof}

\subsection{\texorpdfstring{$\FPn$}{FPn} conditions over \texorpdfstring{$\FF_p$}{Fp}}\label{subsection:FPn over FFp}

Throughout this section, we fix a prime $p$ and work over the ring $\FF_p$ with the family $\mathcal{P}$ of all finite $p$-subgroups of groups in $\mathcal{F}$.

\begin{Lemma}\cite[Lemma 5.3]{HambletonPamukYalcin-EquivariantCWComplexesAndTheOrbitCategory}  For any finite subgroup $H \in \mathcal{P}$ and $\HeckeP$-module $M$, $D_HM$ extends to a cohomological Mackey functor.
\end{Lemma}

\begin{Prop}\label{prop:mackey HPFPn iff OPFPn}
 $G$ is $\HeckeP\negthinspace\FP_n$ over $\FF_p$ if and only if $G$ is $\OP\negthinspace\FP_n$ over $\FF_p$.
\end{Prop}
The proof is basically that of Proposition \ref{prop:MFFPn implies OFFPn} combined with the lemma above.
\begin{proof}
We know already from Proposition \ref{prop:MFFPn implies HFFPn} and Corollary \ref{cor:OFFPn iff MFFPn} that $\OP\negthinspace\FP_n$ implies $\HeckeP\negthinspace\FP_n$.

Let $M_\lambda$, for $\lambda \in \Lambda$, be a directed system of $\OP$-modules with colimit zero.  Using the notation of Proposition \ref{prop:MFFPn implies OFFPn} there is an exact sequence of directed systems for each $i \ge 0$,
\[ 0 \longrightarrow C^iM_\lambda \longrightarrow DC^iM_\lambda \longrightarrow C^{i+1}M_\lambda \longrightarrow 0, \]
each of which has colimit zero.  Moreover, $DC^iM_\lambda$ extends to a cohomological Mackey functor so using the Bieri--Eckmann criterion (Theorem \ref{theorem:C bieri-eckmann criterion}), if $m \le n$ then for all $i \ge 0$,
\[
\varinjlim_{\Lambda} H^m_{\OP} (G, DC^iM_\lambda) = 0.
\]
Thus,
\begin{align*}
 \varinjlim H^m_{\OP}(G, M_\lambda) &= \varinjlim  H^m_{\OP}(G, C^0M_\lambda) \\
&= \varinjlim  H^{m-1}_{\OP}(G, C^1M_\lambda) \\
&= \cdots \\
&=  \varinjlim H^0_{\OP}(G, C^mM_\lambda) \\
&= 0.
\end{align*}
Where the final zero is because $G$ is $\OP\negthinspace\FP_0$ (by \cite[Proposition 4.2]{LearyNucinkis-GroupsActingPrimePowerOrder} and Theorem \ref{theorem:HFFPn iff FFPn}).
\end{proof}

\begin{Cor}\label{cor:mackey HFFPn complete description}
$G$ is $\HFFP_n$ over $\FF_p$ if and only if $\mathcal{P}$ contains finitely many conjugacy classes, and $WH$ is $\FP_n$ over $\FF_p$ for all $H \in \mathcal{P}$.
\end{Cor}
\begin{proof}
$G$ is $\HFFP_n$ if and only if $G$ is $\OP\negthinspace\FP_n$ by Theorem \ref{theorem:HF HeckeFcd = HeckePcd and HeckeFFPn = HeckePFPn} and Proposition \ref{prop:mackey HPFPn iff OPFPn}.  Now use that $G$ is $\OP\negthinspace\FP_n$ if and only if $\mathcal{P}$ contains finitely many conjugacy classes, and $WH$ is $\FP_n$ for all $H \in \mathcal{P}$ (Corollary \ref{cor:uFPn equivalent conditions}).
\end{proof}

\begin{Prop}
 If $G$ is virtually torsion-free then the conditions virtually $\FP$ over $\FF_p$ and $\HeckeF\FP$ over $\FF_p$ are equivalent.
\end{Prop}
\begin{proof}
 If $G$ is virtually $\FP$ over $\FF_p$ then $G$ has finitely many conjugacy classes of finite $p$-subgroups \cite[IX.(13.2)]{Brown}. A result of Hamilton gives that for any finite $p$-subgroup $H$ of $G$, $WH$ is virtually $\FP$ over $\FF_p$, in particular $WH$ is $\FP_\infty$ over $\FF_p$ \cite[Theorem 7]{Hamilton-WhenIsGroupCohomologyFinitary}.  Finally, \cite[Proposition 34]{Hamilton-WhenIsGroupCohomologyFinitary} gives that $G$ acts properly on a finite-dimensional $\FF_p$-acyclic space, thus in particular $\HFcd_{\FF_p} G < \infty$.  The other direction is obvious.
\end{proof}

In \cite{LearyNucinkis-GroupsActingPrimePowerOrder} it is conjectured that, if $\mathcal{F} = \Fin$, $G$ is $\FFP_\infty$ if and only if $G$ is $\FP_\infty$ and has finitely many conjugacy classes of finite $p$-subgroups for all primes $p$.  One could generalise this and ask:

\begin{Question}\label{question:HFFPn iff HFFP0 and FPn}
 Let $\mathcal{F} = \Fin$ and $n \in \NN \cup \{ \infty \}$.  
 \begin{enumerate}
  \item If $G$ is $\FP_n$ over $\ZZ$ with finitely many conjugacy classes of finite $p$-subgroups for all primes $p$, then is $G$ of type $\HFFP_n$ over $\ZZ$?
  \item Fixing a prime $p$, if $G$ is $\FP_n$ over $\FF_p$ with finitely many conjugacy classes of finite $p$-subgroups then is $G$ of type $\HFFP_n$ over $\FF_p$?
 \end{enumerate}
\end{Question}

A problem with finding a counterexample to Question \ref{question:HFFPn iff HFFP0 and FPn}(2) is that if $G$ admits a cocompact action on a finite-dimensional $\FF_p$-acyclic space $X$ then, via Smith theory, $X^P$ is $\FF_p$-acyclic for any finite $p$-subgroup $P$ and thus $WP$ is $\FP_n$ over $\FF_p$.  For this reason one cannot use the examples of Leary and Nucinkis in \cite{LearyNucinkis-SomeGroupsOfTypeVF}, their construction requires actions of finite groups on finite dimensional $\FF_p$-acyclic flag complexes with fixed point sets that are not $\FF_p$-acyclic.

\chapter{Gorenstein cohomology and \texorpdfstring{$\mathfrak{F}$}{F}-cohomology}\label{chapter:G}
This chapter contains material that has appeared in:
\begin{itemize}
 \item On the Gorenstein and $\mathfrak{F}$-cohomological dimensions (2013, to appear Bull. Lond. Math. Soc.) \cite{Me-GorensteinAndF}.
\end{itemize}

We study the Gorenstein cohomological dimension $\Gcd G$ and prove the following result.

\theoremstyle{plain}\newtheorem*{CustomThmH*}{Theorem \ref{thm:Fcd finite then Fcd=Gcd}}
\begin{CustomThmH*}
If $\Fincd G < \infty$ then $\Fincd G = \Gcd G$.
\end{CustomThmH*}

The proof is via the construction in Theorem \ref{theorem:FAM LES} of a long exact sequence relating the $\mathfrak{F}$-cohomology, the complete $\mathfrak{F}$-cohomology, and a new cohomology theory we call the $\mathfrak{F}_G$-cohomology.  The construction is analogous to the construction of the long exact sequence of Avramov--Martsinkovsky  relating the group cohomology, complete cohomology, and Gorenstein cohomology \cite[\S 7]{AvramovMartsinkovsky-AbsoluteRelativeAndTateCohomology-FiniteGorenstein}\cite[Theorem 3.11]{AsadollahiBahlekehSalarian-HierachyCohomologicalDimensions}.

In Section \ref{section:Group Extensions} we use Theorem \ref{thm:Fcd finite then Fcd=Gcd} and subadditivity of the Gorenstein cohomological dimension to study the behaviour of the $\mathfrak{F}$-cohomological dimension under group extensions.

\theoremstyle{plain}\newtheorem*{CustomThmI*}{Corollary \ref{cor:Fcd subadditive if finite}}
\begin{CustomThmI*}
Given a short exact sequence of groups
\[ 1 \longrightarrow N \longrightarrow G \longrightarrow Q \longrightarrow 1,  \]
if $\Fincd G < \infty$ then $ \Fincd G \le \Fincd N + \Fincd Q $.
\end{CustomThmI*}

Finally, in Section \ref{section:cdQ} we use the Avramov--Martsinkovsky long exact sequence  to prove the following.

\theoremstyle{plain}\newtheorem*{CustomThmJ*}{Proposition \ref{prop:cdQG finite then cdQG le GcdG}}
\begin{CustomThmJ*}
 If $\Gcd G < \infty$ and $\cd_{\QQ} G < \infty$ then $\cd_{\QQ} G \le \Gcd G$.
\end{CustomThmJ*}

\section{Preliminaries}\label{section:preliminaries}

\subsection{Complete resolutions and complete cohomology}\label{subsection:complete resolutions}

A \emph{weak complete resolution} of a module $M$ is an acyclic complex $T_*$ of projective modules which coincides with an ordinary projective resolution $P_*$ of $M$ in sufficiently high degree.  The degree in which the two coincide is called the \emph{coincidence index}.  A weak complete resolution is called a \emph{strong complete resolution} if $\Hom_{RG} (T_*, Q)$ is acyclic for every projective module $Q$.  We avoid the term ``complete resolution'' since some authors use it to refer to a weak complete resolution and others to a strong complete resolution.  \index{Weak complete resolution}\index{Strong complete resolution}

\begin{Prop}\cite[Proposition 2.8]{AsadollahiBahlekehSalarian-HierachyCohomologicalDimensions}
 A group $G$ admits a strong complete resolution if and only if $\Gcd G < \infty$.
\end{Prop}

The advantage of strong complete resolutions is that given strong complete resolutions $T_*$ and $S_*$ of modules $M$ and $N$, any module homomorphism $M \to N$ lifts to a morphism of strong complete resolutions $T_* \to S_*$ \cite[Lemma 2.4]{CornickKropholler-OnCompleteResolutions}.  Thus they can be used to define a cohomology theory: given a strong complete resolution $T_*$ of $M$ we define \index{Complete Ext functor $\widehat{\Ext}^*_{RG}(M,-)$}
\[
\widehat{\Ext}^*_{RG}(M,-) \cong H^*\Hom_{RG}(T_*, -).
\]
We also set $\widehat{H}^*(G, -) = \widehat{\Ext}^*_{RG}(R, -)$. \index{Complete cohomology $\widehat{H}^*(G, -)$}  This coincides with the complete cohomology of Mislin \cite{Mislin-TateViaSatellites}, Vogel \cite{Goichot-HomologieDeTateVogel}, and Benson--Carlson \cite{BensonCarlson-ProductsInNegativeCohomology} (see \cite[Theorem 1.2]{CornickKropholler-OnCompleteResolutions} for a proof).  Recall that the complete cohomology is itself a generalisation of the Farrell--Tate Cohomology, defined only for groups with finite virtual cohomological dimension \cite[\S X]{Brown}.

Even weak complete resolutions do not always exist, for example a free Abelian group of infinite rank cannot admit a weak complete resolution \cite[Corollary 2.10]{MislinTalelli-FreelyProperlyOnFDHomotopySpheres}.  It is conjectured by Dembegioti and Talelli that a $\ZZ G$-module admits a weak complete resolution if and only if it admits a strong complete resolution \cite[Conjecture B]{DembegiotiTalelli-ANoteOnCompleteResolutions}.

\subsection{\texorpdfstring{$\mathfrak{F}$}{F}-cohomology}
This section contains two technical lemmas we will need later.  

If $M$ is any $RG$-module and $F_i = R \Delta^i$ is the standard $\mathfrak{F}$-split resolution of $R$ \cite[p.342]{Nucinkis-EasyAlgebraicCharacterisationOfUniversalProperGSpaces}, then $F_* \otimes_{R} M$ is an $\mathfrak{F}$-split $\mathfrak{F}$-projective resolution of $M$.  Thus we've shown:

\begin{Lemma}\label{lemma:Fsplit resolutions exist}
 $\mathfrak{F}$-split $\mathfrak{F}$-projective resolutions exist for all $RG$-modules $M$.
\end{Lemma}

There is also a version of the Horseshoe lemma.
\begin{Lemma}[Horseshoe lemma]\label{lemma:Fcohomology horseshoe}
 If 
\[0 \longrightarrow A \longrightarrow B \longrightarrow C \longrightarrow 0 \]
is an $\mathfrak{F}$-split short exact sequence and $P_*$ and $Q_*$ are $\mathfrak{F}$-split $\mathfrak{F}$-projective resolutions of $A$ and $C$ respectively then there is an $\mathfrak{F}$-split $\mathfrak{F}$-projective resolution $S_*$ of $B$ such that $S_i = P_i \oplus Q_i$ and there is an $\mathfrak{F}$-split short exact sequence of augmented complexes 
\[
0 \longrightarrow \tilde P_* \longrightarrow \tilde S_* \longrightarrow \tilde Q_* \longrightarrow 0.
\]
\end{Lemma}
The proof is simlar to \cite[Lemma 8.2.1]{EnochsJenda-RelativeHomologicalAlgebra1}.
\begin{proof}
 First build the diagram below as in, for example, \cite[Proposition 6.24]{Rotman-HomologicalAlgebra} where it is shown to commute and have exact rows and columns.  Here, $K_A$, $K_B$ and $K_C$ are the kernels of the maps $P_0 \longrightarrow A$, $P_0 \oplus Q_0 \longrightarrow B$, and $Q_0 \longrightarrow C$ respectively.
 \[
  \xymatrix{
  & 0 \ar[d] & 0 \ar[d] & 0 \ar[d] & \\
  0 \ar[r] & K_A \ar[r] \ar[d] & K_B \ar[d] \ar[r] & K_C \ar[r] \ar[d] & 0 \\
  0 \ar[r] & P_0 \ar[r] \ar[d] & P_0 \oplus Q_0 \ar[d] \ar[r] & Q_0 \ar[r] \ar[d] & 0 \\
  0 \ar[r] & A \ar[r] \ar[d] & B \ar[d] \ar[r] & C \ar[r] \ar[d] & 0 \\
  & 0 & 0 & 0 & 
  }
 \]

 Since $P_0$ and $Q_0$ are both $\mathfrak{F}$-projective, $P_0 \oplus Q_0$ is $\mathfrak{F}$-projective, and since the middle row is split, it is $\mathfrak{F}$-split.
 
 Let $\Delta$ be the $G$-set $\coprod_{H \in \Fin}G/H$ and apply $-\otimes R\Delta$ to the commutative diagram to obtain a new commutative diagram with exact left column, right column, bottom row, and central row. 
 \[
 (P_0 \oplus Q_0) \otimes R\Delta \longrightarrow  B \otimes R\Delta
 \]
is surjective is because the tensor product is right exact, and an application of the 5-Lemma \cite[Lemma 2.72]{Rotman-HomologicalAlgebra} shows
\[
K_B \otimes R\Delta \longrightarrow (P_0 \oplus Q_0)\otimes R\Delta
 \]
is injective.  Hence the central column of our new commutative diagram is exact.  The $3\times 3$-Lemma provides that the top row is exact too \cite[Ex 2.32]{Rotman-HomologicalAlgebra}, thus all rows and columns of the first commutative diagram are $\mathfrak{F}$-split.  

Now repeat this process, but starting with the $\mathfrak{F}$-split short exact sequence
\[
 0 \longrightarrow K_A \longrightarrow K_B \longrightarrow K_C \longrightarrow 0.
\]
\end{proof}

\subsection{Complete \texorpdfstring{$\mathfrak{F}$}{F}-cohomology}\label{subsection:complete F cohomology}

In \cite{Nucinkis-CohomologyRelativeGSet}, Nucinkis constructs a complete $\mathfrak{F}$-cohomology, we give a brief outline here.  An \emph{$\mathfrak{F}$-complete resolution} $T_*$ of $M$ is an acyclic $\mathfrak{F}$-split complex of $\mathfrak{F}$-projectives which coincides with an $\mathfrak{F}$-split $\mathfrak{F}$-projective resolution of $M$ in high enough dimensions.  \index{F-complete resolution@$\mathfrak{F}$-complete resolution}
An \emph{$\mathfrak{F}$-strong} $\mathfrak{F}$-complete resolution $T_*$ has $\Hom_{RG}(T_*, Q)$ exact for all $\mathfrak{F}$-projectives $Q$.  \index{F-strong resolution@$\mathfrak{F}$-strong resolution}
Given such a $T_*$ we define \index{F-complete Ext@$\mathfrak{F}$-complete Ext functor $\widehat{\mathfrak{F}\negthinspace\Ext}_{RG}^*(M,-)$}\index{F-complete cohomology@$\mathfrak{F}$-complete cohomology $\FhatH^*(G, -)$}
\[ \widehat{\mathfrak{F}\negthinspace\Ext}_{RG}^*(M,-) = H^* \Hom_{RG}(T_*, -) \]
\[ \FhatH^*(G, -) = \widehat{\mathfrak{F}\negthinspace\Ext}_{RG}^*(R,-). \]
Nucinkis also describes a Mislin style construction and a Benson--Carlson construction of complete $\mathfrak{F}$-cohomology defined for all groups, proves they are equivalent, and proves that whenever there exists an $\mathfrak{F}$-complete resolution they agree with the definition above.

\subsection{Gorenstein cohomology}\label{subsection:Gorenstein cohomology}
The Gorenstein cohomology is, like the $\mathfrak{F}$-cohomology, a special case of the relative homology of Mac Lane \cite[\S IX]{MacLane-Homology} and Eilenberg--Moore \cite{EilenbergMoore-FoundationsOfRelativeHomologicalAlgebra}.

Recall that a module is Gorenstein projective if it is a cokernel in a strong complete resolution.  
An acyclic complex $C_*$ of Gorenstein projective modules is \emph{G-proper} if $\Hom_{R G}(Q, C_*)$ is exact for every Gorenstein projective $Q$.\index{G-proper chain complex}  The class of G-proper short exact sequences is allowable in the sense of Mac Lane \cite[\S IX.4]{MacLane-Homology}.  The projectives objects with respect to G-proper short exact sequences are exactly the Gorenstein projectives (for the definition of a projective object with respect to a class of short exact sequences see \cite[p.261]{MacLane-Homology}).  For $M$ and $N$ any $RG$-modules, we define \index{Gorenstein Ext functor $\GExt_{RG}^*(M, -)$}\index{Gorenstein cohomology $\GH^*(G, -)$}
\[ \GExt_{RG}^*(M, N) = H^*\Hom_{RG} (P_*, N) \]
\[ \GH^*(G, N) = \GExt_{RG}^*(R, N) \]
where $P_*$ is a G-proper resolution of $M$ by Gorenstein projectives.  

The usual method of producing a ``Gorenstein projective dimension'' of a module $M$ in this setting would be to look at the shortest length of a G-proper resolution of $M$ by Gorenstein projectives.  A priori this could be larger than the Gorenstein projective dimension defined in the introduction, where the G-proper condition is not required.  Fortunately there is the following theorem of Holm:

\begin{Theorem}\label{theorem:Holms theorem on proper vs non proper res}\cite[Theorem 2.10]{Holm-GorensteinHomologicalDimensions}
 If $M$ has finite Gorenstein projective dimension then $M$ admits a G-proper Gorenstein projective resolution of length $\Gpd M$.
\end{Theorem}

Generalising an argument of Avramov and Martsinkovsky in \cite[\S 7]{AvramovMartsinkovsky-AbsoluteRelativeAndTateCohomology-FiniteGorenstein} Asadollahi, Bahlekeh, and Salarian construct a long exact sequence:

\begin{Theorem}[Avramov--Martsinkovsky long exact sequence]\label{theorem:AM LES}\cite[Theorem 3.11]{AsadollahiBahlekehSalarian-HierachyCohomologicalDimensions}\index{Avramov--Martsinkovsky long exact sequence}
For a group $G$ with $\Gcd G < \infty$, there is a long exact sequence of cohomology functors
\[ 0 \longrightarrow \GH^1(G, -) \longrightarrow H^1(G, -) \longrightarrow \cdots \]
\[ \cdots \longrightarrow \GH^n(G, -) \longrightarrow H^n(G, -) \longrightarrow \widehat{H}^n(G, -) \longrightarrow \GH^{n+1}(G, -) \longrightarrow \cdots\]
\end{Theorem}

The construction relies on the complete cohomology being calculable via a complete resolution, hence the requirement that $\Gcd G < \infty$.

We will need the following lemma later:

\begin{Lemma}\label{lemma:Gproper res of R is Fsplit}
 Any G-proper resolution of $R$ is $\mathfrak{F}$-split.
\end{Lemma}
\begin{proof}
 If $P_*$ is a G-proper resolution of $R$ then since $R[G/H]$ is a Gorenstein projective \cite[Lemma 2.21]{AsadollahiBahlekehSalarian-HierachyCohomologicalDimensions},
\[ \Hom_{RG}(R[G/H], P_*) \cong \Hom_{RH}(R, P_*) \cong P_*^H \]
 is exact, thus by the argument of Proposition \ref{prop:HF HFFPn implies FFPn} $P_*$ is $\mathfrak{F}$-split.
\end{proof}

\section{\texorpdfstring{$\mathfrak{F}_G$}{F\_G}-cohomology}\label{section:FGcohomology}

\subsection{Construction}

We define another special case of relative homology, which we call the $\mathfrak{F}_G$-cohomology.  It enables us to build an Avramov--Martsinkovsky long exact sequence of cohomology functors containing $\FinH^*(G, -)$ and $\FhatH^*(G, -)$.  

We define an \emph{$\mathfrak{F}_G$-projective} to be the cokernel in a $\mathfrak{F}$-complete $\mathfrak{F}$-strong resolution and say a complex $C_*$ of $RG$-modules is \emph{$\mathfrak{F}_G$-proper} if $\Hom_{RG}(Q, C_*)$ is exact for any $\mathfrak{F}_G$-projective $Q$.\index{FG-projective@$\mathfrak{F}_G$-projective}\index{FG proper chain@$\mathfrak{F}_G$-proper chain complex}  The $\mathfrak{F}_G$-proper short exact sequences form an allowable class in the sense of Mac Lane, whose projective objects are the $\mathfrak{F}_G$-projectives---to check the class of $\mathfrak{F}_G$-proper short exact sequences is allowable we need only check that given a $\mathfrak{F}_G$-proper short exact sequence, any isomorphic short exact sequence is $\mathfrak{F}_G$-proper and that for any $RG$-module $A$ the short exact sequences
\[ 0 \longrightarrow A \stackrel{\id}{\longrightarrow} A \longrightarrow 0 \longrightarrow 0 \]
and
\[ 0 \longrightarrow 0 \longrightarrow A \stackrel{\id}{\longrightarrow} A \longrightarrow 0 \]
are $\mathfrak{F}_G$-proper.

We don't know if the class of $\mathfrak{F}_G$-projectives is precovering (see \cite[\S 8]{EnochsJenda-RelativeHomologicalAlgebra1}), so we don't know if there always exists an $\mathfrak{F}_G$-proper $\mathfrak{F}_G$-projective resolution.  However, if $A$ and $B$ admit $\mathfrak{F}_G$-proper $\mathfrak{F}_G$-resolutions $P_*$ and $Q_*$, respectively, then any map $A \longrightarrow B$ induces a map of resolutions $P_* \longrightarrow Q_*$ which is unique up to chain homotopy equivalence \cite[IX.4.3]{MacLane-Homology} and we have a slightly weaker form of the Horseshoe lemma.

\begin{Lemma}[Horseshoe lemma]\label{lemma:weak horseshoe for FG}
 Suppose 
 \[ 0 \longrightarrow A \longrightarrow B \longrightarrow C \longrightarrow 0 \]
 is a $\mathfrak{F}_G$-proper short exact sequence of $RG$-modules and both $A$ and $C$ admit $\mathfrak{F}_G$-proper $\mathfrak{F}_G$-projective resolutions $P_*$ and $Q_*$ then there is an $\mathfrak{F}_G$-proper resolution $S_*$ of $B$ such that $S_i = P_i \oplus Q_i$ and there is an $\mathfrak{F}_G$-proper short exact sequence of augmented complexes 
\[0 
\longrightarrow \tilde P_* \longrightarrow \tilde S_* \longrightarrow \tilde Q_* \longrightarrow 0.
\]
\end{Lemma}
The proof is similar to that of \cite[8.2.1]{EnochsJenda-RelativeHomologicalAlgebra1} and Lemma \ref{lemma:Fcohomology horseshoe}.
\begin{proof}
 First build the same commutative diagram as in the proof of Lemma \ref{lemma:Fcohomology horseshoe}.  Since $P_0$ and $Q_0$ are both $\mathfrak{F}_G$-projective, $P_0 \oplus Q_0$ is $\mathfrak{F}_G$-projective, and since the middle row is split, it is $\mathfrak{F}_G$-proper.
 
 Let $T$ be an $\mathfrak{F}_G$-projective and apply $\Hom_{RG}(T, -)$ to obtain a new commutative diagram with exact left column, right column, bottom row, and central row.  That 
 \[
\Hom_{RG}(T, P_0 \oplus Q_0) \longrightarrow \Hom_{RG}(T, B)
 \]
is surjective is an application of the 5-Lemma \cite[Lemma 2.72]{Rotman-HomologicalAlgebra}, and another application of the same lemma shows
\[
\Hom_{RG}(T,K_B) \longrightarrow \Hom_{RG}(T, P_0 \oplus Q_0)
 \]
is injective.  Hence the central column of our new commutative diagram is exact, and an application of the $3\times 3$-Lemma shows the top row is exact \cite[Ex 2.32]{Rotman-HomologicalAlgebra}, thus the original commutative diagram is $\mathfrak{F}$-split.  The rest of the proof is the same as that of Lemma \ref{lemma:Fcohomology horseshoe}.
\end{proof}

For any module $M$ which admits an $\mathfrak{F}_G$-proper resolution $P_*$ by $\mathfrak{F}_G$-projectives we define \index{FG Ext@$\mathfrak{F}_G$ Ext functor $\FGExt^*_{RG}(M, -)$}\index{FG cohomology@$\mathfrak{F}_G$ cohomology $\FGH^*(G, -)$}
\[
\FGExt^*_{RG}(M, N) = H^*\Hom_{RG}(P_*, N).
\]
We define also
\[
\FGH^*(G, -) = \FGExt^*_{RG} (R, -) .
\]

The next lemma follows from Lemma \ref{lemma:weak horseshoe for FG}, see \cite[8.2.3]{EnochsJenda-RelativeHomologicalAlgebra1}.
\begin{Lemma}\label{lemma:long exact sequences in FGExt}
Suppose 
 \[ 0 \longrightarrow A \longrightarrow B \longrightarrow C \longrightarrow 0 \]
 is a $\mathfrak{F}_G$-proper short exact sequence of $RG$-modules and both $A$ and $C$ admit $\mathfrak{F}_G$-proper $\mathfrak{F}_G$-projective resolutions, then there is an $\FGExt^*_{RG}(-,M)$ long exact sequence for any $RG$-module $M$.
\end{Lemma}
 
For any $RG$-module $M$ the \emph{$\mathfrak{F}_G$ projective dimension} of $G$ denoted $\FGpd M$ is the minimal length of an $\mathfrak{F}_G$-proper resolution of $M$ by $\mathfrak{F}_G$-projectives.  We set $\FGcd G = \FGpd R$.  Note that these finiteness conditions will not be defined unless $R$ admits an $\mathfrak{F}_G$-proper resolution by $\mathfrak{F}_G$-projectives.  \index{FG projective dimension@$\mathfrak{F}_G$ projective dimension $\FGpd$}\index{FG cohomological dimension@$\mathfrak{F}_G$ cohomological dimension $\FGcd G$}

One could think of $\mathfrak{F}_G$-cohomology as the ``Gorenstein cohomology relative $\mathfrak{F}$''.
\subsection{Technical results}

We need some results for the $\mathfrak{F}_G$-cohomology whose analogues are well-known for Gorenstein cohomology \cite{Holm-GorensteinHomologicalDimensions}.

We say an $RG$-module $M$ admits a \emph{right} resolution by $\mathfrak{F}$-projectives if there exists an exact chain complex
\[ 0 \longrightarrow M \longrightarrow T_{-1} \longrightarrow T_{-2} \longrightarrow \cdots \]
where the $T_i$ are $\mathfrak{F}$-projectives.  $\mathfrak{F}$-strong right resolutions and $\mathfrak{F}$-split right resolutions are defined as for any chain complex.

\begin{Lemma}\label{lemma:FGprog iff Ext vanishes and strong right res}
 An $RG$-module $M$ is $\mathfrak{F}_G$-projective if and only if $M$ satisfies
\begin{equation*}
 \FinExt^i_{RG}(M, Q) \cong 0 \text{ for all $\mathfrak{F}$-projective $Q$ } \tag{$\star$} 
\end{equation*}
for all $i \ge 1$ and $M$ admits a right $\mathfrak{F}$-strong $\mathfrak{F}$-split resolution by $\mathfrak{F}$-projectives.
\end{Lemma}
\begin{proof}
 If $M$ is the cokernel of a $\mathfrak{F}$-strong $\mathfrak{F}$-complete resolution $T_*$ then for all $i \ge 1$ and any $\mathfrak{F}$-projective $Q$,
\[ \FinExt^i_{RG}( M, Q ) \cong H^i\Hom_{RG}(T_*^+, Q) \]
where $T_*^+$ denotes the resolution $T_i^+ = T_i$ if $i \ge 0$ and $T_i^+ = 0$ for $i < 0$.  Then $(\star)$ follows because $T_*$ is $\mathfrak{F}$-strong.  

Conversely given $(\star)$ and an $\mathfrak{F}$-strong right resolution $T_*^-$ then let $T_*^+$ be the standard $\mathfrak{F}$-split resolution for $M$ (Lemma \ref{lemma:Fsplit resolutions exist}), $(\star)$ ensures that $T_*^+$ is $\mathfrak{F}$-strong and splicing together $T_*^+$ and $T_*^-$ gives the required resolution.
\end{proof}

\begin{Lemma}\label{lemma:FExt(Gproj, finite Fpd) = 0}
 If $\Finpd N < \infty$ and $M$ is $\mathfrak{F}_G$-projective then $\FinExt^i_{RG}(M, N) = 0$ for all $ i \ge 1$.
\end{Lemma}
\begin{proof}
Let $P_* \longrightarrow N$ be a $\mathfrak{F}$-split $\mathfrak{F}$-projective resolution then by a standard dimension shifting argument 
\[  \FinExt^i(M, N) \cong \FinExt^{i+j}(M, K_j)  \]
where $K_j$ is the $j^\text{th}$ syzygy of $P_*$.  Since $K_j$ is $\mathfrak{F}$-projective for $j \ge n$ the result follows from Lemma \ref{lemma:FGprog iff Ext vanishes and strong right res}.
\end{proof}

\begin{Prop}\label{prop:syzygy Fproj then FGproper}
Let $A$ be any $RG$-module and $P_* \longrightarrow A$ a length $n$ $\mathfrak{F}$-split resolution of $A$ with $P_i$ $\mathfrak{F}$-projective for $i \ge 1$, then $P_*$ is $\mathfrak{F}_G$-proper.
\end{Prop}
\begin{proof}
The case $n= 0$ is obvious.  If $n = 1$ then for any $\mathfrak{F}_G$-projective $Q$, there is a long exact sequence
\[ 0 \longrightarrow \Hom_{RG}(Q, P_1) \longrightarrow \Hom_{RG}(Q, P_0) \longrightarrow \Hom_{RG}(Q, A) \] 
\[ \longrightarrow \FinExt^1_{RG}(Q, P_1) \longrightarrow \cdots \] 
but $\FinExt^1_{RG}(Q, P_1) = 0$ by Lemma \ref{lemma:FExt(Gproj, finite Fpd) = 0}.  

Assume $n \ge 2$ and let $K_*$ be the syzygies of $P_*$, then there is an $\mathfrak{F}$-split resolution
\[ 0 \longrightarrow P_n \longrightarrow \cdots \longrightarrow P_{i+1} \longrightarrow K_i \longrightarrow 0 \] 
so $\Finpd K_i < \infty$ for all $i \ge 0$.  Thus every short exact sequence
\[ 0 \longrightarrow K_i \longrightarrow P_i \longrightarrow K_{i-1} \]
is $\mathfrak{F}_G$-proper by Lemma \ref{lemma:FExt(Gproj, finite Fpd) = 0}, so $P_*$ is $\mathfrak{F}_G$-proper.
\end{proof}

\begin{Lemma}[Comparison Lemma]\label{lemma:GF comparison lemma}
 Let $A$ and $B$ be two $RG$-modules with $\mathfrak{F}$-strong $\mathfrak{F}$-split right resolutions by $\mathfrak{F}$-projectives called $S^*$ and $T^*$ respectively, then any map $f: A \to B$ lifts to a map $f_*$ of complexes as shown below:
\[
\xymatrix{
0 \ar[r] & A \ar[r] \ar^f[d] & S^{1} \ar[r] \ar^{f_1}[d] & S^2 \ar[r] \ar^{f_2}[d] & \cdots \\
0 \ar[r] & B \ar[r] & T^{1} \ar[r]  & T^{2} \ar[r] & \cdots 
}
\]
The map of complexes is unique up to chain homotopy and if $f$ is $\mathfrak{F}$-split then so is $f_*$.
\end{Lemma}
\begin{proof}
 The lemma without the $\mathfrak{F}$-splitting comes from dualising \cite[p.169]{EnochsJenda-RelativeHomologicalAlgebra1}, see also \cite[Proposition 1.8]{Holm-GorensteinHomologicalDimensions}. 

 Assume $f$ is $\mathfrak{F}$-split and consider the map of complexes restricted to $R H$ for some finite subgroup $H$ of $G$.  Let $\iota_*^T$ and $\iota_*^S$ denote the splittings of the top and bottom rows and $s_*$ the splitting of $f_*$, constructed only up to degree $i-1$.  The base case of the induction, when $i = 0$, holds because $f$ is $\mathfrak{F}$-split.
\[
\xymatrix@+15pt{
\cdots \ar^{\partial_{i-2}^S}@/^/[r] & S^{i-1} \ar^{\iota_{i-2}^S}@/^/[l] \ar^{\partial_{i-1}^S}@/^/[r] \ar^{f_{i-1}}@/^/[d] & S^i \ar^{\iota_{i-1}^S}@/^/[l] \ar^{\partial_{i}^S}@/^/[r] \ar^{f_{i}}@/^/[d] & \ar^{\iota_{i}^S}@/^/[l] \cdots \\
\cdots \ar^{\partial_{i-2}^T}@/^/[r] & T^{i-1} \ar^{\iota_{i-2}^T}@/^/[l] \ar^{\partial_{i-1}^T}@/^/[r] \ar^{s_{i-1}}@/^/[u] & T^i \ar^{\iota_{i-1}^T}@/^/[l] \ar^{\partial_{i}^T}@/^/[r] & \ar^{\iota_{i}^T}@/^/[l] \cdots 
}
\]
Let $s_{i} = \partial_{i-1}^S \circ s_{i-1} \circ \iota_{i-1}^T$.  Then,
\begin{align*}
 f_i \circ s_i &= f_i \circ \partial_{i-1}^S \circ s_{i-1} \circ \iota_{i-1}^T \\
 &= \partial_{i-1}^T \circ f_{i-1} \circ s_{i-1} \circ \iota^T_{i-1} \\
 &= \partial_{i-1}^T \circ \iota_{i-1}^T \\
 &= \id_{T^i}
\end{align*}
where the second equality is the commutativity condition coming from the fact that $f_*$ is a chain map.
\end{proof}

\subsection{An Avramov--Martsinkovsky long exact sequence in \texorpdfstring{$\mathfrak{F}$}{F}-coho\-mology}

\begin{Theorem}\label{theorem:FAM LES}
Given an $\mathfrak{F}$-strong $\mathfrak{F}$-complete resolution of $R$ there is a long exact sequence
\[ 0 \longrightarrow \mathfrak{F}_G H^1(G, -) \longrightarrow \cdots  \]
\[ \cdots \longrightarrow \FhatH^{n-1}(G, -) \longrightarrow \FGH^n(G, -) \longrightarrow \FinH^n(G, -)  \]
\[ \longrightarrow \FhatH^n(G, -) \longrightarrow \FGH^{n+1}(G, -) \longrightarrow \cdots. \]
\end{Theorem}
\begin{proof}
 We follow the proof in \cite[\S 3]{AsadollahiBahlekehSalarian-HierachyCohomologicalDimensions}.  Let $T_*$ be an $\mathfrak{F}$-strong $\mathfrak{F}$-complete resolution coinciding with an $\mathfrak{F}$-projective $\mathfrak{F}$-split resolution $P_*$ in degrees $n$ and above.  We may choose $\theta_* : T_* \longrightarrow P_*$ to be $\mathfrak{F}$-split by Lemma \ref{lemma:GF comparison lemma} and without loss of generality we may also assume that $\theta_i$ is surjective for all $i$.

 Truncating at position $0$ and adding cokernels gives the bottom two rows of the diagram below, the row above is the row of kernels.  Note that the map $A \to R$ is necessarily surjective since the maps $T_0 \to P_0$ and $P_0 \to R$ are surjective.

\[ 
\xymatrix{
\cdots \ar[r] & 0   \ar[r] \ar[d] & K_{n-1} \ar[r] \ar[d] & \cdots \ar[r] & K_0 \ar[r] \ar[d] & K   \ar[r] \ar[d] & 0 \\
\cdots \ar[r] & T_n \ar[r] \ar[d] & T_{n-1} \ar[r] \ar[d] & \cdots \ar[r] & T_0 \ar[r] \ar[d] & A   \ar[r] \ar[d] & 0 \\
\cdots \ar[r] & P_n \ar[r]        & P_{n-1} \ar[r]        & \cdots \ar[r] & P_0 \ar[r]        & R \ar[r]        & 0
}
\]
We make some observations about the diagram:  Firstly, since the module $A$ is the cokernel of a $\mathfrak{F}$-strong $\mathfrak{F}$-complete resolution, $A$ is $\mathfrak{F}_G$-projective.  Secondly, in degree $i \ge 0$ the columns are $\mathfrak{F}$-split and the $P_i$ are $\mathfrak{F}$-projective, thus the $K_i$ are $\mathfrak{F}$-projective for all $i \ge 0$.  Thirdly the far right vertical short exact sequence is $\mathfrak{F}$-split since the degree $0$ column and the rows are $\mathfrak{F}$-split. Finally the top row is exact and $\mathfrak{F}$-split since the other two rows are.  

Apply the functor $\Hom_{RG} (-, M)$ for an arbitrary $R G$-module $M$ and take homology.  This gives a long exact sequence
\[ 
\cdots \longrightarrow  \FinH^i (G, M) \longrightarrow \FhatH^i (G, M) \longrightarrow H^i\Hom_{RG}(K_*, M) \longrightarrow \cdots. 
\]
We can simplify the right-hand term:
\begin{align*}
H^i\Hom_{RG}(K_*, M) &\cong \FGExt^i_{RG}(K, M) \\
&\cong \FGH^{i+1}(G, M)
\end{align*}
 where the first isomorphism is because, by Proposition \ref{prop:syzygy Fproj then FGproper}, the top row is $\mathfrak{F}_G$-proper.  For the second isomorphism note that the short exact sequence
\[ 0 \longrightarrow K \longrightarrow A \longrightarrow R \longrightarrow 0 \]
is $\mathfrak{F}_G$-proper by Proposition \ref{prop:syzygy Fproj then FGproper}, so 
\[ 0 \longrightarrow K_{n-1} \longrightarrow \cdots \longrightarrow K_0 \longrightarrow A \longrightarrow R \longrightarrow 0 \]
is an $\mathfrak{F}_G$-proper $\mathfrak{F}_G$-projective resolution of $R$.  Thus the second isomorphism follows from the short exact sequence and Lemma \ref{lemma:long exact sequences in FGExt}.
\end{proof}

\begin{Cor}
 If there exists an $\mathfrak{F}$-strong $\mathfrak{F}$-complete resolution of $R$ then $\FGcd G < \infty$.
\end{Cor}
\begin{proof}
 In the proof of the theorem we assumed an $\mathfrak{F}$-strong $\mathfrak{F}$-complete resolution of $R$ and built a finite length $\mathfrak{F}_G$-proper resolution of $R$ by $\mathfrak{F}_G$-projectives.
\end{proof}

\begin{Prop}\label{prop:AM seqs diagram}
If the Avramov--Martsinkovsky long exact sequence and the long exact sequence of Theorem \ref{theorem:FAM LES} both exist, then there is a commutative diagram:
\[
\xymatrix{
 \cdots \ar[r] & \FhatH^{n-1} \ar^{\gamma_{n-1}}[d] \ar[r] & \FGH^n \ar[r] \ar^{\alpha_n}[d] & \FinH^n \ar[r] \ar^{\beta_n}[d] & \FhatH^n \ar^{\gamma_n}[d] \ar[r] & \FGH^{n+1} \ar[r] \ar^{\alpha_{n+1}}[d] & \cdots \\
 \cdots \ar[r] &  \widehat{H}^{n-1} \ar[r] & \GH^n \ar[r] \ar_{\eta_n}[ru] & H^n \ar[r] & \widehat{H}^n \ar[r]   & \GH^{n+1} \ar[r] \ar[ru] & \cdots \\
}
\]
Where for conciseness we have written $H^n$ for $H^n(G, -)$ etc.
\end{Prop}
\begin{proof}
The Avramov--Martsinkovsky long exact sequence is constructed analogously to in the proof of Theorem \ref{theorem:AM LES}, we give a quick sketch below as we will need the notation.  Take a strong complete resolution $T^\prime_*$ of $R$ coinciding with a projective resolution $P_*^\prime$ in high dimensions and let $A^\prime$ be the zeroth cokernel of $T_*^\prime$.  Thus $A^\prime$ is Gorenstein projective.  Again, the map $T_*^\prime \to P_*^\prime$ is assumed surjective and the kernel $K_*^\prime$ is a projective resolution of $K^\prime$, the kernel of the map $A^\prime \longrightarrow R$.  Applying $\Hom_{RG}(-, M)$, for some $RG$-module $M$, to the short exact sequence of complexes 
\[0 \longrightarrow K_*^\prime \longrightarrow T_*^\prime \longrightarrow P_*^\prime \longrightarrow 0  \] 
gives the Avramov--Martsinkovsky long exact sequence.

Let $T_*$, $P_*$, $K_*$, $K$ and $A$ be as defined in the proof of Theorem \ref{theorem:AM LES}.  There is a commutative diagram of chain complexes
\[ 
\xymatrix{
0 \ar[r] & {K}_* \ar[r]  & {T}_* \ar[r] & {P}_* \ar[r] & 0\\
0 \ar[r] & {K}_*^\prime \ar[r] \ar^\alpha[u] & {T}_*^\prime \ar[r] \ar^\gamma[u] & {P}_*^\prime \ar[r] \ar^\beta[u] & 0 
}\]
where the maps $\beta$ exists by the comparison theorem for projective resolutions and $\gamma$ exists by the comparison theorem for strong complete resolutions \cite[Lemma 2.4]{CornickKropholler-OnCompleteResolutions}.   The map $\alpha$ is the induced map on the kernels.  Applying the functor $\Hom_{RG}(-, M)$ for some $RG$-module $M$, and taking homology, the maps $\alpha$, $\beta$ and $\gamma$ induce the maps $\alpha_*$, $\beta_*$ and $\gamma_*$.

Finally we construct the map $\eta_n: \GH^n(G, -) \longrightarrow \FinH^n(G, -)$.  Let $B_*$ be a G-proper Gorenstein projective resolution and recall $P_*$ is an $\mathfrak{F}$-split resolution by $\mathfrak{F}$-projectives.  Then $B_*$ is $\mathfrak{F}$-split (Lemma \ref{lemma:Gproper res of R is Fsplit}) so there is a chain map  $P_* \to B_*$ inducing $\eta_*$ on cohomology.

Commutativity is obvious for the diagram with the maps $\eta_i$ removed, leaving us with two relations to prove.  Let 
\[
\varepsilon_n^G : \GH^n(G, -) \longrightarrow H^n(G, -)
\]
denote the map from the commutative diagram.  This is the map induced by comparison of a resolution of Gorenstein projectives and ordinary projectives \cite[3.2,3.11]{AsadollahiBahlekehSalarian-HierachyCohomologicalDimensions}.  We get $\beta_* \circ \eta_* = \varepsilon^G_* $, since all the maps are induced by comparison of resolutions, and such maps are unique up to chain homotopy equivalence.

The final commutativity relation, that $\eta_* \circ \alpha_* = \varepsilon_*^{\mathfrak{F}_G} $, is the most difficult to show.  Here 
\[\varepsilon^{\mathfrak{F}_G}_n: \FGH^n(G, -) \longrightarrow \FinH^n(G, -)\]
denotes the map from the commutative diagram, it is induced by comparison of resolutions.

Here is a commutative diagram showing the resolutions involved:
\[
\xymatrix@-15pt{
& 0 \ar[rr] & & K \ar[rr] & & A \ar[rr] & & R \ar[rr] & & 0 \\
0 \ar[rr] & & K_* \ar[rr] \ar[ru] & \ar[u] & T_* \ar[ru] \ar[rr] & \ar[u] & P_* \ar[ru] \ar[rr] & \ar[u] & 0 & \\
& 0 \ar@{-}[r] & \ar[r] & K^\prime \ar@{-}[r] \ar@{-}[u] & \ar[r] & A^\prime \ar@{-}[r] \ar@{-}[u] & \ar[r] & R \ar@{-}[u] \ar[rr] & & 0 \\
0 \ar[rr] & & K_*^\prime \ar[rr] \ar[ru] \ar[uu] & & T_*^\prime \ar[rr] \ar[ru] \ar[uu] & & P_*^\prime \ar[rr] \ar[ru] \ar[uu] & & 0
}
\]

Let $L_*$ be the chain complex defined by $L_i = K_{i-1}$ for all $i \ge 1$ and $L_0 = A$, with boundary map at $i =1$ the composition of the maps $K_0 \to K$ and $K \to A$.  Thus $L_*$ is acyclic except at degree zero where $H_0L_* = R$.  Similarly, let $L_*^\prime$ denote chain complex with $L_i^\prime = K_{i-1}^\prime$ for all $i \ge 1$ and $L_0^\prime = A^\prime$ augmented by $A^\prime$, so $L_*^\prime$ is acyclic except at degree zero where $H_0L_*^\prime = R$.  Note that $L_*$ is an $\mathfrak{F}_G$-proper resolution of $R$ by Proposition \ref{prop:syzygy Fproj then FGproper} and $L_*^\prime$ is a G-proper resolution of $R$ by the Gorenstein cohomology version of the same proposition.

Recall that the maps $\varepsilon_*^{\mathfrak{F}_G}$ and $\eta_*$ are induced by comparison of resolutions:  $\varepsilon^{\mathfrak{F}_G}_*$ is induced by a map $P_* \to L_*$ and $\eta_*$ is induced by a map $P_* \to L_*^\prime$.  The map 
\[\FGExt^i_{RG}(K, -) \longrightarrow \GExt^i_{RG}(K^\prime, -)\]
is induced by $\alpha: K_*^\prime \longrightarrow K_*$.  Thus the map 
\[\alpha_* : \FGH^n(G, -) \longrightarrow \GH^n(G, -)\]
is induced by $L_*^\prime \longrightarrow L_*$.  The diagram below is the one we must show commutes.  
\[
\xymatrix{
\FGH^n(G, -) \cong H^n\Hom_{RG}(L_* , -) \ar^{\alpha_n}[d] \ar^{\varepsilon_n^{\mathfrak{F}_G}}[r] & \FinH^n(G, -) \cong H^n \Hom_{RG}(P_* , -) \\
\GH^n(G, -) \cong H^n\Hom_{RG}(L^\prime_* , -) \ar_{\eta_n}[ur] & \\
}
\]
Since the composition $P_*$ to $L^\prime_*$ to $L_*$ is a map of resolutions from $P_*$ to $L_*$, and such maps are unique up to chain homotopy equivalence, this completes the proof.
\end{proof}
\begin{Cor}
Given an $\mathfrak{F}$-strong $\mathfrak{F}$-complete resolution of $R$, $\Gcd G = n < \infty$ implies $\FinH^i(G, -)$ injects into $\FhatH^i(G, -)$ for all $i \ge n+1$.
\end{Cor}
\begin{proof}
$\Gcd G < \infty$ implies the Avramov--Martsinkovsky long exact sequence exists (Theorem \ref{theorem:AM LES}). Consider the the commutative diagram of Proposition \ref{prop:AM seqs diagram}.  The map
\[ \FGH^i(G, -) \longrightarrow \FinH^i(G, -) \]
factors as $\eta_i \circ \alpha_i = 0$, so since $\GH^{i}(G, -) = 0$ for all $i \ge n+1$, $\FinH^i(G, -)$ injects into $\FhatH^i(G, -)$ for all $i \ge n+1$.  
\end{proof}

\begin{Theorem}\label{thm:Fcd finite then Fcd=Gcd}
 If $\Fincd G < \infty$ then $\Fincd G = \Gcd G$.
\end{Theorem}
\begin{proof}~
 We know already that $\Gcd G \le \Fincd G$ (see Section \ref{section:I Gorenstein}).  If $\Fincd G < \infty$ then it is trivially true that $\mathfrak{F}$ admits an $\mathfrak{F}$-strong $\mathfrak{F}$-complete resolution, thus $\FinH^i(G, -)$ injects into $\FhatH^i(G, -)$ for all $i \ge \Gcd G + 1$, but $\FhatH^i(G, -)$ is always zero since $\Fincd G < \infty$ \cite[4.1(i)]{Kropholler-OnGroupsOfTypeFP_infty}.
\end{proof}

\begin{Example}\label{example:HF} 
Let $R = \ZZ$ for this example.  Kropholler introduced the class $\HF$ of hierarchically decomposable groups in \cite{Kropholler-OnGroupsOfTypeFP_infty} as the smallest class of groups such that if there exists a finite-dimensional contractible $G$-CW complex with stabilisers in $\HF$ then $G \in \HF$.  Let $\HF_b$ denote the subclass of $\HF$ containing groups with a bound on the orders of their finite subgroups.

The $\ZZ G$-module $B(G, \ZZ)$ of bounded functions from $G$ to $\ZZ$ was first studied in \cite{KrophollerTalelli-PropertyOfFundamentalGroupsOfGraphOfFiniteGroups}, Kropholler and Mislin proved that if $G$ is $\HF$ with a bound on lengths of chains of finite subgroups and $\pd_{\ZZ G} B(G, \ZZ ) < \infty$ then $\OFincd G < \infty$, in particular $\Fincd G < \infty$ \cite[Theorem B]{KrophollerMislin-GroupsOnFinDimSpacesWithFinStab}.  If $\Gcd G < \infty$ then $\pd_{\ZZ G}B(G, \ZZ ) < \infty$ \cite[2.10]{AsadollahiBahlekehSalarian-HierachyCohomologicalDimensions}\cite[Theorem C]{CornickKropholler-HomologicalFinitenessConditions}.  Thus if $G \in \HF_b$ then $\Gcd G = \Fincd G$.
\end{Example}

\section{Group extensions}\label{section:Group Extensions}

Recall Theorem \ref{theorem:Fcd=HFcd}, that $\Fincd G =\HFincd G $ for all groups.  The invariant $\HFincd G$ was studied by Degrijse in \cite{Degrijse-ProperActionsAndMackeyFunctors} where he proves the following (though stated for $\HFincd G$ not $\Fincd G$):  
\begin{Theorem}\cite[Theorem B]{Degrijse-ProperActionsAndMackeyFunctors}
 Let
\[ 
1 \longrightarrow N \longrightarrow G \longrightarrow Q \longrightarrow 1  
\]
be a short exact sequence of groups such that every finite index overgroup of $N$ in $G$ has a bound on the orders of the finite subgroups not contained in $N$.  If $\Fincd G < \infty$ then $ \Fincd G \le \Fincd N + \Fincd Q $.
\end{Theorem}

Since Gorenstein cohomological dimension is subadditive under extensions \cite[Remark 2.9(2)]{BahlekehDembegiotiTalelli-GorensteinDimensionAndProperActions}, an application of Theorem \ref{thm:Fcd finite then Fcd=Gcd} removes the condition on the orders of finite subgroups:
\begin{Cor}\label{cor:Fcd subadditive if finite}
 Given a short exact sequence of groups
\[ 1 \longrightarrow N \longrightarrow G \longrightarrow Q \longrightarrow 1  \]
if $\Fincd G < \infty$ then $ \Fincd G \le \Fincd N + \Fincd Q $.
\end{Cor}

For further discussion on the behaviour of $\Fincd G$ (equivalently $\HFincd G$) under group extensions and other standard constructions see Section \ref{subsection:mackey HF closure properties}.

\section{Rational cohomological dimension}\label{section:cdQ}

For this section, let $R = \ZZ$.  Gandini has shown that for groups in $\HF$, $\cd_{\QQ} G \le \Gcd G$ \cite[Remark 4.14]{Gandini-CohomologicalInvariants} and this is the only result we are aware of relating $\cd_{\QQ}G$ and $\Gcd G$.  In Proposition \ref{prop:cdQG finite then cdQG le GcdG} we show that $\cd_{\QQ}G \le \Gcd G$ for all groups with $\cd_{\QQ} G < \infty$.  Recall there are examples of torsion-free groups with $\cd_{\QQ} G < \cd_{\ZZ} G$ \cite[Example 8.5.8]{Davis} and $\Gcd G = \cd_{\ZZ} G$ whenever $\cd_{\ZZ} G < \infty$ \cite[Corollary 2.9]{AsadollahiBahlekehSalarian-HierachyCohomologicalDimensions}, so we cannot hope for equality of $\cd_{\QQ} G$ and $\Gcd G $ in general.

\begin{Question}\label{question:G exist groups Gcd finite cdQ infinite}
 Are there groups $G$ with $\Gcd G < \infty$ but $\cd_{\QQ} G = \infty$?
\end{Question}

Recall from Section \ref{section:I Gorenstein} that $\silp RG$ denotes the supremum of the injective lengths (injective dimensions) of all projective $RG$-modules and $\spli RG$ denotes the infimum of the projective lengths (projective dimensions) of all injective $RG$-modules.

\begin{Lemma}
 For any group $G$, $\silp \QQ G \le \silp \ZZ G$.
\end{Lemma}
\begin{proof}
 Using \cite[Theorem 4.4]{Emmanouil-OnCertainCohomologicalInvariantsOfGroups}, $\silp \QQ G = \spli \QQ G$ and $\silp \ZZ G = \spli \ZZ G$.  Combining with \cite[Lemma 6.4]{GedrichGruenberg-CompleteCohomologicalFunctors} that $\spli \QQ G \le \spli \ZZ G$ gives the result.
\end{proof}

\begin{Lemma}\label{lemma:complete tensored with Q}
 If $\Gcd G < \infty$ then for any $\QQ G$-module $M$ there is a natural isomorphism
\[ 
\widehat{H}^*(G, M) \otimes \QQ \cong \widehat{\Ext}^*_{\QQ G}(\QQ,M).
\]
\end{Lemma}
\begin{proof}
 Let $T_*$ be a strong complete resolution of $\ZZ$ by $\ZZ G$-modules, then $T_* \otimes \QQ$ is a strong complete resolution of $\QQ$ by $\QQ G$-modules.  By an obvious generalisation of \cite[Lemma 2.2]{MislinTalelli-FreelyProperlyOnFDHomotopySpheres}, if $\silp \QQ G \le \infty$ then any complete $\QQ G$-module resolution is a strong complete $\QQ G$-module resolution, so since $\silp \QQ G < \silp \ZZ G \le \infty$, $T_* \otimes \QQ$ is a strong complete resolution.  This gives a chain of isomorphisms for any $\QQ G$-module $M$:
\begin{align*}
 \widehat{H}^*(G, M) \otimes \QQ &\cong H^* \Hom_{\ZZ G}(T_*, M) \otimes \QQ \\
 &\cong H^* \Hom_{\QQ G}(T_* \otimes \QQ, M)  \\
 &\cong \widehat{\Ext}^*_{\QQ G}(\QQ , M).
\end{align*}
\end{proof}

\begin{Prop}\label{prop:cdQG finite then cdQG le GcdG}
 If $\cd_{\QQ} G < \infty$ then $\cd_{\QQ} G \le \Gcd G$.
\end{Prop}
\begin{proof}
There is nothing to show if $\Gcd G = \infty$ so assume that $\Gcd G < \infty$.  Since $\QQ$ is flat over $\ZZ$, tensoring the Avramov--Martsinkovsky long exact sequence with $\QQ$ preserves exactness.  Combining this with Lemma \ref{lemma:complete tensored with Q} and the well-known fact that for any $\QQ G$-module $M$ there is a natural isomorphism \cite[p.2]{Bieri-HomDimOfDiscreteGroups} 
 \[H^*(G, M) \otimes \QQ \cong \Ext^*_{\QQ G}(\QQ, M)\]
 gives the long exact sequence
\[ 
\cdots \longrightarrow \GH^i(G, M) \otimes \QQ \longrightarrow \Ext^i_{\QQ G}(\QQ, M) \longrightarrow \widehat{\Ext}^i_{\QQ G}(\QQ,M) \longrightarrow \cdots.
\]
Since $\cd_{\QQ} G < \infty$, we have that $\widehat{\Ext}^i_{\QQ G}(\QQ,M) = 0$ \cite[4.1(i)]{Kropholler-OnGroupsOfTypeFP_infty}.  Thus there is an isomorphism for all $i$,
\[
\GH^i(G, M) \otimes \QQ \cong \Ext^i_{\QQ G}(\QQ, M)
\]
and the result follows.
\end{proof}

\chapter{Bredon duality groups}\label{chapter:BD}
This chapter contains material that has appeared in:
\begin{itemize}
 \item Bredon--Poincar\'e duality groups (2013, to appear J.~Group Theory) \cite{Me-BredonPoincareDuality}.
\end{itemize}

In this chapter we study Bredon duality and Bredon--Poincar\'e duality groups.  Recall that a \emph{Bredon duality group over $R$} is a group $G$ of type $\OFinFP$ over $R$ such that for every finite subgroup $H$ of $G$ there is an integer $d_H$ with 
\[ 
H^i(WH, R [WH]) = \left\{ \begin{array}{l l} \text{$R$-flat} & \text{if $i = d_H$,} \\ 0 & \text{else.} \end{array} \right. 
\]
Furthermore, $G$ is said to be \emph{Bredon--Poincar\'e duality} over $R$ if for all finite subgroups $H$, 
\[
H^{d_H}(WH, R [WH]) = R.
\]

In Section \ref{section:BD examples} we give several sources of examples of both Bredon duality and Bredon--Poincar\'e duality groups, including the example below of Jonathan Block and Schmuel Weinberger, suggested to us by Jim Davis.

\theoremstyle{plain}\newtheorem*{CustomThmO}{Theorem \ref{theorem:BPD groups with no manifold model}}
\begin{CustomThmO}
  There exist examples of Bredon--Poincar\'e duality groups over $\ZZ$, such that $WH$ is finitely presented for all finite subgroups $H$ but $G$ doesn't admit a cocompact manifold model $M$ for $\EFin G$.
\end{CustomThmO}

This is a counterexample to a possible generalisation of Wall's conjecture, which asks if a finitely presented Poincar\'e duality group admits a cocompact manifold model for $\E G$:  Let $G$ be Bredon--Poincar\'e duality over $\ZZ$, such that $WH$ is finitely presented for all finite subgroups $H$, does $G$ admit a cocompact manifold model $M$ for $\EFin G$?  

Section \ref{section:BD low dimensions} contains an analysis of Bredon duality and Bredon--Poincar\'e duality groups of low dimension and Section \ref{section:BD extensions} looks at when these properties are preserved under group extensions.  

Recall that given a Bredon duality group $G$ of dimension $n$ we write $\mathcal{V}(G)$ for the set
\[
\mathcal{V}(G) = \{ d_F : F \text{ a non-trivial finite subgroup of }G\} \subseteq \{0, \ldots, n \}.
\]
In Example \ref{example:duality arbitrary V} we build Bredon duality groups with arbitrary $\mathcal{V}(G)$ and in Section \ref{section:reflection groups} we build Bredon--Poincar\'e duality groups with a large selection of $\mathcal{V}(G)$, although we are unable to produce arbitrary $\mathcal{V}(G)$.

One might hope to give a definition of Bredon--Poincar\'e duality groups in terms of Bredon cohomology only, we show in Section \ref{section:wrong notion of duality} that the na\"ive idea of asking that a group be $\OFinFP$ with 
\[
H^i_{\OFin}(G, R[?,-]_{\OFin}) \cong \left\{ \begin{array}{l l} \uR & \text{if $i = n$,} \\ 0 & \text{else,} \end{array} \right.
\]
is not the correct definition, namely we show in Theorem \ref{theorem:bredon wrong duality} that any such group is necessarily a torsion-free Poincar\'e duality group over $R$.
 
\section{Preliminary observations}

Recall that a group $G$ is \emph{$R$-torsion-free} if the order of every finite subgroup of $G$ is invertible in $R$, equivalently the order of every finite order element is invertible in $R$ (see page \pageref{Rtorsionfree}).

Recall that a Bredon duality group is said to be dimension $n$ if $\OFincd G = n$.

\begin{Lemma}\label{lemma:uDn over Z then uDn over R}
 If $G$ is Bredon duality of dimension $n$ over $\ZZ$ then $G$ is Bredon duality of dimension $n$ over any ring $R$, with the same values of $d_H$ for all finite subgroups $H$.
\end{Lemma}
\begin{proof}
 Since $G$ is $\OFinFP$ over $\ZZ$, $G$ is $\OFinFP$ over $R$ (Lemma \ref{lemma:ucdR less ucdZZ and OFFPn over ZZ implies OFFPn over R}).  Also because $G$ is $\OFinFP$ over $\ZZ$, $WH$ is $\FP_\infty$ over $\ZZ$ for all finite subgroups $H$ (Corollary \ref{cor:uFPn equivalent conditions}) and we may apply \cite[Corollary 3.6]{Bieri-HomDimOfDiscreteGroups} to get a short exact sequence
\[0 \longrightarrow H^q(WH, \ZZ [WH]) \otimes_\ZZ R \longrightarrow H^q(WH, R \otimes_\ZZ \ZZ [WH]) \]
\[\longrightarrow \Tor_1^\ZZ(H^{q+1}(WH, \ZZ [WH]) ,R) \longrightarrow 0.\]
$H^{q+1}(WH, \ZZ [WH])$ is $\ZZ$-flat for all $q$ giving an isomorphism
\[H^q(WH, \ZZ [WH]) \otimes_\ZZ R \cong H^q(WH, R [WH]).\]
Observing that if an Abelian group $M$ is $\ZZ$-flat then $M \otimes R$ is $R$-flat completes the proof.
\end{proof}

\begin{Lemma}\label{Lemma:observations}
 If $G$ is $R$-torsion-free and Bredon duality of dimension $n$ over $R$ then $d_H = \cd_R WH$ and $d_{1} \le n$.
\end{Lemma}

To prove the Lemma we need the following proposition, an analogue of \cite[VIII.6.7]{Brown} for arbitrary rings $R$ and proved in exactly the same way.
\begin{Prop}\label{prop:FP over Q then cdQ = max H^n}
 If $G$ is $\FP$ over $R$ then 
 \[
 \cd_R G = \max \{ n  :  H^n(G, R  G) \neq 0\}.
 \]
\end{Prop}

\begin{proof}[Proof of Lemma \ref{Lemma:observations}]
If $G$ is $R$-torsion-free then for any finite subgroup $H$, 
\[\cd_R N_GH \le \cd_R G \le \OFincd_R G\]
and $N_GH$ is $\FP_\infty$ over $R$ by Corollary \ref{cor:uFPn equivalent conditions} and Lemma \ref{lemma:no R torsion then cdR less ucd R etc}.  The short exact sequence 
\[
 1 \longrightarrow H \longrightarrow N_G H \longrightarrow WH \longrightarrow 1
\]
and Lemma \ref{lemma:Finite by Q implies HG is HQ} implies that
\[H^i(N_G H, R[N_G H]) \cong H^i(WH, R[WH]).\]
Thus Proposition \ref{prop:FP over Q then cdQ = max H^n} shows $d_H = \cd_R N_GH = \cd_R WH$.  Finally, $d_{1} \le n$ because $\cd_R G \le \OFincd_R G$ (Lemma \ref{lemma:no R torsion then cdR less ucd R etc}). 
\end{proof}

In the proposition below $\Fincd G$ denotes the $\mathfrak{F}$-cohomological dimension (see Section \ref{section:I Fcohomology}) and $\Gcd G$ denotes the Gorenstein cohomological dimension (see Section \ref{section:I Gorenstein}).

This proposition implies that if $G$ is Bredon--Poincar\'e duality over $R$ then $\Gcd G = \Fincd G = d_1$ and if $G$ is also virtually torsion-free then $\vcd G = d_1$ also.
\begin{Prop}\label{prop:FP and Ofincd finite then Gcd=Fcd=sup}
 If $G$ is $\FP_\infty$ with $\OFincd G < \infty$ then 
 \[
    \Gcd G = \Fincd G = \sup \{ n : H^n(G, RG) \neq 0 \}  ,
 \]
 and if $G$ is also virtually torsion-free then $\vcd G = \Gcd G$ also.
\end{Prop}
\begin{proof}
 This proof uses an argument due to Degrijse and Mart\'{\i}nez-P\'erez in \cite{DegrijseMartinezPerez-DimensionInvariants}.  By \cite[Theorem 2.20]{Holm-GorensteinHomologicalDimensions} the Gorenstein cohomological dimension can be characterised as
 \[ \Gcd G = \sup \{ n : H^n(G, P) \neq 0 \text{ for $P$ any projective $RG$-module} \}.  \] 
 As $G$ is $\FP_\infty$ we need only check when $P = RG$.  Since $\Fincd G \le \OFincd G < \infty$, we can conclude that $\Fincd G = \Gcd G$ (Theorem \ref{thm:Fcd finite then Fcd=Gcd}) and finally for virtually torsion-free groups $\Fincd G = \vcd G$ \cite[Theorem 5.1]{MartinezPerezNucinkis-MackeyFunctorsForInfiniteGroups}.
\end{proof}

\section{Examples}\label{section:BD examples}
In this section we provide several sources of examples of Bredon duality and Bredon--Poincar\'e duality groups, showing that these properties are not too rare.

\subsection{Smooth actions on manifolds}\label{subsection:duality smooth actions on manifolds}

Recall from the introduction that if $G$ has a cocompact manifold model $M$ for $\EFin G$ such that $M^H$ is a submanifold for all finite subgroups $H$ then $G$ is Bredon--Poincar\'e duality.  The following lemma gives a condition which guarantees that $M^H$ is a submanifold of $M$.
\begin{Lemma}\label{Lemma:proper and locally linear => fixed points are submanifolds}\cite[10.1 p.177]{Davis} If $G$ is a discrete group acting properly and locally linearly on a manifold $M$ then the fixed points subsets of finite subgroups of $G$ are submanifolds of $M$.
\end{Lemma}
Locally linear is a technical condition, the definition of which can be found in \cite[Definition 10.1.1]{Davis}, for our purposes it is enough to know that if $M$ is a smooth manifold and $G$ acts by diffeomorphisms then the action is locally linear.  The locally linear condition is necessary however---in \cite{DavisLeary-DiscreteGroupActionsOnAsphericalManifolds} examples are given of virtually torsion-free groups acting as a discrete cocompact group of isometries of a CAT(0) manifold which are not Bredon duality.

\begin{Example}
Let $p$ be a prime, let $K$ be a cyclic group of order $p$, and let $G$ be the wreath product
 \[G = \ZZ \wr K = \left( \bigoplus_{K} \ZZ \right) \rtimes K.\]
$G$ acts properly and by diffeomorphisms on $\RR^p$: The copies of $\ZZ$ act by translation along the axes, and the $K$ permutes the axes.  The action is cocompact with fundamental domain the quotient of the $p$-torus by the action of $K$.  The finite subgroup $K$ is a representative of the only conjugacy class of finite subgroups in $G$, and has fixed point set the line $\{(\lambda,\cdots,\lambda)  :  \lambda \in \RR \}$.  
If $z = (z_1, \ldots, z_p) \in \ZZ^p$ then 
\[
 (\RR^p)^{K^z} = \{ (\lambda + z_1, \ldots, \lambda + z_p)  :  \lambda \in \RR \}.
\]
Hence $\RR^p$ is a model for $\EFin G$ and, invoking Lemma \ref{Lemma:proper and locally linear => fixed points are submanifolds},  $G$ is a Bredon--Poincar\'e duality group of dimension $p$. 

We can explicitly calculate the Weyl group $WK$: let $L$ denote the copy of $\ZZ$ inside $\ZZ^p$ generated by $(1,1,\ldots,1)$, then the normaliser $N_GK$ is $L \rtimes K$ and thus the Weyl group $WK$ is isomorphic to $\ZZ$.  Since $K$ is a representative of the only non-trivial conjugacy class of finite subgroups it provides the only element of $\mathcal{V}(G)$, thus  $\mathcal{V}(G) = \{1\}$.
\end{Example}

\begin{Example}\label{example:duality Zn antipodal}
Fixing positive integers $m \le n$, if $G = \ZZ^n \rtimes C_2$ where $C_2$, the cyclic group of order 2, acts as the antipodal map on $\ZZ^{n-m} \le \ZZ^n$ then 
\[
N_GC_2 = C_GC_2 = \{ g \in G  :  gz = zg \}.  
\]
But this is exactly the fixed points of the action of $C_2$ on $G$, hence $N_G C_2 = \ZZ^m \rtimes C_2$ and
\[ 
H^i(W_GC_2, R[W_G C_2]) \cong \left\{ \begin{array}{l l} R & \text{ if } i = m, \\ 0 & \text{ else.} \end{array}\right. 
\]
$G$ embeds as a discrete subgroup of $\text{Isom}(\RR^n) = \RR^n \rtimes GL_n(\RR)$ and acts properly and cocompactly on $\RR^n$.  It follows that $G$ is $\OFinFP$ and $\OFincd G = n$ so $G$ is Bredon--Poincar\'e duality of dimension $n$ over any ring $R$ with $\mathcal{V} = \{m\}$.
\end{Example}

\begin{Example}\label{example:duality Zn antipodal generalised}
Similarly to the previous example we can take 
\[
G = \ZZ^n \rtimes \bigoplus_{i=1}^n C_2
\]
where the $j^\text{th}$ copy of $C_2$ acts antipodally on the $j^\text{th}$ copy of $\ZZ$ in $\ZZ^n$.  Note that $G$ is isomorphic to $(D_\infty)^n$ where $D_\infty$ denotes the infinite dihedral group.  As before $G$ embeds as a discrete subgroup of $\text{Isom}(\RR^n) = \RR^n \rtimes GL_n(\RR)$ and acts properly and cocompactly on $\RR^n$.  Thus $G$ is $\OFinFP$ and $\OFincd G = n$, so $G$ is Bredon--Poincar\'e duality of dimension $n$ over any ring $R$ with $\mathcal{V}(G) = \{0, \ldots, n\}$.

More generally, we could take a subgroup $\bigoplus_{i=1}^m C_2 \longhookrightarrow \bigoplus_{i=1}^nC_2$ and form the semi-direct product of $\ZZ^n$ with this subgroup.  Although this gives us a range of possible values for $\mathcal{V}(G)$ it is impossible to produce a full range of values.  Consider the case $m=2$, so we have a group 
\[
G = \ZZ^n \rtimes (A \times B) 
\]
where $A \cong B \cong C_2$, and both $A$ and $B$ act either trivially or antipodally on each coordinate of $\ZZ^n$.  We can describe the normaliser $N_GA$ by an element $(a_1, \ldots, a_n) \in \{0,1\}^n$, so $A$ acts trivially on the $i^\text{th}$ copy of $\ZZ$ if $a_i = 1$ and acts antipodally otherwise.  Thus, 
\[
N_GA = \left( \bigoplus_{i=1}^n \left\{ \begin{array}{l l} \ZZ & \text{ if $a_i = 1$ } \\ 0 & \text{ else.} \end{array} \right\} \right)  \rtimes (A \times B).  
\]
Similarly we can describe $N_G B$ by an element $(b_1, \ldots, b_n) \in \{ 0, 1\}^n$.  One calculates that the normaliser $N_G(A \times B)$ is described by the element 
\[
(a_1\wedge b_1, \ldots, a_n \wedge b_n) 
\]
where $\wedge$ denotes the boolean AND function.  This is because the $i^\text{th}$ copy of $\ZZ$ is normalised by $A \times B$ if and only if it is normalised by $A$ and also by $B$.

If $C$ denotes the subgroup of $A \times B$ generated by the element $(1,1)$ then the normaliser of $N_G C$ is described by the element
\[
 ( \neg ( a_1 \oplus b_1)  , \ldots, \neg (  a_n \oplus  b_n))
\]
where $\oplus$ denotes the boolean XOR function, and $\neg$ the unary negation operator.  

Now, using the above it can be shown that, for example, a Bredon--Poincar\'e duality group of dimension $4$ with the form
\[
G = \ZZ^4 \rtimes \bigoplus_{i = 1}^m C_2 
\]
cannot have $\mathcal{V}(G) = \{1, 3\}$.  Assume that such a $G$ exists, clearly $m \ge 2$, let $A$ and $B$ denote two of the $C_2$ summands of $\oplus_{1=1}^mC_2$.  Without loss of generality we can assume that $A$ and $B$ don't have the same action on $\ZZ^4$.  If $d_A = d_B = 1$ then by the description of the normaliser of $A \times B$ above, $d_{A \times B} = 0$, a contradiction.  If $d_A = d_B = 3$ then in order for $A$ and $B$ not to have the same action on $\ZZ^4$, we must have (up to some reordering of the coordinates)
\[
(a_1, \ldots, a_4) = (1,1,1,0)  
\]
\[
(b_1, \ldots, b_4) = (0,1,1,1).
\]
So $d_{A \times B} = 2$, a contradiction.  Finally, if $d_A = 1$ and $d_B = 3$ then let $C$ be the subgroup of $A \times B$ generated by $(1,1)$.  There are two possibilities, up to reordering of the coordinates, either 
\[
(a_1, \ldots, a_4) = (1,1,1,0) 
\]
\[
(b_1, \ldots, b_4) = (1,0,0,0) 
\]
or
\[
(a_1, \ldots, a_4) = (1,1,1,0) 
\]
\[
(b_1, \ldots, b_4) = (0,0,0,1).
\]
In the first case, $d_C =2$, and in the second case $d_{A \times B} = 0$, both contradictions.
\end{Example}

\begin{Example}\label{example:FarbWeinberger}In \cite[Theorem 6.1]{FarbWeinberger-IsometriesRigidityAndUniversalCovers}, Farb and Weinberger construct a Bredon--Poincar\'e duality group $G$ arising from a proper cocompact action on $\RR^n$ by diffeomorphisms, however $G$ is not virtually torsion-free. 
\end{Example}

\begin{Remark}\label{remark:duality Cp restrictions on fixed points}{\bf Restrictions on the dimensions of the fixed point sets.}
Suppose $G$ is a group acting smoothly on an $m$-dimensional manifold $M$, and suppose furthermore that $G$ contains a finite cyclic subgroup $C_p$ fixing a point $x \in M$.  There is an induced linear action of $C_p$ on the tangent space $ T_x M \cong \RR^m$, equivalently a representation of $C_p$ into the orthogonal group $O(m)$.  We can use this to give some small restrictions on the possible dimensions of the submanifold $M^{C_p}$, and hence on the values of $d_{C_p}$.

A representation of $C_p$ in $O(m)$ is simply a matrix $N$ with $N^p = 1$.  Using the Jordan--Chevalley decomposition (expressing a matrix as the product of its semi-simple and nilpotent parts), we see that $N$ is semi-simple, so viewing $N$ as a matrix over $\CC$ it is diagonalisable.  However, since $N^p = 1$ and the characteristic polynomial has coefficients in $\RR$, all the eigenvalues come in pairs $\omega,\omega^{-1}$, where $\omega$ is a $p^\text{th}$ root of unity.  Thus $N$ is conjugate via complex matrices to 
\[
\begin{pmatrix}
 \omega_1 & & & & & \\
  & \omega_1^{-1} & & & & \\
  &  & \ddots & & & \\
  &  &  & \omega_{\frac{m}{2}} & & \\
  &  &  &  & \omega_{\frac{m}{2}}^{-1} \\
\end{pmatrix}
\text{or}
\begin{pmatrix}
 \omega_1 & & & & & & \\
  & \omega_1^{-1} & & & & &  \\
  &  & \ddots & & & &  \\
  &  &  & \omega_{\frac{m-1}{2}} & & &  \\
  &  &  &  & \omega_{\frac{m-1}{2}}^{-1} & \\
  & & & & & \pm 1
\end{pmatrix}
\]
depending on whether $m$ is even or odd.  The blank space in the matrices should be filled with zeros.  Note that the $\pm 1$ term can only be a $-1$ if $p=2$.  The matrix
\[
  \begin{pmatrix}
    \omega & 0\\
    0 & \omega^{-1} 
  \end{pmatrix}
 \]
is conjugate via complex matrices to 
\[
 R_\theta = 
  \begin{pmatrix}
    \cos \theta & - \sin \theta \\
    \sin\theta & \cos \theta
  \end{pmatrix}
\]
Thus $N$ is conjugate via complex matrices to $R_{\theta_1} \oplus \cdots \oplus R_{\theta_{m/2}}$ or $R_{\theta_1} \oplus \cdots \oplus R_{\theta_{(m-1)/2}} \oplus (\pm 1)$, and by \cite[5.11]{Zhang-MatrixTheory}, they are conjugate via real matrices as well.  Hence the dimensions of the fixed point sets are the same.  Noting that the rotation matrix $R_\theta$ fixes only the origin when $\theta \neq 0$, we conclude that for $p\neq 2$, the fixed point set $M^{C_p}$ must be even dimensional if $m$ is even, and odd dimensional otherwise.

Consider the case that $G$ is a Bredon--Poincar\'e duality group, arising from a smooth cocompact action on an $m$-dimensional manifold $M$, and $C_p$ for $p \neq 2$ is some finite subgroup of $G$.  Then $d_{C_p}$ is exactly the dimension of the submanifold $M^{C_p}$, and by the discussion above $d_{C_p}$ is odd dimensional if $m$ is odd dimensional, even dimensional otherwise.  As demonstrated by Example \ref{example:duality Zn antipodal}, there are no restrictions when $p=2$.
\end{Remark}

\subsection{A counterexample to the generalised \texorpdfstring{PD$^\text{n}$}{PDn} conjecture}\label{section:counterexample to generalised PDn}

Let $G$ be Bredon--Poincar\'e duality over $\ZZ$, such that $WH$ is finitely presented for all finite subgroups $H$.  One might ask if $G$ admits a cocompact manifold model $M$ for $\EFin G$.  This is generalisation of the famous PD$^\text{n}$-conjecture, due to Wall \cite{Wall-HomologicalGroupTheory}.  This example is due to Jonathan Block and Schmuel Weinberger and was suggested to us by Jim Davis.

\begin{Theorem}\label{theorem:BPD groups with no manifold model}
 There exist examples of Bredon--Poincar\'e duality groups over $\ZZ$, such that $WH$ is finitely presented for all finite subgroups $H$, but there doesn't exist a cocompact manifold model $M$ for $\EFin G$.
\end{Theorem}

Combining Theorems 1.5 and 1.8 of \cite{BlockWeinberger-GeneralizedNielsenRealizationProblem} gives the following example.

\begin{Theorem}[Block--Weinberger]\label{theorem:block--weinberger}
There exists a short exact sequence of groups 
\[ 
1 \longrightarrow K \longrightarrow G \longrightarrow Q \longrightarrow 1
\]
with $Q$ finite, such that
\begin{enumerate}
 \item All torsion in $G$ is contained in $K$.
 \item There exists a cocompact manifold model for $\EFin K$.
 \item $\gdFin G < \infty$.
 \item There exists no manifold model for $\EFin G$.
\end{enumerate} 
\end{Theorem}

\begin{proof}[Proof of Theorem \ref{theorem:BPD groups with no manifold model}]
Let $G$ be one of the groups constructed by Block and Weinberger in the theorem above.  Since $K$ has a cocompact model for $\EFin K$ it has finitely many conjugacy classes of finite subgroups hence $G$ has finitely many conjugacy classes of finite subgroups, since all torsion in $G$ is contained in $K$.  Let $H$ be a finite subgroup of $G$, so $H$ is necessarily a subgroup of $K$ and the normaliser $N_KH$ is finite index in $N_GH$.  Since there is a cocompact model for $\EFin K$, the normaliser $N_KH$ is $\FP_\infty$ and finitely presented \cite[Theorem 0.1]{LuckMeintrup-UniversalSpaceGrpActionsCompactIsotropy} hence $N_GH$ and $W_GH$ are $\FP_\infty$ and finitely presented too \cite[VIII.5.1]{Brown}\cite[2.2.5]{Robinson}.  Using Corollary \ref{cor:uFPn equivalent conditions}, $G$ is of type $\OFinFP$.  

Finally, using \cite[III.(6.5)]{Brown}, there is a chain of isomorphisms for all natural numbers $i$,
\begin{align*}
H^i(W_GH, R[W_GH]) &\cong H^i(N_GH, R[N_GH]) \\
&\cong H^i(N_KH, R[N_KH]) \\
&\cong H^i(W_KH, R[W_KH])
\end{align*}
proving that the Weyl groups of finite subgroups have the correct cohomology.
\end{proof}

\begin{Remark}
 Although it doesn't appear in the statements of \cite[Theorems 1.5, 1.8]{BlockWeinberger-GeneralizedNielsenRealizationProblem}, Block and Weinberger do prove that there is a cocompact model for $\EFin G$, in their notation this is the space $\widetilde{X}$.
\end{Remark}

\subsection{Actions on \texorpdfstring{$R$}{R}-homology manifolds}\label{subsection:actions on homology manifolds}

Following \cite{DicksLeary-SubgroupsOfCoxeterGroups} we define an \emph{$R$-homology $n$-manifold} to be a locally finite simplicial complex $M$ such that the link $\sigma$ of every $i$-simplex of $M$ satisfies \index{R-homology manifold@$R$-homology manifold}
\[
H_j(\sigma, R) = \left\{ \begin{array}{l l} R & \text{ if } j = n - i - 1 \text{ or } j = 0, \\ 0 & \text{ else,} \end{array} \right. 
\]
for all $i$ such that $n - i - 1 \ge 0$ and the link is empty if $n - i - 1 < 0$.  In particular $M$ is an $n$-dimensional simplicial complex.  $M$ is called \emph{orientable} if we can choose an orientation for each $n$-simplex which is consistent along the $(n-1)$-simplices and we say that $M$ is \emph{$R$-orientable} if either $M$ is orientable or if $R$ has characteristic $2$.\index{R-orientable@$R$-orientable space}

A topological space $X$ is called $R$-acyclic if the reduced homology $\tilde{H}_*(X, R)$ is trivial. \index{R-acyclic@$R$-acyclic space}

\begin{Theorem}\label{theorem:proper actions on Racyclic Rhomology then BPD}
 If $G$ is a group acting properly and cocompactly on an $R$-acyclic $R$-orientable $R$-homology $n$-manifold $M$ then
\[ 
H^i(G, RG) \cong \left\{ \begin{array}{l l} R & \text{ if } i = n, \\ 0 & \text{ else.}\end{array} \right. 
\]
\end{Theorem}
\begin{proof}
 By \cite[Lemma F.2.2]{Davis} $ H^i(G, RG) \cong H^i_c(M, R) $, where $H^i_c$ denotes cohomology with compact supports.  By Poincar\'e duality for $R$-orientable $R$-homology manifolds (see for example \cite[Theorem 5]{DicksLeary-SubgroupsOfCoxeterGroups}), there is a duality isomorphism $H^i_c(M, R) \cong H_{n-i}(M, R) $.  Finally, since $M$ is assumed $R$-acyclic,
\[
H_{n-i}(M, R) \cong \left\{ \begin{array}{l l} R & \text{ if } i = n, \\ 0 & \text{ else.}\end{array} \right.
\]
\end{proof}

\begin{Example}\label{example:DicksLeary example PD over R not Z}
In \cite[Example 3]{DicksLeary-SubgroupsOfCoxeterGroups}, Dicks and Leary construct a group which is Poincar\'e duality over $R$, arising from an action on an $R$-orientable $R$-acyclic $R$-homology manifold, but which is not Poincar\'e duality over $\ZZ$.  Here $R$ may be any ring for which a fixed prime $q$ is invertible, for example $R = \FF_p$ for $p \neq q$ or $R = \QQ$.
\end{Example}

\begin{Cor}\label{cor:actions on Rorient Rhom manifolds gives BPD}
 Let $G$ be a group which admits a cocompact model $X$ for $\EFin G$ such that for every finite subgroup $H$ of $G$, $X^H$ is an $R$-orientable $R$-acyclic $R$-homology manifold, then $G$ is Bredon--Poincar\'e duality over $R$.
\end{Cor}

\begin{Remark}
 In the case $R = \ZZ$ we can drop the condition that $M$ be orientable since this is implied by being acyclic.  This is because if $M$ is acyclic then $\pi_1(M)$ is perfect, thus $\pi_1(M)$ has no normal subgroups of prime index, in particular $M$ has no index $2$ subgroups.  But if $M$ were non-orientable then the existence of an orientable double cover (see for example \cite[p.234]{Hatcher}) would imply that $\pi_1(M)$ has a subgroup of index $2$.
\end{Remark}

Let $p$ be a prime and $\FF_p$ the field of $p$ elements.  A consequence of Smith theory is the following theorem, for background on Smith theory see \cite[\S III]{Bredon-IntroductionToCompactTransformationGroups}

\begin{Theorem}[{\cite[\S 5 Theorem 2.2]{Borel-SeminarOnTransformationGroups}\cite[10.4.3]{Davis}}]\label{thm:smith}
 If $G$ is a finite $p$-group acting properly on an $\FF_p$-homology manifold $M$ then the fixed point set $M^G$ is also an $\FF_p$-homology manifold.  If $p \neq 2$ then $M^G$ has even codimension in $M$.
\end{Theorem}

\begin{Cor}[Actions on homology manifolds]\label{cor:actions on hom manifolds}~
\begin{enumerate}
 \item Let $G$ have an $n$-dimensional $\FF_p$-homology manifold model $M$ for $\EFin G$.  If $H$ is a finite $p$-subgroup of $G$ then $M^H$ is an $\FF_p$-homology manifold.  In particular if all finite subgroups of $G$ are $p$-groups and $M$ is cocompact then $G$ is Bredon--Poincar\'e duality over $\FF_p$.  If $ p \neq 2$ and $H$ is a finite $p$-subgroup of $G$ then $n - d_H$ is even.
 \item Let $G$ have an $n$-dimensional $\ZZ$-homology manifold model $M$ for $\EFin G$.  If $p \neq 2$ is a prime and $H$ is a finite $p$-subgroup of $G$ such that $M^H$ is a $\ZZ$-homology manifold then $n - \dim M^H$ is even.
\end{enumerate}
\end{Cor}

\begin{Remark}
 Given a group $G$ with subgroup $H$ which is not of prime power order, looking at the Sylow $p$-subgroups can give further restrictions.  For example if $P_i$ for $i \in I$ is a set of Sylow $p$-subgroups of $H$, one for each prime $p$, then $G$ is generated by the $P_i$ \cite[Ex. 4.10]{Rotman-Groups}.  Thus if $G$ acts on an $R$-homology manifold then the fixed points of $H$ are exactly the intersection of the fixed points of the $P_i$.  
\end{Remark}

\subsection{One relator groups}

The following lemma is adapted from \cite[5.2]{BieriEckmann-HomologicalDualityGeneralizingPoincareDuality}.
\begin{Lemma}\label{lemma:tech lemma for one-relator}
 If $G$ is $\FP_2$ with $\cd G = 2$ and $H^1(G, \ZZ G) = 0$ then $G$ is a duality group.
\end{Lemma}
\begin{proof}
We must show that $H^2(G, \ZZ G)$ is a flat $\ZZ$-module.  Consider the short exact sequence of $\ZZ G$ modules
\[
0 \longrightarrow \ZZ G \stackrel{\times p}{\longrightarrow} \ZZ G \longrightarrow \FF_pG \longrightarrow 0.
\]
This yields a long exact sequence
\[
\cdots \longrightarrow H^1(G, \FF_p G) \longrightarrow H^2(G, \ZZ G) \stackrel{\times p}{\longrightarrow} H^2(G, \ZZ G) \longrightarrow \cdots.
\]
By \cite[Corollary 3.6]{Bieri-HomDimOfDiscreteGroups}, $H^1(G, \FF_p G) \cong  H^1(G, \ZZ G) \otimes_\ZZ \FF_p = 0$.  Hence the map $H^2(G, \ZZ G) \stackrel{\times p}{\longrightarrow} H^2(G, \ZZ G)$ must have zero kernel for all $p$, in other words $H^2(G, \ZZ G)$ is torsion-free, but the torsion-free $\ZZ$-modules are exactly the flat $\ZZ$-modules.
\end{proof}

Let $G$ be a one-relator group---a group admitting a presentation of finitely many generators and one relator (see \cite[\S 5]{LyndonSchupp} for background), then:
\begin{enumerate}
 \item $G$ is $\OFinFP$ and $\OFincd_\ZZ G = 2$ \cite[4.12]{Luck-SurveyOnClassifyingSpaces}.
 \item $G$ contains a torsion-free subgroup $Q$ of finite index \cite{FischerKarrassSolitar-OneRelatorGroupsHavingElementsOfFiniteOrder}.
  \item The normaliser of every non-trivial finite subgroup $F$ is finite.  In fact, every such $F$ is subconjugate to a finite cyclic self-normalising subgroup $C$ of $G$, and furthermore the normaliser $N_GF$ is subconjugate to $C$ \cite[II.5.17,II.5.19]{LyndonSchupp}.  Thus, 
 \[H^i(N_GF, \ZZ [N_G F]) = \left\{\begin{array}{l l}0 & \text{ if $i > 0$,} \\ \ZZ & \text{ if $i = 0$.} \end{array} \right.\]
\end{enumerate}
Assume further that $H^1(G, \ZZ G) = 0$.

If $\cd_{\ZZ} Q \le 1$ then $Q$ is either trivial or a finitely generated free group so $G$ is either finite or virtually finitely generated free.  Thus $G$ is Bredon duality over $\ZZ$ by \ref{lemma:duality D0}, \ref{remark:duality uD1}, and \ref{prop:duality uD1 Rtf}.  Assume therefore that $\cd_{\ZZ} Q = 2$.  Being finite index in $G$, $Q$ is also $\FP_2$ and $H^1(Q, \ZZ Q) = H^1(G, \ZZ G) = 0$ \cite[III.(6.5)]{Brown}, thus by Lemma \ref{lemma:tech lemma for one-relator} $Q$ is a duality group and $G$ is virtual duality.  Combining with (iii) above, $G$ is also Bredon duality of dimension $2$.

 In summary:
 
\begin{Prop}
 If $G$ is a one relator group with $H^1(G, \ZZ G) = 0$ then $G$ is Bredon duality over any ring $R$.
\end{Prop}

\begin{Remark}
If $G$ is a one relator group with $H^1(G, \ZZ G) \neq 0$ then, since $G$ is $\OFinFP_0$, $G$ has bounded orders of finite subgroups by Proposition \ref{prop:uFP_0 iff finitely many conj classes of finite subgroups}.  By a result of Linnell, $G$ admits a decomposition as the fundamental group of a finite graph of groups with finite edge groups and vertex groups $G_v$ satisfying $H^1(G_v, \ZZ G_v) = 0$ \cite{Linnell-OnAccessibilityOfGroups}.  These vertex groups are subgroups of virtually torsion-free groups so in particular virtually torsion-free with $\OFincd_\ZZ G_v \le 2$.  Lemma \ref{lemma:duality FP2 split is FP2} below gives that the vertex groups are $\FP_2$ and Lemma \ref{lemma:tech lemma for one-relator} shows that the edge groups are virtually duality.  
\end{Remark}

\begin{Lemma}\label{lemma:duality FP2 split is FP2}
 Let $G$ be a group which splits as a finite graph of groups with finite edge groups $G_e$, indexed by $E$, and vertex groups $G_v$, indexed by $V$.  Then if $G$ is $\FP_2$, so are the vertex groups $G_v$.
\end{Lemma}
\begin{proof}
Fix a vertex group $G_v$.  Let $M_\lambda$, for $\lambda \in \Lambda$, be a directed system of $\ZZ G_v$ modules with $\varinjlim M_\lambda = 0$.  To use the Bieri--Eckmann criterion \cite[Theorem 1.3]{Bieri-HomDimOfDiscreteGroups}, we must show that $\varinjlim H^i(G_v, M_\lambda) = 0$ for $i = 1,2$.
 
 The Mayer--Vietoris sequence associated to the graph of groups is 
 \[ 
 \cdots \longrightarrow H^i(G, -) \longrightarrow \bigoplus_{u \in V} H^i(G_u, -) \longrightarrow \bigoplus_{e \in E} H^i(G_e, -) \longrightarrow \cdots. 
 \]
 Now $\varinjlim M_\lambda = 0$, so $\varinjlim \Ind_{\ZZ G_v}^{\ZZ G} M_\lambda = 0$ as well.  Evaluating the Mayer--Vietoris sequence at $\Ind_{\ZZ G_v}^{\ZZ G} M_\lambda$, taking the limit, and using the Bieri--Eckmann criterion, implies 
 \[\varinjlim_{\Lambda} \bigoplus_{u \in V} H^i(G_u, \Ind_{\ZZ G_v}^{\ZZ G} M_\lambda ) = 0.\]
In particular $\varinjlim H^i(G_v, \Ind_{\ZZ G_v}^{\ZZ G} M_\lambda ) = 0$ and because, as $\ZZ G_v$-module, $M_\lambda$ is a direct summand of $\Ind_{\ZZ G_v}^{\ZZ G} M_\lambda$ \cite[VII.5.1]{Brown}, this implies $\varinjlim H^i(G_v, M_\lambda ) = 0$.  
\end{proof}

\subsection{Discrete subgroups of Lie groups}
If $L$ is a Lie group with finitely many path components, $K$ a maximal compact subgroup and $G$ a discrete subgroup then $L/K$ is a model for $\EFin G$.  The space $L/K$ is a manifold and the action of $G$ on $L/K$ is smooth so the fixed point subsets of finite groups are submanifolds of $L/K$, using Lemma \ref{Lemma:proper and locally linear => fixed points are submanifolds}.  If we assume that the action is cocompact then $G$ is seen to be of type $\OFinFP$, $\OFincd G = \dim L/K$ and $G$ is a Bredon--Poincar\'e duality group.  See \cite[Theorem 5.24]{Luck-SurveyOnClassifyingSpaces} for a statement of these results.

\begin{Example}\label{example:disc sub of lie not vtf}
In \cite{Raghunathan-TorsionInCoCompactLatticesInSpin,Raghunathan-CorrigendumTorsionInCoCompactLatticesInSpin}, examples of cocompact lattices in finite covers of the Lie group $\text{Spin}(2,n)$ are given which are not virtually torsion-free. 
\end{Example}

\subsection{Virtually soluble groups}
For $G$ a soluble group the \emph{Hirsch length} $hG$ is the sum of the torsion-free ranks of the factors in an abelian series \cite[p.422]{Robinson}.  Hillman later extended this definition to elementary amenable groups \cite{Hillman-ElementaryAmenableGroupsAnd4Manifolds}. \index{hG@$hG$, Hirsch length}

Much of the observation below appears in \cite[Example 5.6]{MartinezPerez-EulerClassesAndBredonForRestrictedFamilies}.

Torsion-free soluble groups of type $\FP_\infty$ are duality \cite{Kropholler-CohomologicalDimensionOfSolubleGroups}.  We combine this with \cite{MartinezPerezNucinkis-VirtuallySolubleGroupsOfTypeFPinfty}, that virtually soluble groups of type $\FP_\infty$ are $\OFinFP$ with $\OFincd G = hG$, and deduce that if $G$ is a virtually soluble duality group (equivalently virtually soluble of type $\FP_\infty$) then $G$ virtually duality of type $\OFinFP$ with $\OFincd G = hG$.  We claim $G$ is also Bredon duality, so we must check the cohomology condition on the Weyl groups.  Since $G$ is $\OFinFP$, the normalisers $N_G F$ of any finite subgroup $F$ of $G$ are $\FP_\infty$ (Corollary \ref{cor:uFPn equivalent conditions}). Subgroups of virtually-soluble groups are virtually-soluble \cite[5.1.1]{Robinson}, so the normalisers $N_GF$ are virtually-soluble $\FP_\infty$ and hence virtually duality, and so the Weyl groups satisfy the required condition on cohomology.

Additionally, if $G$ is a virtually soluble Poincar\'e duality group then we claim $G$ is Bredon--Poincar\'e duality.  By \cite[Theorem 9.23]{Bieri-HomDimOfDiscreteGroups}, $G$ is virtually-polycyclic.  Subgroups of virtually-polycyclic groups are virtually-polycyclic \cite[p.52]{Robinson}, so $N_GF$ is polycyclic $\FP_\infty$ for all finite subgroups $F$ and, since polycylic groups are Poincar\'e duality,
\[
H^{d_F}(N_G F, \ZZ [N_G F]) = \ZZ.
\]

\begin{Prop}\label{prop:bredon equivalent to soluble virtual duality}The following conditions on a virtually-soluble group $G$ are equivalent:
\begin{enumerate}
 \item $G$ is $\FP_\infty$.
 \item $G$ is virtually duality.
 \item $G$ is virtually torsion-free and $\vcd G = hG < \infty$.
 \item $G$ is Bredon duality.
\end{enumerate}
Additionally, if $G$ is Bredon duality then $G$ is virtually Poincar\'e duality if and only if $G$ is virtually-polycyclic if and only if $G$ is Bredon--Poincar\'e duality.
\end{Prop}
\begin{proof}
 The equivalence of the first three is \cite{Kropholler-CohomologicalDimensionOfSolubleGroups} and \cite{Kropholler-OnGroupsOfTypeFP_infty}.  The rest is the discussion above.
\end{proof}

\subsection{Elementary amenable groups}

If $G$ is an elementary amenable group of type $\FP_\infty$ then $G$ is virtually soluble \cite[p.4]{KMN-CohomologicalFinitenessConditionsForElementaryAmenable}, in particular Bredon duality over $\ZZ$ of dimension $hG$.  The converse, that every elementary amenable Bredon duality group is $\FP_\infty$, is obvious.

If $G$ is elementary amenable $\FP_\infty$ then the condition $H^n(G, \ZZ G) \cong \ZZ$ implies that $G$ is Bredon--Poincar\'e duality, so for all finite subgroups,
\[
H^{d_F}(N_GF, \ZZ [N_GF]) \cong \ZZ.
\]
A natural question is whether 
\[
H^{d_F}(N_G F, \ZZ [N_G F]) = \ZZ 
\]
can ever occur for an elementary amenable, or indeed a soluble Bredon duality, but not Bredon--Poincar\'e duality group.  An example of this behaviour is given below.

\begin{Example}
We construct a finite index extension of the Baumslag--Solitar group $BS(1,p)$, for $p$ a prime.
\[
BS(1,p) = \langle x,y  :  y^{-1}xy = x^p\rangle 
\]
This has a normal series \cite[p.60]{LennoxRobinson},
\[
1 \unlhd \langle x \rangle \unlhd \langle \langle x \rangle \rangle \unlhd  BS(1,p),
\]
where $\langle \langle x \rangle \rangle$ denotes the normal closure of $x$.  The quotients of this normal series are $\langle x \rangle / 1 \cong \ZZ$, $\langle \langle x \rangle \rangle / \langle x \rangle \cong C_{p^\infty}$ and $BS(1,p)/\langle \langle x \rangle \rangle \cong \ZZ$, where $C_{p^\infty}$ denotes the Pr\"ufer group (see \cite[p.94]{Robinson}).  Clearly $BS(1,p)$ is finitely generated torsion-free soluble with $hBS(1,p) = 2$, but not polycyclic, since $C_{p^\infty}$ does not have the maximal condition on subgroups \cite[5.4.12]{Robinson}, thus $BS(1,p)$ is not Poincar\'e duality.  Also since $BS(1,p)$ is an HNN extension of $\langle x \rangle \cong \ZZ$ it has cohomological dimension $2$ \cite[Proposition 6.12]{Bieri-HomDimOfDiscreteGroups} and thus $\cd BS(1,p) = hBS(1,p)$.  By Proposition \ref{prop:bredon equivalent to soluble virtual duality}, $BS(1,p)$ is a Bredon duality group.

Recall that elements of $BS(1,p)$ can be put in a normal form: $y^i x^k y^{-j} $ where $i,j\ge 0$ and if $i,j > 0$ then $n \nmid k$.  Consider the automorphism $\varphi$ of $BS(1,p)$, sending $x \mapsto x^{-1}$ and $y \mapsto y$, an automorphism since it is its own inverse and because the relation $ y^{-1}xy = x^p $ in $BS(1,p)$ implies the relation $y^{-1}x^{-1}y = x^{-p} $.  Let  $y^i x^k, y^{-j} $ be an element in normal form, then
\[
\varphi :y^i x^k y^{-j} \longmapsto y^i x^{-k} y^{-j}.
\]
So the only fixed points of $\varphi$ are in the subgroup $\langle y \rangle \cong \ZZ$.  Form the extension 
\[
1 \longrightarrow BS(1,p) \longrightarrow G \longrightarrow C_2 \longrightarrow 1 
\]
where $C_2$ acts by the automorphism $\varphi$.  The property of being soluble is extension closed \cite[5.1.1]{Robinson}, so $G$ is soluble virtual duality and Bredon duality by Proposition \ref{prop:bredon equivalent to soluble virtual duality}.  The normaliser 
\[
N_G C_2 = C_G C_2 = \{ g \in G  :  gz = zg \text{ for the generator }z \in C_2 \}
\]
is exactly $\langle y \rangle \times C_2 \cong \ZZ \times C_2$.  Thus $W_GC_2 \cong \ZZ$ and $H^1(W_GC_2, \ZZ[W_GC_2]) \cong \ZZ$.  Since $BS(1,p)$ is finite index in $G$, by \cite[III.(6.5)]{Brown}
\[
H^2(G, \ZZ G) \cong H^2(BS(1,p), \ZZ [BS(1,p)]).
\]
However since $BS(1,p)$ is not Poincar\'e duality, $ H^n(BS(1,p), \ZZ [BS(1,p)])$ is $\ZZ$-flat but not isomorphic to $\ZZ$.
\end{Example}

\begin{Remark}
Baues \cite{Baues-Infrasolvmanifolds} and Dekimpe \cite{Dekimpe-PolycyclicGroupAdmitsNILAffine} proved independently that any virtually polycyclic group $G$ can be realised as a NIL affine crystallographic group---$G$ acts properly, cocompactly, and by diffeomorphisms on a simply connected nilpotent Lie group of dimension $hG$.  Any connected, simply connected nilpotent Lie group is diffeomorphic to some Euclidean space \cite[\S I.16]{Knapp} and hence contractible, so any elementary amenable Bredon--Poincar\'e duality group has a cocompact manifold model for $\EFin G$.
\end{Remark}

\section{Finite extensions of right-angled Coxeter groups}\label{section:reflection groups}\index{Reflection group trick}

Recall Corollary \ref{cor:actions on hom manifolds}, that if $G$ has a cocompact $n$-dimensional $\ZZ$-homology manifold model $M$ for $\EFin G$ such that all fixed point sets $M^H$ are $\ZZ$-homology manifolds, and if $H$ is a finite $p$-subgroup of $G$ with $p \neq 2$ then $n - d_H$ is even.  In this section we construct Bredon--Poincar\'e duality groups $G$ over $\ZZ$ of arbitrary dimension such that, for any fixed prime $p \neq 2$:
\begin{enumerate}
 \item All of the finite subgroups of $G$ are $p$-groups.
 \item $\mathcal{V}(G)$ is any set with $n-d_H$ even for all finite subgroups $H$.
\end{enumerate}

The method of constructing these examples was recommended to us by Ian Leary and utilises methods from \cite[\S 2]{DavisLeary-DiscreteGroupActionsOnAsphericalManifolds} and \cite[\S 11]{Davis}.  We write ``$\Gamma$'' instead of ``$W$'', as is used in \cite{DavisLeary-DiscreteGroupActionsOnAsphericalManifolds}, to denote a Coxeter group so the notation can't be confused with our use of $W_GH$ for the Weyl group.

Let $M$ be an compact contractible $n$-manifold with boundary $\partial M$, such that $\partial M$ is triangulated as a flag complex.  Let $G$ be a group acting on $M$ such that the induced action on the boundary is by simplicial automorphisms.  Let $\Gamma$ be the right angled Coxeter group corresponding to the flag complex $\partial M$, the group $G$ acts by automorphisms on $\Gamma$ and we can form the semi-direct product $\Gamma \rtimes G$.  Moreover there is a space $\mathcal{U} = \mathcal{U}(M, \partial M, G)$ such that: \index{Davis complex $\mathcal{U}$}
\begin{enumerate}
 \item $\mathcal{U}$ is a contractible $n$-manifold without boundary.
 \item $\Gamma \rtimes G$ acts properly and cocompactly on $\mathcal{U}$.
 \item For any finite subgroup $H$ of $G$, we have $\mathcal{U}^H = \mathcal{U}(M^H, (\partial M)^H, W_GH)$, in particular $\dim \mathcal{U}^H = \dim M^H$.
 \item $W_{\Gamma \rtimes G}H = \Gamma^H \rtimes W_GH$, where $\Gamma^H$ is the right-angled Coxeter group associated to the flag complex $(\partial M)^H$.
 \item If $M$ is the cone on a finite complex then $\mathcal{U}$ has a CAT(0) cubical structure ($1$-connected cubical complex whose links are simplicial flag complexes) such that the action of $\Gamma \rtimes G$ is by isometries. \index{CAT(0) cubical complex} 
\end{enumerate}

Every Coxeter group contains a finite-index torsion-free subgroup \cite[Corollary D.1.4]{Davis}, let $\Gamma^\prime$ denote such a subgroup of $\Gamma$ and assume that $\Gamma^\prime$ is normal.  Then $\Gamma^\prime \rtimes G$ is finite index in $\Gamma \rtimes G$ and so acts properly and cocompactly on $\mathcal{U}$ also.

\begin{Lemma}\label{lemma:arbVG construction is uEG}
 If $M$ is the cone on a finite complex then $\mathcal{U}$ is a cocompact model for $\EFin (\Gamma \rtimes G)$ and for $\EFin (\Gamma^\prime \rtimes G)$.  In particular, $\Gamma^\prime \rtimes G$ is of type $\OFinFP$.
\end{Lemma}

\begin{proof}
 A CAT(0) cubical complex has a CAT(0) metric \cite[Remark 2.1]{Wise-FromRichesToRaags} and any contractible CAT(0) space on which a group acts properly is a model for the classifying space for proper actions \cite[Corollary II.2.8]{BridsonHaefliger} (see also \cite[Theorem 4.6]{Luck-SurveyOnClassifyingSpaces}).  
\end{proof}

\begin{Lemma}\label{lemma:arbVG construction, conj of finite subgroups}
 Let $M$ be the cone on the finite complex $\partial M$.  If $K$ is a finite subgroup of $\Gamma^\prime \rtimes G$ then $K$ is subconjugate in $\Gamma \rtimes G$ to $G$.
\end{Lemma}
\begin{proof}
Since, by Lemma \ref{lemma:arbVG construction is uEG}, $\mathcal{U}$ is a model for $\EFin (\Gamma \rtimes G)$, the finite subgroup $K$ necessarily fixes a vertex $v$ of $\mathcal{U}$ and hence is a subgroup of the stabiliser of $v$. 

Recall from \cite[\S 5]{Davis} and Section \ref{section:bredon interesting examples} that 
\[
 \mathcal{U} = \Gamma \times M / \sim
\]
where the identification is along $\Gamma \times \partial M$ only and the action of $\Gamma \rtimes G$ on $\mathcal{U}$ is given by 
\[
(\gamma^\prime, g) \cdot (\gamma, m) = (\gamma^\prime\gamma, gm).
\]

A fundamental domain for the $\Gamma$-action is the copy of $M$ inside $\mathcal{U}$ given by $1 \times M$ and as such the stabiliser of any vertex is conjugate via an element of $\Gamma$ to the stabiliser of a vertex $v^\prime$ in $1 \times M$.  Finally, the only elements from $\Gamma^\prime \rtimes G$ stabilising $v^\prime \in 1\times M$ are contained in $G$ ($\Gamma^\prime$ moves $M$ about $\mathcal{U}$ freely, whereas $G$ stabilises $M$ setwise). 
\end{proof}

\begin{Theorem}\label{theorem:arbVG general setup}
Let $G$ be a finite group with real representation $\rho : G \longhookrightarrow \text{GL}_n \RR$ and, for all subgroups $H$ of $G$, let $d_H$ denote the dimension of the subspace of $\RR^n$ fixed by $H$.  Then there exists a Bredon--Poincar\'e duality group $\Gamma^\prime \rtimes G$ of dimension $n$ such that
\[
\mathcal{V} (\Gamma^\prime \rtimes G) = \{ d_H  :  H \le G\}.
\]
\end{Theorem}
\begin{proof}
 Restrict $\rho$ to an action on $(D^n, S^{n-1})$ and choose a $G$-equivariant flag triangulation of $S^{n-1}$ (use, for example, \cite{Illman-SmoothEquiTriangulationsForGFinite}).  We obtain a Coxeter group $\Gamma$, normal finite-index torsion-free subgroup $\Gamma^\prime$, and space $\mathcal{U}$.  Lemma \ref{lemma:arbVG construction is uEG} gives that $\Gamma^\prime \rtimes G$ is of type $\OFinFP$.  

 Since $\Gamma^\prime \rtimes G$ has an $n$-dimensional model for $\EFin \Gamma^\prime \rtimes G$ we have that $\gdFin \Gamma^\prime \rtimes G \le n$, by Lemma \ref{Lemma:observations}(1) $\cd_\QQ \Gamma^\prime \rtimes G = d_{1}$, and by Theorem \ref{theorem:proper actions on Racyclic Rhomology then BPD} and Corollary \ref{cor:actions on Rorient Rhom manifolds gives BPD} $d_{1} = n$.  Using the chain of inequalities 
 \[
 n = \cd_\QQ \Gamma^\prime \rtimes G \le \OFincd \Gamma^\prime \rtimes G \le \gdFin \Gamma^\prime \rtimes G \le n,
 \]
 shows that $\OFincd \Gamma^\prime \rtimes G = n$.  It remains only to check the condition on the cohomology of the Weyl groups of the finite subgroups.
 
 For any finite subgroup $H$ of $G$, the Weyl group $W_{\Gamma^\prime \rtimes G} H$ acts properly and cocompactly on $\mathcal{U}(M^H, (\partial M)^H, W_G H)$ which is a contractible $d_H$-manifold without boundary.  By Theorem \ref{theorem:proper actions on Racyclic Rhomology then BPD},
 \[
 H^n(W_{\Gamma^\prime \rtimes G} H, \ZZ[W_{\Gamma^\prime \rtimes G} H]) = \left\{ \begin{array}{l l}\ZZ & \text{ if } i = d_H, \\ 0 & \text{else.} \end{array} \right.  
 \]
 
 If $K$ is any finite subgroup of $\Gamma^\prime \rtimes G$ then, by Lemma \ref{lemma:arbVG construction, conj of finite subgroups}, $K$ is conjugate in $\Gamma \rtimes G$ to some $H \le G$.  In particular the normalisers of $H$ and $K$ in $\Gamma \rtimes G$ are isomorphic.  Also, since $\Gamma^\prime \rtimes G$ is finite index in $\Gamma \rtimes G$, the normaliser $N_{\Gamma^\prime \rtimes G}K$ is finite index in $N_{\Gamma \rtimes G}K$, thus:
\begin{align*} 
H^n(N_{\Gamma^\prime \rtimes G}K, \ZZ[N_{\Gamma^\prime \rtimes G}K]) 
&\cong H^n(N_{\Gamma \rtimes G}K, \ZZ[N_{\Gamma \rtimes G}K]) \\
&\cong H^n(N_{\Gamma \rtimes G}H, \ZZ[N_{\Gamma \rtimes G}H]) \\
&\cong H^n(N_{\Gamma^\prime \rtimes G}H, \ZZ[N_{\Gamma^\prime \rtimes G}H]).
\end{align*}
From the short exact sequence
\[
1 \longrightarrow K \longrightarrow N_{\Gamma^\prime \rtimes G} K \longrightarrow W_{\Gamma^\prime \rtimes G} K \longrightarrow 1
\]
and Lemma \ref{lemma:Finite by Q implies HG is HQ}, 
\[
H^n(N_{\Gamma^\prime \rtimes G}H, \ZZ[N_{\Gamma^\prime \rtimes G}H]) \cong H^n(W_{\Gamma^\prime \rtimes G}K, \ZZ[W_{\Gamma^\prime \rtimes G}K]).
\]
Thus,
\[
H^n(W_{\Gamma^\prime \rtimes G} K, \ZZ[W_{\Gamma^\prime \rtimes G} K]) = \left\{ \begin{array}{l l}\ZZ & \text{ if } i = d_H, \\ 0 & \text{else.} \end{array} \right.  
\]
\end{proof}

\begin{Example}\label{example:arb VG for pgroups}
We construct a group using Theorem \ref{theorem:arbVG general setup} with the properties mentioned at the beginning of this section.  It will be of the form $\Gamma^\prime \rtimes C_{p^m}$, where $C_{p^m}$ is the cyclic group of order $p^m$.

For $i$ between $1$ and $m$ let $w_i$ be any collection of positive integers and let $n = \sum_i 2w_i$. 
If $c$ is a generator of the cyclic group $C_{p^m}$, then $C_{p^m}$ embeds into the orthogonal group $O(n)$ via the real representation
 \begin{align*}
 \rho: C_{p^m} &\longhookrightarrow O(n) \\
 c &\longmapsto (R_{2 \pi/p})^{\oplus w_1} \oplus (R_{2 \pi/p^2})^{\oplus  w_2} \oplus \cdots \oplus (R_{2 \pi/p^m})^{\oplus w_m}
 \end{align*}
where $R_\theta$ is the the $2$-dimensional rotation matrix of angle $\theta$.  The image is in $O(n)$ since we chose $n$ such that $2w_1 + \cdots 2w_n = n$.
 
If $i$ is some integer between $1$ and $m$ then there is a unique subgroup $C_{p^{m-i+1}}$ of $C_{p^m}$ with generator $c^{p^{i}}$, in fact this enumerates all subgroups of $C_{p^m}$ except the trivial subgroup.  Under $\rho$, this generator maps to
\[
\rho : c^{p^{i}} \longmapsto R_0^{\oplus w_1} \oplus \cdots \oplus R_0^{\oplus w_i} \oplus \left( R_{p^{i} 2 \pi / p^{i+1} } \right)^{\oplus w_{i+1}} \oplus \cdots \oplus \left( R_{p^{i} 2 \pi / p^{m} } \right)^{\oplus w_n}.
\]
In other words, the fixed point set corresponding to $C_{p^{m-i+1}}$ is $\RR^{2w_1 + \cdots + 2w_i}$.  Thus the set of dimensions of the fixed point subspaces of non-trivial finite subgroups of $C_{p^m}$ are
 \[
 \{2w_1, 2(w_1 + w_2), \ldots, 2(w_1 + w_2 + \ldots + w_{m-1}) \}.
 \]
Applying Theorem \ref{theorem:arbVG general setup} gives a group $\Gamma^\prime \rtimes C_{p^m}$ of type $\OFinFP$ with 
\[\OFincd G = n = \sum^m_{i=1} 2w_i\]
and such that 
\[
\mathcal{V}( \Gamma^\prime \rtimes C_{p^m} ) = \{2w_1, 2(w_1 + w_2), \ldots, 2(w_1 + w_2 + \ldots + w_{m-1}) \}.
\]
Since there were no restrictions on the integers $w_i$, using this technique we can build an even dimensional Bredon--Poincar\'e duality group with any $\mathcal{V}(G)$, as long as all the integers $d_H$ are even.

The case $n$ is odd reduces to the case $n$ is even.  Proposition \ref{prop:duality direct products of BPD groups} shows that if a group $G$ is Bredon--Poincar\'e duality then taking the direct product with $\ZZ$ gives a Bredon--Poincar\'e duality group $G \times \ZZ$ where
\[
\mathcal{V}(G \times \ZZ) \cong \{ v + 1  :  v \in \mathcal{V}(G)\}.
\]
Thus we can build a group with odd $n_H$ and $\mathcal{V}$ containing only odd elements by building a group with even $n_H$ and then taking a direct product with $\ZZ$.
\end{Example}

\section{Low dimensions}\label{section:BD low dimensions}

This section is devoted to the study of Bredon duality groups and Bredon--Poincar\'e duality groups of low dimension.  We completely classify those of dimension $0$ in Lemma \ref{lemma:duality D0}.  We partially classify those of dimension 1---see Propositions \ref{prop:duality uD1 Rtf} and \ref{prop:duality uPD1 Rtf}, and Question \ref{question:duality D1 not Rtf}.  There is a discussion of the dimension $2$ case.

Recall that a group $G$ is duality of dimension $0$ over $R$ if and only if $\lvert G \rvert$ is finite and invertible in $R$, and any such group is necessarily Poincar\'e duality \cite[Proposition 9.17(a)]{Bieri-HomDimOfDiscreteGroups}.

\begin{Lemma}\label{lemma:duality D0}
 $G$ is Bredon duality of dimension $0$ over $R$ if and only if $\vert G \vert$ is finite.  Any such group is necessarily Bredon--Poincar\'e duality.
\end{Lemma}
\begin{proof}
 By \cite[13.2.11]{Geoghegan},
\[
H^0(G, R G) = \left\{ \begin{array}{l l} R & \text{if $\vert G \vert$ is finite,} \\ 0 & \text{else.} \end{array}\right.
\]
Hence if $G$ is Bredon duality of dimension $0$ then $G$ is finite and moreover $G$ is Bredon--Poincar\'e duality.
 
Conversely, if $G$ is finite then $\OFincd_R G = 0$ and $G$ is $\OFinFP_\infty$ over $R$ (Propositions \ref{prop:ucdG=0 over R iff G finite} and \ref{prop:uFP_0 iff finitely many conj classes of finite subgroups}).  Finally the Weyl groups of any finite subgroup will be finite so by  \cite[13.2.11,13.3.1]{Geoghegan},
\[
H^n(WH, R [WH]) = \left\{ \begin{array}{c c} R & \text{ if } n = 0, \\ 0 & \text{ if } n > 0. \end{array}\right.
\]
Thus $G$ is Bredon--Poincar\'e duality of dimension $0$.
\end{proof}

The duality groups of dimension $1$ over $R$ are exactly the groups of type $\FP_1$ over $R$ (equivalently finitely generated groups \cite[Proposition 2.1]{Bieri-HomDimOfDiscreteGroups}) with $\cd_R G = 1$ \cite[Proposition 9.17(b)]{Bieri-HomDimOfDiscreteGroups}.

\begin{Prop}\label{prop:duality uD1 Rtf}
If $G$ is infinite $R$-torsion free, then the following are equivalent:
\begin{enumerate}
 \item $G$ is Bredon duality over $R$, of dimension $1$.
 \item $G$ is finitely generated and virtually-free.
 \item $G$ is virtually duality over $R$, of dimension $1$. 
\end{enumerate}
\end{Prop}
\begin{proof}
 That $2 \Rightarrow 3$ is \cite[Proposition 9.17(b)]{Bieri-HomDimOfDiscreteGroups}.  For $3 \Rightarrow 2$, let $G$ be virtually duality over $R$ of dimension $1$, then $\cd_R G \le 1$ so by \cite{Dunwoody-AccessabilityAndGroupsOfCohomologicalDimensionOne} $G$ acts properly on a tree.  Since $G$ is assumed finitely generated, $G$ is virtually-free \cite[Theorem 3.3]{Antolin-CayleyGraphsOfVirtuallyFreeGroups}.

 For $1 \Rightarrow 2$, if $G$ is Bredon duality over $R$ of dimension $1$, then $G$ is automatically finitely generated and $\OFincd_R G = 1$.  By Lemma \ref{lemma:no R torsion then cdR less ucd R etc} $\cd_R G = 1$ so, as above, by \cite{Dunwoody-AccessabilityAndGroupsOfCohomologicalDimensionOne} and \cite[Theorem 3.3]{Antolin-CayleyGraphsOfVirtuallyFreeGroups}, $G$ is virtually-free.

 For $2 \Rightarrow 1$, if $G$ is virtually finitely generated free then $G$ acts properly and cocompactly on a tree \cite[Theorem 3.3]{Antolin-CayleyGraphsOfVirtuallyFreeGroups}, so $G$ is $\OFinFP$ over $R$ with $\OFincd_R G = 1$.  As $G$ is $\OFinFP$, for any finite subgroup $K$, the normaliser $N_G K$ is finitely generated.  Subgroups of virtually-free groups are virtually-free, so $N_G K$ is virtually finitely generated free, in particular a virtual duality group \cite[Proposition 9.17(b)]{Bieri-HomDimOfDiscreteGroups}, so
\[
H^i(WK, \ZZ[WK]) = H^i(N_GK, \ZZ[N_G K]) = \left\{ \begin{array}{l l} \text{$\ZZ$-flat} & \text{ for }i = d_K, \\ 0 & \text{ else,} \end{array} \right. 
\]
where $d_K = 0$ or $1$.  Thus $G$ is Bredon duality over $\ZZ$ and hence also over $R$.
\end{proof}

\begin{Remark}\label{remark:duality uD1}
 The only place that the condition $G$ be $R$-torsion-free was used was in the implication $1 \Rightarrow 2$, the problem for groups which are not $R$-torsion-free is that the condition $\OFincd_R G \le 1$ is not known to imply that $G$ acts properly on a tree.  
 
 If we take $R=\ZZ$ then $\OFincd_\ZZ G \le 1$ implies $G$ acts properly on a tree by a result of Dunwoody \cite{Dunwoody-AccessabilityAndGroupsOfCohomologicalDimensionOne}.  Thus over $\ZZ$, $G$ is Bredon duality of dimension $1$ if and only $G$ is finitely generated virtually free, if and only if $G$ is virtually duality of dimension $1$.
\end{Remark}

\begin{Question}\label{question:duality D1 not Rtf}
 What characterises Bredon duality groups of dimension $1$ over $R$?
\end{Question}

\begin{Prop}\label{prop:duality uPD1 Rtf}
 If $G$ is infinite then the following are equivalent:
\begin{enumerate}
 \item $G$ is Bredon--Poincar\'e duality over $R$, of dimension $1$.
 \item $G$ is virtually infinite cyclic.
 \item $G$ is virtually Poincar\'e duality over $R$, of dimension $1$. 
\end{enumerate}
\end{Prop}
\begin{proof}
  The equivalence follows from the fact that for $G$ a finitely generated group, $G$ is virtually infinite cyclic if and only if $H^1(G, RG) \cong R$ \cite[13.5.5,13.5.9]{Geoghegan}.
\end{proof}

In dimension $2$ we can only classify Bredon--Poincar\'e duality groups over $\ZZ$.
The following result appears in \cite[Example 5.7]{MartinezPerez-EulerClassesAndBredonForRestrictedFamilies}, but a proof is not given there.  Recall that a surface group is the fundamental group of a compact surface without boundary. \index{Surface group}

\begin{Lemma}\label{lemma:duality vsurface is uPD2}
 If $G$ is virtually a surface group then $G$ is Bredon--Poincar\'e duality.
\end{Lemma}

\begin{proof}
As $G$ is a virtual surface group, $G$ has finite index subgroup $H$ with $H$ the fundamental group of some compact surface without boundary.  Firstly, assume $H = \pi_1(S_g)$ where $S_g$ is the orientable surface of genus $g$.  If $g = 0$ then $S_g$ is the $2$-sphere and $G$ is a finite group, thus $G$ is Bredon--Poincar\'e duality by Lemma \ref{lemma:duality D0}.  If $g > 0$ then by \cite[Lemma 4.4(b)]{Mislin-ClassifyingSpacesProperActionsMappingClassGroups} $G$ is $\OFinFP$ over $\ZZ$ with $\OFincd_\ZZ G \le 2$.  

We now treat the cases $g = 1$ and $g>1$ separately. If $g > 1$ then, in the same lemma, Mislin shows that the upper half-plane is a model for $\EFin G$ with $G$ acting by hyperbolic isometries.  Thus \cite[\S 10.1]{Davis} gives that the fixed point sets are all submanifolds, hence $G$ is Bredon--Poincar\'e duality of dimension $2$.  If $g=1$ then by \cite[Lemma 4.3]{Mislin-ClassifyingSpacesProperActionsMappingClassGroups}, $G$ acts by affine maps on $\RR^2$ so again $\RR^2$ is a model for $\EFin G$ whose fixed point sets are submanifolds, and thus $G$ is Bredon--Poincar\'e duality of dimension $2$.  

Now we treat the non-orientable case, so $H = \pi_1(T_k)$ where $T_k$ is a closed non-orientable surface of genus $k$.  In particular $T_k$ has Euler characteristic $\chi(T_k) = 2-k$.  $H$ has an index $2$ subgroup $H^\prime$ isomorphic to the fundamental group of the closed orientable surface of Euler characteristic $2\chi(T_k)$, thus $H^\prime = \pi_1(S_{k-1})$.  If $k = 1$ then $H = \ZZ/2$ and $G$ is a finite group, thus Bredon--Poincar\'e duality by Lemma \ref{lemma:duality D0}.  
Assume then that $k > 1$, we are now back in the situation above where $G$ is virtually $\pi_1(S_g) $ for $g > 0$ and as such $G$ is Bredon--Poincar\'e duality of dimension $n$, by the previous part of the proof.
\end{proof}

\begin{Prop}\label{lemma:duality vPD2ZZ iff uPD2ZZ}
The following conditions are equivalent:
\begin{enumerate}
 \item $G$ is virtually Poincar\'e duality of dimension $2$ over $\ZZ$.
 \item $G$ is virtually surface.
 \item $G$ is Bredon--Poincar\'e duality of dimension $2$ over $\ZZ$.
\end{enumerate}
\end{Prop}
\begin{proof}
 That $1 \Leftrightarrow 2$ is \cite{Eckmann-PoincareDualityGroupsOfDimensionTwoAreSurfaceGroups} and that  $2 \Rightarrow 3$ is Lemma \ref{lemma:duality vsurface is uPD2}.  The implication $3 \Rightarrow 1$ is provided by \cite[Theorem 0.1]{Bowditch-PlanarGroupsAndTheSeifertConjecture} which states that any $\FP_2$ group with $H^2(G, \QQ G) = \QQ$ is a virtual surface group and hence a virtual Poincar\'e duality group.  If $G$ is Bredon--Poincar\'e duality of dimension $2$ then $H^i(G, \QQ G) = H^i(G, \ZZ G) \otimes \QQ = \QQ$ (see proof of Lemma \ref{Lemma:observations}(1)) and $G$ is $\FP_2$ so we may apply the aforementioned theorem.
\end{proof}

The above proposition doesn't extend from Poincar\'e duality to just duality, as demonstrated by \cite[p.163]{Schneebeli-VirtualPropertiesAndGroupExtensions} where an example, based on Higman's group, is given of a Bredon duality group of dimension $2$ over $\ZZ$ which is not virtual duality.  This example is extension of a finite group by a torsion-free duality group of dimension $2$.  Schneebeli proves that the group is not virtually torsion-free, that it is Bredon duality follows from Proposition \ref{prop:finite-by-uDn is uDn}.

\begin{Question}\label{question:BD characterise BPD dim 2}
  Is there an easy characterisation of Bredon duality, or Bredon--Poincar\'e duality groups, of dimension $2$ over $R$?
\end{Question}

\section{Extensions}\label{section:BD extensions}

In the classical case, extensions of duality groups by duality groups are always duality \cite[9.10]{Bieri-HomDimOfDiscreteGroups}.  In the Bredon case the situation is more complex, for example semi-direct products of torsion-free groups by finite groups may not even be $\OFinFP_0$ \cite{LearyNucinkis-SomeGroupsOfTypeVF}.  Davis and Leary build examples of finite index extensions of Poincar\'e duality groups which are not Bredon duality, although they are $\OFinFP_\infty$ \cite[Theorem 2]{DavisLeary-DiscreteGroupActionsOnAsphericalManifolds}, and examples of virtual duality groups which are not of type $\OFinFP_\infty$ \cite[Theorem 1]{DavisLeary-DiscreteGroupActionsOnAsphericalManifolds}.  In \cite{FarrellLafont-FiniteAutomorphismsOfNegCurvedPDGroups}, Farrell and Lafont give examples of prime index extensions of $\delta$-hyperbolic Poincar\'e duality groups which are not Bredon--Poincar\'e duality.  In \cite[\S 5]{MartinezPerez-EulerClassesAndBredonForRestrictedFamilies},
 Mart\'{\i}nez-P\'erez considers $p$-power extensions of duality groups over fields of characteristic $p$, showing that if $Q$ is a $p$-group and $G$ is Poincar\'e duality of dimension $n$ over a field of characteristic $p$ then then $G \rtimes Q$ is Bredon--Poincar\'e duality of dimension $n$.  These results do not extend from Poincar\'e duality groups to duality groups however \cite[\S 6]{MartinezPerez-EulerClassesAndBredonForRestrictedFamilies}.

We study direct products of Bredon duality groups and extensions of the form finite-by-Bredon duality.

\subsection{Direct products}
\begin{Lemma}\label{lemma:bredon G uFP and H uFP then GxH is uFP}
For all groups $G_1$ and $G_2$,
\begin{enumerate}
 \item If $G_1$ and $G_2$ are $\OFinFP$ over $R$ then $G_1 \times G_2$ is $\OFinFP$ over $R$.
 \item $\OFincd_R G_1 \times G_2 \le \OFincd_R G_1 + \OFincd_R G_2$.
\end{enumerate}
\end{Lemma}
\begin{proof}
That $\OFincd_R G_1 \times G_2 \le \OFincd_R G_1 + \OFincd_R G_2$ is a special case of \cite[3.62]{Fluch}, where Fluch proves that given projective resolutions $P_*$ of $\uR$ by $\OFin$-modules for $G_1$ and $Q_*$ of $\uR$ by $\OFin$-modules for $G_2$, the total complex of the tensor product double complex is a projective resolution of $\uR$ by projective $\OFin$-modules for $G_1 \times G_2$.  So to prove that $G_1 \times G_2$ is $\OFinFP$ it is sufficient to show that if $P_*$ and $Q_*$ are finite type resolutions, then so is the total complex, but this follows from \cite[3.52]{Fluch}.
\end{proof}

\begin{Lemma}\label{lemma:finite subgroups of direct product are fi in projections}
 If $L$ is a finite subgroup of $G_1 \times G_2$ then the normaliser $N_{G_1 \times G_2} L$ is finite index in $N_{G_1} \pi_1L \times N_{G_2} \pi_2 L$, where $\pi_1$ and $\pi_2$ are the projection maps from $G_1 \times G_2$ onto the factors $G_1$ and $G_2$.
\end{Lemma}
\begin{proof}
 It's straightforward to check that
 \[
 N_{G_1 \times G_2} L \le N_{G_1} \pi_1L \times N_{G_2} \pi_2 L.
 \]
 
 Next, observe that $N_{G_1} \pi_1L \times N_{G_2} \pi_2 L$ acts by conjugation on $\pi_1 L \times \pi_2 L$ and the setwise stabiliser of $L \le (\pi_1 L \times \pi_2 L)$ is exactly $N_{G_1 \times G_2}L$.  Since $\pi_1 L \times \pi_2 L $ is finite, any stabiliser of a subset is necessarily finite-index (via the orbit-stabiliser theorem), thus $ N_{G_1 \times G_2}L$ is finite index in $ N_{G_1} \pi_1L \times N_{G_2} \pi_2 L$.
 \end{proof}

\begin{Prop}\label{prop:duality direct products of BPD groups}
 If $G_1$ and $G_2$ are Bredon duality (resp. Bredon--Poincar\'e duality), then $G = G_1 \times G_2$ is Bredon duality (resp. Bredon--Poincar\'e duality).  Furthermore,
\begin{align*}
\mathcal{V}(G_1 \times G_2) = \left\{ v_1 + v_2 : v_i \in \mathcal{V}(G_i) \right\} \cup \left\{ v_1 + d_1(G_2) : v_1 \in \mathcal{V}(G_1) \right\} \\ \cup \left\{ d_1(G_1) + v_2 : v_2 \in \mathcal{V}(G_2) \right\}.
\end{align*}
\end{Prop}
\begin{proof}
 By Lemma \ref{lemma:bredon G uFP and H uFP then GxH is uFP}, $G_1 \times G_2$ is $\OFinFP$.  If $L$ is some finite subgroup of $G$, then, via Lemma \ref{lemma:finite subgroups of direct product are fi in projections}, the normaliser $N_G L$ is finite index in $N_{G_1} \pi_1L \times N_{G_2} \pi_2 L$ so an application of Shapiro's Lemma \cite[III.(6.5) p.73]{Brown} gives that for all $i$,
\[
H^i(N_{G} L, R[N_{G} L] ) \cong H^i( N_{G_1} \pi_1L \times N_{G_2} \pi_2 L, R[N_{G_1} \pi_1L \times N_{G_2} \pi_2 L]).
\]
Noting the isomorphism of $RG$-modules 
\[
R[N_{G_1} \pi_1L \times N_{G_2} \pi_2 L] \cong R[N_{G_1} \pi_1L] \otimes R[N_{G_2} \pi_2 L], 
\]
the K\"unneth formula for group cohomology (see \cite[p.109]{Brown}) is:
\[\xymatrix@-7pt{
0 \ar[d] \\
 \bigoplus_{i + j = k} \big( H^i(N_{G_1} \pi_1L, R[N_{G_1} \pi_1L]) \otimes H^j(N_{G_1} \pi_1L, R[N_{G_1} \pi_1L]) \big) \ar[d]\\
 H^k(N_{G_1} \pi_1L \times N_{G_2} \pi_2L, R[N_{G_1} \pi_1L \times N_{G_2} \pi_2 L]) \ar[d] \\
 \makebox[.8\textwidth]{$\bigoplus_{i + j = k+1} \Tor_1^R( H^i(N_{G_1} \pi_1L,R[N_{G_1} \pi_1L] ), H^j(N_{G_2} \pi_2L, R[N_{G_2} \pi_2 L]) )$} \ar[d] \\
0}\]
Here we are using that the $R[N_{G_i} \pi_iL]$ are $R$-free.  Since $H^i(N_{G_1} \pi_1L,R[N_{G_1} \pi_1L] )$ is $R$-flat the $\Tor_1$ term is zero.  Hence the central term is non-zero only when $i = d_{\pi_1L}$ and $j = d_{\pi_2L}$, in which case it is $R$-flat.  Furthermore, $d_L = d_{\pi_1 L} + d_{\pi_2 L}$.

If $G_1$ and $G_2$ are Bredon--Poincar\'e duality then the central term in this case is $R$.

Since if $L$ is non-trivial one of $\pi_1L$ and $\pi_2L$ must be non-trivial, the argument above implies that
\begin{align*}
\mathcal{V}(G_1 \times G_2) \subseteq \left\{ v_1 + v_2 : v_i \in \mathcal{V}(G_i) \right\} \cup \left\{ v_1 + d_1(G_2) : v_1 \in \mathcal{V}(G_1) \right\} \\ \cup \left\{ d_1(G_1) + v_2 : v_2 \in \mathcal{V}(G_2) \right\}.
\end{align*}
For the other inclusion let 
\begin{align*}
 v \in \left\{ v_1 + v_2 : v_i \in \mathcal{V}(G_i) \right\} \cup \left\{ v_1 + d_1(G_2) : v_1 \in \mathcal{V}(G) \right\} \\ \cup \left\{ d_1(G_1) + v_2 : v_2 \in \mathcal{V}(G) \right\}.
\end{align*}
Thus there exist finite finite subgroups $L_1$ of $G_1$ and $L_2$ of $G_2$ such that $d_{L_1} = v_1$, $d_{L_2} = v_2$, and one of the $L_i$ is non-trivial.  Using the K\"unneth formula again, one calculates that $d_{L_1 \times L_2} = v$.
\end{proof}

\subsection{Finite-by-duality groups}

Throughout this section, $F$, $G$ and $Q$ will denote groups in a short exact sequence
\[1 \longrightarrow F \longrightarrow G \stackrel{\pi}{\longrightarrow} Q \longrightarrow 1,\]
where $F$ is finite.  This section builds up to the proof of Proposition \ref{prop:finite-by-uDn is uDn} that if $Q$ is Bredon duality of dimension $n$ over $R$, then $G$ is also.

\begin{Lemma}\label{lemma:Finite by Q implies HG is HQ}
$ H^i(G, R G) \cong H^i(Q, R Q ) $ for all $i$.
\end{Lemma}
\begin{proof}

 The Lyndon--Hochschild--Serre spectral sequence associated to the extension is \cite[VII\S 6]{Brown}
\[
H^p(Q, H^q(F, R G) ) \underset{p}{\Rightarrow} H^{p+q}(G, R G).
\]
$R G$ is projective as a $R F$-module so by \cite[Proposition 5.3, Lemma 5.7]{Bieri-HomDimOfDiscreteGroups}, 
\[H^q(F, R G) \cong H^q(F, R F) \otimes_{R F} R G = \left\{ \begin{array}{l l} R \otimes_{R F} R G = R Q & \text{ if } q = 0, \\ 0 & \text{ else.} \end{array} \right. \]
 The spectral sequence collapses to $H^i(G, R G) \cong H^i(Q, R Q )$.
\end{proof}

\begin{Lemma}\label{lemma:Finite-by-uFP0 is uFP0}
 If $Q$ is $\OFinFP_0$, then $G$ is $\OFinFP_0$.
\end{Lemma}
\begin{proof}
Let $B_i$ for $i = 0, \ldots, n$ be a collection of conjugacy class representatives of all finite subgroups in $Q$.  For each $i$, let $B_i^j$ be a collection of conjugacy class representatives of finite subgroups in $G$ which project onto $B_i$.  Since $F$ is finite $\pi^{-1}(B_i)$ is finite and there are only finitely many $j$ for each $i$, we claim that these $B_i^j$ are conjugacy class representatives for all finite subgroups in $G$.

Let $K$ be some finite subgroup of $G$, we need to check it is conjugate to some $B_i^j$.  $A= \pi(K)$ is conjugate to $B_i$, let $q \in Q$ be such that $q^{-1}Aq = B_i$ and let $g \in G$ be such that $\pi(g) = q$. 

$\pi(g^{-1} K g) = q^{-1} A q = B_i$ so $g^{-1}Kg$ is conjugate to some $B_i^j$ and hence $K$ is conjugate to some $B_i^j$.  Since we have already observed that for each $i$, the set $\{B_i^j\}_j$ is finite, $G$ has finitely many conjugacy classes of finite subgroups.
\end{proof}

\begin{Lemma}\label{lemma:pi G to Q with finite kernel then NGK is finite index in etc}
If $K$ is a finite subgroup of $G$ then $N_GK$ is finite index in $N_G (\pi^{-1}\circ\pi(K))$.
\end{Lemma}
\begin{proof}
 $N_GK$ is a subgroup of $N_G (\pi^{-1}\circ\pi(K))$ since if $g^{-1}Kg = K$ then 
\[
\left(\pi^{-1} \circ \pi(g)\right) \left(\pi^{-1} \circ \pi(K)\right) \left(\pi^{-1} \circ \pi(g)\right)^{-1} = \pi^{-1} \circ \pi(K),
\]
but $g \in \pi^{-1} \circ \pi(g)$ so $g\left( \pi^{-1} \circ \pi(K) \right) g^{-1} = \pi^{-1} \circ \pi(K)$.

Consider the action of $N_G(\pi^{-1} \circ \pi(K))$ on $\pi^{-1} \circ \pi (K)$ by conjugation, the setwise stabiliser of $K$ is exactly $N_GK$.  Since $\pi^{-1} \circ \pi(K)$ is finite, any stabiliser is finite index via the orbit-stabiliser theorem.  We conclude that $N_GK$ is finite index in $N_G(\pi^{-1} \circ \pi(K))$.
\end{proof}

\begin{Lemma}\label{lemma:pi G to Q and L subgroup of Q - structure of preimage of NQL}
 If $L$ is a subgroup of $Q$ then $N_G\pi^{-1}(L) = \pi^{-1} N_Q L$.
\end{Lemma}
\begin{proof}
If $g \in N_G\pi^{-1}(L)$ then $g^{-1}\pi^{-1}(L)g = \pi^{-1}(L)$ so applying $\pi$ gives that $\pi(g)^{-1}L\pi(g) = L$.  Thus $\pi(g) \in N_QL$, equivalently $g \in \pi^{-1}N_QL$.  

Conversely if $g \in \pi^{-1} (N_Q L)$ then $\pi(g)^{-1}L\pi(g) = L$ so 
\[
\left( \pi^{-1}\circ\pi(g) \right)^{-1} \pi^{-1}(L) \left( \pi^{-1}\circ\pi(g) \right) = \pi^{-1}(L).
\]
Since $g \in \pi^{-1}\circ\pi(g)$, we have that $g^{-1}\pi^{-1}(L)g = \pi^{-1}(L)$.
\end{proof}

\begin{Prop}\label{prop:finite-by-uDn is uDn}
$Q$ is Bredon duality of dimension $n$ over $R$ if and only if $G$ is Bredon duality of dimension $n$ over $R$.   Moreover, $\mathcal{V}(G) = \mathcal{V}(Q)$.
\end{Prop}

\begin{proof}
Assume that $Q$ is Bredon duality of dimension $n$ of $R$.  Let $K$ be a finite subgroup of $G$. We combine Lemma \ref{lemma:pi G to Q with finite kernel then NGK is finite index in etc} and Lemma \ref{lemma:pi G to Q and L 
subgroup of Q - structure of preimage of NQL} to see that $N_GK$ is finite index in $N_G (\pi^{-1} \circ \pi(K)) = \pi^{-1}\left( N_Q \pi(K) \right)$.  Hence 
\begin{align*}
H^i\left(W_GK , R[W_GK]\right) 
& \cong H^i\left(N_GK , R[N_GK]\right) \\
&\cong H^i\left(\pi^{-1} \left( N_Q \pi(K) \right), R\left[\pi^{-1}\left( N_Q\pi(K)\right) \right]\right) \\
&\cong H^i\left( N_Q \pi(K) , R\left[ N_Q\pi(K) \right]\right) \\
&\cong H^i\left( W_Q \pi(K) , R\left[ W_Q\pi(K) \right]\right)
\end{align*}
where the first isomorphism is from the short exact sequence
\[
 1 \longrightarrow K \longrightarrow N_GK \longrightarrow W_GK \longrightarrow 1
\]
and Lemma \ref{lemma:Finite by Q implies HG is HQ}, the fourth isomorphism is from the same lemma and a similar short exact sequence containing $N_QK$, and the third isomorphism follows from Lemma \ref{lemma:Finite by Q implies HG is HQ} and the short exact sequence
\[
1 \longrightarrow F \longrightarrow \pi^{-1}\left( N_Q\pi(K)\right) \longrightarrow N_Q\pi(K) \longrightarrow 1.
\]
Since $Q$ is Bredon duality of dimension $n$ this gives the condition on the cohomology of the Weyl groups.

$G$ is $\OFinFP_0$ by Lemma \ref{lemma:Finite-by-uFP0 is uFP0}, and $\OFincd G = \OFincd Q = n$ by \cite[Theorem 5.5]{Nucinkis-OnDimensionsInBredonCohomology}.  So by Corollary \ref{cor:uFPn equivalent conditions}, it remains to show that the Weyl groups of the finite subgroups are $\FP_\infty$.  For any finite subgroup $K$ of $G$, the short exact sequence above and \cite[Proposition 1.4]{Bieri-HomDimOfDiscreteGroups} gives that $\pi^{-1}\left( N_Q\pi(K)\right)$ is $\FP_\infty$.  But, as discussed at the beginning of the proof, $N_GK$ is finite index in $N_G (\pi^{-1} \circ \pi(K)) = \pi^{-1}\left( N_Q \pi(K) \right)$, so $N_G K$ is $\FP_\infty$ also. 

For the converse, assume that $G$ is Bredon duality of dimension $n$ over $R$.  Let $K$ be a finite subgroup of $Q$ then
\begin{align*}
 H^i(W_QK, R[W_QK]) &\cong H^i(N_QK, R[N_QK]) \\
 &\cong H^i( \pi^{-1}(N_QK), R[\pi^{-1}(N_QK)] ) \\
 &\cong H^i( N_G \pi^{-1} K, R[N_G\pi^{-1}K]) \\
 &\cong H^i( W_G \pi^{-1} K, R[W_G\pi^{-1}K]), 
\end{align*}
where the first isomorphism is from the short exact sequence
\[
 1 \longrightarrow K \longrightarrow N_QK \longrightarrow W_QK \longrightarrow 1
\]
and Lemma \ref{lemma:Finite by Q implies HG is HQ}, the fourth isomorphism is from the same lemma and a similar short exact sequence containing $N_G\pi^{-1}K$, the second isomorphism follows from Lemma \ref{lemma:Finite by Q implies HG is HQ} and the short exact sequence
\[
1 \longrightarrow F \longrightarrow \pi^{-1}\left( N_QK\right) \longrightarrow N_QK \longrightarrow 1,
\]
and the third isomorphism is from Lemma \ref{lemma:pi G to Q and L subgroup of Q - structure of preimage of NQL}.

Since $G$ is Bredon duality of dimension $n$ this gives the condition on the cohomology of the Weyl groups.  Finally, since $G$ is $\OFinFP_\infty$, thus also $Q$ is $\OFinFP_\infty$.
\end{proof}

\section{Graphs of groups}
An amalgamated free product of two duality groups of dimension $n$ over a duality group of dimension $n-1$ is duality of dimension $n$, similarly an HNN extension of a duality group of dimension $n$ relative to a duality group of dimension $n-1$ is duality of dimension $n$ \cite[Proposition 9.15]{Bieri-HomDimOfDiscreteGroups}.  
Unfortunately we know of no such result for Bredon--Poincar\'e duality groups:  the problem is how to obtain the correct condition on the cohomology of the Weyl groups of the finite subgroups.  However by putting some restrictions on the graph of groups, we can obtain some useful examples.  For instance using graphs of groups of Bredon duality groups we will be able to build Bredon duality groups $G$ with arbitrary $\mathcal{V}(G)$.  

Throughout this section, $G$ is the fundamental group of a finite graph of groups.  Let $T = (V, E)$ denote the associated Bass--Serre tree, we denote by $G_v$ the stabiliser of the vertex $v \in V$ and we denote by $G_e$ the stabiliser of the edge $e \in E$.  See \cite{Serre} for the necessary background on Bass--Serre trees and graphs of groups.

We need some preliminary results, showing that a graph of groups is $\OFinFP$ if all groups involved are $\OFinFP$.  See \cite{Serre} for the necessary background on Bass--Serre trees and graphs of groups.

\begin{Lemma}\cite[Lemma 3.2]{GandiniNucinkis-SomeH1FGroupsAndAConjectureOfKrophollerAndMislin}\label{lemma: duality mayervietoris of graph of groups}
There is an exact sequence, arising from the Bass--Serre tree.
 \[ \cdots \longrightarrow H^i_{\OFin} \left(G, -\right) \longrightarrow \bigoplus_{v \in V} H^i_{\OFin} \left( G_v, \Res^{G}_{G_v} - \right) \]
\[ \longrightarrow \bigoplus_{e \in E}H^i_{\OFin} \left(G_e, \Res^{G}_{G_e} - \right) \longrightarrow \cdots  \]
\end{Lemma}

\begin{Lemma}\label{lemma:uFPn for graph of groups}
 If all vertex groups $G_v$ are of type $\OFinFP_n$ and all edge groups $G_e$ are of type $\OFinFP_{n-1}$ over $R$ then $G$ is of type $\OFinFP_n$ over $R$.
\end{Lemma}

\begin{proof}
 Let $M_\lambda$, for $\lambda \in \Lambda$, be a directed system of $\OFin$-modules with colimit zero.  For any subgroup $H$ of $G$, the directed system $\Res^{G}_{H} M_\lambda$ also has colimit zero.  The long exact sequence of Lemma \ref{lemma: duality mayervietoris of graph of groups}, and the exactness of colimits gives that for all $i$, there is an exact sequence 
\[
\cdots  \longrightarrow \varinjlim_{\lambda \in \Lambda} H^{i-1}_{\OFin} \left(G, M_\lambda \right) \longrightarrow \bigoplus_{v \in V} \varinjlim_{\lambda \in \Lambda} H^i_{\OFin} \left( G_v, \Res^{G}_{G_v} M_\lambda \right) 
\]
\[
\longrightarrow  \bigoplus_{e \in E} \varinjlim_{\lambda \in \Lambda} H^i_{\OFin} \left(G_e, \Res^{G}_{G_e} M_\lambda \right) \longrightarrow \cdots.
\]
If $i \le n$ then by the Bieri--Eckmann criterion (Theorem \ref{theorem:C bieri-eckmann criterion}), the left and right hand terms vanish, thus the central term vanishes.  Another application of the Bieri--Eckmann criterion gives that $G$ is $\OFinFP_n$.
\end{proof}

\begin{Lemma}\label{lemma:ucd for graph of groups}
 If $\OFincd_R G_v \le n$ for all vertex groups $G_v$ and $ \OFincd_R G_e \le n-1$ for all edge groups $G_e$ then $\OFincd_R G \le n$.
\end{Lemma}
\begin{proof}
Use the long exact sequence of Lemma \ref{lemma: duality mayervietoris of graph of groups}.
\end{proof}

\begin{Lemma}\label{lemma:cohomology for graph of groups}
 If there exists a positive integer $n$ such that:
 \begin{enumerate}
  \item For every $v \in V$, $H^i(G_v, RG_v)$ is $R$-flat if $i = n$ and $0$ otherwise.
  \item For every $e \in E$, $H^i(G_e, RG_e)$ is $R$-flat if $i = n-1$ and $0$ otherwise.
 \end{enumerate}
 Then $H^i(G, RG)$ is $R$-flat if $i = n$ and $0$ else.
\end{Lemma}
\begin{proof}
The Mayer--Vietoris sequence associated to the graph of groups is
\[
\cdots \longrightarrow H^q(G, RG) \longrightarrow \bigoplus_{v \in V} H^q\left(G_v, RG\right) \longrightarrow \bigoplus_{e \in E} H^q\left(G_e, RG \right) \longrightarrow \cdots.
\]
$H^q\left(G_v, RG\right) = H^q(G_v, RG_v) \otimes_{RG_v} RG$ by \cite[Proposition 5.4]{Bieri-HomDimOfDiscreteGroups} so we have 
\[
H^q(G, RG) = 0 \text{ for } q \neq n,
\]
and a short exact sequence
\[
0 \longrightarrow \bigoplus_{e \in E} H^{n-1}(G_e, RG_e) \otimes_{RG_e} RG \longrightarrow H^n(G, RG) 
\]
\nopagebreak
\[
\longrightarrow \bigoplus_{v \in V} H^n(G_v, RG_v)\otimes_{RG_v} RG \longrightarrow 0.
\]
Finally, extensions of flat modules by flat modules are flat (use, for example, the long exact sequence associated to $\Tor^{RG}_*$).
\end{proof}

\begin{Remark}\label{remark:duality lemma:cohomology for graph of groups}
 In the above, if $H^n(G, RG_v) \cong R$ and $H^{n-1}(G_e, RG_e) \cong R$ for all vertex and edge groups then $H^n(G, RG) $ will not be isomorphic to $R$.
\end{Remark}

\begin{Lemma}\label{lemma:duality normaliser in tree of groups}
If $K$ is a subgroup of the vertex group $G_v$ and $K$ is not subconjugate to any edge group then $N_G K = N_{G_v} K$. 
\end{Lemma}
\begin{proof}
Let $T$ be the Bass--Serre tree, then the normaliser $N_G K$ fixes $T^K$ setwise, but $T^K$ is the single vertex $v$ (if $w \neq v$ was also fixed by $K$ then $K$ would fix all edges on the path from $v$ to $w$, but it is assumed that $K$ is not subconjugate to any edge stabiliser).  Thus, $N_G K \le G_v$.
\end{proof}

\begin{Example}
Let $S_n$ denote the star graph of $n+1$ vertices---a single central vertex $v_0$, and a single edge connecting every other vertex $v_i$ to the central vertex.  Let $G$ be the fundamental group of a graph of groups on $S_n$, where the central vertex group $G_0$ is torsion-free duality of dimension $n$, the edge groups are torsion-free duality of dimension $n-1$ and the remaining vertex groups $G_i$ are Bredon duality of dimension $n$ with $n_1 = n$.

By Lemmas \ref{lemma:uFPn for graph of groups} and \ref{lemma:ucd for graph of groups}, $G$ is $\OFinFP$ of dimension $n$, so to prove it is Bredon duality it suffices to check the cohomology of the Weyl groups of the finite subgroups.  Any non-trivial finite subgroup is subconjugate to a unique vertex group $G_i$, and cannot be subconjugate to an edge group since they are assumed torsion-free.  If $K$ is a subgroup of $G_i$ then by Lemma \ref{lemma:duality normaliser in tree of groups}, $H^i(N_G K, R[N_GK]) \cong H^i(N_{G_i} K, R[N_{G_i}K])$ and the condition follows as $G_i$ was assumed to be Bredon duality.  Finally, for the trivial subgroup we must calculate $H^i(G, RG)$, which is Lemma \ref{lemma:cohomology for graph of groups}.

$\mathcal{V}(G)$ is easily calculable too, 
\[
\mathcal{V}(G) = \{ v : v \in \mathcal{V}(G_i) \text{ for some }i \in \{1, \ldots, n\}  \}.
\]
\end{Example}

\begin{Example}[A Bredon duality group with prescribed $\mathcal{V}(G)$]\label{example:duality arbitrary V} We specialise the above example.
 Let $\mathcal{V} = \{v_1, \ldots, v_t\} \subset \{0, 1, \ldots, n\}$ be given.  Choosing $G_i = \ZZ^n \rtimes \ZZ_2$ as in Example \ref{example:duality Zn antipodal} so that $\mathcal{V}(G_i) = v_i$, let $G_0 = \ZZ^n$, let the edge groups be $\ZZ^{n-1}$, and choose injections $\ZZ^{n-1} \to \ZZ^n$ and $\ZZ^{n-1} \to \ZZ^n \rtimes \ZZ_2$ from the edge groups into the vertex groups.  Then form the graph of groups as in the previous example to get, for $G$ the fundamental group of the graph of groups,
\[
\mathcal{V}(G) = \{v_1, \ldots, v_t\}.
\]
\end{Example}

Because of Remark \ref{remark:duality lemma:cohomology for graph of groups} the groups constructed in the example above will not be Bredon--Poincar\'e duality groups. 

\section{The wrong notion of Bredon duality}\label{section:wrong notion of duality}

This section grew out of an investigation into which groups were $\OFinFP$ over some ring $R$ with
\[
H^i_{\OFin}(G, R[ -, ? ]_{\OFin}) \cong \left\{ \begin{array}{l l} \uR(?) & \text{if $i = n$,} \\ 0 & \text{else.} \end{array} \right. 
\]
One might hope that this na\"ive definition would give a duality similar to Poincar\'e duality, we show this is not the case.  Namely we prove in Theorem \ref{theorem:bredon wrong duality} that the only groups satisfying this property are torsion-free, and hence torsion-free Poincar\'e duality groups over $R$.  We need a couple of technical results before we can prove the theorem.

Recall from Section \ref{section:bredon fg projectives and duality} that for $M$ a contravariant $\OFin$-module we denote by $M^D$ the dual module 
\[
M^D = \Hom_{\OFin} \left( M(-), R[ -, ?]_{\OFin} \right).
\]
Note that $M^D$ is a covariant $\OFin$-module.  Similarly for $A$ a covariant $\OFin$-module,
\[
A^D = \Hom_{\OFin} \left( A(-), R[?, -]_{\OFin} \right). 
\]

\begin{Lemma}\label{lemma:length n cov res then cdR<=n}
 If there exists a length $n$ resolution of the constant covariant module $\uR$ by projective covariant $\OFin$-modules then $G$ is $R$-torsion free and $\cd_R G \le n$.
\end{Lemma}

\begin{proof}
Let $P_* \longtwoheadrightarrow \uR$ be a length $n$ projective covariant resolution of $\uR$, evaluating at $G/1$ gives a length $n$ resolution of $R$ by projective $RG$-modules (Propositions \ref{prop:C properties of res ind coind}(1) and \ref{prop:co and contra res properties})(1).  Thus $\cd_R G \le n$ and it follows that $G$ is $R$-torsion free.
\end{proof}

Let $M$ be an $RG$-module and recall from Section \ref{section:bredon res ind coind} that inducing $M$ to a covariant $\OFin$-module gives 
$\Ind_{RG}^{\OFin} M = M \otimes_{RG} R[G/1, -]_{\OFin}$.  The covariant induction functor maps projective modules to projective modules and satisfies the following adjoint isomorphism for any covariant $\OFin$-module $A$ (Propositions \ref{prop:C properties of res ind coind} and \ref{prop:ind coind and res adjoint isos}),
\[
\Hom_{\OFin}( \Ind_{RG}^{\OFin} M, A ) \cong \Hom_{RG}(M, A(G/1) ).
\]

\begin{Lemma}\label{lemma:cdR<=n then length n cov res}
 If $\cd_R G \le n$ then there exists a length $n$ projective covariant resolution of $\uR$.
\end{Lemma}
\begin{proof}
 Let $P_*$ be a length $n$ projective $RG$-module resolution of $R$, then we claim $\Ind_{RG}^{\OFin}P_*$ is a projective covariant resolution of $\uR$.  One can easily check that $\Ind_{RG}^{\OFin} R = \uR$ (Example \ref{example:cov E_1R is uR}) and since $G$ is necessarily $R$-torsion-free, $\Ind_{RG}^{\OFin}P_*$ is exact (Proposition \ref{prop:co and contra ind properties}).
\end{proof}

\begin{Theorem}\label{theorem:bredon wrong duality}
If $G$ is $\OFinFP$ with $\OFincd_R G = n$ and 
\[H^i_{\OFin}(G, R[ -, ? ]_{\OFin}) \cong \left\{ \begin{array}{l l} \uR(?) & \text{if $i = n$,} \\ 0 & \text{else,} \end{array} \right. \]
then $G$ is torsion-free.  Note that in the above, $\uR$ denotes the constant covariant $\OFin$-module.
\end{Theorem}
\begin{proof}
Choose a length $n$ finite type projective $\OFin$-module resolution $ P_*$ of $\uR$ then by the assumption on $H^n_{\OFin}(G, R[ -, ? ]_{\OFin})$, we know that $P_*^D$ is a covariant resolution by finitely generated projectives of $\uR$:
\[
0 \rightarrow  P_0^D(-) \stackrel{\partial_1^D}{\rightarrow} P_1^D(-) \stackrel{\partial_2^D}{\rightarrow} \cdots \stackrel{\partial_n^D}{\rightarrow} P_n^D(-) \rightarrow H^n_{\OFin}(G, R[  ?, - ]_{\OFin}) \cong \uR(-)  \rightarrow 0.
\]
By Lemma \ref{lemma:length n cov res then cdR<=n} $G$ is $R$-torsion-free and $\cd_R G \le n$.  Since $G$ is $\OFinFP_\infty$, $G$ is $\FP_\infty$ (Corollary \ref{cor:uFPn equivalent conditions}) and we may choose a length $n$ finite type projective $RG$-module resolution $Q_*$ of $R$.  Lemma \ref{lemma:cdR<=n then length n cov res} gives that $\Ind_{RG}^{\OFin}Q_* \longtwoheadrightarrow \uR$ is a projective covariant resolution.  

By the $\OFin$-module analogue of the comparison theorem \cite[2.2.6]{Weibel}, the two projective covariant resolutions of $\uR$ are chain homotopy equivalent.  Any additive functor preserves chain homotopy equivalences, so applying the dual functor to both complexes gives a chain homotopy equivalence between
\[
0 \longrightarrow \uR^D \cong 0 \longrightarrow (\Ind_{RG}^{\OFin} Q_0)^D \longrightarrow  \cdots \longrightarrow (\Ind_{RG}^{\OFin} Q_n)^D 
\]
and
\[
0 \longrightarrow \uR^D \cong 0 \longrightarrow P_n^{DD} \longrightarrow  P_{n-1}^{DD} \longrightarrow \cdots \longrightarrow P_0^{DD},
\]
(that $\uR^D \cong 0 $ is Example \ref{example:dual of constant functor}).  Since $\Hom_{\OFin}$ is left exact we know both complexes above are left exact.  Lemma \ref{lemma:double dual is nat iso} gives the commutative diagram below.
\[
\xymatrix{
 0 \ar[r] & P_n^{DD} \ar[r] \ar^\cong[d] &  \cdots \ar[r] & P_1^{DD} \ar[r] \ar^\cong[d] &  P_0^{DD}  \ar^\cong[d]\\
 0 \ar[r] & P_n \ar[r] &  \cdots \ar[r]  & P_1 \ar[r] & P_0
}\]

The lower complex, $P_*$, satisfies $H_0P_* \cong \uR$ and $H_iP_* = 0$ for all $i \neq 0$.  Thus the same is true for the top complex, and also the complex $\Ind_{RG}^{\OFin} Q_*^D$, since this is homotopy equivalent to it.  In particular, there is an epimorphism of $\OFin$-modules, 
\[
\Ind_{RG}^{\OFin} Q_n^D \longtwoheadrightarrow \uR.
\]
The left hand side simplifies, using the adjoint isomorphism
\[
\Ind_{RG}^{\OFin} Q_n^D = \Hom_{\OFin} \left( \Ind_{RG}^{\OFin} Q_n , R[?, -]_{\OFin} \right) \cong \Hom_{RG}(Q_n, R[?, G/1]_{\OFin}) .
\]
Since $\Hom_{RG}(Q_n, R[?, G/1]_{\OFin})(G/H) = 0$ if $H \neq 1$, this module cannot surject onto $\uR$ unless $G$ is torsion-free. 
\end{proof} 

\chapter{Houghton's groups}\label{chapter:H}

This chapter, with the exception of Sections \ref{subsection:houghton cdVCyc LW} and \ref{subsection:houghton cdVCyc calc}, contains material that has appeared in:
\begin{itemize}
  \item Centralisers in Houghton's Groups (2012, to appear Proc. Edinburgh Math. Soc.) \cite{Me-HoughtonsGroups}.
\end{itemize}
Sections \ref{subsection:houghton cdVCyc LW} and \ref{subsection:houghton cdVCyc calc} contain unpublished joint work with Nansen Petrosyan.

Section \ref{section:centralisersOfFiniteSubgroups} contains an analysis of the centralisers of finite subgroups in Houghton's group.  As Corollary \ref{cor:stabOfH is FPn-1 not FPn} we obtain that centralisers of finite subgroups are $\FP_{n-1}$ but not $\FP_n$.  In Section \ref{section:centralisersOfElements} our analysis is extended to arbitrary elements and virtually cyclic subgroups.  Using this information elements in $H_n$ are constructed whose centralisers are $\FP_i$  for any $0 \le i \le n-3$.  

In Section \ref{section:Browns model} the space that Brown constructed in \cite{Brown-FinitenessPropertiesOfGroups}, in order to prove that $H_n$ is $\FP_{n-1}$ but not $\FP_n$, is shown to be a model for $\EFin H_n$, the classifying space for proper actions of $H_n$.  Finally Section \ref{section:finitenessConditionsSatisfied} contains a discussion of Bredon (co)homological finiteness conditions satisfied by Houghton's group, namely we show in Proposition \ref{prop:HisnotQuasiFP0} that $H_n$ is not quasi-$\OFinFP_0$ and in Proposition \ref{prop:cd=gd=n} that the Bredon cohomological dimension with respect to the family of finite subgroups and virtually cyclic subgroups are both equal to $n$.  See Section \ref{subsection:bredon quasi-FPn} for the definition of quasi-$\OFinFP_n$.

Fixing a natural number $n >1$, define \emph{Houghton's group} $H_n$ to be the group of permutations of $S = \mathbb{N} \times \{ 1, \ldots, n\}$ which are ``eventually translations'', ie. for any given permutation $h \in H_n$ there are collections $\{z_1, \ldots, z_n\} \in \NN^n$ and $\{m_1, \ldots, m_n\} \in \ZZ^n$ with 
\begin{equation}\label{eventuallyATranslation}
h(i,x) = (i+m_x, x) \text{ for all }x \in \{1, \ldots, n \} \text{ and all }i \ge z_x.
\end{equation}
\index{Hn@$H_n$, Houghton's group}

Define a map $\phi$ as follows: \index{$\phi$, group homomorphism $H_n \longrightarrow \ZZ^{n-1}$}
\begin{align}\label{def:phi}
 \phi : H_n &\to \{ (m_1, \ldots, m_n) \in \mathbb{Z}^n \: : \: \sum m_i = 0\} \cong \ZZ^{n-1} \\
 \phi : h &\mapsto (m_1, \ldots, m_n).
\end{align}
Its kernel is exactly the permutations which are ``eventually zero'' on $S$, ie. the infinite symmetric group $\Sym_\infty$ (the finite support permutations of the countable set $S$).  

\section{Centralisers of finite subgroups in \texorpdfstring{$\Hn$}{Hn}}\label{section:centralisersOfFiniteSubgroups}
First we recall some properties of group actions on sets, before specialising to Houghton's group.
\begin{Prop}\label{prop:groupActionsOnSets}
 If $G$ is a group acting on a countable set $X$ and $H$ is any subgroup of $G$ then
\begin{enumerate}
 \item If $x$ and $y$ are in the same $G$-orbit then their isotropy subgroups $G_x$ and $G_y$ are $G$-conjugate.
 \item If $g \in C_G(H)$ then $H_{gx} = H_x$ for all $x \in X$.
 \item Partition $X$ into $\{X_a\}_{a=1}^t$, where $t \in \NN \cup \{\infty\}$, via the equivalence relation $x \sim y$ if and only if $H_x$ is $H$-conjugate to $H_y$.  Any two points in the same $H$-orbit will lie in the same partition and any $c \in C_G(H)$ maps $X_a$ onto $X_a$ for all $a$.
 \item Let $G$ act faithfully on $X$, with the property that for all $g \in G$ and $X_a \subseteq X$ as in the previous section, there exists a group element $g_a \in G$ which fixes $X \setminus X_a$ and acts as $g$ does on $X_a$.  Then $C_G(H) = C_1 \times \cdots \times C_t$ where $C_a$ is the subgroup of $C_G(H)$ acting trivially on $X \setminus X_a$.
\end{enumerate}
\end{Prop}
\begin{proof}(1) and (2) are standard results.
\begin{enumerate}
 \setcounter{enumi}{2} 
\item This follows immediately from (1) and (2).
 \item This follows from (3) and our new assumption on $G$:  Let $c \in C_G(H)$ and $c_a$ be the element given by the assumption.  Since the action of $G$ on $X$ is faithful, $c_a$ is necessarily unique.  That the action is faithful also implies $c = c_1 \cdots c_t$ and that any two $c_a$ and $c_b$ commute in $G$ because they act non-trivially only on distinct $X_a$.  Thus we have the necessary isomorphism $C_G(H) \longrightarrow C_1 \times \cdots \times C_t$.
\end{enumerate}
\end{proof}
Let $Q \le H_n$ be a finite subgroup of Houghton's group $H_n$ and $S_Q  = S \setminus S^Q$ the set of points of $S$ which are \emph{not fixed} by $Q$.  $Q$ being finite implies $\phi(Q) = 0$ as any element $q$ with $\phi(q) \neq 0$ necessarily has infinite order.  For every $q \in Q$ there exists $\{z_1, \ldots, z_n \} \in \NN^n$ such that
\[
q(i,x) = (i, x) \text{ if } i \ge z_x.
\]
Taking $z_x^\prime$ to be the maximum of these $z_x$ over all elements in $Q$, then $Q$ must fix the set $\{(i,x)\: : \: i \ge z^\prime_x\}$ and in particular $S_Q \subseteq \{(i,x)\: : \: i < z^\prime_x\}$ is finite.

We need to see that the subgroup $Q \le H_n$ acting on the set $S$ satisfies the conditions of Proposition \ref{prop:groupActionsOnSets}(4).  We give the following lemma in more generality than is needed here, as it will come in useful later on.  That the action is faithful is automatic as an element $h \in H_n$ is uniquely determined by its action on the set $S$.
\begin{Lemma}\label{lemma:prop conditions are satisfied}
 Let $Q \le H_n$ be a subgroup, which is either finite or of the form $F \rtimes \ZZ$ for $F$ a finite subgroup of $H_n$.  Partition $S$ with respect to $Q$  into sets $\{S_a\}_{a = 1}^t$ as in Proposition \ref{prop:groupActionsOnSets}(3) applied to the action of $H_n$ on $S$ and the subgroup $Q$ of $H_n$.  Then the conditions of Proposition \ref{prop:groupActionsOnSets}(4) are satisfied.
\end{Lemma}
\begin{proof}
Fix $a \in \{1, \ldots, t\}$ and let $h_a$ denote the permutation of $S$ which fixes $S \setminus S_a$ and acts as $h$ does on $S_a$.  We wish to show that $h_a$ is an element of $H_n$.

There are only finitely many elements in $Q$ with finite order so as in the argument just before this lemma we may choose integers $z_x$ for $x \in \{1, \ldots, n\}$ such that if $q$ is a finite order element of $Q$ then $q(i,x) = (i,x)$ whenever $i \ge z_x$.  If $Q$ is a finite group then either:
\begin{itemize}
 \item $S_a$ is fixed by $Q$, in which case 
\[
\{(i,x) \: : \: i \ge z_x \: , \: x \in \{1, \ldots, n \} \} \subseteq S_a
\]
 so $h_a(i,x) = h(i,x)$ for all $i \ge z_x$.  In particular for large enough $i$, $h_a$ acts as a translation on $(i,x)$ and is hence an element of $H_n$.
\end{itemize}
Or
\begin{itemize}
 \item $S_a$ is not fixed by $Q$, in which case 
 \[
 S_a \subseteq \{(i,x) \: : \: i < z_x \: , \: x \in \{1, \ldots, n \} \}.
 \]
In particular $S_a$ is finite and $h_a(i,x) = (i,x)$ for all $i \ge z_x$.  Hence $h_a$ is an element of $H_n$.
\end{itemize}

It remains to treat the case where $Q= F \rtimes \ZZ$.  Write $w$ for a generator of $\ZZ$ in $F \rtimes \ZZ$.  By choosing a larger $z_x$ if needed we may assume $w$ acts either trivially or as a translation on $(i,x)$ whenever $i \ge z_x$.  Hence for any $x \in \{1, \ldots, n\}$, the isotropy group in $Q$ of $\{(i,x) \: :\: i \ge z_x \}$ is either $F$ or $Q$.  

If $S_a$ has isotropy group $Q$ or $F$ then for some $x \in \{1, \ldots, n\}$, either
\begin{itemize}
 \item \[
        S_a \cap \{(i,x) : i \ge z_x\} = \{(i,x) : i \ge z_x\}
       \]
 in which case $h_a(i,x) = h(i,x)$ for $i \ge z_x$.  In particular for large enough $i$, $h_a$ acts as a translation on $(i,x)$ and hence is an element of $H_n$.
\end{itemize}
Or
\begin{itemize} 
\item \[
S_a \cap \{(i,x) : i \ge z_x\} = \emptyset
\]
 in which case $h_a(i,x) = (i,x)$ for $i \ge z_x$.  In particular for large enough $i$, $h_a$ fixes $(i,x)$ and hence is an element of $H_n$.
\end{itemize}
If $S_a$ is the set corresponding to an isotropy group not equal to $F$ or $Q$ then 
\[
S_a \subseteq \{ (i,x) \: : \: i \ge z_x \: , \: x \in \{1, \ldots, n\} \}.
\]
So $h_a$ fixes $(i,x)$ for $i \ge z_x$ and hence $h_a$ is an element of $H_n$.
\end{proof}

Partition $S$ into disjoint sets according to the $Q$-conjugacy classes of the stabilisers, as in Proposition \ref{prop:groupActionsOnSets}(3).  The set with isotropy in $Q$ equal to $Q$ is $S^Q$ and since $S_Q$ is finite the partition is finite, thus
\[ 
S = S^Q \cup S_1 \cup \cdots \cup S_t.
\]
Proposition \ref{prop:groupActionsOnSets}(4) gives that  
\[C_{H_n}(Q) = H_n|_{S^Q} \times C_1 \times \ldots \times C_t \]
where each $C_a$ acts only on $S_a$ and leaves $S^Q$ and $S_b$ fixed for $a \neq b$ (where $a, b \in \{1, \ldots, t\}$).  The first element of the direct product decomposition is the subgroup of $C_{H_n}(Q)$ acting only on $S^Q$ and leaving $S \setminus S^Q$ fixed.  This is $H_n|_{S^Q}$ ($H_n$ restricted to $S^Q$) because, as the action of $Q$ on $S^Q$ is trivial, any permutation of $S^Q$ will centralise $Q$.  Choose a bijection $S^Q \to S$ such that for all $x$, $(i,x) \mapsto (i + m_x, x)$ for large enough $i$ and some $m_x \in \ZZ$, this induces an isomorphism between $H_n|_{S^Q}$ and $H_n$.

To give an explicit definition of the group $C_a$ we need three lemmas.

\begin{Lemma}\label{lemma:Ca iso T}
 $C_a$ is isomorphic to the group $T$ of $Q$-set automorphisms of $S_a$.
\end{Lemma}
\begin{proof}
 An element $c \in C_a$ determines a $Q$-set automorphism of $S_a$, giving a map $C_a \to T$.  Since the action of $C_a$ on $S_a$ is faithful this map is injective.  Any $Q$-set automorphism $\alpha$ of $S_a$ may be extended to a $Q$-set automorphism of $S$, where $\alpha$ acts trivially on $S \setminus S_a$.  Since $S_a$ is a finite set, $\alpha$ acts trivially on $(i,x)$ for large enough $i$ and any $x \in \{1, \ldots, n\}$, and hence $\alpha$ is an element of $H_n$.  Finally, since $\alpha$ is a $Q$-set automorphism $q \alpha s = \alpha q s$, equivalently $\alpha^{-1} q \alpha s = s$, for all $s \in S$ and $q \in Q$, showing that $\alpha \in C_a$ and so the map $C_a \to T$ is surjective.
\end{proof}

\begin{Lemma}\label{lemma:S_a is disjoint union of Q/Qa}
 $S_a$ is $Q$-set isomorphic to the disjoint union of $r$ copies of $Q / Q_a$, where $Q_a$ is an isotropy group of $S_a$ and $r = \vert S_a \vert / \vert Q : Q_a \vert$.
\end{Lemma}
\begin{proof}
 $S_a$ is finite and so splits as a disjoint union of finitely many $Q$-orbits.  Choose orbit representatives $\{s_1, \ldots, s_r \} \subset S_a$ for these orbits, these $s_k$ may be chosen to have the same $Q$-stabilisers: If $Q_{s_1} \neq Q_{s_2}$ then there is some $q \in Q$ such that $Q_{qs_2} = qQ_{s_2} q^{-1} = Q_{s_1}$ (the partitions $S_a$ were chosen to have this property by Proposition \ref{prop:groupActionsOnSets}), iterating this procedure we get a set of representatives who all have isotropy group $Q_{s_1}$.  Now set $Q_a = Q_{s_1}$ and note that there are $\vert Q:Q_a \vert $ elements in each of the $Q$-orbits so $ r \vert Q : Q_a \vert = \vert S_a \vert$.
\end{proof}

Recall that if $G$ is any group and $r \ge 1$ is some natural number then the wreath product $ G \wr \Sym_r $ is the semi-direct product
\[ G \wr \Sym_r = \prod_{k = 1}^r G \rtimes \Sym_r \]
where the symmetric group $\Sym_r$ acts by permuting the factors in the direct product.  

Recall also that for any subgroup $H$ of a group $G$, the Weyl group $W_G H$ is defined to be $W_G H = N_G H / H$.
\begin{Lemma}\label{lemma:description of Ca}
The group $C_a$ is isomorphic to the wreath product $ W_QQ_a \wr \Sym_r $, where $Q_a$ is some isotropy group of $S_a$ and $r = \vert S_a \vert / \vert Q : Q_a \vert$.
\end{Lemma}
\begin{proof}
 Using Lemmas \ref{lemma:Ca iso T} and \ref{lemma:S_a is disjoint union of Q/Qa}, $C_a$ is isomorphic to the group of $Q$-set automorphisms of the disjoint union of $r$ copies of $Q/Q_a$.

 To begin, we show the group of automorphisms of the $Q$-set $Q / Q_a$ is isomorphic to $W_QQ_a$.  An automorphism $\alpha : Q/Q_a \to Q/Q_a$ is determined by the image $\alpha(Q_a) = qQ_a$ of the identity coset and such an element determines an automorphism if and only if $q^{-1} Q_a q \le Q_a $, equivalently $q \in N_Q Q_a$.  Since two elements $q_1,q_2 \in Q$ will determine the same automorphism if and only if $q_1 Q_a = q_2Q_a$, the group of $Q$-set automorphisms of $Q/Q_a$ is the Weyl group $W_QQ_a$.

 For the general case, note that if $c \in C_a$ then $c$ permutes the $Q$-orbits $\{Qs_1, \ldots, Qs_r\}$, so there is a map $\pi : C_a \to \Sym_r$.  Assume that the representatives $\{s_1, \ldots, s_r\}$ have been chosen, as in the proof of Lemma \ref{lemma:S_a is disjoint union of Q/Qa}, to have the same $Q$-stabilisers.  The map $\pi$ is split by the map 
\begin{align*}
 \iota:\Sym_r &\to C_a \\
 \sigma &\mapsto \left( \iota(\sigma):  qs_k \mapsto qs_{\sigma(k)} \text{ for all }q \in Q \right).
\end{align*}
Each $\iota(\sigma)$ is a well defined element of $H_n$ since 
\[
qs_k = \tilde{q} s_k \Leftrightarrow \tilde{q}^{-1}q \in Q_{s_k} = Q_{s_{\sigma(k)}} \Leftrightarrow qs_{\sigma(k)} = \tilde{q} s_{\sigma(k)}.
\]

The kernel of the map $\pi$ is exactly the elements of $C_a$ which fix each $Q$-orbit but may permute the elements inside the $Q$-orbits, by the previous part this is exactly $\prod_{k = 1}^r W_QQ_a$.  For any $\sigma \in \Sym_r$, the element $\iota(\sigma)$ acts on $\prod_{k = 1}^r W_QQ_a$ by permuting the factors, so the group $C_a$ is indeed isomorphic to the wreath product.
\end{proof}

\begin{figure}[ht]
 \centering
 \includegraphics[width = 0.4\textwidth]{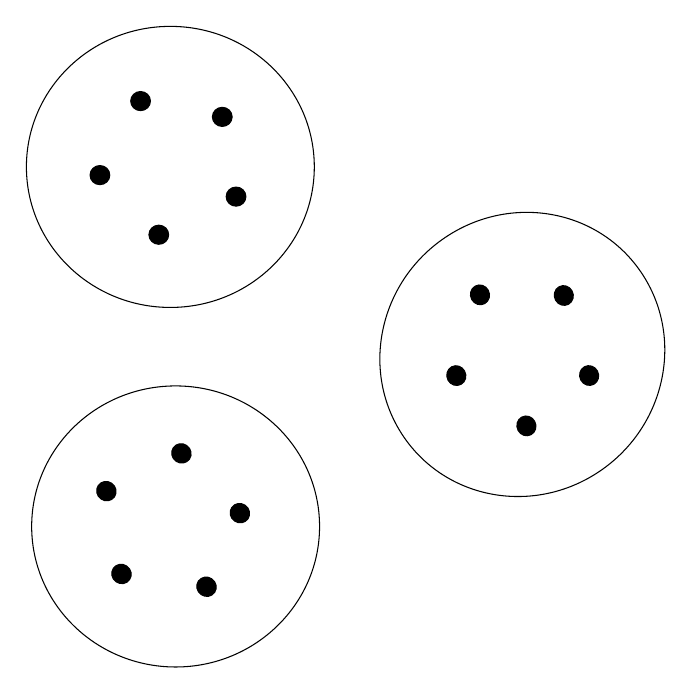}
 \caption[A graphical representation of the set $S_a$ described in Section \ref{section:centralisersOfFiniteSubgroups}.]{A representation of $S_a$.  The large circles are the sets $\{Qs_1 \ldots, Qs_r\}$ (in this figure $r = 3$).  Elements of $\Sym_r$ permute only the large circles, while elements of $\prod_{k = 1}^r W_QQ_a$ leave the large circles fixed and permute only elements inside them.}
\end{figure}  

The centraliser $C_{H_n}Q$ can now be completely described.

\begin{Prop}\label{prop:stabOfH}
 The centraliser $C_{H_n}(Q)$ of any finite subgroup $Q \le H_n$ splits as a direct product 
\[ 
C_{H_n}(Q) \cong  H_n|_{S^Q} \times C_1 \times \cdots \times C_t,
\]
where $H_n|_{S^Q} \cong H_n$ is Houghton's group restricted to $S^Q$ and for all $a \in \{1, \ldots, t\}$,
\[ 
C_a \cong W_Q Q_a \wr \Sym_r  
\]
for $Q_a$ is an isotropy group of $S_a$ and r = $\vert S_a \vert / \vert Q : Q_a \vert$.  In particular $H_n$ is finite index in $C_{H_n}(Q)$.
\end{Prop}

\begin{proof}
We have already proven that 
\[
C_{H_n}(Q) \cong  H_n|_{S^Q} \times C_1 \times \cdots \times C_t 
\]
and Lemma \ref{lemma:description of Ca} gives the required description of $C_a$.
\end{proof}

\begin{Cor}\label{cor:stabOfH is FPn-1 not FPn}
 If $Q$ is a finite subgroup of $H_n$ then the centraliser $C_{H_n}(Q)$ is $\FP_{n-1}$ but not $\FP_n$.
\end{Cor}
\begin{proof}
 $H_n$ is finite index in the centraliser $C_{H_n}(Q)$ by Proposition \ref{prop:stabOfH}.  Appealing to Brown's result \cite[5.1]{Brown-FinitenessPropertiesOfGroups} that $H_n$ is $\FP_{n-1}$ but not $\FP_n$, and that a group is $\FP_n$ if and only if a finite index subgroup is $\FP_n$ \cite[VIII.5.5.1]{Brown} we can deduce $C_{H_n}(Q)$ is $\FP_{n-1}$ but not $\FP_n$. 
\end{proof}

\section{Centralisers of elements in \texorpdfstring{$\Hn$}{Hn}}\label{section:centralisersOfElements}

If $q \in H_n$ is an element of finite order then the subgroup $Q = \langle q \rangle$ is a finite subgroup and the previous section may be used to describe the centraliser $C_{H_n}(q) = C_{H_n}(Q)$.  Thus for an element $q$ of finite order $C_{H_n}(q) \cong C \times H_n$ for some finite group $C$.

If $q \in H_n$ is an element of infinite order and $Q = \langle q \rangle$ then we may apply Proposition \ref{prop:groupActionsOnSets}(3) to split up $S$ into a disjoint collection $\{ S_a \: : \: a \in A \subseteq \NN\} \cup S^Q$ ($S^Q$ is the element of the collection associated to the isotropy group $Q$).  Assume that $S_0$ is the set associated to the trivial isotropy group.  Since $q$ is a translation on $(i,x) \in S = \NN \times \{1, \ldots, n\}$ for large enough $i$ and points acted on by such a translation have trivial isotropy, there are only finitely many elements of $S$ whose isotropy group is neither the trivial group nor $Q$.  Hence $S_a$ is finite for $a \neq 0$ and the set $A$ is finite.  From now on let $A = \{0, \ldots, t\}$.  We now use Lemma \ref{lemma:prop conditions are satisfied} and Proposition \ref{prop:groupActionsOnSets}(4) as in the previous section: $C_{H_n}(Q)$ splits as 
\[
C_{H_n}(Q) \cong C_0 \times C_1 \times \cdots \times C_t \times H_n|_{S^Q}
\]
where $C_a$ acts only on $S_a$ and $H_n|_{S^Q}$ is Houghton's group restricted to $S^Q$.  Unlike in the last section, $H_n|_{S^Q}$ may not be isomorphic to $H_n$.  Let $J \subseteq \{1, \ldots, n\}$ satisfy
\[
 x \in J \text{ if and only if } (i,x) \in S^Q \text{ for all }i \ge z_x \text{, some }z_x \in \NN.
\]
If $x \notin J$ then for large enough $i$, $q$ must act as a non-trivial translation on $(i,x)$, and the set $\left(\NN \times \{x\}  \right) \cap S^Q $ is finite.  Clearly $\vert J \vert \le n-2$, but different elements $q$ may give values $0 \le \vert J \vert \le n-2$.  In the case $\vert J \vert = 0$, $S^Q$ is necessarily finite and so $H_n|_{S^Q}$ is isomorphic to a finite symmetric group on $S^Q$.  It is also possible that $S^Q  = \emptyset$, in which case $H_n|_{S^Q}$ is just the trivial group.  If $\vert J \vert \neq 0$ then the argument proceeds by choosing a bijection
\[S^Q \to \NN \times J\]
such that $(i,x) \mapsto (i + m_x, x)$ for some $m_x \in \ZZ$ whenever $i$ is large enough and $x \in J$.  This set map induces a group isomorphism between $H_n|_{S^Q}$ and $H_{\vert J \vert}$ (Houghton's group on the set $J \times \NN$).   

Lemma \ref{lemma:description of Ca} describes the groups $C_a$ for $a \neq 0$, so it remains only to treat the case $a=0$.  We cannot use the arguments used for $a \neq 0$ here as the set $S_0$ is not finite, in particular Lemma \ref{lemma:Ca iso T} doesn't apply: Every $Q$-set isomorphism of $S_0$ is realised by an element of the infinite support permutation group on $S_0$, but there are $Q$-set isomorphisms of $S_0$ which are not realised by an element of $H_n$.

The next three lemmas are needed to describe $C_0$, this description will use the graph $\Gamma$ which we now describe.  The vertices of $\Gamma$ are those $x \in \{1, \ldots, n\}$ for which $q$ acts non-trivially on infinitely many elements of $ \NN \times \{ x\} $.  Equivalently, the vertices are the elements of $\{1, \ldots, n\} \setminus J$.  There is an edge from $x$ to $y$ in $\Gamma$ if there exists $s \in S_0$ and $N \in \NN$ such that for all $m \ge N$ we have $q^{-m}s \in \NN \times \{ x \} $ and $q^m s \in \NN \times \{ y\}$.  Let $\pi_0 \Gamma$ denote the path components of $\Gamma$, and for any vertex $x$ of $\Gamma$ denote by $[x]$ the element of $\pi_0 \Gamma$ corresponding to that vertex.

Let $z \in \NN$ be some integer such that for all $i \ge z$, $q$ acts trivially or as a translation on $(i, x)$ for all $x \in \{1, \ldots, n\}$.  Fix $z$ for the remainder of this section.

For each path component $[x]$ in $\pi_0 \Gamma$, let $S_0^{[x]}$ denote the smallest $Q$-subset of $S_0$ containing the set $\{(i,y) \: : \: i \ge z , \: y \in [x] \}$.  Note that $(i,y) \notin S_0^{[x]}$ for any $y \notin [x]$ and $i \ge z$, since if $(i,x)$ and $(j,y)$ are two elements of $S_0$ in the same $Q$-orbit with $i \ge z$ and $j \ge z$ then there is an edge between $x$ and $y$ in $\Gamma$:  If $(i,x) = q^k(j,y)$ and $q$ acts as a positive translation on the element $(i,x)$ then let $N=k$ and $s = (i,x)$, similarly for when $q$ acts as a negative translation.  This gives a $Q$-set decomposition of $S_0$ as
\[S_0 = \coprod_{[x] \in \pi_0 \Gamma} S_0^{[x]}, \]
where $\coprod$ denotes disjoint union.

\begin{Lemma}\label{lemma:decomposition of elements in S0x}
Let $[x] \in \pi_0 \Gamma$, if $C_0^{[x]}$ denotes the subgroup of $C_0$ which acts non-trivially only on $S_0^{[x]}$ then there is an isomorphism
\[ C_0 \cong C_0^{[x_1]} \times \cdots \times C_0^{[x_r]}, \]
where $[x_1], [x_2], \ldots, [x_r]$ are all elements of $\pi_0 \Gamma$.
\end{Lemma}
\begin{proof}
If $c \in C_0$ and $[x] \in \pi_0 \Gamma$ then let $c_{[x]}$ denote the permutation of $S$ such that $c_{[x]}$ acts as $c$ does on $S_0^{[x]}$, and acts trivially on $S \setminus S_0^{[x]}$.  We will show that $c_{[x]}$ is an element of $C_0$.  Since the action of $C_0$ on $S_0$ is faithful it follows that the elements $c_{[x]}$ and $c_{[y]}$ commute and 
\[ c = c_{[x_1]} c_{[x_2]} \cdots c_{[x_r]}, \]
which suffices to prove the lemma.  

Let $y \in \{1,\ldots,n\}$.  The element $c_{[x]}$ acts trivially on $(i,y)$ for $i \ge z$ if $y \notin [x]$ and acts as as $c$ does on $(i,y)$ for $i \ge z$ if $y \in [x]$, thus $c_{[x]}$ is an element of $H_n$.  Since $c_{[x]}$ is also a $Q$-set automorphism of $S$, $c_{[x]}$ is a member of $C_0$.
\end{proof}

\begin{Lemma}\label{lemma:uniquness of action on S0x}
 Let $[x] \in \pi_0\Gamma$, let $c \in C_0$, and let $z^\prime \in \NN$ be such that $c$ acts either trivially or as a translation on $(i,x)$ for all $x \in \{1, \ldots, n\}$ and $i \ge z^\prime$.  Then the action of $c$ on some element $(i, x) \in S$ for $i \ge z^\prime$ completely determines the action of $c$ on $S_0^{[x]}$.
\end{Lemma}
\begin{proof}
Firstly, note that knowing the action of $c$ on some element $(i,x)$ for $i \ge z^\prime$ determines the action of $c$ on the set $\{(i,x) \: : \: i \ge z^\prime \}$, since we chose $z^\prime$ in order to have this property.

Let $y \in [x]$ such that there is an edge from $x$ to $y$, so there is a natural number $N$ and element $s \in S_0^{[x]}$ such that $ q^N s = (i,x)$ and $q^{-N} s = (j,y)$ for some natural numbers $i$ and $j$.  By choosing $N$ larger if necessary we can take $i, j \ge z^\prime$.  The action of $c$ on $(j,y)$ is now completely determined by the action on $(i,x)$, since
\[
c(j,y) = c q^{-2N}(i,x) = q^{-2N} c(i,x).
\]
For any $y \in [x]$ there is a path from $x$ to $y$ in $\Gamma$, so we've determined the action of $c$ on the set $X = \{(j,y) \: : \: j \ge z^\prime \: , \: y \in [x]\}$.  If $s \in S_0^{[x]} \setminus X$ then, since $S_0^{[x]} \setminus X$ is finite, there is some integer $m$ with $q^m s = x \in X$.  So $c s = c q^{-m} x = q^{-m} c x $, which completely determines the action of $c$ on $s$.
\end{proof}

\begin{Lemma}\label{lemma:C0x iso Z}
 For any $[x] \in \pi_0\Gamma$, there is an isomorphism
\[
C_0^{[x]} \cong \ZZ.
\]
\end{Lemma}
\begin{proof}
 By Lemma \ref{lemma:uniquness of action on S0x} the action is completely determined by the action on some element $(i,x)$ for large enough $i$, and the action on this element is necessarily by translation by some element $m_x(c)$.  This defines an injective homomorphism $C_0^{[x]} \to \ZZ$, sending $c \mapsto m_x(c)$.  Let $q_{[x]}$ be the element of $C_0^{[x]}$ described in the proof of Lemma \ref{lemma:decomposition of elements in S0x}, $q_{[x]}$ is a non-trivial element of $C_0^{[x]}$ so $C_0^{[x]}$ is mapped isomorphically onto a non-trivial subgroup of $\ZZ$.
\end{proof}

Combining Lemmas \ref{lemma:decomposition of elements in S0x} and \ref{lemma:C0x iso Z} shows $C_0 \cong \ZZ^r$ where $r = \vert \pi_0 \Gamma \vert$.  

Recall that the vertices of $\Gamma$ are indexed by the set $\{1, \ldots, n\} \setminus J$.  Since there are no isolated vertices in $\Gamma$, $\vert \pi_0 \Gamma \vert \le \lfloor (n- \vert J \vert ) / 2 \rfloor$ (where $\lfloor - \rfloor$ denotes the integer floor function).  Recalling that $0 \le \vert J \vert \le n-2$, the set $\{1, \ldots, n\} \setminus J$ is necessarily non-empty so $1 \le \vert \pi_0 \Gamma \vert $, combining these gives
\[
1 \le \vert \pi_0 \Gamma \vert \le \lfloor (n-\vert J \vert)/2 \rfloor.
\]

We can now completely describe the centraliser $C_{H_n}(q)$.

\begin{Theorem}\label{thm:stabOfInfiniteQ}~
\begin{enumerate}
 \item If $q \in H_n$ is an element of finite order then 
 \[C_{H_n}(q) \cong H_n\vert_{S^Q} \times C_1 \times \cdots \times C_t\]
 where $H_n\vert_{S^Q} \cong H_n$ is Houghton's group restricted to $S^Q$ and for all $a \in \{1, \ldots, t\}$, 
 \[ C_a \cong W_QQ_a \wr \Sym_r \]
 for $Q_a$ an isotropy group of $S_a$ and $r = \vert S_a \vert / \vert Q : Q_a \vert$.  In particular $H_n$ is finite index in $C_{H_n}Q$.
 \item If $q \in H_n$ is an element of infinite order then either
\[C_{H_n}(q) \cong H_k \times \ZZ^r \times C_1 \times \cdots \times C_t\]
 or
\[C_{H_n}(q) \cong F \times \ZZ^r \times C_1 \times \cdots \times C_t\]
 where $F$ is some finite symmetric group, $H_k$ is Houghton's group with $0 \le k \le n-2$, and the groups $C_a$ are as in the previous part.   In the first case $1 \le r \le \lfloor (n-k) / 2 \rfloor$, and in the second case $1 \le r \le \lfloor n / 2 \rfloor$.  
\end{enumerate}
\end{Theorem}

 In Corollary \ref{cor:stabOfH is FPn-1 not FPn} it was proved that for an element $q$ of finite order, $C_{H_n}(q)$ is $\FP_{n-1}$ but not $\FP_n$.  The situation is much worse for elements $q$ of infinite order, in which case the centraliser may not even be finitely generated, for example when $n$ is odd and $q$ is the element acting on $S = \NN \times \{1, \ldots, n\}$ as
\[
q :\left\{ \begin{array}{l l}
  (i,x) \mapsto (i+1,x) & \text{ if }x \le (n-1)/2 \\
  (i,x) \mapsto (i-1,x) & \text{ if }(n+1)/2 \le x \le n-1 \text{ and } i \neq 0\\
  (0,x) \mapsto (0,x - ((n-1)/2)) & \text{ if } (n+1)/2 \le x \le n-1 \\
  (i,n) \mapsto (i,n)
   \end{array} \right.
\]
then the only fixed points are on the ray $\NN \times \{n \}$.  The argument leading up to Theorem \ref{thm:stabOfInfiniteQ} shows that the centraliser is a direct product of groups, one of which is Houghton's group $H_1$ which is isomorphic to the infinite symmetric group and hence not finitely generated.  In particular for this $q$, the centraliser $C_{H_n}(q)$ is not even $\FP_1$.  A similar example can easily be constructed when $n$ is even.

All the groups in the direct product decomposition from Theorem \ref{thm:stabOfInfiniteQ} except $H_k$ are $\FP_\infty$, being built by extensions from finite groups and free Abelian groups.  By choosing various infinite order elements $q$, for example by modifying the example of the previous paragraph, the centralisers can be chosen to be $\FP_k$ for $0 \le k \le n-3$.  The upper bound of $n-3$ arises because any infinite order element $q$ must necessarily be ``eventually a translation'' (in the sense of \eqref{eventuallyATranslation}) on  $\NN \times \{x \} $ for \emph{at least} two $x$.  As such the copy of Houghton's group in the centraliser can act on at most $n-2$ rays and is thus at largest $H_{n-2}$, which is $\FP_{n-3}$.  

\begin{Cor}\label{cor:stabOfVCycSubgroups}
 If $Q$ is an infinite virtually cyclic subgroup of $H_n$ then either
\[ C_{H_n}(Q) \cong H_k \times \ZZ^r \times C_1 \times \cdots \times C_t \]
or 
\[ C_{H_n}(Q) \cong F \times \ZZ^r \times C_1 \times \cdots \times C_t \]
where the elements in the decomposition are all as in Theorem \ref{thm:stabOfInfiniteQ}.
\end{Cor}

This corollary can be proved by reducing to the case of Theorem \ref{thm:stabOfInfiniteQ}, but before that we require the following lemma.
\begin{Lemma}
 Every infinite virtually cyclic subgroup $Q$ of $H_n$ is finite-by-$\ZZ$.
\end{Lemma}
\begin{proof}
By \cite[Proposition 4]{JuanPinedaLeary-ClassifyingSpacesForVCsubgrps}, $Q$ is either finite-by-$\ZZ$ or finite-by-$\text{D}_\infty$ where $\text{D}_\infty$ denotes the infinite dihedral group, we show the latter cannot occur.  Assume that there is a short exact sequence of groups
\[
1 \longrightarrow F \longhookrightarrow Q \stackrel{\pi}{\longrightarrow} \text{D}_\infty \longrightarrow 1,
\]
regarding $F$ as a subgroup of $Q$.  Let $a,b$ generate $\text{D}_\infty$, so that 
\[
{\text{D}_\infty = \langle a,b \mid a^2 = b^2 = 1\rangle}.
\]
Let $p,q \in Q$ be lifts of $a, b$, such that $\pi(p) = a$, $\pi(q) = b$, then $p^2 \in F$.  Since $F$ is finite, $p^2$ has finite order and hence $p$ has finite order.  The same argument shows that $q$ has finite order.  The element $pq \in Q$ necessarily has infinite order as $\pi(pq)$ is infinite order in $D_\infty$.

However, since $p$ and $q$ are finite order elements of $H_n$, by the argument at the beginning of Section \ref{section:centralisersOfFiniteSubgroups} they both permute only a finite subset of $S$.  Thus $pq$ permutes a finite subset of $S$ and is of finite order, but this contradicts the previous paragraph.

\end{proof}
\begin{proof}[Proof of Corollary \ref{cor:stabOfVCycSubgroups}]
Using the previous lemma, write $Q$ as $Q = F \rtimes \ZZ$ where $F$ is a finite group.  As $F$ is finite, the set $S_F$ of points not fixed by $F$ is finite (see the argument at the beginning of Section \ref{section:centralisersOfFiniteSubgroups}).  Let $z \in \NN$ be such that for $i \ge z$, $F$ acts trivially on $(i,x)$ for all $x$, and $\ZZ$ acts on $(i,x)$ either trivially or as a translation.  Applying Lemma \ref{lemma:prop conditions are satisfied} and Proposition \ref{prop:groupActionsOnSets}, $S$ splits as a disjoint union 
\[
S = S^Q \cup S_0 \cup S_1 \cup \cdots \cup S_t
\]
where $S^Q$ is the fixed point set, $S_0$ is the set with isotropy group $F$ and the $S_a$ for $1 \le a \le t$ are subsets of $\{(i,x)\: :  \: i \le z \}$, and hence all finite.  By Proposition \ref{prop:groupActionsOnSets}, $C_{H_n}(Q)$ splits as a direct product
\[
C = H_n|_{S^Q} \times C_0 \times C_1 \times \ldots \times C_t 
\]
where $H_n|_{S^Q}$ denotes Houghton's group restricted to $S^Q$.  The argument of Theorem \ref{thm:stabOfInfiniteQ} showing that $H_n|_{S^Q}$ is isomorphic to either a finite symmetric group or to $H_k$ for some $0 \le k \le n-2$ goes through with no change, as does the proof of the structure of the groups $C_a$ for $1 \le a \le t$.  It remains to observe that because every element in $S_0$ is fixed by $F$, any element of $H_n$ centralising $\ZZ$ and fixing $S \setminus S_0$ necessarily also centralises $Q$ and is thus a member of $C_0$.  This reduces us again to the case of Theorem \ref{thm:stabOfInfiniteQ} showing that $C_0 \cong \ZZ^r$ for some natural number $1 \le r \le \lfloor (n-k)/2 \rfloor$, or $1 \le r \le \lfloor n /2 \rfloor$ if $H_n\vert_{S^Q}$ is a finite symmetric group.
\end{proof}

\section{Brown's model for \texorpdfstring{$\EFin\Hn$}{the classifying space for proper actions of Hn}}\label{section:Browns model}

The main result of this section will be Corollary \ref{cor:Browns model is an EunderbarH}, where the construction of Brown \cite{Brown-FinitenessPropertiesOfGroups} used to prove that $H_n$ is $\FP_{n-1}$ but not $\FP_n$ is shown to be a model for $\EFin H_n$.

In this section, maps are written from left to right.

Write $\mathcal{M}$ for the monoid of injective maps $S \to S$ with the property that every permutation is ``eventually a translation'' (in the sense of \eqref{eventuallyATranslation}), and write $T$ for the free monoid generated by $\{t_1, \ldots, t_n\}$ where 
\[
(i,x)t_y = \left\{ \begin{array}{c c}(i+1,x) & \text{ if } x = y, \\ (i,x) & \text{ if } x \neq y. \end{array}\right. 
\]
The elements of $T$ will be called \emph{translations}.  The map $\phi : H_n \to \ZZ^n$, defined in \eqref{def:phi}, extends naturally to a map $\phi : \mathcal{M} \to \ZZ^n$.  
Give  $\mathcal{M}$ a poset structure by setting $\alpha \le \beta$ if $\beta = t\alpha$ for some $t \in T$.  The monoid $\mathcal{M}$ can be given the obvious action on the right by $H_n$, which in turn gives an action of $H_n$ on the poset $(\mathcal{M}, \le)$ since $\beta = t\alpha$ implies $\beta h = t\alpha h$ for all $h \in H_n$.  Let $\vert \mathcal{M} \vert$ be the geometric realisation of this poset, namely simplicies in $\vert \mathcal{M} \vert$ are finite ordered collections of elements in $\mathcal{M}$ with the obvious face maps. 
An element $h \in H_n$ fixes a vertex $\{\alpha\} \in \vert \mathcal{M} \vert$ if and only if $s \alpha h = s \alpha$ for all $s \in S$ if and only if $h$ fixes $S\alpha$, so the stabiliser $(H_n)_\alpha$ may only permute the finite set $S \setminus S\alpha$ and we may deduce:
\begin{Prop}\label{prop:vmvHasFiniteIsotropy}
Stabilisers of simplicies in $\vmv$ are finite.
\end{Prop}
We now build up to the the proof that $\vmv$ is a model for $\EFin H_n$ with a few lemmas.

\begin{Prop}\label{prop:vmvQfinite=>contractible}
If $Q \le H_n$ is a finite group then the fixed point set $\vmv^Q$ is non-empty and contractible.
\end{Prop}
\begin{proof} 
For all $q \in Q$, choose $\{z_0(q), \ldots, z_n(q)\}$ to be an $n$-tuple of natural numbers such that $(i,x)q = (i,x)$ whenever $i \ge z_x(q)$ for all $i$.  $Q$ then fixes all elements $(i,x) \in S$ with $i \ge \max_Q z_x(q)$.  Define a translation $t = t_1^{\max_Q z_1(q)} \cdots t_n^{\max_Q z_n(q)}$, $t \in \mathcal{M}^Q$ so $\{t \}$ is a vertex of $\vmv^Q$ and  $\vmv^Q \neq \emptyset$.  

If $\{m\}, \{n\} \in \vmv^Q$ then let $a, b \in T$ be two translations such that 
\[
\phi(m) - \phi(n) = \phi(b) -  \phi(a)
\]
(recall that for a translation $t$, $\phi(t)$ must be an $n$-tuple of positive numbers).  Thus $\phi(am) = \phi(bn)$, and since $am, bn \in \mathcal{M}$ there exist $n$-tuples $\{z_1, \ldots, z_n\}$ and $\{z_1^\prime, \ldots, z_n^\prime\}$ such that $am$ acts as a translation for all $(i,x) \in S$ with $i \ge z_x$ and $bn$ acts as a translation for all $(i,x) \in S$ with $i \ge z_x^\prime$.  Let 
\[
c = t_1^{\max\{z_1, z_1^\prime\}} \ldots t_n^{\max\{z_n, z_n^\prime\}}
\]
so that $cam = cbn$, further pre-composing $c$ with a large translation (for example that from the first section of this proof) we can assume that $cam = cbn \in \mathcal{M}^Q$, and $\{cam=cbn\} \in \, \vmv^Q$.  This shows the poset $\mathcal{M}^Q$ is directed and hence the simplicial realisation $\vert \mathcal{M}^Q \vert = \vmv^Q$ is contractible. 
\end{proof}
\begin{Prop}\label{prop:vmvQinfinite=>empty}
 If $Q \le H_n$ is an infinite group then $\vmv^Q = \emptyset$.
\end{Prop}
\begin{proof}
 Consider an infinite subgroup $Q \le H_n$ with $\vmv^Q \neq \emptyset$ and choose some vertex ${\{m\} \in \vmv^Q}$.  For any $q \in Q$, since $mq = m$ it must be that $\phi(m)+ \phi(q) = \phi(m)$ and so $\phi(q) = 0$, hence $Q$ is a subgroup of $\Sym_\infty \le H_n$.  Furthermore $Q$ must permute an infinite subset of $S$ (if it permuted just a finite set it would be a finite subgroup).  That $mq = m$ implies that this infinite subset is a subset of $S\setminus Sm$ but this is finite by construction.  So the fixed point subset $\vmv^Q$ for any infinite subgroup $Q$ is empty.
\end{proof}

\begin{Cor}\label{cor:Browns model is an EunderbarH}
 $\vmv$ is a model for $\EFin H_n$.
\end{Cor}
\begin{proof}
Combine Propositions \ref{prop:vmvHasFiniteIsotropy}, \ref{prop:vmvQfinite=>contractible} and \ref{prop:vmvQinfinite=>empty}.
\end{proof}

\section{Finiteness conditions satisfied by \texorpdfstring{$\Hn$}{Hn}} \label{section:finitenessConditionsSatisfied}

Recall from Proposition \ref{prop:uFP_0 iff finitely many conj classes of finite subgroups} that a group $G$ is $\OFinFP_0$ if and only if it has finitely many conjugacy classes of finite subgroups.  $G$ satisfies the weaker quasi-$\OFinFP_0$ condition if and only if it has  finitely many conjugacy classes of subgroups isomorphic to a given finite subgroup (see Section \ref{subsection:bredon quasi-FPn}).
\begin{Prop}\label{prop:HisnotQuasiFP0}
 $H_n$ is not quasi-$\OFinFP_0$.
\end{Prop}
Before the above proposition is proved, we need a lemma.  In the infinite symmetric group $\Sym_\infty$ acting on the set $S$, elements can be represented by products of disjoint cycles.  We use the standard notation for a cycle: $(s_1, s_2, \ldots, s_m)$ represents the element of $\Sym_\infty$ sending $s_i \mapsto s_{i+1}$ for $i < n$ and $s_n \mapsto s_1$.  Any element of finite order in $H_n$ is contained in the infinite symmetric group $\Sym_\infty$ by the argument at the beginning of Section \ref{section:centralisersOfFiniteSubgroups}.  We say two elements of $\Sym_\infty$ have the same \emph{cycle type} if they have the same number of cycles of length $m$ for each $m \in \NN$.

\begin{Lemma}\label{lemma:cycle type in H}
If $q$ is a finite order element of $H_n$ and $h$ is an arbitrary element of $H_n$, then $hqh^{-1}$ is the permutation given in the disjoint cycle notation by applying $h$ to each element in each disjoint cycle of $q$.  In particular, if $q$ is represented by the single cycle $(s_1, \ldots s_m)$, then $hqh^{-1}$ is represented by $(hs_1,\ldots,hs_m)$.

Furthermore, two finite order elements of $H_n$ are conjugate if and only if they have the same cycle type.
\end{Lemma}
\begin{proof} The proof of the first part is analogous to \cite[Lemma 3.4]{Rotman-Groups}.  Let $q$ be an element of finite order and $h$ an arbitrary element of $H_n$.  If $q$ fixes $s \in S$ then $hqh^{-1}$ fixes $hs$.  If $q(i) = j$, $h(i) = k$ and $h(j) = l$, for $i,j,k,l \in S$, then $hqh^{-1}(k) = l$ exactly as required.

By the above, conjugate elements have the same cycle type.  For the converse, notice any two finite order elements with the same cycle type necessarily lie in $\Sym_r$ for some $r \in \NN$ so by \cite[Theorem 3.5]{Rotman-Groups} they are conjugated by an element of $\Sym_r$.
\end{proof}

\begin{proof}[Proof of Proposition \ref{prop:HisnotQuasiFP0}]
Choosing a collection of elements $q_i$ for each $i \in \NN_{\ge 1}$, so that $q_i$ has $i$ disjoint 2-cycles gives a collection of isomorphic subgroups which are all non-conjugate by Lemma \ref{lemma:cycle type in H}.
\end{proof}

\begin{Prop}\label{prop:cd=gd=n} $\OFincd H_n = \gdFin H_n = n$.
\end{Prop}
\begin{proof}As described in the introduction, $H_n$ can be written as 
\[\Sym_\infty \longhookrightarrow H_n \longtwoheadrightarrow \ZZ^{n-1}.\]
Now, $\gdFin \ZZ^{n-1} = \gd \ZZ^{n-1} = n-1$ and $\gdFin \Sym_\infty = 1$ by \cite[Theorem 4.3]{LueckWeiermann-ClassifyingspaceForVCYC}, as it is the colimit of its finite subgroups each of which have proper geometric dimension $0$, and the directed category over which the colimit is taken has homotopy dimension $1$ \cite[Lemma 4.2]{LueckWeiermann-ClassifyingspaceForVCYC}.  $\ZZ^{n-1}$ is torsion free and so has a bound of $1$ on the orders of its finite subgroups and we deduce from \cite[Theorem 3.1]{Luck-TypeOfTheClassifyingSpace} that $\gdFin H_n \le n-1 + 1 = n$.

To deduce the other bound, we use an argument due to Gandini \cite{Gandini-BoundingTheHomologicalFinitenessLength}.  Assume that $\OFincd H_n \le n-1$.  By \cite[Theorem 2]{BradyLearyNucinkis-AlgAndGeoDimGroupsWithTorsion} we have
\[\cd_\QQ H_n \le \OFincd H_n = n-1. \]
In \cite[Theorem 5.1]{Brown-FinitenessPropertiesOfGroups}, it is proved that $H_n$ is $\FP_{n-1}$ (but not $\FP_n$), combining this with \cite[Proposition 1]{LearyNucinkis-BoundingOrdersOfFiniteSubgroups} we deduce that there is a bound on the orders of the finite subgroups of $H_n$, but this is obviously a contradiction, thus 
\[n \le \OFincd H_n \le \gdFin H_n \le n. \]
\end{proof}

The remainder of this section is devoted to proving the following.  

\begin{Theorem}\label{theorem:houghton cdVCyc Hn = n}
 $\OVCyccd H_n = n$.
\end{Theorem}

The proof is based on a pushout of L\"uck and Weiermann \cite{LueckWeiermann-ClassifyingspaceForVCYC}, described below.

\subsection{The pushout of L\"uck and Weiermann}\label{subsection:houghton cdVCyc LW}

For any group $G$, we say two infinite virtually cyclic subgroups $K$ and $K^\prime$ of $G$ are \emph{commensurate}, written $K \sim K^\prime$, if $\lvert K \cap K^\prime \rvert = \infty$.  Commensurability is an equivalence relation and we write $[\VCyc \setminus \Fin]$ for the set of equivalence classes.  The \emph{normaliser} of an equivalence class $[K]$ is defined to be the stabiliser of the action of $G$ on $[\VCyc \setminus \Fin]$ by conjugation, namely
\[ N_G[K] = \{ x \in G : K^x \sim K \}. \]
Associated to each infinite virtually cyclic subgroup $K$ we define the subfamily $\VCyc[K]$ of $\VCyc$ by
\[\VCyc[K] = \{ L \in \VCyc \setminus \Fin : L \sim K \} \cup \left( \Fin \cap K \right). \]

If $X$ is a right $G$-space and $Y$ is a left $G$-space then we denote by $X \times_G Y$ the \emph{twisted product} of $X$ and $Y$, defined to be the quotient space of $X \times Y$ under the action $g \cdot (x,y) = (xg^{-1}, gy)$ \cite[\S II.2]{Bredon-IntroductionToCompactTransformationGroups}.  If $Y$ is a left $H$-space for some subgroup $H$ of $G$ then $G \times_H Y$ is a left $G$-space via the usual left action of $G$ on itself.
\begin{Prop}[{\cite[Proposition I.(4.3)]{Dieck-TransformationGroups}}]\label{prop:adjoint iso for twisted product}
Let $H$ be a subgroup of $G$, let $Y$ be a left $H$-space, and let $X$ be a left $G$-space.  There is an adjoint isomorphism
\[
[G \times_H Y, X]_G = [Y, X]_H,
\]
where $[Z, X]_G$ denotes the set of $G$-homotopy classes of $G$-equivariant maps between two $G$-spaces $Z$ and $X$. 
\end{Prop}

\begin{Theorem}\cite[Theorem 2.3, Remark 2.5]{LueckWeiermann-ClassifyingspaceForVCYC}
Let $I$ denote a complete set of representatives of the $G$-orbits in $[\VCyc \setminus \Fin]$.  Choosing arbitrary $N_G[K]$-CW-models for $E_{\Fin} N_G[K]$ and $E_{\VCyc[K]} N_G[K]$ and an arbitrary $G$-CW-model for $E_{\Fin}G$, the cellular $G$-pushout described below may be constructed with the maps $i$ and $f_{[K]}$ equivariant cellular maps, and  either with $i$ an inclusion of $G$-CW-complexes or with every $f_{[K]}$ an inclusion of $N_G[K]$-CW-complexes and $i$ a cellular $G$-map.
 \[
 \xymatrix{
  \coprod\limits_{[K] \in I} G \times_{N_G[K]} E_{\Fin} N_G[K] \ar^-i[r] \ar^-{\coprod_{[K]\in I} \id_G \times_{N_G[K]}f_{[K]}}[d] & E_{\Fin}G \ar[d] \\
  \coprod\limits_{[K] \in I} G \times_{N_G[K]} E_{\VCyc[K]} N_G[K] \ar[r] & X \\  
 } 
\]  
Moreover the space $X$ defined by the pushout is a model for $\E_{\VCyc}\negthinspace G$.
\end{Theorem}

We can describe explicitly the $G$-homotopy classes of the maps $i$ and $f_{[K]}$ in the pushout above:  By restricting the $G$-action, any model for $E_{\Fin} G$ is a model for $E_{\Fin} N_G[K]$ so there is a $N_G[K]$-map $E_{\Fin} N_G[K] \to E_{\Fin} G$, and using the adjoint isomorphism of Proposition \ref{prop:adjoint iso for twisted product} there is a $G$-map $ G \times_{N_G[K]} E_{\Fin} N_G[K] \to E_{\Fin}G$.  The coproduct of these maps, one for each $[K] \in I$, is the map $i$.  Since $E_{\Fin}N_G[K]$ is an $N_G[K]$-space with finite isotropy, it is a priori an $N_G[K]$-space with isotropy in $\VCyc[K]$, there is a map $E_{\Fin}N_G[K] \to E_{\VCyc[K]} N_G[K] $.  This is the map $f_{[K]}$.

This pushout gives a long exact sequence in Bredon cohomology \cite[Lemma 13.7]{Lueck}:
 \[
 \cdots \longrightarrow H^i_{\mathcal{O}_{\VCyc}}(G, -) \longrightarrow 
 \]
 \nopagebreak
 \[
 \left( \prod\limits_{[K] \in I} H^i_{\mathcal{O}_{\VCyc[K]}}( N_{G}[K], \Res^{\mathcal{O}_{\VCyc}G}_{\mathcal{O}_{\VCyc[K]}N_G[K]} -) \right) \oplus 
H^i_{\OFin}(G, \Res^{\mathcal{O}_{\VCyc}G}_{\mathcal{O}_{\Fin}G} -)  
\]
\nopagebreak
\[
\longrightarrow \prod\limits_{[K] \in I} H^{i}_{\OFin} (N_{G}[K], \Res^{\mathcal{O}_{\VCyc}G}_{\mathcal{O}_{\Fin}N_G[K]} -) \longrightarrow \cdots.
\]
For brevity we will usually omit the restriction maps from now on.

Given an infinite virtually cyclic subgroup $K$ of $G$, let $\pi^K: N_GK \to WK$ denote the projection map and for any $\mathcal{O}_{\VCyc[K]}N_GK$-module let $\pi^K_* M$ denote the $\OFin WK$-module given by
\[
 \pi^K_* M : WK / L \mapsto M(N_GK / \pi^{-1}(L)),
\]
for any finite subgroup $L$ of $WK$.

\begin{Lemma}\cite[Lemma 4.2]{DegrijsePetrosyan-GeometricDimensionForVCYC}
 If $K$ is an infinite virtually cyclic subgroup and $N_G[K] = N_GK$ then there is an isomorphism,
  \[
  H^i_{\mathcal{O}_{\VCyc[K]}}( N_G K, -) \cong H^i_{\OFin}( W_{G}K, \pi^K_* -).
  \]
\end{Lemma}
 
Combining this lemma with the long exact sequence gives the following.

\begin{Prop}\label{prop:houghton usable LES for uucd}
If every $G$-orbit in $[\VCyc \setminus \Fin]$ contains a $K$ such that $N_{G}[K] = N_GK$ then, letting $A$ be a set of representatives with that property, there is a long exact sequence:
\[
 \cdots \longrightarrow H^i_{\mathcal{O}_{\VCyc}}(G, -) \longrightarrow \left( \prod\limits_{K \in A} H^i_{\OFin}( W_{G}K, \pi^K_* -) \right) \oplus H^i_{\OFin}(G, -)  
\]
\[
 \longrightarrow \prod\limits_{K \in A} H^{i}_{\OFin} (N_{G}K, -) \longrightarrow \cdots.
\]
\end{Prop}

\subsection{Calculation of \texorpdfstring{$\OVCyccd H_n$}{cdVCyc Hn}}\label{subsection:houghton cdVCyc calc}
 
\begin{Lemma}\label{lemma:houghton cdVCyc exists L comm K with nice normaliser}
 For every infinite virtually cyclic subgroup $K$ of $H_n$ there exists $L$ commensurate to $K$ with 
 \[ N_G[K] = N_G[L] = N_GL. \]
 Moreover, we may assume $L \cong \ZZ$.
\end{Lemma}
\begin{proof}
 Firstly replace $K$ with a finite-index subgroup isomorphic to $\ZZ$ and let $k$ generate $K$.  Consider the action of $K$ on $S$.  There are finitely many finite $K$-orbits (see Section \ref{section:centralisersOfElements}), so for large enough $m$, the subgroup $ \langle k^m \rangle$ acts semifreely (freely away from the fixed point set), let $L = \langle k^m \rangle$.
 
 If $S^L \neq \emptyset$ then let $s \in S^L$ and pick any $n \in N_{H_n}[L]$, so the stabiliser of $ns$ is $L^n$.  Since $L$ acts semi-freely $L^n = L$ or $L^n = 1$ but since $L^n \sim L$ this forces $L^n = L$, thus $N_{H_n}[L] = N_{H_n}L$. 
 
 If $S^L = \emptyset$ then let $n \in N_{H_n}[L]$, let $l \in L$, and let $i \in \{1, \ldots, n\}$.  So for $x$ large enough,
 \[ l : (i, x) \mapsto (i, x+t_l) \]
 \[ n : (i, x) \mapsto (i, x+t_n). \]
 For some $t_l, t_n \in \ZZ$ with $t_l \neq 0$ ($L$ acts non-trivially everywhere),
 \[ n^{-1}ln : (i,x) \mapsto (i, x +  t_n  + t_l - t_n ) = (i, x + t_l)\]
 so $n^{-1}ln$ acts as $L$ does on all $(i,x)$ for large enough $x$, in particular for all but finitely many elements of $S$.  Since all orbits are infinite this means $n^{-1}ln$ acts as $L$ does on all of $S$.  Hence $n^{-1}ln = l$, in particular $n \in N_GL$.
\end{proof}

Let $A$ denote a set of representatives of $H_n$-orbits in $[\VCyc \setminus \Fin]$ such that for all $L \in A$, we have $N_G[L] = N_GL$ and $L \cong \ZZ$.

\begin{Lemma}\label{lemma:Hn exists K with NK equal CK}
 For any $K \in A$ we have $N_{H_n}K \cong C_{H_n}K$.
\end{Lemma}
\begin{proof}
 Recall that there is a short exact sequence
 \[ 1 \longrightarrow C_{H_n}K \longrightarrow N_{H_n}K \longrightarrow Q \longrightarrow 1 \]
where $Q$ is a subgroup of $\Aut(K)$ \cite[1.6.13]{Robinson}. 

 Let $n \in N_{H_n}K$ and choose some $k \in K$ generating $K$.  Assume that $K$ acts non-trivially on the $i^\text{th}$ ray, so for $x$ large enough,
 \[ k : (i, x) \mapsto (i, x+t_k) \]
 \[ n : (i, x) \mapsto (i, x+t_n) \]
 for some $t_k, t_n \in \ZZ$ with $t_k \neq 0$.  Let $a \in \ZZ$ be such that $n^{-1}kn = k^a$, then
 \[ n^{-1}kn : (i, x) \mapsto (i, x+t_k)  \]
 but
 \[ k^a : (i, x) \mapsto (i, x+at_k).  \]
 Thus $a = 1$ and $n$ acts as the trivial automorphism on $K$, thus $Q = 1$, proving the lemma.
\end{proof}

\begin{Lemma}\label{lemma:cd_Fin NHK le n-1}
If $K \in A$ then $\OFincd N_{H_n}K \le n-1$ and $\OFincd W_{H_n}K \le n-2$.
\end{Lemma}
\begin{proof}
Recall from Corollary \ref{cor:stabOfVCycSubgroups} that 
\[
C_{H_n}K \cong H_k \times \ZZ^r \times F,
\]
where $F$ is a finite group, $0 \le k \le n-2$, and $1 \le r \le \lfloor (n-k)/2\rfloor$.  Thus,
\begin{align*}
 \OFincd N_{H_n}K &= \OFincd C_{H_n}K \\
  &= \OFincd (H_k \times \ZZ^r \times F) \\
  &\le \OFincd H_k + \OFincd \ZZ^r + \OFincd F \\
  &=k+r \\
  &\le \max_{1 \le k \le n-2} \left( k + \lfloor (n-k)/2 \rfloor \right),
\end{align*}
where we've used Proposition \ref{prop:cd=gd=n} that $\OFincd H_k = k$, Lemma \ref{lemma:Hn exists K with NK equal CK}, and Lemma \ref{lemma:bredon G uFP and H uFP then GxH is uFP}.  We claim that $\max_{1 \le k \le n-2} k + \lfloor (n-k)/2 \rfloor  = n-1$, indeed we can always achieve $n-1$ by choosing $k = n-2$ and $k + \lfloor (n-k)/2 \rfloor$ is an increasing function of $k$.
 
Examining the proof of Corollary \ref{cor:stabOfVCycSubgroups} that $C_{H_n}K \cong H_k \times \ZZ^r $, we see that $K$ is a subgroup of $\ZZ^r$, so 
\[
W_GH \cong H_k \times \ZZ^{r-1} \times F^\prime,
\]
for some finite subgroup $F^\prime$, which gives the second inequality.
\end{proof}

\begin{proof}[Proof of \ref{theorem:houghton cdVCyc Hn = n}]
Via Lemma \ref{lemma:houghton cdVCyc exists L comm K with nice normaliser} we have the long exact sequence of Proposition \ref{prop:houghton usable LES for uucd},
\[ \cdots \longrightarrow \prod\limits_{K \in A} H^{i-1}_{\OFin} (N_{H_n}K, -) \longrightarrow H^i_{\mathcal{O}_{\VCyc}}(H_n, -) \]
\[\longrightarrow \left( \prod\limits_{K \in A} H^i_{\OFin} (W_{H_n}K, \pi^K_* -) \right) \oplus H^i_{\OFin}(H_n, -) \longrightarrow \cdots. \]

Let $i = n+1$ then using Lemma \ref{lemma:cd_Fin NHK le n-1}, the left and right hand terms vanish.  Thus the central term vanishes proving $\OVCyccd H_n \le n$.

Let $i = n$ then, using Lemma \ref{lemma:cd_Fin NHK le n-1} again, there is a long exact sequence which terminates as
\[ \cdots \longrightarrow H^n_{\mathcal{O}_{\VCyc}}(H_n, -) \longrightarrow H^n_{\OFin}(H_n, -) \longrightarrow 0. \] 

Let $M$ be an $\OFin$-module such that $H^n_{\OFin}(H_n, M) \neq 0$ then we may extend $M$ to an $\mathcal{O}_{\VCyc}$-module by setting $M(G/K) = 0$ for all virtually cyclic subgroups $K$ and thus $H^n_{\mathcal{O}_{\VCyc}}(H_n, M) \neq 0$.  In particular, $\OVCyccd H_n \ge n$.
\end{proof}

\backmatter
\bibliographystyle{amsalpha}
\bibliography{thesis.bib}

\cleardoublepage

\printindex

\end{document}